\documentclass[10pt]{amsart}
\usepackage{amssymb,amsmath,amsfonts,amsthm,graphics,mathrsfs,amscd,hyperref}
\usepackage[makeroom]{cancel}
\usepackage[hmargin=1in,vmargin=1in]{geometry}
\usepackage[all,cmtip]{xy}

\theoremstyle{plain}

\newtheorem{definition}[equation]{Definition}

\newtheorem{lemma}[equation]{Lemma}
\newtheorem{proposition}[equation]{Proposition}
\newtheorem{theorem}[equation]{Theorem}

\newtheorem{notation}{Notation}
\newtheorem{remark}[equation]{Remark}

\numberwithin{equation}{subsection}
\newtheorem{conjecture}{Conjecture}

\title{On deformations of foliated complex analytic structures III : Unfoldings}
\author{Chunghoon Kim}

\begin{document}

\maketitle
\begin{abstract}
We study unfoldings of foliated complex analytic structures defined by locally free subsheaves of cotangent sheaves. We prove theorem of existence for unfoldings of foliated complex analytic structures defined by locally free subsheaves of cotangent sheaves for one-parameter families. We prove theorem of completeness for unfoldings of foliated complex analytic structures defined by locally free subsheaves of cotangent sheaves. We also present unfoldings of foliated nonsingular varieties defined by locally free subsheaves of cotangent sheaves in the language of functors of Artin rings which is the algebraic version of unfoldings of foliated complex analytic structures. We identify first-order unfoldings and obstructions.
\end{abstract}

\tableofcontents

Before proceeding the discussion on unfoldings of singular holomorphic foliations on compact complex manifolds defined by locally free subsheaves of cotangent sheaves, the author would like to remark that he could not succeed in proving `theorem of existence for unfoldings' in the complete form (see Conjecture \ref{conjecture1}). We only prove theorem of existence for unfoldings for one-parameter families (see Theorem \ref{u56}). The author had difficulties in dealing with multi-parameters. However the complete form of theorem of existence for unfoldings is fundamental for the general framework of deformation theory in the spirit of Kodaira-Spencer.  By the same reason, given an arbitrary small extension $0\to (t)\to \tilde{A} \to A \to 0$ in the language of functors of Artin rings, the author had difficulties of identifying obstructions of liftings of unfoldings of $\left(X, \mathcal{N}_{\mathcal{F}_0}^* \right)$ over $A$ to $\tilde{A}$ with elements in the second cohomology group $\mathbb{H}^2\left(X, \tilde{\mathcal{N}}_{\mathcal{F}_0}^{*\bullet} \right)$, where $\left(X, \mathcal{N}_{\mathcal{F}_0}^*\right)$ is a compact foliated complex manifold defined by a locally free susbsheaf $\mathcal{N}_{\mathcal{F}_0}^*$ of $\Omega_X^1$, and the complex $\tilde{\mathcal{N}}_{\mathcal{F}_0}^{*\bullet}$ is the unfolded dual leaf complex associated to $(X, \mathcal{N}_{\mathcal{F}_0}^*)$ (see Appendix \ref{AppendixA}). That's why the author confine our discussion to small extensions of the form $0\to \left(t^m\right)\to k[t]/\left(t^{m+1}\right)\to k[t]/\left(t^m\right)\to 0$ involving only one parameter when we identify obstructions. On the other hand, we prove theorem of completeness for unfoldings in the complete form (see Theorem \ref{uc1}), which also implies the rigidity of unfoldings of $\mathbb{H}^1\left( X, \tilde{\mathcal{N}}_{\mathcal{F}_0}^{*\bullet} \right)=0$ (for the direct proof, see Part IV-2), which is stronger than rigidity of infinitesimal unfoldings over Artin rings, but the author could not prove the rigidity of infinitesimal unfoldings of $(X, \mathcal{N}_{\mathcal{F}_0}^* )$ with $\mathbb{H}^1\left( X, \tilde{\mathcal{N}}_{\mathcal{F}_0}^{*\bullet} \right)=0$ over Artin rings by the method (employed from infinitesimal deformations of nonsingular varieties, and Poisson deformations, and foliated deformations over Artin rings) which is based on the transitive action of first-order deformations on the liftings, which does not seem to necessarily hold for unfoldings for arbitrary small extensions. However infinitesimal rigidity of $(X, \mathcal{N}_{\mathcal{F}_0}^*)$ of $\mathbb{H}^1\left(X, \tilde{\mathcal{N}}_{\mathcal{F}_0}^{*\bullet}  \right)=0$ over Artin rings should be true in the light of theorem of completeness for unfoldings. That is one of the reasons that  the author believes that the general framework of deformation theory in complex algebraic geometry might hold for unfoldings (i.e. the complete form of theorem of existence for unfoldings holds for compact foliated complex manifold $\left( M, \mathcal{N}_{\mathcal{F}_0}^* \right)$ with $\mathcal{N}_{\mathcal{F}_0}^*$ locally free and $\mathbb{H}^2\left( M, \tilde{\mathcal{N}}_{\mathcal{F}_0}^{*\bullet} \right)=0$), and he is missing something. On the other hand, there is the concept of Kuranishi families over (possibly) singular bases in deformation theory, and unobstructedness is equivalent to the smoothness of the bases. However, in the case of unfoldings, the author is not confident with how to define the concept of  unfolded families over singular bases which generalize the concept of unfolded analytic families over smooth bases that we study in the main body of Part III. In Appendix \ref{b2} we only define infinitesimal unfoldings over local Artin rings instead of general parameter spaces. Is theory of  unfoldings different since parameter spaces are also involved unlike deformations of complex structures, holomorphic Poisson structures, and foliated complex analytic structures? 

\section{Unfoldings of foliated complex analytic structures}

\begin{definition}[compare \cite{Suw83} p.185 and Part I] \label{u1}
Suppose that given a domain $B\subset \mathbb{C}^m$, there is a set $\left\{\left(M_t, \mathcal{N}_{\mathcal{F}_t}^*\right)|t\in  B \right\}$ of $n$-dimensional compact $($singularly$)$ foliated complex manifolds $\left(M_t, \mathcal{N}_{\mathcal{F}_t}^*\right)$ of codimension $q$, so that for each $t\in B$, we have an exact sequence
\begin{align}
0\to \mathcal{N}_{\mathcal{F}_t}^*\to \Omega_{M_t}^1\to \Omega_{M_t}^1/\mathcal{N}_{\mathcal{F}_t}^*\to 0
\end{align}
where $\Omega_{M_t}^1/\mathcal{N}_{\mathcal{F}_t}^*$ is a torsion-free $\mathcal{O}_{M_t}$-module. We say that $\left\{\left(M_t, \mathcal{N}_{\mathcal{F}_t}^*\right)|t\in B \right\}$ is a family of unfoldings or unfolded $($complex$)$ analytic family if there exist a complex manifold $\mathcal{M}$ and a holomorphic map $\pi:\mathcal{M}\to B$ such that
\begin{enumerate}
\item we have an exact sequence
\begin{align*}
0\to \tilde{\mathcal{N}}_\mathcal{F}^*\to \Omega_\mathcal{M}^1\to \Omega_\mathcal{M}^1/\tilde{\mathcal{N}}_{\mathcal{F}}^*\to 0
\end{align*}
such that $\tilde{\mathcal{N}}_{\mathcal{F}}^*$ is a coherent subsheaf of the cotangent sheaf $\Omega_\mathcal{M}^1$ on $\mathcal{M}$.
\item the rank of the Jacobian of $\pi$ is equal to $m$ at every point of $\mathcal{M}$.
\item from the exact sequence $0\to \pi^* \Omega_B^1\to \Omega_\mathcal{M}^1\xrightarrow{\alpha} \Omega_{\mathcal{M}/B}^1 \to 0$, let us denote the image of $\tilde{\mathcal{N}}_{\mathcal{F}}^*$ under $\alpha$ by $\mathcal{N}_{\mathcal{F}}^*$. Then $\alpha$ is isomorphic on $\tilde{\mathcal{N}}_{\mathcal{F}}^*$, i.e $ \tilde{\mathcal{N}}_{\mathcal{F}}^* \cong \alpha\left( \tilde{\mathcal{N}}_{\mathcal{F}}^*  \right)= \mathcal{\mathcal{N}}_{\mathcal{F}}^*$ and the exact sequence
\begin{align*}
0\to \mathcal{N}_{\mathcal{F}}^*\to \Omega_{\mathcal{M}/B}^1\to \Omega_{\mathcal{M}/B}^1/\mathcal{N}_{\mathcal{F}}^*\to 0
\end{align*}
defines a foliated $($complex$)$ analytic family of $\left\{\left(M_t, \mathcal{N}_{\mathcal{F}_t}^* \right)\right\}$ in terms of cotangent sheaves in the sense of  \textnormal{Part I}, so that $\Omega_{\mathcal{M}/B}^1/\mathcal{N}_{\mathcal{F}}^*$ is flat over $B$.
\item $\tilde{\mathcal{N}}_{\mathcal{F}}^*$ satisfies the integrability condition. In other words,  let $\tilde{S}$ be the locus of points of $\mathcal{M}$ such that $\frac{\Omega_{\mathcal{M}}^1}{ \tilde{\mathcal{N}}_\mathcal{F}^*}$ is not locally free.  Then we have $d \tilde{\mathcal{N}}_{\mathcal{F}}^*\subset \tilde{\mathcal{N}}_{\mathcal{F}}^* \bigwedge \Omega_\mathcal{M}^1$ at every point of $\mathcal{M}- \tilde{S}$.  
\end{enumerate}
We will denote the unfolded analytic family by $\left(\mathcal{M}, \tilde{\mathcal{N}}_{\mathcal{F}}^*, B, \pi \right)$ and the associated foliated analytic family by $\left(\mathcal{M}, \mathcal{N}_\mathcal{F}^*, B, \pi \right)$ in terms of cotangent sheaves.
\end{definition}

\begin{remark}
Let $\left( \mathcal{M}, \tilde{\mathcal{N}}_{\mathcal{F}}^*, B, \pi \right)$ with $\tilde{\mathcal{N}}_\mathcal{F}^*$ be locally free be an unfolded analytic family. Let $\Delta$ be an open set of $B$. Then the restriction $\left( \mathcal{M}_\Delta =\pi^{-1}\left(\Delta \right) , \tilde{\mathcal{N}}_\mathcal{F}^*|_\Delta,  \Delta, \pi |_\Delta   \right)$ is also an unfolded complex analytic family. We will denote the family by $\left( \mathcal{M}_\Delta, \tilde{\mathcal{N}}_{\mathcal{F}_\Delta}^*, \Delta, \pi \right)$.
\end{remark}

\begin{remark}
From $\tilde{\mathcal{N}}_\mathcal{F}^*\cong \mathcal{N}_{\mathcal{F}}^*$, we have an exact sequence $0\to \pi^* \Omega_B^1\to \frac{ \Omega_\mathcal{M}^1}{ \tilde{\mathcal{N}}_\mathcal{F}^*} \xrightarrow{\alpha} \frac{\Omega_{\mathcal{M}/B}^1}{\mathcal{N}_{\mathcal{F}}^*}\to 0$. Since $\frac{\Omega_{\mathcal{M}/B}^1}{\mathcal{N}_{\mathcal{F}}^*}$ is flat over $B$, and $\pi^*\Omega_B^1$ is locally free, we see that $\frac{ \Omega_\mathcal{M}^1}{ \tilde{\mathcal{N}}_\mathcal{F}^*} $ is flat over $B$. Flatness of $\Omega_\mathcal{M}^1/\tilde{\mathcal{N}}_\mathcal{F}^*$ implies that $\tilde{\mathcal{N}}_\mathcal{F}^*$ is flat over $B$ and we have an exact sequence
\begin{align*}
0\to \tilde{\mathcal{N}}_{\mathcal{F}}^*|_{M_t}\to\Omega_\mathcal{M}^1|_{M_t}\to \frac{\Omega_\mathcal{M}^1}{\tilde{\mathcal{N}}_\mathcal{F}^*}|_{M_t}\to 0
\end{align*}
From the exact sequence $0\to \pi^*\Omega_B^1|_{M_t}\to \Omega_\mathcal{M}^1|_{M_t}\xrightarrow{\alpha_t} \Omega_{M_t}^1\to 0$, the image of $\tilde{\mathcal{N}}_{\mathcal{F}}^*|_{M_t}$ under $\alpha_t$ is $\mathcal{N}_{\mathcal{F}_t}^*$, and $\tilde{\mathcal{N}}_\mathcal{F}^*|_{M_t}\cong \mathcal{N}_{\mathcal{F}_t}^*$.  Assume that $\mathcal{N}_{\mathcal{F}_{t_0}}^*$ is a locally free subsheaf of $\Omega_{M_{t_0}}^1$ for some $t_0\in B$. Then since $\tilde{\mathcal{N}}_\mathcal{F}^*|_{M_{t_0}}\cong \mathcal{N}_{\mathcal{F}_{t_0}}^*$ is locally free, we see that $\tilde{\mathcal{N}}_\mathcal{F}^*$ is locally free in some neighborhood of $t_0$. In the following we will assume that $\tilde{\mathcal{N}}_\mathcal{F}^*$ is a locally free $\mathcal{O}_\mathcal{M}$-module and so $\mathcal{N}_{\mathcal{F}}^*$ is a locally free $\mathcal{O}_\mathcal{M}$-module.
\end{remark}

\begin{remark}
Let $\tilde{S}$ be the locus of points of $\mathcal{M}$ such that $\frac{\Omega_{\mathcal{M}}^1}{\tilde{\mathcal{N}}_\mathcal{F}^*}$ is not locally free. From the exact sequence $0 \to \pi^* \Omega_B^1 \to \frac{\Omega_\mathcal{M}^1}{\tilde{\mathcal{N}}_\mathcal{F}^*}\to \frac{\Omega_{\mathcal{M}/B}^1}{\mathcal{N}_\mathcal{F}^*}\to 0$, we see that $\tilde{S}\subset S=\bigcup_{t\in B} S_t$ where $S_t= \textnormal{Sing}\left( \mathcal{F}_t \right)$. Then the condition $(4)$ $d \tilde{\mathcal{N}}_{\mathcal{F}}^*\subset \tilde{\mathcal{N}}_{\mathcal{F}}^* \bigwedge \Omega_\mathcal{M}^1$ at every point of $\mathcal{M}- \tilde{S}$ implies that in particular  $d \tilde{\mathcal{N}}_{\mathcal{F}}^*\subset \tilde{\mathcal{N}}_{\mathcal{F}}^* \bigwedge \Omega_\mathcal{M}^1$ at every point of $\mathcal{M}- \bigcup_{t\in B} S_t$.
\end{remark}

\begin{remark}\label{u3}
From \textnormal{Definition \ref{u1}}, we can choose a system of local complex coordinate $\{z_i,...,z_j,...\},z_j:p\to z_j(p)$, and coordinate polydisks $\mathcal{U}_j$ with respect to $z_j$, satisfying the following conditions\footnote{for the detail, see \cite{Kod05} p.60}$:$
\begin{enumerate}
\item $z_j(p)= \left(z_j^1(p),...,z_j^n(p), t_j^1,...,t_j^m\right),\left(t_j^1,...,t_j^m\right)=\pi(p)$;
\item $\mathcal{U}= \left\{\mathcal{U}_j|j=1,2,... \right\}$ is locally finite.
\end{enumerate}
Then $\left\{p\to \left(z_j^1(p),...,z_j^n(p) \right)|\mathcal{U}_j\cap M_t\ne \emptyset \right\}$ gives a system of local complex coordinates on $M_t$. In terms of these coordinates, $\pi$ is the projection given by $\pi: \left(z_j^1,...,z_j^n, t_1,...,t_m \right)\to \left(t_i,...,t_m\right)$. For $j,k$ with $\mathcal{U}_j\cap \mathcal{U}_k\ne \emptyset$, we denote the coordinate transformation from $z_k$ to $z_j$ by
\begin{align*}
f_{jk}: \left(z_k^1,...,z_k^n,t_j \right) \to \left(z_j^1,...,z_j^n, t_j \right)=f_{jk}\left(z_k^1,...,z_k^n,t_j \right)
\end{align*}
Thus $f_{jk}$ is given by
\begin{align*}
z_j^\alpha=f_{jk}^\alpha \left(z_k^1,...,z_k^n, t_k^1,...,t_k^m \right),\alpha=1,...,n.
\end{align*}

On the other hand, we may assume that the locally free sheaf $\tilde{\mathcal{N}}_\mathcal{F}^*$ is trivialized in terms of the covering $\mathcal{U}$ so that
\begin{align*}
\Gamma\left(\mathcal{U}_j, \tilde{\mathcal{N}}_{\mathcal{F}}^*\right)&\cong \bigoplus^q \Gamma\left(\mathcal{U}, \mathcal{O}_\mathcal{M}\right)\\
                \omega_j^\alpha(z_j,t)&\mapsto e_j^\alpha= \left(0,..., \overbrace{1}^{\alpha-\textnormal{th}},...,0 \right), \alpha=1,...,q.
\end{align*}
which is generated by holomorphic forms $\omega_j^\alpha(z_j,t)\in \Gamma\left(\mathcal{U}_j, \tilde{\mathcal{N}}_\mathcal{F}^* \right),\alpha=1,...,q$ of the form
\begin{align}\label{u2}
\omega_j^\alpha(z_j,t):=w_j^\alpha(z_j,t)+v_j^\alpha(z_j,t)=\sum_{\beta=1}^n w_j^{\alpha\beta}(z_j,t_j)dz_j^\beta+\sum_{\eta=1}^m v_j^{\alpha\eta}(z_j,t_j) dt_j^\eta,\,\,\,\,\alpha=1,...,q
\end{align}
For $(x,a)\in \mathcal{U}_j-S$ where $S=\bigcup_{t\in B}\textnormal{Sing}(\mathcal{F}_t)$, there exists a neighborhood $U_{(x,a)}$ of $(x,a)$ such that $U_{(x,a)}\subset \mathcal{U}_j-S$, and for $\alpha=1,...,q$,
\begin{align*}
d\omega_j^\alpha(z_j,t)=\sum_{\beta=1}^q c_{j_{(x,a)}}^{\alpha\beta}(z_j,t)\wedge \omega_j^\beta(z_j,t)\,\,\,\,\,\,\,\,\,\textnormal{for some}\,\,\,\,\,c_{j_{(x,a)}}^{\alpha\beta}(z_j,t)\in \Gamma(\mathcal{U}_{(x,a)},\Omega_\mathcal{M}^1)
\end{align*}
Then $\omega_j^1(z_j,t)\wedge \cdots \wedge \omega_j^q(z_j,t)\wedge d\omega_j^\alpha(z_j,t)=0$ on $\mathcal{U}_j\cap \left(M_t-S_t\right)$, so that it vanishes on $\mathcal{U}_j\cap M_t$ for each $t\in B$ and so on $\mathcal{U}_j$. Hence we have
\begin{align}\label{u8}
\omega_j^1(z_j,t)\wedge \cdots \wedge \omega_j^q(z_j,t)\wedge d\omega_j^\alpha(z_j,t)=0,\,\,\,\,\,\alpha=1,...,q\,\,\,\textnormal{on}\,\,\,\mathcal{U}_j.
\end{align}
In particular, by considering the terms of $dz_j^{\alpha_1}\wedge \cdots \wedge dz_j^{\alpha_{q+2}}$ not containing $dt_j^1,..., dt_j^m$, we see that
\begin{align*}
w_j^1(z_j,t)\wedge \cdots \wedge w_j^q(z_j,t)\wedge d_z w_j^\alpha(z_j,t)=0,\,\,\,\,\,\alpha=1,...,q\,\,\,\textnormal{on}\,\,\,\mathcal{U}_j.
\end{align*}
where $d_z$ is the relative differential on $\Omega_{\mathcal{M}/B}^1$ over $B$.

For $j,k$ with $\mathcal{U}_j\cap \mathcal{U}_k \ne \emptyset$, we have invertible $q\times q$ matrices $H_{jk}:=\left(h_{jk}^{\alpha\beta}(z_j,t)\right)$ with components $h_{jk}^\alpha(z_k,t)\in \Gamma\left(\mathcal{U}_j\cap \mathcal{U}_k, \mathcal{O}_\mathcal{M}\right)$ such that
\begin{center}
$\left[\begin{matrix}
\omega_j^1(z_j,t)\\
\omega_j^2(z_j,t)\\
\cdot\\
\cdot\\
\cdot\\
\omega_j^q(z_j,t)
\end{matrix}\right]=
\left[\begin{matrix}
h_{jk}^{11}(z_k,t)& \cdots & h_{jk}^{1q}(z_k,t)\\
h_{jk}^{21}(z_k,t) &\cdots & h_{jk}^{2q}(z_k,t)\\
\cdot  &\cdots & \cdot \\
\cdot & \cdots &\cdot\\
\cdot & \cdots &\cdot\\
h_{jk}^{q1}(z_k,t) &\cdots & h_{jk}^{qq}(z_k,t)
\end{matrix}\right]
\left[\begin{matrix}
\omega_k^1(z_k,t)\\
\omega_k^2(z_k,t)\\
\cdot\\
\cdot\\
\cdot\\
\omega_k^q(z_k,t)
\end{matrix}\right]$
\end{center}
and for $i,j,k\in \mathcal{U}_i\cap \mathcal{U}_j \cap \mathcal{U}_k\ne \emptyset$, we have
\begin{align*}
H_{ik}=H_{ij}\cdot H_{jk}\,\,\,\,\,\textnormal{on}\,\,\,\mathcal{U}_i\cap \mathcal{U}_j\cap\mathcal{U}_k.
\end{align*}
More specifically, $\omega_j^\alpha(z_j,t)=\sum_{\beta=1}^q h_{jk}^{\alpha\beta}\omega_k^\beta(z_k,t)$ implies that
\begin{align}
\sum_{\gamma=1}^n w_j^{\alpha\gamma}(z_j,t)dz_j^\gamma+\sum_{\eta=1}^m v_j^{\alpha\eta}(z_j,t)dt_j^\eta=\sum_{\beta=1}^q h_{jk}^{\alpha\beta}(z_k,t)\left(\sum_{\xi=1}^n w_k^{\beta\xi}(z_k,t)dz_k^\xi+\sum_{\eta=1}^m v_k^{\beta \eta}(z_k,t)dt_k^\eta \right)
\end{align}
Then we have
\begin{align*}
&\sum_{\gamma=1}^n\sum_{\xi=1}^n w_j^{\alpha\gamma}(f_{jk}(z_k,t),t)\frac{\partial f_{jk}^\gamma(z_k,t)}{\partial z_k^\xi}dz_k^\xi +\sum_{\gamma=1}^n \sum_{\eta=1}^m w_j^{\alpha\gamma}(f_{jk}(z_k,t),t)\frac{\partial f_{jk}^\gamma(z_k,t)}{\partial t_k^\eta} dt_k^\eta+\sum_{\eta=1}^m v_j^{\alpha\eta}(f_{jk}(z_k,t),t) dt_k^\eta\\
&=\sum_{\beta=1}^q\sum_{\xi=1}^n h_{jk}^{\alpha\beta}(z_k,t) w_k^{\beta\xi}(z_k,t) dz_k^\xi +\sum_{\beta=1}^q \sum_{\eta=1}^m  h_{jk}^{\alpha\beta}(z_k,t) v_k^{\beta\eta}(z_k,t)dt_k^\eta
\end{align*}
By considering the coefficients of $dz_k^\xi$, and $dt^\eta:=dt_k^\eta= dt_j^\eta$, we have
\begin{align}
\sum_{\gamma=1}^n w_j^{\alpha\gamma}(f_{jk}(z_k,t),t)\frac{\partial f_{jk}^\gamma(z_k,t)}{\partial z_k^\xi} &= \sum_{\beta=1}^q h_{jk}^{\alpha\beta}(z_k,t) w_k^{\beta\xi}(z_k,t)\\
\sum_{\gamma=1}^n\sum_{\eta=1}^m w_j^{\alpha\gamma}(f_{jk}(z_k,t),t)\frac{\partial f_{jk}^\gamma(z_k,t)}{\partial t^\eta}dt^\eta + \sum_{\eta=1}^m v_j^{\alpha\eta}(f_{jk}(z_k,t),t)dt^\eta & = \sum_{\beta=1}^q \sum_{\eta=1}^m h_{jk}^{\alpha\beta}(z_k,t) v_k^{\beta\eta}(z_k,t)dt^\eta,\,\,\,\,\,\,\,\,\,\alpha=1,...,q. \label{u6}
\end{align} 
Let $\frac{\partial}{\partial t}=\sum_{\eta=1}^m c_\eta\frac{\partial}{\partial t_\eta}$. By applying $i_{\frac{\partial}{\partial t}}$ to $(\ref{u6})$, we have
\begin{align}\label{u7}
\sum_{\gamma=1}^n \sum_{\eta=1}^m w_j^{\alpha\gamma}(f_{jk}(z_k,t),t)c_\eta\frac{\partial f_{jk}^{\gamma}(z_k,t)}{\partial t^\eta }+\sum_{\eta=1}^m c_\eta v_j^{\alpha\eta}(f_{jk}(z_k,t),t)=\sum_{\beta=1}^q \sum_{\eta=1}^m h_{jk}^{\alpha\beta}c_\eta v_k^{\beta\eta}(z_k,t)
\end{align}

\end{remark}

\subsection{Unfolded dual leaf complex controlling unfoldings of singular holomorphic foliations}\

Let $\left( \mathcal{M}, \tilde{\mathcal{N}}_\mathcal{F}^* , B, \pi \right)$ be an unfolded complex analytic family with $\tilde{\mathcal{N}}_\mathcal{F}^*$ locally free as in Definition \ref{u1}, so that $\mathcal{N}_{\mathcal{F}_t}^*$ defines a (singular) holomorphic foliation on a compact complex manifold $M_t$ with $\mathcal{N}_{\mathcal{F}_t}^*$ locally free. Then we have the unfolded dual leaf complex on $M_t$ associated to $\mathcal{N}_{\mathcal{F}_t}^*$ (see Appendix \ref{AppendixA})
{\small{\begin{align*}
\tilde{\mathcal{N}}_{\mathcal{F}_t}^{*\bullet}: \Theta_{M_t}  \xrightarrow{G_0^t} \mathscr{H}om_{\mathcal{O}_{M_t}}\left(\mathcal{N}_{\mathcal{F}_t}^*, \frac{\Omega_{M_t}^1}{\mathcal{N}_{\mathcal{F}_t}^*} \right) \bigoplus & \mathscr{H}om_{\mathcal{O}_{M_t}}\left(\mathcal{N}_{\mathcal{F}_t}^*, \mathcal{O}_{M_t} \right) \xrightarrow{G_1^t} \mathscr{H}om_{\mathcal{O}_{M_t}}\left(\mathcal{N}_{\mathcal{F}_0}, \tilde{\mathcal{S}}_t^2 \right) \bigoplus \mathscr{H}om_{\mathcal{O}_{M_t}}\left( \mathcal{N}_{\mathcal{F}_t}, \tilde{\mathcal{S}}_t^1 \right) \\
&\xrightarrow{G_2^t} \mathscr{H}om_{\mathcal{O}_{M_t}}\left( \mathcal{N}_{\mathcal{F}_t}^*, \tilde{\mathcal{S}}_t^3 \right) \bigoplus \mathscr{H}om_{\mathcal{O}_{M_t}}\left( \mathcal{N}_{\mathcal{F}_t}^*, \tilde{\mathcal{S}}_t^2 \right) \xrightarrow{G_3^t} \cdots 
\end{align*}}}
We will denote the $i$-th cohomology group by $\mathbb{H}^i\left( M_t , \tilde{\mathcal{N}}_{\mathcal{F}_t}^{* \bullet} \right)$. We can compute $\mathbb{H}^i\left(M_t, \tilde{\mathcal{N}}_{\mathcal{F}_t}^{*\bullet} \right)$  by the following \v Cech resolution of $\mathcal{N}_{\mathcal{F}_t}^{*\bullet}$. Here $\delta$ is the \v Cech map and $\mathcal{U}_t= \mathcal{U}\cap M_t=\left\{ U_j^t :=\mathcal{U}_j \cap M_t | j=1,2,... \right\}$ is an open covering of $M_t$:

{\small{\begin{center}
$\begin{CD}
\cdots\\
@AG_2^t AA \\
C^0\left(\mathcal{U}_t, \left( \left( \mathcal{N}_{\mathcal{F}_t}^* \right)^*\otimes  \tilde{\mathcal{S}}_t^2 \right) \bigoplus \left( \left( \mathcal{N}_{\mathcal{F}_0}^* \right)^* \otimes \tilde{\mathcal{S}}^1 \right) \right) @>-\delta>> \cdots \\
@AG_1^t AA @AG_1^t AA \\
C^0\left( \mathcal{U}_t, \left(\left(\mathcal{N}_{\mathcal{F}_t}^*\right)^* \otimes \frac{\Omega_{M_t}^1}{\mathcal{N}_{\mathcal{F}_t}^*}\right) \bigoplus \left( \mathcal{N}_{\mathcal{F}_t}^* \right)^*\right) @>\delta>> C^1\left( \mathcal{U}_t, \left(\left(\mathcal{N}_{\mathcal{F}_t}^*\right)^* \otimes \frac{\Omega_{M_t}^1}{\mathcal{N}_{\mathcal{F}_t}^*}\right) \bigoplus \left( \mathcal{N}_{\mathcal{F}_t}^* \right)^*\right)  @>-\delta>> \cdots \\
@AG_0^tAA @AG_0^tAA @AG_0^tAA \\
C^0\left(\mathcal{U}_t, \Theta_{M_t} \right) @>-\delta>> C^1\left( \mathcal{U}_t, \Theta_{M_t} \right) @>\delta>> C^2\left( \mathcal{U}_t, \Theta_{M_t} \right) @>-\delta >> \cdots
\end{CD}$
\end{center}}}

We will  relate the first cohomology group $\mathbb{H}^1\left( M_t, \tilde{\mathcal{N}}_{\mathcal{F}_t}^{*\bullet} \right)$ to infinitesimal unfoldings of $\pi^{-1}(t)=\left(M_t, \mathcal{N}_{\mathcal{F}_t}^* \right)$ in the unfolded analytic family $\left(\mathcal{M}, \tilde{\mathcal{N}}_\mathcal{F}^*, \pi, B \right)$. 

\subsection{Infinitesimal unfoldings}\

Let $\left( \mathcal{M},\tilde{\mathcal{N}}_{\mathcal{F}}^*, B, \pi \right)$ be an unfolded analytic family with $\tilde{\mathcal{N}}_\mathcal{F}^*$ locally free as in Definition \ref{u1}. We keep the notations in Remark \ref{u3}.
For $\frac{\partial}{\partial t} = \sum_{\eta=1}^m c_\eta \frac{\partial}{\partial t_\eta} \in T_{t}B$ with $c_\eta \in\mathbb{C}$, we define an element $\beta_j(t)$ in $\Gamma\left(U_j^t,\mathscr{H}om_{\mathcal{O}_{M_t}}\left( \mathcal{N}_{\mathcal{F}_t}^*, \frac{\Omega_{M_t}^1}{\mathcal{N}_{\mathcal{F}_t}^*} \right)\right) $ as in Part I in the following way:
\begin{align*}
\beta_j(t):\Gamma\left(U_j^t,\mathcal{N}_{\mathcal{F}_t}^*\right)&\to \Gamma\left(U_j^t, \frac{\Omega_{M_t}^1}{\mathcal{N}_{\mathcal{F}_t}^*}\right)\\
 w_j^\alpha(z_j,t)&\mapsto -\overline{ \frac{\partial w_j^\alpha(z_j,t)}{\partial t}}= -\sum_{\eta=1}^m c_\eta \overline{\frac{\partial w_j^\alpha(z_j,t)}{\partial t_\eta}}:=  -\overline{\sum_{\eta=1}^m \sum_{\beta=1}^n c_\eta \frac{\partial w_j^{\alpha\beta}(z_j,t)}{\partial t_\eta}}dz_j^\beta,\,\,\,\,\,\,\,\alpha=1,...,q
\end{align*}
where $\overline{\frac{\partial w_j^\alpha(z_j,t)}{\partial t}}$ denotes the image of $\frac{\partial w_j^\alpha(z_j,t)}{\partial t}$ in the natural map $\Omega_{M_t}^1\to \Omega_{M_t}^1/\mathcal{N}_{\mathcal{F}_t}^*$, and linearly extends to $\Gamma\left(U_j^t, \mathcal{N}_{\mathcal{F}_t}^* \right)$. We also define $\tilde{\beta}_j(t)\in \Gamma\left(U_j^t, \mathscr{H}om_{\mathcal{O}_{M_t}}\left( \mathcal{N}_{\mathcal{F}_t}^*, \Omega_{M_t}^1 \right) \right)$ in the following way:
\begin{align*}
\tilde{\beta}_j(t):\Gamma\left(U_j^t, \mathcal{N}_{\mathcal{F}_t}^*\right)&\to \Gamma\left(U_j^t,\Omega_{M_t}^1  \right) \\
                                                         w_j^\alpha(z_j,t)&\mapsto -\frac{\partial w_j^\alpha(z_j,t)}{\partial t} = - \sum_{\eta=1}^m \sum_{\beta=1}^n c_\eta \frac{\partial w_j^{\alpha\beta}(z_j,t)}{\partial t_\eta}
\end{align*}
and linearly extends to $\Gamma\left( U_j^t, \mathcal{N}_{\mathcal{F}_t}^* \right)$. We define an element $\gamma_j(t)$ in $\Gamma\left(U_j^t, \mathscr{H}om_{\mathcal{O}_{M_t}}\left(\mathcal{N}_{\mathcal{F}_t}^*, \mathcal{O}_{M_t} \right)\right)$ in the following way:
\begin{align*}
\gamma_j(t):\Gamma\left(U_j^t, \mathcal{N}_{\mathcal{F}_t}^*\right)&\to \Gamma\left(U_j^t, \mathcal{O}_{M_t}\right)\\
w_j^\alpha(z_j,t) &\mapsto -i_{\frac{\partial}{\partial t}}\left(v_j^\alpha(z_j,t) \right)= -i_{\frac{\partial}{\partial t}}\left(\sum_{\eta=1}^m v_j^{\alpha\eta}(z_j,t_j)dt_j^\eta \right)= -\sum_{\eta=1}^m c_\eta v_j^{\alpha\eta}(z_j,t)
\end{align*}
and we define
\begin{align*}
\theta_{jk}(t):= \sum_{\gamma=1}^n \frac{\partial f_{jk}^\gamma(z_k,t)}{\partial t}\frac{\partial}{\partial z_j^\gamma }=\sum_{\eta=1}^m \sum_{\gamma=1}^n c_\eta \frac{\partial f_{jk}^\gamma(z_k,t)}{\partial t^\eta }\frac{\partial}{\partial z_j^\gamma} \in \Gamma\left(U_{ij}^t, \Theta_{M_t}\right)
\end{align*}

Then we have 
\begin{proposition} \label{u115}
{\small{\begin{align*}
\left(\left\{\beta_j(t) \right\},\left\{\gamma_j(t) \right\}, \left\{\theta_{jk}(t) \right\} \right)\in C^0\left(\mathcal{U}_t, \mathscr{H}om_{\mathcal{O}_{M_t}}\left( \mathcal{N}_{\mathcal{F}_t}^*, \frac{\Omega_{M_t}^1}{\mathcal{N}_{\mathcal{F}_t}^*} \right)\right)\bigoplus C^0\left( \mathcal{U}_t , \mathscr{H}om_{\mathcal{O}_{M_t}}\left( \mathcal{N}_{\mathcal{F}_t}^*, \mathcal{O}_{M_t} \right) \right)\bigoplus C^1\left(\mathcal{U}_t, \Theta_{M_t}\right)
\end{align*}}}
defines a $1$-cocycle in the above \v Cech resolution of $\tilde{\mathcal{N}}_{\mathcal{F}_t}^{*\bullet}$ and call its cohomology class in $\mathbb{H}^1\left( M_t, \tilde{\mathcal{N}}_{\mathcal{F}_t}^{*\bullet} \right)$ the infinitesimal unfoldings along $\frac{\partial}{\partial t}= \sum_{\eta=1}^m c_\eta \frac{\partial}{\partial t^\eta}$. This expression is independent of the choice of local coordinates.
\end{proposition}

\begin{proof}
First we note that $\delta\left( \left\{ \theta_{jk}(t) \right\} \right)=0$ (see \cite{Kod05} p.201). From Part I we showed that $E_1(\beta_j(t))=0, \delta\left( \{\beta_j(t) \} \right) + E_0\left( \left\{\theta_{jk}(t)\right\}\right)=0$. It remains to show that 
\begin{align}\label{u4}
\left(\gamma_k(t)- \gamma_j(t)\right)(w_j^\alpha)+ i_{\theta_{jk}(t)}(w_j^\alpha)=0\,\,\,\,\,\,\,\,\textnormal{for $\alpha=1,...,q$}
\end{align}
 \begin{align}\label{u5}
 - w_j^1 \wedge \cdots \wedge w_j^q \wedge \tilde{\beta}_j(t) +\sum_{\beta=1}^q (-1)^{q-\beta} \gamma_j(t) w_j^1 \wedge \cdots \wedge \widehat{w}_j^\beta \wedge \cdots \wedge w_j^q \wedge dw_j^\alpha + w_j^1 \wedge \cdots \wedge w_j^q \wedge d\gamma_j^\alpha(t)=0
 \end{align}
for $\alpha=1,...,q$. Let us show $(\ref{u4})$. In fact, since $w_j^\alpha= \sum_{\beta=1}^q h_{jk}^{\alpha\beta} w_k^\beta$, $(\ref{u7})$ implies that
\begin{align*}
\left(\gamma_k(t)- \gamma_j(t)\right)(w_j^\alpha)+ i_{\theta_{jk}(t)}(w_j^\alpha) = - \sum_{\beta=1}^q h_{jk}^{\alpha\beta} \left( \sum_{\eta=1}^m c_\eta v_k^{\beta \eta} \right)+ \sum_{\eta=1}^m c_\eta v_j^{\alpha \eta} + \sum_{\eta=1}^m \sum_{\gamma=1}^q c_\eta \frac{\partial f_{jk}^\gamma}{\partial t^\eta} w_j^{\alpha \gamma}=0
\end{align*}
Let us show $(\ref{u5})$. In fact, from $(\ref{u8})$, we have
{\small{\begin{align*}
\omega_j^1(z_j,t)\wedge \cdots \wedge \overbrace{\left(\sum_{\beta=1}^n w_j^{\gamma\beta}(z_j,t)dz_j^\beta+\sum_{\eta=1}^m v_j^{\gamma\eta}(z_j,t) dt^\eta \right)}^{\gamma-\textnormal{th}}\wedge \cdots \wedge \omega_j^1(z_j,t)\wedge d\left(\sum_{\beta=1}^n w_j^{\alpha\beta}(z_j,t)dz_j^\beta+\sum_{\eta=1}^m v_j^{\alpha\eta}(z_j,t) dt^\eta\right)=0
\end{align*}}}
Let us collect the terms involving $dz_j^{\gamma_1}\wedge \cdots dz_j^{\gamma_{q+2}}\wedge dt^\eta$. Then we have
\begin{align*}
&\sum_{\gamma=1}^q w_j^1(z_j,t)\wedge \cdots \wedge \left( \sum_{\eta=1}^m v_j^{\gamma \eta}(z_j,t)dt^\eta \right) \wedge \cdots \wedge w_j^q(z_j,t)\wedge d_t w_j^\alpha(z_j,t)\\
&+w_j^1(z_j,t)\wedge \cdots \wedge w_j^q(z_j,t)\wedge \left( \sum_{\beta=1}^n \sum_{\eta=1}^m \frac{\partial w_j^{\alpha\beta}(z_j,t)}{\partial t^\eta}dt_\eta\wedge dz_j^\beta + \sum_{\eta=1}^m\sum_{\gamma=1}^n \frac{\partial v_j^{\alpha\eta}}{\partial z_j^\gamma} dz_j^\gamma \wedge dt_j^\eta\right) = 0
\end{align*}
Here $d_t$ is the differential on $M_t$. Then we have
{\Small{\begin{align*}
\sum_{\gamma=1}^q \sum_{\eta=1}^m (-1)^{q-\gamma} v_j^{\gamma\eta} w_j^1\wedge \cdots \wedge \hat{w}_j^\gamma \wedge \cdots \wedge w_j^q\wedge d_t w_j^\alpha \wedge dt^\eta+  \sum_{\eta=1}^m w_j^1\wedge \cdots \wedge w_j^q\wedge d_t v_j^{\alpha\eta} \wedge dt^\eta-w_j^1\wedge \cdots \wedge w_j^q\wedge \left( \sum_{\beta=1}^n \sum_{\eta=1}^m \frac{\partial w_j^{\alpha\beta}}{\partial t^\eta } dz_j^\beta\right) \wedge d t^\eta =0
\end{align*}}}
Let $\frac{\partial}{\partial t}=\sum_{\eta=1}^m c_\eta\frac{\partial}{\partial t^\eta }$. Then by applying $i_{\frac{\partial}{\partial t}}$, we have
{\tiny{\begin{align*}
0&=i_{\frac{\partial}{\partial t}}\left(\sum_{\gamma=1}^q\sum_{\eta=1}^m (-1)^{q-\gamma} v_j^{\gamma\eta} w_j^1\wedge \cdots \wedge \hat{w}_j^\gamma \wedge \cdots \wedge w_j^q\wedge d_t w_j^\alpha \wedge dt^\eta+ \sum_{\eta=1}^m w_j^1\wedge \cdots \wedge w_j^q\wedge d_t v_j^{\alpha\eta} \wedge dt^\eta-w_j^1\wedge \cdots \wedge w_j^q\wedge \left( \sum_{\beta=1}^n \sum_{\eta=1}^m \frac{\partial w_j^{\alpha\beta}}{\partial t^\eta} dz_j^\beta\right) \wedge dt^\eta \right)\\
&=\sum_{\gamma=1}^q (-1)^{q-\gamma}\left(\sum_{\eta=1}^m c_\eta v_j^{\gamma\eta} \right)  w_j^1\wedge \cdots \wedge \hat{w}_j^\gamma \wedge \cdots \wedge w_j^q\wedge d_t w_j^\alpha+ w_j^1\wedge \cdots \wedge w_j^q\wedge d_t  \left( \sum_{\eta=1}^mc_\eta v_j^{\alpha\eta} \right) -w_j^1\wedge \cdots \wedge w_j^q\wedge \left( \sum_{\beta=1}^n \sum_{\eta=1}^m c_\eta \frac{\partial w_j^{\alpha\beta}}{\partial t^\eta} dz_j^\beta\right)
\end{align*}}}
This implies $(\ref{u5})$. Next we show that $\left( \left\{\beta_j(t) \right\}, \left\{ \gamma_j(t) \right\}, \left\{  \theta_{jk}(t)\right\} \right)$ is independent of choice of system of local coordinates. We can show that the infinitesimal deformation does not change under the refinement of the open covering (see \cite{Kod05} p.190). Since we can choose a common refinement for two systems of local coordinates, it is sufficient to show that given two local coordinates $x_j=(z_j,t)$ and $u_j=(y_j,t)$ on each $\mathcal{U}_j$, the infinitesimal unfolding $\left(\left\{\eta_{jk}(t) \right\}, \left\{ \chi_j(t) \right\} , \left\{\pi_j(t) \right\} \right)$ with respect to $\{u_j\}$ coincides with $\left( \left\{\theta_{jk} (t) \right\} , \left\{ \beta_j (t) \right\}  ,\left\{ \gamma_j(t)\right\} \right)$ with respect to $\{x_j\}$. Let $\tilde{\mathcal{N}}_\mathcal{F}^*$ be generated by $W_j^\alpha (y_j,t)  + V_j^\alpha (y_j,t) = \sum_{\gamma=1}^n W_j^{\alpha \gamma}(y_j,t)dy_j^\gamma + \sum_{\eta=1}^m V_j^{\alpha \eta}(y_j,t) dt^\eta,\alpha=1,...,q  $ in terms of coordinates $u_j=(y_j,t)$ on $\mathcal{U}_j$. Let $(y_k,t)\to (y_j,t)=(e_{jk}(y_k,t),t)$ be the coordinate transformation of $\{u_j\}$ on $\mathcal{U}_j\cap \mathcal{U}_k\ne \emptyset$. Then for $\frac{\partial}{\partial t}= \sum_{\eta=1}^m c_\eta \frac{\partial}{\partial t^\eta}$ with $c_\eta\in \mathbb{C}$,  we have $\eta_{jk}(t)=\sum_{\alpha=1}^n \frac{\partial e_{jk}^\alpha(y_k,t)}{\partial t}\frac{\partial}{\partial y_j^\alpha}, y_k= e_{kj}(y_j,t)$, and 
\begin{align*}
\chi_j (t):\Gamma\left(U_j^t,  \mathcal{N}_{\mathcal{F}_t}^* \right) \to \Gamma\left( U_j^t, \frac{\Omega_{M_t}^1}{\mathcal{N}_{\mathcal{F}_t}^*} \right), \,\,\,\,\,& W_j^\alpha \mapsto - \overline{\frac{\partial W_j^\alpha}{\partial t}}:= - \overline{\sum_{\gamma=1}^n \frac{\partial W_j^{\alpha \gamma}(y_j,t)}{\partial t} \frac{\partial}{\partial y_j^\gamma}}\\
\pi_j(t):\Gamma\left(U_j^t,  \mathcal{N}_{\mathcal{F}_t}^* \right) \to \Gamma\left( U_j^t, \mathcal{O}_{M_t} \right),\,\,\,\,\,& W_j^\alpha \mapsto - i_{\frac{\partial}{\partial t}}\left(V_j^\alpha \right):= - \sum_{\eta=1}^m c_j^\eta  V_j^{\alpha \eta}(y_i,t)
\end{align*}
We show that $\left(\left\{ \theta_{jk}(t) \right\}, \left\{ \beta_j(t)\right\}, \left\{ \gamma_j(t) \right\} \right)$ is cohomologous to $\left(\left\{ \eta_{jk}(t)\right\}, \left\{ \chi_j(t)\right\}  , \left\{ \pi_j(t) \right\}\right)$. Let $y_j^\alpha=g_j^\alpha\left(z_j^1,..., z_j^n, t \right)$ for $\alpha=1,...,n$, define the coordinate transformation from $x_j=(z_j,t)$ to $u_j=(y_j,t)$. Since $\left\{w_j^\alpha+ v_j^\alpha\right\}$ and $\left\{W_j^\alpha + V_j^\alpha \right\},\alpha=1,...,n$ are two bases of $\Gamma\left(\mathcal{U}_j, \tilde{\mathcal{N}}_\mathcal{F}^* \right)$, we can write by setting $g_j(z_j,t)=\left( g_j^1(z_j,t),..., g_j^n(z_j,t) \right)$,
\begin{align}
W_j^\alpha(y_j&,t) + V_j^\alpha(y_i,t)= \sum_{\beta=1}^q b_j^{\alpha\beta}(z_j,t) \left( w_j^\beta(z_j,t) + V_j^\beta(z_j,t)\right) \notag \\
\Longrightarrow &\sum_{\gamma=1}^n W_j^{\alpha\gamma}(y_j,t)dy_j^\gamma + \sum_{\eta=1}^m V_j^{\alpha \eta} (y_j,t) dt^\eta =\sum_{\beta=1}^q b_j^{\alpha\beta}(z_j,t)\left(\sum_{\xi=1}^n w_j^{\beta \xi}(z_j,t)dz_j^\xi + \sum_{\eta=1}^m v_j^{\beta \eta}(z_j,t)dt^\eta \right) \notag \\
\Longrightarrow &\sum_{\gamma=1}^n W_j^{\alpha\gamma}(y_j,t) d_zy_j^\gamma= \sum_{\beta=1}^q \sum_{\xi=1}^n b_j^{\alpha\beta}(z_j,t) w_j^{\beta \xi}(z_j,t)dz_z z_j^\xi,  \notag\\
&\sum_{\gamma=1}^n W_j^{\alpha \gamma}(g_j(z_j,t),t) \frac{\partial g_j^\gamma(z_j,t)}{\partial t^\eta} dt^\eta+\sum_{\eta=1}^m V_j^{\alpha \eta}(g_j(z_j,t),t)dt^\eta= \sum_{\beta=1}^q \sum_{\eta=1}^m b_j^{\alpha\beta}(z_j,t) v_j^{\beta \eta}\left( z_j,t\right)dt^\eta \label{u9}
\end{align}
 for some $b_j^{\alpha \gamma}(z_j,t)\in \Gamma\left(\mathcal{U}_j, \mathcal{O}_{\mathcal{M}} \right)$. We set $\theta_j(t)= \sum_{\alpha=1}^n \frac{\partial g_j^\alpha(z_j,t)}{\partial t} \frac{\partial}{\partial y_j^\alpha}= \sum_{\alpha=1}^n \sum_{\eta=1}^m c_\eta \frac{\partial g_j^{\alpha}(z_j,t)}{\partial t^\eta} \frac{\partial}{\partial y_j^\alpha}, y_j^\alpha= g_j^\alpha(z_j,t)$. Then we have $-\delta\left(\left\{\theta_j(t) \right\} \right)=\left\{\eta_{jk}(t)- \theta_{jk}(t) \right\}$ (for the detail, see \cite{Kod05} p.191-192), and we showed $\overline{\mathcal{L}_{\theta_{j}(t)}\left(W_j^\alpha \right)} = \xi_j(t)\left(W_j^\alpha \right)- \beta_j(t)\left( W_j^\alpha\right)$ in Part I. We claim that $i_{\theta_j(t)}\left(W_j^\alpha \right) = \pi_j(t)\left(W_j^\alpha \right) - \gamma_j(t)\left(W_j^\alpha \right)$. In fact, since $W_j^\alpha=\sum_{\beta=1}^q b_j^{\alpha\beta}(z_j,t) w_j^\beta(z_j,t)$ on $M_t$, we have
 \begin{align*}
 \gamma_j(t)\left( W_j^\alpha \right)- \pi_j(t)\left( W_j^\alpha\right) + i_{\theta_j(t)}\left(W_j^\alpha \right)= -\sum_{\beta=1}^q \sum_{\eta=1}^m b_j^{\alpha\beta} c_\eta v_j^{\beta \eta} + \sum_{\eta=1}^m c_\eta V_j^{\alpha \eta} + \sum_{\gamma=1}^n \sum_{\eta=1}^m c_\eta W_j^{\alpha \gamma}\frac{\partial g_j^\gamma}{\partial t^\eta} =0
 \end{align*}
which follows from applying $i_{\frac{\partial}{\partial t}}$ to $(\ref{u9})$. Hence we see that $\left(\left\{ \theta_{jk}(t) \right\}, \left\{ \beta_j(t)\right\} , \left\{\gamma_j(t) \right\}\right)$ is cohomologous to $\left(\left\{ \eta_{jk}(t)\right\}, \left\{ \chi_j (t)\right\} , \left\{ \pi_j(t) \right\} \right)$.

\end{proof} 

\begin{definition}[unfolded Kodaira-Spencer map]\label{u10}
Let $\left(\mathcal{M}, \tilde{\mathcal{N}}_\mathcal{F}^*, B, \pi\right)$ with $\tilde{\mathcal{N}}_{\mathcal{F}}^*$ locally free be an unfolded complex analytic family of unfoldings of $\left(M_t, \mathcal{N}_{\mathcal{F}_t}^*\right)=\pi^{-1}(t),t\in B$, where $B$ is a domain of $\mathbb{C}^m$. We keep the notations in the proof of \textnormal{Proposition \ref{u115}}. For a tangent vector $\frac{\partial}{\partial t}=\sum_{\eta=1}^m c_\eta \frac{\partial}{\partial t^\eta}, c_\eta \in \mathbb{C}$, of $B$, we put 
\begin{align*}
\frac{\partial \tilde{\mathcal{N}}_{\mathcal{F}_t}^*}{\partial t}:=&\left\{\left(  w_j^\alpha(z_j,t) \mapsto  - \overline{\frac{\partial w_j^\alpha(z_j,t)}{\partial t} } =- \overline{\sum_{\eta=1}^m c_\eta \frac{\partial w_j^\alpha(z_j,t)}{\partial t^\eta}}     \right)_{\alpha=1,...,q} \right\} \bigoplus \left\{ \left( w_j^\alpha(z_j,t) \mapsto - \sum_{\eta=1}^m c_\eta v_j^{\alpha \eta}(z_j,t) \right)_{\alpha=1,...,q}   \right\} \\
& \in C^0\left( \mathcal{U}_t, \mathscr{H}om_{\mathcal{O}_{M_t}}\left(\mathcal{N}_{\mathcal{F}_t}^*, \frac{\Omega_{M_t}^1}{\mathcal{N}_{\mathcal{F}_t}^*} \right) \right) \bigoplus C^0\left( \mathcal{U}_t,  \mathscr{H}om_{\mathcal{O}_{M_t}} \left( \mathcal{N}_{\mathcal{F}_t}^*, \mathcal{O}_{M_t} \right) \right)
\end{align*}

Then the unfolded Kodaira-Spencer map at $t$ is defined to be a $\mathbb{C}$-linear map
{\Small{\begin{align*}
\tilde{\varphi}_t&:T_t(B)\to \mathbb{H}^1\left(M_t, \tilde{\mathcal{N}}_{\mathcal{F}_t}^{*\bullet} \right)\\
                \frac{\partial}{\partial t}&\mapsto \frac{\partial \left(M_t, \tilde{\mathcal{N}}_{\mathcal{F}_t}^* \right)}{\partial t}\\
                &:=\left(\frac{\partial M_t}{\partial t}= \left\{  \sum_{\alpha=1}^n \frac{\partial f_{jk}^\alpha(z_k,t)}{\partial t}\frac{\partial}{\partial z_j^\alpha} \right\}, \frac{\partial \mathcal{N}_{\mathcal{F}_t}^*}{\partial t}= \left\{ \left(w_j^\alpha(z_j,t)\mapsto -\overline{\frac{\partial w_j^\alpha(z_j,t)}{\partial t}}\right)_{\alpha=1,...,q} \right\} , \left\{ \left( w_j^\alpha(z_j,t) \mapsto - \sum_{\eta=1}^m c_\eta v_j^{\alpha \eta}(z_j,t) \right)_{\alpha=1,...,q}   \right\} \right) \\
                 &\in C^1\left(\mathcal{U}_t, \Theta_{M_t} \right) \bigoplus C^0\left( \mathcal{U}_t, \mathscr{H}om_{\mathcal{O}_{M_t}}\left(\mathcal{N}_{\mathcal{F}_t}^*, \frac{\Omega_{M_t}^1}{\mathcal{N}_{\mathcal{F}_t}^*} \right) \right) \bigoplus C^0\left( \mathcal{U}_t,  \mathscr{H}om_{\mathcal{O}_{M_t}} \left( \mathcal{N}_{\mathcal{F}_t}^*, \mathcal{O}_{M_t} \right) \right)
\end{align*}}}
where $\rho_t:T_t(B) \to H^1\left( M_t, \Theta_{M_t}\right), \frac{\partial}{\partial t}\mapsto \frac{\partial M_t}{\partial t}$ is the Kodaira-Spencer map at $t$ of the underlying complex analytic family $\left( \mathcal{M}, B, \pi \right)$ $($see \cite{Kod05} \textnormal{p.201}$)$, and $\varphi_t:T_t(B)\to \mathbb{H}^1\left(M_t, \mathcal{N}_{\mathcal{F}_t}^{*\bullet} \right), \frac{\partial}{\partial t} \mapsto \frac{\partial\left(M_t, \mathcal{N}_{\mathcal{F}_t}^* \right)}{\partial t}$ is the foliated Kodaira-Spencer map at $t$ of the underlying foliated analytic family $\left( \mathcal{M}, \mathcal{N}_\mathcal{F}^*, B, \pi \right)$ in terms of cotangent sheaves $($see \textnormal{Part I}$)$.
\end{definition}

\section{Theorem of existence for unfoldings of foliated complex analytic structures for one-parameter families}

\subsection{Preliminaries}\

Let $\left(\mathcal{M}, \tilde{\mathcal{N}}_\mathcal{F}^* , B, \pi \right)$ with $\tilde{\mathcal{N}}_\mathcal{F}^*$ locally free be an unfolded analytic family, where $B$ is a domain of $\mathbb{C}^m$ with the coordinate $\left(t_1,..., t_m\right)$ containing the origin $0$ as in Definition \ref{u1}. Define $|t|=\max_\lambda |t_\lambda|$ for $t=\left(t_1,...,t_m\right)\in \mathbb{C}^m$, and let $\Delta= \left\{t\in \mathbb{C}^m | |t|<r \right\}$ the polydisk  of radius $r>0$. If we take a sufficiently  small $\Delta \subset B$, then $\left(\mathcal{M}_\Delta, \tilde{\mathcal{N}}_{\mathcal{F}_\Delta}^* \right)=\pi^{-1}\left(\Delta\right)$ is represented in the form
\begin{align*}
\left(\mathcal{M}_\Delta, \tilde{\mathcal{N}}_{\mathcal{F}_\Delta}^* \right)=\bigcup_j \left(U_j\times \Delta, \tilde{\mathcal{N}}_\mathcal{F}^*|_{U_j\times \Delta}\right)
\end{align*}
We denote a point of $U_j$ by $\xi_j=\left(\xi_j^1,...,\xi_j^n\right)$ and $\Gamma\left(U_j\times \Delta ,  \tilde{\mathcal{N}}_\mathcal{F}^* \right)$ is generated by $\omega_j^1\left(\xi_j,t\right),...,\omega_j^q\left(\xi_j,t\right)$, where
\begin{align*}
\omega_j^\alpha\left(\xi_j,t\right)=\sum_{\beta=1}^n w_j^{\alpha\beta}\left(\xi_j, t\right) d\xi_j^\beta + \sum_{\gamma=1}^m v_j^{\alpha \gamma}\left( \xi_j,t \right) dt_\gamma
\end{align*}
for some $w_j^{\alpha\beta}(\xi_j,t)\in \Gamma \left( U_j \times \Delta , \mathcal{O}_\mathcal{M} \right)$ and some $v_j^{\alpha \gamma}(\xi_j,t) \in \Gamma\left( U_j\times \Delta, \mathcal{O}_{\mathcal{M}} \right)$, and satisfies
\begin{align}\label{oui1}
\omega_j^1(\xi_j,t)\wedge \cdots \wedge \omega_j^q(\xi_j,t)\wedge d\omega _j^\alpha(\xi_j,t)=0,\,\,\,\,\,\,\,\alpha=1,...,q,\,\,\,\,\,\textnormal{on}\,\,\,\,\, U_j \times \Delta.
\end{align}
For simplicity, we assume that $U_j=\left\{\xi_j\in \mathbb{C}^m | |\xi_j|<1\right\}$ where $|\xi_j|=\max_a |\xi_j^a|$. $(\xi_j,t)\in U_j\times \Delta$ and $(\xi_k, t)\in U_k\times \Delta$ are the same point on $\mathcal{M}_\Delta$ if $\xi_j^\alpha=f_{jk}^\alpha(\xi_k,t),\alpha=1,...,n$ where $f_{jk}(\xi_k,t)$ is a holomorphic map of $\xi_k^1,..., \xi_k^n, t_1,...,t_m$, defined on $(U_k\times \Delta) \cap (U_j\times \Delta)$, and we have the following relation
\begin{align}\label{oui9}
\omega_j^\alpha(\xi_j,t)=\sum_{\beta=1}^q h_{jk}^{\alpha\beta}(\xi_k,t) \omega_k^\beta(\xi_k,t),\,\,\,\,\,\alpha=1,...,q
\end{align}
for a holomorphic function $h_{jk}^{\alpha\beta}(\xi_k,t)$ on $\left(U_k\times \Delta\right) \cap \left( U_j\times \Delta\right)$. More precisely
\begin{align}
\sum_{\gamma=1}^n w_j^{\alpha\gamma}(f_{jk}(\xi_k,t),t)\frac{\partial f_{jk}^\gamma(\xi_k,t)}{\partial z_k^\xi} &= \sum_{\beta=1}^q h_{jk}^{\alpha\beta}(\xi_k,t) w_k^{\beta\xi}(\xi_k,t)\\
\sum_{\gamma=1}^n\sum_{\eta=1}^m w_j^{\alpha\gamma}(f_{jk}(\xi_k,t),t)\frac{\partial f_{jk}^\gamma(\xi_k,t)}{\partial t_\eta}dt_\eta + \sum_{\eta=1}^m v_j^{\alpha\eta}(f_{jk}(\xi_k,t),t)dt_\eta & = \sum_{\beta=1}^q \sum_{\eta=1}^m h_{jk}^{\alpha\beta}(\xi_k,t) v_k^{\beta\eta}(\xi_k,t)dt_\eta,\,\,\,\,\,\,\,\,\,\alpha=1,...,q. \label{oui3}
\end{align} 
and on $\left(U_i\times \Delta \right)\cap \left(U_j \times \Delta \right) \cap \left( U_k \times \Delta \right)$, we have
\begin{align}
h_{ik}^{\alpha \gamma}(\xi_k,t)= \sum_{\beta=1}^q h_{ij}^{\alpha\beta}\left(f_{jk}(\xi_k,t),t\right) h_{jk}^{\beta \gamma}(\xi_k,t),\,\,\,\,\,\alpha,\gamma=1,...,q
\end{align}
By \cite{Kod05} Theorem 2.3, $M_t$ is diffeomorphic to $M_0=\pi^{-1}(0)$ as differentiable manifolds for each $t\in \Delta$. We put $M:=M_0$. By \cite{Kod05} Theorem 2.5, if we take a sufficiently small $\Delta$, there is a diffeomorphism $\Psi$ of $M\times \Delta$ onto $\mathcal{M}_\Delta$ as differentiable manifolds such that $\pi\circ \Psi$ is the projection $M\times \Delta \to \Delta$. Let $z=\left(z^1,..., z^n \right)$ be local coordinates of $M=M_0$. Then we have $\pi\circ \Psi(z,t)=t, t\in \Delta$. For $\Psi(z,t)\in U_j\times \Delta$, put
\begin{align*}
\Psi\left( z, t\right) = \left(\xi_j^1(z,t), ..., \xi_j^n(z,t), t_1,..., t_m \right)
\end{align*}
Then each component $\xi_j^\alpha=\xi_j^\alpha(z,t),\alpha=1,...,n$ is a $C^\infty$ function. If we identify $\mathcal{M}_\Delta =\Psi\left(M\times \Delta \right)$ with $M\times \Delta$ via $\Psi$, $\left(\mathcal{M}_\Delta,  \tilde{\mathcal{N}}_{\mathcal{F}_\Delta}^* \right)$ is considered as a complex manifold with the complex structure defined on the $C^\infty$ manifold $M\times \Delta$ by the system of local complex coordinates on $U_j\times \Delta$
\begin{align}\label{ne1}
\left\{ \left(\xi_j,t \right)| j=1 ,2, 3,...   \right\},\,\,\,\,\,(\xi_j,t)=\left( \xi_j^1(z,t), ..., \xi_j^n(z,t), t_1,..., t_m \right)
\end{align}
and unfolding $\tilde{\mathcal{N}}_{\mathcal{F}_\Delta}^*$ on $U_j\times \Delta$ with respect to the coordinates $(\ref{ne1})$ generated by $\omega_j^\alpha\left(\xi_j,t \right)= w_j^\alpha(\xi_j,t)+ v_j^\alpha(\xi_j,t)= \sum_{\beta=1}^n w_j^{\alpha\beta}(\xi_j,t)d\xi_j^\beta + \sum_{\gamma=1}^m v_j^{\alpha \gamma}(\xi_j,t) dt_\gamma ,\alpha=1,...,q$. 
We note that since $\left(z^1,..., z^n\right)$ and $\left(\xi_j^1(z,0), ..., \xi_j^n(z,0)\right)$ are local complex coordinates on $M=M_0$, $\xi_j^\alpha(z,0)$ are holomorphic functions of $z^1,...,z^n,\alpha=1,...,n$. We also note that if we take $\Delta$ sufficiently small, we have
\begin{align}\label{ne2}
\det\left( \frac{\partial \xi_j^\alpha(z,t)}{\partial z^\lambda} \right)_{\alpha, \lambda=1,...,n} \ne 0\,\,\,\,\,\, \, \textnormal{for}\,\,\,t\in \Delta
\end{align}
We put $\mathscr{U}_j=\Psi^{-1}\left(U_j\times \Delta \right)$. Then $\mathscr{U}_j \subset M\times \Delta$ is the domain of $\xi_j^{\alpha}(z,t)$. From $(\ref{ne2})$, we can define a $(0,1)$-form $\varphi_j^\lambda(z,t)= \sum_{v=1}^n \varphi_{jv}^\lambda(z,t) d\bar{z}^v$ for each $\lambda=1,...,n $ such that $\bar{\partial}\xi_j^\alpha =  \sum_{\lambda, v =1}^n \varphi_{jv}^\lambda(z,t)\frac{\partial \xi_j^\alpha}{\partial z^\lambda}d\bar{z}^v$. Then on $\mathscr{U}_j \cap \mathscr{U}_k$, we have $\sum_{\lambda=1}^n \varphi_j^\lambda(z,t)\frac{\partial}{\partial z^\lambda} = \sum_{\lambda=1}^n \varphi_k^\lambda(z,t) \frac{\partial}{\partial z^\lambda}$ (see \cite{Kod05} p.262). Then if for $(z,t)\in \mathscr{U}_j$, we define $\varphi(z,t):=\sum_{\lambda=1}^n \varphi_j^\lambda(z,t)\frac{\partial}{\partial z^\lambda}$, then $\varphi(t):=\varphi(z,t)\in A^{0,1}\left(M, \Theta_M \right)$ for $t\in \Delta$, and satisfies
\begin{align}\label{ne3}
\varphi(0)=0,\,\,\,\,\,\,\,\bar{\partial}\varphi(t)- \frac{1}{2}\left[ \varphi(t), \varphi(t) \right]=0
\end{align}
(see \cite{Kod05} p.263, p.265), We also point out that
\begin{theorem}\label{ne4}
If we take a sufficiently small polydisk $\Delta$ as above, then for $t\in \Delta$, a local $C^\infty$ function $f$ on $M$ is holomorphic with respect to the complex structure $M_t$ if and only if
\begin{align*}
\left( \bar{\partial} - \bar{\partial} \varphi(t)\right)f=0
\end{align*} 
\end{theorem}
\begin{proof}
See \cite{Kod05} Theorem 5.3 p.263
\end{proof}

We note that in the view of $(\ref{ne2})$ we let $\left(\psi_j^{\alpha\lambda}(z,t)\right)_{\alpha,\lambda=1,...,n}$ be the invertible matrix of $\left( \frac{\partial \xi_j^\alpha(z,t)}{\partial z_\lambda}\right)_{\alpha,\lambda=1,...,n}$, so that we have
\begin{align}\label{ca1}
\sum_{\gamma=1}^n \frac{\partial \xi_j^\eta(z,t)}{\partial z^\gamma} \psi_j^{\gamma\beta}(z_j,t)=\delta_\eta^\beta,\,\,\,\,\,\,\,\,\,\,\,\,\,\,\,\,\,\,\,\,\,\, \sum_{\gamma=1}^n \psi_j^{\eta \gamma}(z,t)\frac{\partial \xi_j^{\gamma}(z,t)}{\partial z^\beta} = \delta_\eta^\beta
\end{align}

We note that $\omega_j^\alpha(\xi_j,t)= w_j^\alpha(\xi_j,t)+ v_j^\alpha(\xi_j,t)=\sum_{\beta=1}^n w_j^{\alpha\beta}(\xi_j,t)d\xi_j^\beta+\sum_{\gamma=1}^r v_j^{\alpha\gamma}(\xi_j,t)dt_\gamma$ implies that
\begin{align} \label{ue1}
&\omega_j^\alpha(\xi_j,t)=w_j^\alpha(\xi_j,t)+v_j^\alpha(\xi_j,t)\\
&=\sum_{\beta,\eta=1}^n w_j^{\alpha\beta}\left(\xi_j(z,t),t \right)\frac{\partial \xi_j^\beta}{\partial z^\eta}dz^\eta +\sum_{\gamma=1}^r\left(v_j^{\alpha\gamma}\left(\xi_j(z,t),t \right)+\sum_{\beta=1}^n w_j^{\alpha\beta}\left(\xi_j(z,t),t \right)\frac{\partial \xi_j^\beta(z,t)}{\partial t_\gamma} \right) dt_\gamma + \sum_{\beta,\eta=1}^n w_j^{\alpha\beta}\left(\xi_j(z,t),t \right)\frac{\partial \xi_j^\beta}{\partial \bar{z}^\eta}d\bar{z}^\eta \notag \\
&=\sum_{\eta=1}^n P_j^{\alpha\eta}(z,t) dz^\eta +\sum_{\gamma=1}^r Q_j^{\alpha \gamma}(z,t)dt_\gamma+\cdots. \notag
\end{align}
From the first term $\sum_{\beta,\eta=1}^n P_j^{\alpha\eta}(z,t)dz^\eta$ in $(\ref{ue1})$, and from $(\ref{ca1})$, we can recover $w_j^\alpha(\xi_j,t)$ since we have $w_j^{\alpha\eta}=\sum_{\beta,\gamma=1}^n P_j^{\alpha \gamma}\psi_j^{\gamma\eta}$, and from the second term $\sum_{\gamma=1}^r Q_j^{\alpha \gamma}(z,t)dt_\gamma$ and $(\ref{ca1})$, we can recover $v_j^\alpha(\xi_j,t)$ since $v_j^{\alpha \gamma}= Q^{\alpha\gamma}(z,t)-\sum_{\eta,\beta=1}^n P^{\alpha\eta}(z,t)\psi_j^{\eta\beta}\frac{\partial \xi_j^\beta}{\partial t_\gamma}$. Keeping this in mind, we prove
\begin{theorem}
If we take a sufficiently small polydisk $\Delta$ as above, then for $t\in \Delta$, a $C^\infty$ 1-form $w(z,t)+v(z,t)=\sum_{\alpha=1}^n w^\alpha(z,t) dz^\alpha + \sum_{\gamma=1}^r v^\gamma(z,t)dt_\gamma$ defines a holomorphic $1$-form  $\sum_{\alpha,\beta=1}^n w^\alpha \psi_j^{\alpha\beta}d\xi_j^\beta+ \sum_{\gamma=1}^r\left(v^\gamma-\sum_{\alpha,\beta=1}^n w^\alpha\psi_j^{\alpha\beta}\frac{\partial \xi_j^\beta}{\partial t_\gamma} \right)dt_\gamma$ with respect to the complex structure on $M\times \Delta$ defined by $\varphi(t)$ if and only if 
\begin{align}
\bar{\partial} w(z,t)- \mathcal{L}_{\varphi}\left(w(z,t)\right)&=0 \label{ue3} \\
\bar{\partial} v^\gamma (z,t) - \left[\varphi,  v^\gamma (z,t) \right]  &=  i_{\frac{\partial \varphi}{\partial t_\gamma}}\left(w(z,t)\right),\,\,\,\,\,\,\,\gamma=1,...,r \label{ue4}
\end{align}
\end{theorem}

For the formula on bracket and Lie derivative, see Part I.

\begin{proof}
We have shown in Part I that
\begin{align}\label{ue2}
\bar{\partial} \left(\sum_{\alpha=1}^n w^\alpha \psi_j^{\alpha\beta} \right)-\left[ \sum_{\lambda, v=1}^n \varphi_{jv}^\lambda d\bar{z}^v \frac{\partial}{\partial z^\lambda}, \sum_{\alpha=1}^n w^\alpha \psi_j^{\alpha\beta}\right]=0 \iff \bar{\partial} w- L_\varphi(w)=0
\end{align}
Then by Theorem \ref{ne4}, $(\ref{ue2})$ implies that $(\ref{ue3})$ is equivalent to that $\sum_{\alpha, \beta=1}^n w^\alpha \psi_j^{\alpha\beta} d\xi_j^\beta$ is holomorphic with respect to $\varphi(t)$. On the other hand, by Theorem \ref{ne4} and $(\ref{ue2})$, $\sum_{\gamma=1}^r\left(v^\gamma-\sum_{\alpha,\beta=1}^n w^\alpha\psi_j^{\alpha\beta}\frac{\partial \xi_j^\beta}{\partial t_\gamma} \right)dt_\gamma$ is holomorphic with respect to $\varphi(t)$ if and only if
\begin{align}
0&= \bar{\partial}\left( v^\gamma - \sum_{\alpha, \beta=1}^n w^\alpha \psi_j^{\alpha\beta}\frac{\partial \xi_j^\beta}{\partial t_\gamma} \right)- \left[ \varphi,  v^\gamma - \sum_{\alpha, \beta=1}^n w^\alpha \psi_j^{\alpha\beta}\frac{\partial \xi_j^\beta}{\partial t_\gamma}  \right]  \label{ue5} \\
&= \bar{\partial} v^\gamma - \left[ \varphi, v^\gamma \right] - \sum_{\alpha, \beta=1}^n w^\alpha \psi_j^{\alpha\beta} \left( \bar{\partial}\left( \frac{ \partial \xi_j^\beta}{\partial t_\gamma}   \right)  - \left[ \varphi, \frac{\partial \xi_j^\beta}{\partial t_\gamma} \right]  \right) \notag
\end{align}
Since $\bar{\partial} \xi_j^\beta- \left[ \varphi, \xi_j^\beta\right] =0$,  by taking the derivative with respect to $t_\gamma$, we have
\begin{align}\label{ue6}
\bar{\partial}\left( \frac{\partial \xi_j^\beta}{\partial t_\gamma} \right) - \left[ \frac{\partial \varphi}{\partial t_\gamma}, \xi_j^\beta \right] - \left[ \varphi,  \frac{\partial \xi_j^\beta}{\partial t_\gamma} \right] = 0
\end{align}
Then from $(\ref{ue6})$ and $(\ref{ca1})$, we see that $(\ref{ue5})$ is equivalent to
\begin{align*}
0&=\bar{\partial} v^\gamma - \left[ \varphi, v^\gamma \right] - \sum_{\alpha, \beta , \lambda, v=1}^n w^\alpha \psi_j^{\alpha\beta}\frac{\partial \xi_j^\beta}{\partial z^\lambda} \frac{\partial \varphi_{jv}^\lambda}{\partial t_\gamma}d\bar{z}^v = \bar{\partial} v^\gamma - \left[ \varphi, v^\gamma \right] - \sum_{\alpha=1}^n w^\alpha \frac{\partial \varphi_{jv}^\alpha}{\partial t_\gamma} d\bar{z}^v  \\
&= \bar{\partial} v^\gamma - \left[ \varphi, v^\gamma \right] - i_{\sum_{\lambda, v=1}^n \frac{\partial \varphi_{jv}^\lambda}{\partial t_\gamma} d\bar{z}^v \frac{\partial}{\partial z^\lambda}}\left(\sum_{\alpha=1}^n w^\alpha dz^\alpha \right) =  \bar{\partial} v^\gamma - \left[ \varphi, v^\gamma \right]- i_{\frac{\partial \varphi}{\partial t_\gamma}}\left( w\right)
\end{align*}
which implies $(\ref{ue4})$.

\end{proof}

With this preparation, we shall prove theorem of existence for unfoldings of foliated analytic structures for one-parameter families.

\subsection{Theorem of existence for unfoldings of foliated complex analytic structures for one-parameter families}

\begin{theorem}[Theorem of existence for unfoldings of foliated complex analytic structures for one-parameter families (compare \cite{MK06} Theorem for one-parameter families p.55)]\label{u56}
Let $\left(M, \mathcal{N}_{\mathcal{F}_0}^*\right)$ be a compact foliated complex manifold with $\mathcal{N}_{\mathcal{F}_0}^*$ being a locally free subsheaf of $\Omega_M^1$ such that $\mathbb{H}^2\left(M, \tilde{\mathcal{N}}_{\mathcal{F}_0}^{*\bullet} \right)=0$.  Let $\phi\in \mathbb{H}^1 \left(M, \tilde{\mathcal{N}}_{\mathcal{F}_0}^{* \bullet} \right)$. Then there exists an one-parameter unfolded analytic family $(\mathcal{M}, \tilde{\mathcal{N}}_\mathcal{F}^*, B, \pi )$ where $0\in B \subset \mathbb{C}$, and $\pi^{-1}(0)$ induces $\left(M, \mathcal{N}_{\mathcal{F}_0}^*\right)$, such that the image of the unfolded Kodaira-Spencer map at $0$ is generated by $\phi$. In other words, any cohomology class $\phi\in \mathbb{H}^1 \left(M, \mathcal{N}_{\mathcal{F}_0}^{ * \bullet} \right)$ can be integrated to make an one-parameter unfolded complex analytic family.
\end{theorem}

\begin{proof}
We may assume the following:
\begin{enumerate}
\item $M$ is covered by a finite number of coordinate neighborhoods $U_i(i\in I)$ with a system of coordinates $z_i=\left(z_i^1,...,z_i^n\right)$ such that $U_i=\left\{z_i\in \mathbb{C}^n|\max_{1\leq \alpha \leq n} | z_i^\alpha| <1 \right\}$ \label{us2}
\item $\mathcal{N}_{\mathcal{F}_0}^*$ is generated by $\left\{w_{0i}^1,...,w_{0i}^q \right\}$ with $w_{0i}^\alpha\in \Gamma\left(U_i, \mathcal{N}_{\mathcal{F}_0}^*\right)$ and $w_{0j}^\alpha=\sum_{\beta=1}^p h_{0ij}^{\alpha\beta} w_{0j}^\beta$ on $U_i\cap U_j$ for $h_{0ij}^{\alpha\beta}(z_j)\in \Gamma\left( U_{ij}, \mathcal{O}_M \right),\alpha=1,...,q$.
\item $w_{0i}^1\wedge \cdots \wedge w_{0i}^q\wedge dw_{0i}^\alpha=0$ for $\alpha=1,...,q$. 
\item $z_i$ coincides with $z_j$ if and only if $z_i=f_{ij}(z_j)$.
\item for $x\in U_i-S$, where $S=\textnormal{Sing}\left(\mathcal{F}_0\right)$, there is a neighborhood $U_{i_x} \subset U_i-S$ of $x$ such that $dw_{0i}^\alpha =\sum_{\beta=1}^q a_{0i_x}^{\alpha\beta}  \wedge w_{0i}^\beta$ on $U_{i_x}$ for some $a_{0i_x}^{\alpha\beta}\in \Gamma\left(U_{i_x}, \Omega_M^1 \right)$.\label{us3}
\end{enumerate}
Then we note that from $\sum_{\gamma=1}^q h_{0ik}^{\alpha\gamma}w_{0k}^\gamma=w_{0i}^\alpha=\sum_{\beta=1}^q h_{0ij}^{\alpha\beta} w_{0j}^\beta=\sum_{\beta,\gamma=1}^q h_{0ij}^{\alpha\beta} h_{0jk}^{\beta\gamma} w_{0k}^\gamma$, we have $
h_{0ik}^{\alpha\gamma}=\sum_{\beta=1}^q h_{0ij}^{\alpha\beta} h_{0jk}^{\beta\gamma}$.

Let $B=\{t\in \mathbb{C}||t|<\epsilon\}$ with a sufficiently small $\epsilon >0 $. Our purpose is to construct

\begin{enumerate}
\item a $C^\infty (0,1)$-form $\varphi(t)=\varphi_1t+ \varphi_2t^2 +\varphi_3t^3+\cdots$ with coefficients in $\Theta_M$ depending holomorphically on $t$,
\item $C^\infty$-differential form $w_i^\alpha(z_i,t)$ on $U_i\times B$ of the form $w_i^\alpha(z_i,t)+ v_i^\alpha(z_i,t)dt=\sum_{\beta=1}^q w_i^{\alpha\beta}(z_i,t)dz_i^\beta+ v_i^\alpha(z_i,t)dt$ depending holomorphically on $t$,
\item $q\times q$ matrices $H_{ij}:=\left( h_{ij}^{\alpha\beta}(z_j,t) \right)$ with components $C^\infty$ functions $h_{ij}^{\alpha\beta}(z_j,t)$ on $U_{ij}\times B$ depending holomorphically on $t$
\item $C^\infty$ differential form $a_{i_x}^{\alpha\beta}(z_i,t)+ b_{i_x}^{\alpha\beta}(z_i,t)dt$ on $U_{i_x}\times B_{i_x}\subset U_i\times B$ for $x\in U_{i_x}\subset U_i-S$ where $a_{i_x}^{\alpha\beta}(z_i,t)$ is of the form $\sum_{\gamma=1}^n a_{i_x\gamma}^{\alpha\beta}(z_i,t)dz_i^\gamma$ on $U_{i_x}\times B_{i_x}$, and $S=\textnormal{Sing}\left( \mathcal{F}_0 \right)$ such that
\end{enumerate}
\begin{align}\label{ue10}
\varphi(0)=0,\,\,\,\,\,\,\,\,\, w_i^\alpha(z_i,t)=w_i^\alpha,\,\,\,\,\,\,\,\,\, h_{ij}^{\alpha\beta}(z_j,0)=h_{ij}^{\alpha\beta}(z_j),\,\,\,\,\,\,\,\,a_{i_x}^{\alpha\beta}(z_i,0)=a_{0i_x}^{\alpha\beta}(z_i)\,\,\,\text{on}\,\,\, U_{i_x}\times B_{i_x}
\end{align}
\begin{align}
\bar{\partial}\varphi &-\frac{1}{2}\left[\varphi,\varphi\right]=0 \label{u11}\\
\bar{\partial} w_i^\alpha(z_i,t)&- \mathcal{L}_{\varphi} \left(w_i^\alpha(z_i,t)\right)=0\label{u13}\\
\bar{\partial} v_i^{\alpha}(z_i,t) &-\left[ \varphi, v_i^{\alpha}(z_i,t)\right]=i_{\frac{\partial \varphi}{\partial t}}\left(w_i^{\alpha}(z_i,t) \right) \label{uut15}\\
\bar{\partial} h_{ij}^{\alpha\beta}(z_j,t)&-\left[\varphi, h_{ij}^{\alpha\beta}(z_j,t) \right]=0 \label{uut16}\\
w_i^\alpha(z_i,t)+v_i^\alpha(z_i,t)dt&=\sum_{\beta=1}^q h_{ij}^{\alpha\beta}(z_i,t)\left(w_j^\beta(z_j,t)+v_j^\beta(z_j,t)dt \right)  \label{uut17} \\
\left(w_i^1(z_i,t)+v_i^1(z_i,t)dt\right)\wedge \cdots  \wedge &\left(w_i^q(z_i,t)+v_i^q(z_i,t)dt\right)\wedge \partial_{M\times B} \left( w_i^\alpha(z_i,t)+v_i^\alpha(z_i,t) dt\right)=0 \label{uut18}\\
H_{ij}H_{jk}=H_{ik}\iff h_{ik}^{\alpha\gamma}&(z_k,t)=\sum_{\beta=1}^q h_{ij}^{\alpha\beta}\left(f_{ik}(z_k),t \right) h_{jk}^{\beta\gamma}(z_k,t)\,\,\,\,\,\,\,\text{on}\,\,\,U_{ijk}\times B \label{uut19}\\
\partial_{M\times B}\left(w_i^\alpha(z_i,t)+v_i^\alpha(z_i,t)dt \right)=\sum_{\beta=1}^q &\left(a_{i_x}^{\alpha\beta}(z_i,t)+b_{i_x}^{\alpha\beta}(z_i,t)dt\right)\wedge \left(w_i^\beta(z_i,t)+v_i^\beta(z_i,t)dt\right)\,\,\,\,\,\textnormal{on}\,\,\,U_{i_x}\times B_{i_x} \label{u12}
\end{align}

\subsection{Existence of formal solutions}\

We construct solutions of (\ref{u11})-(\ref{u12}) which are formal power series in $t$. We recall from Part I
\begin{notation}\label{te27}
We write the power series expansion of a holomorphic function $P(t)$ in $t_1,..,t_m$ defined on a neighborhood of the origin $0$ in the form$:$ $P(t)=P(0)+ P_1(t)+ \cdots + P_\mu(t)+ \cdots $, where each $P_\mu(t)$ denotes a homogenous polynomial of degree $\mu$ in $t_1,..., t_m$. We set $P^\mu(t):=P(0)+ P_1(t)+ \cdots + P_\mu(t)$. For any power series $P(t)$ and $Q(t)$ in $t=(t_1,..., t_m)$, we indicate $P(t)\equiv_\mu Q(t)$ contains no term of degree $\leq \mu$ in $t$. Then $P^\mu=0$ means $P\equiv_\mu 0$.
\end{notation}
From Notation \ref{te27}, we write $\varphi=\sum_{\mu=0}^\infty \varphi_\mu, h_{ij}^{\alpha\beta \mu}=\sum_{\mu=0}^\infty h_{ij|\mu}^{\alpha\beta}, w_i^\alpha=\sum_{\mu=0}^\infty w_{i|\mu}^\alpha,  v_i^\alpha= \sum_{\mu=0}^\infty v_{i|\mu}^\alpha$ and $a_{i_x}^{\alpha\beta}= \sum_{\mu=0}^\infty a_{i_x|\mu}^{\alpha\beta}$ and $b_{i_x}^{\alpha\beta}=\sum_{\mu=0}^{\infty} b_{i_x|\mu}^{\alpha\beta}$. In view of $(\ref{ue10})$, we set $\varphi_0=0, w_{i|0}^\alpha=w_{0i}^\alpha(z_i)$, and $h_{h|0}^{\alpha\beta}(z_j,t)=h_{0ij}^{\alpha\beta}(z_j)$ and $a_{i_x|0}^{\alpha\beta}= a_{0i_x}^{\alpha\beta}$. Then $(\ref{u13})-(\ref{u12})$ are equivalent to the following system of congruences
\begin{align}
\bar{\partial} \varphi^\mu - \frac{1}{2}&\left[\varphi^\mu, \varphi^\mu \right]\equiv_\mu 0 \label{u14}\\
\bar{\partial}\left(w_i^{\alpha\mu}(z_i,t)\right)&- \mathcal{L}_{\varphi^\mu}\left(w_i^{\alpha\mu}(z_i,t)\right)\equiv_\mu 0\label{u15}\\
\bar{\partial} v_i^{\alpha(\mu-1)}(z_i,t)&-\left[\varphi^\mu, v_i^{\alpha(\mu-1)}(z_i,t) \right]\equiv_{\mu-1} i_{\frac{\partial \varphi^\mu}{\partial t}}(w_i^{\alpha \mu}(z_i,t)) \label{ucc1}\\
\bar{\partial} h_{ij}^{\alpha\beta\mu}-&\left[\varphi^\mu, h_{ij}^{\alpha\beta\mu} \right]\equiv_\mu 0 \label{u16}\\
w_i^{\alpha\mu}+v_i^{\alpha(\mu-1)}dt\equiv_\mu &\sum_{\beta=1}^q h_{ij}^{\alpha\beta\mu}\left( w_j^{\beta\mu}+v_j^{\beta(\mu-1)}dt \right) \label{u17}\\
\left(w_i^{1\mu}+v_i^{1(\mu-1)} dt\right)\wedge \cdots \wedge& \left(w_i^{q\mu}+ v_i^{q(\mu-1)}dt \right)\wedge \partial_{M\times B} \left(w_i^{\alpha\mu} + v_i^{\alpha(\mu-1)}dt\right)\equiv_\mu 0 \label{u18}\\
h_{ik}^{\alpha\gamma\mu}(z_k,t)\equiv_\mu &\sum_{\beta=1}^q h_{ij}^{\alpha\beta\mu}(f_{jk},t)h_{jk}^{\beta\gamma\mu}(z_k,t) \label{ucc5}\\
\partial_{M\times B} \left(w_i^{\alpha\mu} +v_i^{\alpha(\mu-1)} dt \right)\equiv_\mu &\sum_{\beta=1}^q \left(a_{i_x}^{\alpha\beta\mu}(z_i,t)+b_{i_x}^{\alpha\beta(\mu-1)}(z_i,t)dt\right) \wedge \left( w_i^{\beta\mu}+ v_i^{\beta (\mu-1) } dt \right) \label{u19}
\end{align}

We construct solutions of $(\ref{u14})_\mu-(\ref{u19})_\mu$ by induction on $\mu$. We suppose that $\varphi^{\mu-1}, w_i^{\alpha(\mu-1)}, v_i^{\alpha(\mu-2)}$ and $h_{ij}^{\alpha\beta(\mu-1)}, a_{i_x}^{\alpha\beta(\mu-1)}$ and $b_{i_x}^{\alpha\beta(\mu-2)}$  satisfying $(\ref{u14})_{\mu-1}-(\ref{u19})_{\mu-1}$ are already determined. 

We define homogenous polynomials

\begin{align}
\xi_\mu&\equiv_\mu \bar{\partial} \varphi^{\mu-1}-\frac{1}{2} \left[ \varphi^{\mu-1},\varphi^{\mu-1} \right]\label{u20}
\end{align}
\begin{align}
-A_{i|\mu}^\alpha & \equiv_\mu \bar{\partial} w_i^{\alpha(\mu-1)}- L_{\varphi^{\mu-1}}\left(w_i^{\alpha(\mu-1)}\right)\label{u21}
\end{align}
\begin{align}
-\psi_{i|\mu-1}^{\alpha} &\equiv_{\mu-1} \bar{\partial}v_i^{\alpha (\mu-2)}- \left[\varphi^{\mu-1}, v_i^{\alpha(\mu-2)} \right] -i_{\frac{\partial \varphi^{\mu-1} }{\partial t}}\left(w_i^{\alpha(\mu-1)} \right)\label{u22}
\end{align}
\begin{align}
-\sum_{\eta=1}^q B_{ij|\mu}^{\alpha\eta} h_{0ij}^{\eta\beta} &\equiv_\mu \bar{\partial} h_{ij}^{\alpha\beta(\mu-1)}-\left[ \varphi^{\mu-1}, h_{ij}^{\alpha\beta(\mu-1)} \right]
\end{align}
We define 
\begin{align}
C_{ij|\mu}^\alpha+ N_{ij|\mu-1}^\alpha dt&\equiv_\mu \left( w_i^{\alpha(\mu-1)}+ v_i^{\alpha(\mu-2)} dt \right) -\sum_{\beta=1}^q h_{ij}^{\alpha\beta(\mu-1)}\left(w_j^{\beta(\mu-1)}+v_j^{\beta(\mu-2)} dt \right) \notag
\end{align}
so that 
\begin{align}
C_{ij|\mu}^\alpha &\equiv_\mu w_i^{\alpha(\mu-1)} - \sum_{\beta=1}^q h_{ij}^{\alpha\beta(\mu-1)} w_j^{\beta(\mu-1)}
\end{align}
\begin{align}
N_{ij|\mu-1}^\alpha &\equiv_{\mu-1} v_i^{\alpha(\mu-2)} - \sum_{\beta=1}^q h_{ij}^{\alpha\beta(\mu-1)} v_j^{\beta(\mu-2)} 
\end{align}
\begin{align}
\sum_{\beta=1}^q D_{ijk|\mu}^{\alpha\beta} h_{0ik}^{\beta\gamma}&\equiv_\mu  h_{ik}^{\alpha\gamma(\mu-1)}(z_k,t)-\sum_{\beta=1}^q h_{ij}^{\alpha\beta(\mu-1)}(f_{jk},t) h_{jk}^{\beta\gamma(\mu-1)}(z_k,t)
\end{align}
We define
\begin{align}
F_{i|\mu}^\alpha=E_{i|\mu}^\alpha+  M_{i|\mu-1}^{\alpha} \wedge dt &\equiv_\mu \left(w_i^{1(\mu-1)}+v_i^{1(\mu-2)} dt \right) \wedge \cdots \wedge \left( w_i^{q(\mu-1)}+v_i^{q(\mu-2)}dt \right) \wedge \partial_{M\times B} \left( w_i^{\alpha(\mu-1)}+ v_i^{\alpha(\mu-2)}dt \right) \notag
\end{align}
Then we have
\begin{align}
E_{i|\mu}^\alpha \equiv_\mu w_i^{1(\mu-1)} \wedge \cdots \wedge w_i^{q(\mu-1)}\wedge \partial w_i^{\alpha (\mu-1)} \label{utt2}
\end{align}
\begin{align}
 &M_{i|\mu-1}^\alpha \wedge dt\\
 &\equiv_{\mu}  \sum_{\beta=1}^q w_i^{1(\mu-1)}\wedge \cdots \wedge \overbrace{v_i^{\beta(\mu-2)} dt}^{\beta-\textnormal{th}} \wedge \cdots \wedge w_i^{q(\mu-1)}\wedge \partial w_i^{\alpha(\mu-1)}                   + w_i^{1(\mu-1)}\wedge \cdots \wedge w_i^{q(\mu-1)} \wedge \partial v_i^{\alpha(\mu-2)}\wedge dt \notag\\
   &\,\,\,\,\,\,- w_i^{1(\mu-1)} \wedge \cdots \wedge w_i^{q(\mu-1)}\wedge \frac{\partial w_i^{\alpha(\mu-1)}}{\partial t} \wedge dt  \notag \\
                     &\equiv_{\mu} \left(\sum_{\beta=1}^q (-1)^{q-\beta} w_i^{1(\mu-1)}\wedge \cdots \wedge \overbrace{v_i^{\beta(\mu-2)} }^{\beta-\textnormal{th}} \wedge \cdots \wedge w_i^{q(\mu-1)}\wedge \partial w_i^{\alpha(\mu-1)}     + w_i^{1(\mu-1)}\wedge \cdots \wedge w_i^{q(\mu-1)} \wedge  \partial v_i^{\alpha(\mu-2)}   \right)\wedge dt   \notag \\
                     &\,\,\,\,\,\,\,\,\,\,- w_i^{1(\mu-1)}\wedge \cdots \wedge w_i^{q(\mu-1)} \wedge \frac{\partial w_i^{\alpha(\mu-1)}}{\partial t} \wedge dt, \notag
\end{align}
so that
\begin{align}                     
 M_{i|\mu-1}^\alpha &\equiv_{\mu-1} \sum_{\beta=1}^q (-1)^{q-\beta}v_i^{\beta(\mu-2)} w_i^{1(\mu-1)}\wedge \cdots \wedge \widehat{w_i^{\beta(\mu-1)} } \wedge \cdots \wedge w_i^{q(\mu-1)}\wedge \partial w_i^{\alpha(\mu-1)}     + w_i^{1(\mu-1)}\wedge \cdots \wedge w_i^{q(\mu-1)} \wedge \partial v_i^{\alpha(\mu-2)} \label{utt3} \\
 &\,\,\,\,\,\,\,\,\,- w_i^{1(\mu-1)}\wedge \cdots \wedge w_i^{q(\mu-1)} \wedge \frac{\partial w_i^{\alpha(\mu-1)}}{\partial t} \notag
 \end{align}
 We define
 \begin{align}                   
L_{i_x|\mu}^\alpha +  P_{i_x|\mu-1}^{\alpha}\wedge dt \equiv_\mu \partial_{M\times B}\left(w_i^{\alpha(\mu-1)}+v_i^{\alpha(\mu-2)} dt \right)-\sum_{\beta=1}^q \left( a_{i_x}^{\alpha\beta(\mu-1)}+ b_{i_x}^{\alpha\beta(\mu-2)}dt\right) \wedge \left(w_i^{\beta(\mu-1)}+ v_i^{\beta(\mu-2)}dt \right), \notag
\end{align}
so that
\begin{align}
L_{i_x|\mu}^\alpha &\equiv_\mu \partial w_i^{\alpha(\mu-1)} - \sum_{\beta=1}^q a_{i_x}^{\alpha\beta(\mu-1)} \wedge w_i^{\beta(\mu-1)}\,\,\,\,\,\textnormal{on}\,\,\, U_{i_x}
\end{align}
\begin{align}
P_{i_x|\mu-1}^\alpha &\equiv_{\mu-1}-\frac{\partial w_i^{\alpha(\mu-1)}}{\partial t} +\partial v_i^{\alpha(\mu-2)}- \sum_{\beta=1}^q a_{i_x}^{\alpha\beta(\mu-1)} v_i^{\beta(\mu-2)} + \sum_{\beta=1}^q b_{i_x}^{\alpha\beta(\mu-2)} w_i^{\beta(\mu-1)}\,\,\,\,\,\textnormal{on}\,\,\, U_{i_x} \label{uut25}
\end{align}

\begin{lemma}  \label{ue15}
We have the following equalities$:$
\begin{align}
\bar{\partial} \xi_\mu&=0 \label{u29}
\end{align}
\begin{align}
\sum_{\beta,\eta=1}^q h_{0ij}^{\alpha\beta} D_{jkl|\mu}^{\beta\eta}h_{0ji}^{\eta\xi}-D_{ijk|\mu}^{\alpha\xi}&+D_{ijl|\mu}^{\alpha\xi}- D_{ijk|\mu}^{\alpha\xi}=0 \label{u30} 
\end{align}
\begin{align}
\bar{\partial} A_{i|\mu}^\alpha&=L_{\xi_\mu}\left(w_{0i}^\alpha\right) \label{u31} 
\end{align}
\begin{align}
\bar{\partial} \psi_{i|\mu-1}^{\alpha\gamma}&=i_{\frac{\partial \xi_\mu}{\partial t} }\left(w_{0i}^\alpha\right)\label{u32}
\end{align}
\begin{align}
\sum_{\eta,\gamma=1}^q h_{0ij}^{\alpha\eta} B_{jk|\mu}^{\eta\gamma} h_{0ji}^{\gamma\xi}- &B_{ij|\mu}^{\alpha \xi} + B_{ij|\mu}^{\alpha \xi} = \bar{\partial} D_{ijk|\mu}^\xi \label{u33}
\end{align}
\begin{align}
\sum_{\gamma=1}^q h_{0ij}^{\alpha\gamma} C_{jk|\mu}^\gamma-  C_{ik|\mu}^\alpha& + C_{ij|\mu}^\alpha = \sum_{\gamma=1}^q D_{ijk|\mu}^{\alpha\gamma} w_{0i}^\gamma \label{u34} 
\end{align}
\begin{align}
\sum_{\gamma=1}^q h_{0ij}^{\alpha\gamma} N_{jk|\mu-1}^\gamma-  N_{ik|\mu-1}^\alpha& + N_{ij|\mu-1}^\alpha = 0 \label{u35}
\end{align}
\begin{align}
\sum_{\beta=1}^q h_{0ij}^{\alpha\beta} A_{j|\mu}^\beta - A_{i|\mu}^\alpha & = \bar{\partial} C_{ij|\mu}^\alpha - \sum_{\eta=1}^q B_{ij|\mu}^{\alpha\eta} w_{0i}^\eta \label{u36}
\end{align}
\begin{align}
\sum_{\beta=1}^q h_{0ij}^{\alpha\beta} \psi_{j|\mu-1}^\beta - & \psi_{i|\mu-1}^\alpha  = \bar{\partial} N_{ij|\mu-1}^\alpha \label{u37} 
\end{align}
\begin{align}
\bar{\partial} E_{i|\mu}^\alpha= \sum_{\beta=1}^q (-1)^\beta w_{0i}^\beta \wedge \cdots \wedge  \overbrace{ A_{i|\mu}^\beta}^{\beta-\textnormal{th}} &\wedge \cdots \wedge  w_{0i}^q \wedge dw_{0i}^\alpha +(-1)^q w_{0i}^1 \wedge \cdots \wedge w_{0i}^q \wedge \partial A_{i|\mu}^\alpha  \label{u38} 
\end{align}
\begin{align}
\bar{\partial} M_{i|\mu-1}^\alpha = (-1)^q \sum_{\beta=1}^q w_{0i}^1\wedge \cdots \wedge &\overbrace{\psi_{i|\mu-1}^\beta}^{\beta-\textnormal{th}} \wedge \cdots w_{0i}^q \wedge dw_{0i}^\alpha +  (-1)^q w_{0i}^1\wedge \cdots \wedge w_{0i}^q \wedge d\psi_{i|\mu-1}^\alpha  +  (-1)^q w_{0i}^1\wedge \cdots \wedge w_{0i}^q \wedge \frac{\partial A_{i|\mu}^\alpha}{\partial t}       \label{u39} 
\end{align}
\begin{align}
E_{i|\mu}^\alpha- \sum_{\beta=1}^q h_{0ij}^{\alpha\beta} \det\left( h_{0ij}^{\gamma\eta} \right) E_{j|\mu}^\beta &= \left(  \sum_{\beta=1}^q w_{0i}^1\wedge \cdots \wedge \overbrace{C_{ij|\mu}^\beta}^{\beta-\textnormal{th}} \wedge \cdots \wedge w_{0i}^q  \right) \wedge dw_{0i}^\alpha + w_{0i}^1 \wedge \cdots \wedge w_{0i}^q \wedge  \partial C_{ij|\mu}^\alpha \label{u40}
\end{align}
\begin{align}
M_{i|\mu-1}^\alpha - \sum_{\beta=1}^q h_{0ij}^{\alpha\beta} \det\left( h_{0ij}^{\gamma\eta} \right) M_{j|\mu-1}^\beta &= \sum_{\beta=1}^q (-1)^{q-\beta} N_{ij|\mu-1}^\beta w_{0i}^1\wedge \cdots \wedge  \widehat{w_{0i}^\beta} \wedge \cdots \wedge w_{0i}^q \wedge dw_{0i}^\alpha + w_{0i}^1\wedge \cdots \wedge w_{0i}^q \wedge \partial N_{ij|\mu-1}^\alpha  \label{u41} \\
& - w_{0i}^1 \wedge \cdots \wedge w_{0i}^q \wedge \frac{\partial C_{ij|\mu}^\alpha}{\partial t} \notag
\end{align}
\begin{align}
E_{i|\mu}^\alpha = w_{0i}^1 \wedge \cdots \wedge w_{0i}^q \wedge L_{i_x|\mu}^\alpha& \,\,\,\,\,\,\,\text{on}\,\,\,\,\,\, U_{i_x}\subset U_i-S \label{u42}
\end{align}
\begin{align}
M_{i|\mu-1}^\alpha= w_{0i}^1\wedge \cdots \wedge w_{0i}^q \wedge P_{i_x|\mu-1}^\alpha & \,\,\,\,\,\,\,\,\,\text{on}\,\,\,\,\,\, U_{i_x} \subset U_i- S \label{u43}
\end{align}

\end{lemma}

\begin{proof}
We note that $(\ref{u29}),(\ref{u30}),(\ref{u31}),(\ref{u33}),(\ref{u34}),(\ref{u36}),(\ref{u38}),(\ref{u40})$ and $(\ref{u42})$ are proved in Part I. We prove $(\ref{u32})$. In fact,
\begin{align*}
\bar{\partial} \psi_{i|\mu-1}^\alpha & \equiv_{\mu-1} \bar{\partial} \left[\varphi^{\mu-1}, v_i^{\alpha(\mu-2)}\right] + \bar{\partial} i_{\frac{\partial \varphi^{\mu-1}}{\partial t}}\left( w_i^{\alpha(\mu-1)} \right) \\
                                                     & \equiv_{\mu-1} \left[ \bar{\partial} \varphi^{\mu-1}, v_i^{\alpha(\mu-2)}  \right] -\left[ \varphi^{\mu-1}, \bar{\partial} v_i^{\alpha(\mu-2)}  \right] + i_{\frac{\partial}{\partial t}\left( \xi_\mu + \frac{1}{2}\left[ \varphi^{\mu-1}, \varphi^{\mu-1} \right] \right)}\left( w_i^{\alpha(\mu-1)} \right) + i_{\frac{\partial \varphi^{\mu-1}}{\partial t}} \left(  \mathcal{L}_{\varphi^{\mu-1}}\left( w_i^{\alpha(\mu-1)} \right) - A_{i|\mu}^\alpha \right)\\
                                                     &\equiv_{\mu-1} \left[ \xi_\mu + \frac{1}{2}\left[\varphi^{\mu-1}, \varphi^{\mu-1} \right] , v_i^{\alpha(\mu-2)}  \right] -\left[\varphi^{\mu-1}, \left[ \varphi^{\mu-1}, v_i^{\alpha(\mu-2)} \right]  + i_{\frac{\partial \varphi^{\mu-1}}{\partial t}}\left( w_i^{\alpha(\mu-1)} \right) - \psi_{i|\mu-1}^\alpha \right]\\
                                                     &+ i_{\frac{\partial \xi_\mu}{\partial t}}\left(  w_i^{\alpha(\mu-1)} \right) + i_{\left[ \varphi^{\mu-1}, \frac{\partial \varphi^{\mu-1}}{\partial t}  \right]}\left( w_i^{\alpha(\mu-1)} \right) +i_{\frac{\partial \varphi^{\mu-1}}{\partial t}}\left( \mathcal{L}_{\varphi^{\mu-1}}\left( w_i^{\alpha(\mu-1)} \right) \right) \\
                                                      &\equiv_{\mu-1} i_{\frac{\partial \xi_\mu}{\partial t}}\left( w_{0i}^\alpha \right) - \left[ \varphi^{\mu-1}, i_{\frac{\partial \varphi^{\mu-1}}{\partial t}}\left( w_i^{\alpha(\mu-1)} \right)   \right] + i_{\left[ \varphi^{\mu-1}, \frac{\partial \varphi^{\mu-1}}{\partial t}  \right]}\left( w_i^{\alpha(\mu-1)} \right) +i_{\frac{\partial \varphi^{\mu-1}}{\partial t}}\left(  \mathcal{L}_{\varphi^{\mu-1}}\left( w_i^{\alpha(\mu-1)} \right) \right)\\
                                                      &\equiv_{\mu-1} i_{\frac{\partial \xi_\mu}{\partial t}}\left( w_{0i}^\alpha \right)
\end{align*}

We prove $(\ref{u35})$. In fact,
{\Small{\begin{align*}
&C_{ik|\mu}^\alpha+ N_{ik|\mu}^\alpha dt \equiv_\mu \left( w_i^{\alpha(\mu-1)} + v_i^{\alpha(\mu-2)}dt    \right) - \sum_{\beta=1}^q h_{ik}^{\alpha\beta(\mu-1)}\left( w_k^{\beta(\mu-1)} + v_k^{\beta(\mu-2)}dt  \right)\\
&\equiv_\mu \sum_{\beta=1}^q h_{ij}^{\alpha\beta(\mu-1)} \left( w_j^{\beta(\mu-1)} + v_j^{\beta(\mu-2)}dt \right) + C_{ij|\mu}^\alpha + N_{ij|\mu-1}^\alpha dt - \sum_{\beta=1}^q \left(\sum_{\gamma=1}^q D_{ijk|\mu}^{\alpha \gamma} h_{0ik}^{\gamma \beta} + \sum_{\gamma=1}^q h_{ij}^{\alpha \gamma(\mu-1)} h_{jk}^{\gamma \beta(\mu-1)}\right)\left(w_k^{\beta(\mu-1)} + v_k^{\beta(\mu-2)} dt \right)\\
&\equiv_\mu \sum_{\beta=1}^q h_{ij}^{\alpha\beta(\mu-1)} \left( w_j^{\beta(\mu-1)} + v_j^{\beta(\mu-2)}dt \right) + C_{ij|\mu}^\alpha + N_{ij|\mu-1}^\alpha dt  - \sum_{\gamma=1}^q D_{ijk|\mu}^{\alpha \gamma} w_{0i}^\gamma - \sum_{\gamma=1}^q h_{ij}^{\alpha \gamma(\mu-1)}\left( w_j^{\gamma(\mu-1)} + v_j^{\gamma(\mu-2)}dt - C_{jk|\mu}^\gamma - N_{ij|\mu-1}^{\gamma}dt        \right)\\
&\equiv_\mu C_{ij|\mu}^\alpha + N_{ij|\mu-1}^\alpha dt - \sum_{\gamma=1}^q D_{ijk|\mu}^{\alpha \gamma} w_{0i}^\gamma + \sum_{\gamma=1}^q h_{0ij}^{\alpha \gamma} C_{jk|\mu}^\gamma  +   \sum_{\gamma=1}^q h_{0ij}^{\alpha \gamma} N_{ij|\mu-1}^\gamma dt
\end{align*}}}
By considering the coefficient of $dt$, we have $(\ref{u35})$. We prove $(\ref{u37})$. In fact,
{\Small{\begin{align*}
- \psi_{i|\mu-1}^\alpha & \equiv_{\mu-1} \bar{\partial} v_i^{\alpha(\mu-2)} - \left[ \varphi^{\mu-1}, v_i^{\alpha(\mu-2)}  \right] - i_{\frac{\partial \varphi^{\mu-1}}{\partial t}}\left( w_i^{\alpha(\mu-1)}\right) \\
    &\equiv_{\mu-1} \bar{\partial}\left( N_{ij|\mu-1}^\alpha + \sum_{\beta=1}^q h_{ij}^{\alpha\beta(\mu-1)} v_j^{\beta(\mu-2)}  \right) -\left[\varphi^{\mu-1}, N_{ij|\mu-1}^\alpha + \sum_{\beta=1}^q h_{ij}^{\alpha\beta(\mu-1)} v_j^{\beta(\mu-2)}    \right] - i_{\frac{\partial \varphi^{\mu-1}}{\partial t}}\left( C_{ij|\mu}^\alpha + \sum_{\beta=1}^q h_{ij}^{\alpha\beta(\mu-1)} w_j^{\beta(\mu-1)} \right)\\
    &\equiv_{\mu-1} \bar{\partial} N_{ij|\mu-1}^\alpha + \sum_{\beta=1}^q h_{ij}^{\alpha\beta(\mu-1)}\left( \bar{\partial} v_j^{\beta(\mu-2)} -\left[  \varphi^{\mu-1}, v_j^{\beta(\mu-2)} \right] - i_{\frac{\partial \varphi^{\mu-1}}{\partial t}}\left( w_j^{\beta(\mu-1)} \right)  \right) + \sum_{\beta=1}^q \left( \bar{\partial} h_{ij}^{\alpha\beta(\mu-1)}- \left[ \varphi^{\mu-1}, h_{ij}^{\alpha\beta(\mu-1)}\right]   \right) v_j^{\beta(\mu-2)} \\
    &\equiv_{\mu-1} \bar{\partial} N_{ij|\mu-1}^\alpha - \sum_{\beta=1}^q h_{0ij}^{\alpha\beta} \psi_{j|\mu-1}^\beta - \sum_{\beta,\eta=1}^q B_{ij|\mu}^{\alpha\eta} h_{0ij}^{\eta\beta} w_{0j}^\beta \equiv_{\mu-1} \bar{\partial} N_{ij|\mu-1}^\alpha - \sum_{\beta=1}^q h_{0ij}^{\alpha\beta} \psi_{j|\mu-1}^\beta 
\end{align*}}}
We prove $(\ref{u39})$. In fact,
\begin{align*}
\bar{\partial} M_{i|\mu-1}^\alpha &\equiv_{\mu-1} \sum_{\beta=1}^q (-1)^{q-\beta}  \bar{\partial} v_i^{\beta(\mu-2)}  \wedge w_i^{1(\mu-1)} \wedge \cdots \wedge \widehat{w_i^{\beta(\mu-1)}} \wedge \cdots \wedge w_i^{q(\mu-1)} \wedge  \partial w_i^{\alpha(\mu-1)} \\
                            &+ \sum_{\beta=1}^q \sum_{\gamma=1}^{\beta-1} (-1)^{q-\beta} (-1)^{\gamma-1} v_i^{\beta(\mu-1)} w_i^{1(\mu-1)}\wedge \cdots \wedge \bar{\partial} w_i^{\gamma(\mu-1)} \wedge\cdots \wedge \widehat{w_i^{\beta(\mu-1)}} \wedge \cdots \wedge w_i^{q(\mu-1)} \wedge \partial w_i^{\alpha(\mu-1)} \\
                            &+ \sum_{\beta=1}^q \sum_{\gamma=\beta+1}^q (-1)^{q-\beta} (-1)^\gamma v_i^{\beta (\mu-1)} w_i^{1(\mu-1)}\wedge \cdots \wedge \widehat{w_i^{\beta(\mu-1)}} \wedge \cdots \wedge \bar{\partial} w_i^{\gamma(\mu-1)} \wedge \cdots \wedge w_i^{q(\mu-1)} \wedge \partial w_i^{\alpha(\mu-1)}\\
                            &+\sum_{\beta=1}^q (-1)^{q-\beta}(-1)^{q-2} v_i^{\beta(\mu-2)} w_i^{1(\mu-1)}\wedge \cdots \wedge \widehat{w_i^{\beta(\mu-1)}} \wedge \cdots \wedge w_i^{q(\mu-1)}\wedge \partial \left(\bar{\partial} w_i^{\alpha(\mu-1)} \right)\\
                            &+\sum_{\beta=1}^q (-1)^{\beta-1} w_i^{1(\mu-1)} \wedge \cdots \wedge \overbrace{ \bar{\partial} w_i^{\beta(\mu-1)} }^{\beta-\textnormal{th}} \wedge \cdots \wedge w_i^{q(\mu-1)} \wedge \left( \partial v_i^{\alpha(\mu-2) } - \frac{\partial w_i^{\alpha(\mu-1)}}{\partial t} \right) \\
                            &+ (-1)^q w_i^{1(\mu-1)} \wedge \cdots \wedge w_i^{q(\mu-1)} \wedge \left( - \partial \bar{\partial} v_i^{\alpha(\mu-2)} - \frac{\partial}{\partial t}\left(\bar{\partial} w_i^{\alpha(\mu-1)} \right)\right)
\end{align*}
Then we have
{\small{\begin{align*}
\bar{\partial} M_{i|\mu-1}^\alpha &\equiv_{\mu-1} \sum_{\beta=1}^q (-1)^{q-\beta} \left(\left[ \varphi^{\mu-1}, v_i^{\beta(\mu-2)}\right] + i_{\frac{\partial \varphi^{\mu-1}}{\partial t}}\left( w_i^{\beta(\mu-1)} \right) - \psi_{i|\mu-1}^\beta     \right) \wedge w_i^{1(\mu-1)} \wedge \cdots \wedge \widehat{w_i^{\beta(\mu-1)}} \wedge \cdots \wedge w_i^{q(\mu-1)} \wedge \partial w_i^{\alpha(\mu-1)} \\
                            &+ \sum_{\beta=1}^q \sum_{\gamma=1}^{\beta-1} (-1)^{q-\beta} (-1)^{\gamma-1} v_i^{\beta(\mu-1)} w_i^{1(\mu-1)}\wedge \cdots \wedge \left( \mathcal{L}_{\varphi^{\mu-1}}\left( w_i^{\gamma(\mu-1)} \right) - A_{i|\mu}^\gamma \right) \wedge\cdots \wedge \widehat{w_i^{\beta(\mu-1)}} \wedge \cdots \wedge w_i^{q(\mu-1)} \wedge \partial w_i^{\alpha(\mu-1)} \\
                            &+ \sum_{\beta=1}^q \sum_{\gamma=\beta+1}^q (-1)^{q-\beta} (-1)^\gamma v_i^{\beta (\mu-1)} w_i^{1(\mu-1)}\wedge \cdots \wedge \widehat{w_i^{\beta(\mu-1)}} \wedge \cdots \wedge \left( \mathcal{L}_{\varphi^{\mu-1}}\left( w_i^{\gamma(\mu-1)} \right) - A_{i|\mu}^\gamma \right)\wedge \cdots \wedge w_i^{q(\mu-1)} \wedge  \partial w_i^{\alpha(\mu-1)}\\
                            &+\sum_{\beta=1}^q (-1)^{q-\beta}(-1)^{q} v_i^{\beta(\mu-2)} w_i^{1(\mu-1)}\wedge \cdots \wedge \widehat{w_i^{\beta(\mu-1)}} \wedge \cdots \wedge w_i^{q(\mu-1)}\wedge  \partial \left( \mathcal{L}_{\varphi^{\mu-1}}\left( w_i^{\alpha(\mu-1)} \right) - A_{i|\mu}^\alpha  \right)\\
                            &+\sum_{\beta=1}^q (-1)^{\beta-1} w_i^{1(\mu-1)} \wedge \cdots \wedge \overbrace{ \left( \mathcal{L}_{\varphi^{\mu-1}}\left( w_i^{\beta(\mu-1)} \right) - A_{i|\mu}^\beta \right) }^{\beta-\textnormal{th}} \wedge \cdots \wedge w_i^{q(\mu-1)} \wedge \left( \partial v_i^{\alpha(\mu-2) } - \frac{\partial w_i^{\alpha(\mu-1)}}{\partial t} \right) \\
                            &+ (-1)^q w_i^{1(\mu-1)} \wedge \cdots \wedge w_i^{q(\mu-1)} \wedge \left( - \partial \left(\left[ \varphi^{\mu-1}, v_i^{\alpha(\mu-2)}\right] + i_{\frac{\partial \varphi^{\mu-1}}{\partial t}}\left( w_i^{\alpha(\mu-1)} \right) - \psi_{i|\mu-1}^\alpha     \right) - \frac{\partial}{\partial t}\left(\mathcal{L}_{\varphi^{\mu-1}}\left( w_i^{\alpha(\mu-1)} \right) - A_{i|\mu}^\alpha \right)\right)
\end{align*}}}
Then we have
{\small{\begin{align}\label{ue11}
&\bar{\partial} M_{i|\mu-1}^\alpha\\
&\equiv_{\mu-1}  \sum_{\beta=1}^q (-1)^{q} w_{0i}^1\wedge \cdots \wedge \overbrace{\psi_{i|\mu-1}^\beta}^{\beta-\textnormal{th}} \wedge \cdots w_{0i}^q \wedge dw_{0i}^\alpha +(-1)^q w_{0i}^1\wedge \cdots \wedge w_{0i}^q \wedge \partial \psi_{i|\mu-1}^\alpha +(-1)^q w_{0i}^1\wedge \cdots \wedge w_{0i}^q \wedge \frac{\partial A_{i|\mu}^\alpha}{\partial t} \notag\\
&+ \sum_{\beta=1}^q (-1)^{q-\beta} \left[ \varphi^{\mu-1}, v_i^{\beta(\mu-2)} \right] \wedge w_i^{1(\mu-1)} \wedge \cdots \wedge \widehat{w_i^{\beta(\mu-1)}} \wedge \cdots \wedge w_i^{q(\mu-1)} \wedge \partial w_i^{\alpha(\mu-1)} \notag \\
                            &+ \sum_{\beta=1}^q \sum_{\gamma=1}^{\beta-1} (-1)^{q-\beta} (-1)^{\gamma-1} v_i^{\beta(\mu-1)} w_i^{1(\mu-1)}\wedge \cdots \wedge \mathcal{L}_{\varphi^{\mu-1}}\left( w_i^{\gamma(\mu-1)} \right) \wedge\cdots \wedge \widehat{w_i^{\beta(\mu-1)}} \wedge \cdots \wedge w_i^{q(\mu-1)} \wedge \partial w_i^{\alpha(\mu-1)}\notag \\
                                                        &+ \sum_{\beta=1}^q \sum_{\gamma=\beta+1}^q (-1)^{q-\beta} (-1)^\gamma v_i^{\beta (\mu-1)} w_i^{1(\mu-1)}\wedge \cdots \wedge \widehat{w_i^{\beta(\mu-1)}} \wedge \cdots \wedge \mathcal{L}_{\varphi^{\mu-1}}\left( w_i^{\gamma(\mu-1)} \right) \wedge \cdots \wedge w_i^{q(\mu-1)} \wedge \partial w_i^{\alpha(\mu-1)} \notag\\
                                                                                    &+\sum_{\beta=1}^q (-1)^{q-\beta}(-1)^{q} v_i^{\beta(\mu-2)} w_i^{1(\mu-1)}\wedge \cdots \wedge \widehat{w_i^{\beta(\mu-1)}} \wedge \cdots \wedge w_i^{q(\mu-1)}\wedge \partial \left(\mathcal{L}_{\varphi^{\mu-1}}\left( w_i^{\alpha(\mu-1)} \right)   \right) \notag \\
        & + \sum_{\beta=1}^q (-1)^{\beta-1} w_i^{1(\mu-1)} \wedge \cdots \wedge \mathcal{L}_{\varphi^{\mu-1}}\left( w_i^{\beta(\mu-1)} \right) \wedge \cdots \wedge w_i^{q(\mu-1)}   \wedge \partial v_i^{\alpha(\mu-2)}- (-1)^q w_i^{1(\mu-1)}\wedge \cdots \wedge w_i^{q(\mu-1)}\wedge \partial \mathcal{L}_{\varphi^{\mu-1}}\left( v_i^{\alpha(\mu-2)} \right) \notag     \\  
        &- \sum_{\beta=1}^q (-1)^{\beta-1}w_i^{1(\mu-1)}\wedge \cdots \wedge \mathcal{L}_{\varphi^{\mu-1}}\left( w_i^{\beta(\mu-1)} \right) \wedge \cdots \wedge w_i^{q(\mu-1)}\wedge \frac{\partial w_i^{\alpha(\mu-1)}}{\partial t} -(-1)^q w_i^{1(\mu-1)} \wedge \cdots \wedge w_i^{q(\mu-1)}\wedge \mathcal{L}_{\varphi^{\mu-1}}\left( \frac{\partial w_i^{\alpha(\mu-1)}}{\partial t} \right)  \notag \\                                                                          
       &+\sum_{\beta=1}^q (-1)^{q-\beta} i_{\frac{\partial \varphi^{\mu-1}}{\partial t}}\left( w_i^{\beta(\mu-1)}\right) \wedge w_i^{1(\mu-1)} \wedge \cdots \wedge \widehat{w_i^{\beta(\mu-1)}} \wedge \cdots \wedge w_i^{q(\mu-1)} \wedge \partial w_i^{\alpha(\mu-1)} \notag \\
       &-(-1)^q w_i^{1(\mu-1)}\wedge \cdots w_i^{q(\mu-1)}\wedge \left(  \partial i_{\frac{\partial \varphi^{\mu-1}}{\partial t}}\left( w_i^{\alpha(\mu-1)} \right) + \mathcal{L}_{\frac{\partial \varphi^{\mu-1}}{\partial t}}\left(w_i^{\alpha(\mu-1)} \right) \right) \notag
\end{align}}}

On the other hand, we note that
{\small{\begin{align}\label{ue12}
&0\equiv_{\mu-1} \mathcal{L}_{\varphi^{\mu-1}}\left( M_{i|\mu-1}^\alpha \right) \\
&=\sum_{\beta=1}^q (-1)^{q-\beta} \mathcal{L}_{\varphi^{\mu-1}}\left(v_i^{\beta(\mu-2)} w_i^{1(\mu-1)}\wedge \cdots \wedge \widehat{w_i^{\beta(\mu-1)}} \wedge \cdots \wedge w_i^{q(\mu-1)} \wedge  \partial w_i^{\alpha(\mu-1)} \right)+ \mathcal{L}_{\varphi^{\mu-1}}\left(w_i^{1(\mu-1)}\wedge \cdots \wedge w_i^{q(\mu-1)}\wedge \partial v_i^{\alpha(\mu-2)} \right) \notag \\
&- \mathcal{L}_{\varphi^{\mu-1}}\left(w_i^{1(\mu-1)}\wedge \cdots \wedge w_i^{q(\mu-1)}\wedge \frac{\partial w_i^{\alpha(\mu-1)}}{\partial t} \right) \notag \\
&=\sum_{\beta=1}^q (-1)^{q-\beta} \left[ \varphi^{\mu-1}, v_i^{\beta(\mu-2)} \right]  \wedge w_i^{1(\mu-1)}\wedge \cdots \wedge \widehat{w_i^{\beta(\mu-1)}} \wedge \cdots \wedge w_i^{q(\mu-1)} \wedge \partial w_i^{\alpha(\mu-1)}  \notag \\
                            &+ \sum_{\beta=1}^q \sum_{\gamma=1}^{\beta-1} (-1)^{q-\beta} (-1)^{\gamma-1} v_i^{\beta(\mu-1)} w_i^{1(\mu-1)}\wedge \cdots \wedge \mathcal{L}_{\varphi^{\mu-1}}\left( w_i^{\gamma(\mu-1)} \right) \wedge\cdots \wedge \widehat{w_i^{\beta(\mu-1)}} \wedge \cdots \wedge w_i^{q(\mu-1)} \wedge \partial w_i^{\alpha(\mu-1)} \notag \\
                                                        &+ \sum_{\beta=1}^q \sum_{\gamma=\beta+1}^q (-1)^{q-\beta} (-1)^\gamma v_i^{\beta (\mu-1)} w_i^{1(\mu-1)}\wedge \cdots \wedge \widehat{w_i^{\beta(\mu-1)}} \wedge \cdots \wedge \mathcal{L}_{\varphi^{\mu-1}}\left( w_i^{\gamma(\mu-1)} \right) \wedge \cdots \wedge w_i^{q(\mu-1)} \wedge \partial w_i^{\alpha(\mu-1)} \notag \\
                                                                                    &+\sum_{\beta=1}^q (-1)^{q-\beta}(-1)^{q-2} v_i^{\beta(\mu-2)}\wedge \cdots \wedge \widehat{w_i^{\beta(\mu-1)}} \wedge \cdots \wedge w_i^{q(\mu-1)}\wedge \partial \left( \mathcal{L}_{\varphi^{\mu-1}}\left( w_i^{\alpha(\mu-1)} \right)   \right) \notag \\
&+  \sum_{\beta=1}^q (-1)^{\beta-1} w_i^{1(\mu-1)} \wedge \cdots \wedge \mathcal{L}_{\varphi^{\mu-1}}\left( w_i^{\beta(\mu-1)} \right) \wedge \cdots \wedge w_i^{q(\mu-1)} \wedge v_i^{\alpha(\mu-2)}-(-1)^q w_i^{1(\mu-1)} \wedge \cdots \wedge w_i^{q(\mu-1)}\wedge \partial \mathcal{L}_{\varphi^{\mu-1}}\left(v_i^{\alpha(\mu-2)} \right) \notag \\
&- \sum_{\beta=1}^q (-1)^{\beta-1} w_i^{1(\mu-1)} \wedge \cdots \wedge \mathcal{L}_{\varphi^{\mu-1}} \left( w_i^{\beta(\mu-1)} \right) \wedge \cdots \wedge w_i^{q(\mu-1)}\wedge \frac{\partial w_i^{\alpha(\mu-1)}}{\partial t} - (-1)^q w_i^{1(\mu-1)}\wedge \cdots \wedge w_i^{q(\mu-1)}\wedge \mathcal{L}_{\varphi^{\mu-1}}\left( \frac{\partial w_i^{\alpha(\mu-1)}}{\partial t} \right)             \notag                                                                       
\end{align}}}
and we also note that
\begin{align}\label{ue13}
&0\equiv_{\mu-1} i_{\frac{\partial \varphi^{\mu-1}}{\partial t}}\left( E_{i|\mu}^\alpha \right)= i_{\frac{\partial \varphi^{\mu-1}}{\partial t}}\left( w_i^{1(\mu-1)}\wedge \cdots \wedge w_i^{q(\mu-1)} \wedge \partial w_i^{\alpha(\mu-1)} \right)\\
&=\sum_{\beta=1}^q w_i^{1(\mu-1)} \wedge \cdots \wedge i_{\frac{\partial \varphi^{\mu-1}}{\partial t}} \left( w_i^{\beta(\mu-1)} \right) \wedge \cdots \wedge w_i^{q(\mu-1)} \wedge \partial w_i^{\alpha(\mu-1)} + w_i^{1(\mu-1)} \wedge \cdots \wedge w_i^{q(\mu-1)} \wedge i_{\frac{\partial \varphi^{\mu-1}}{\partial  t}} \left( \partial w_i^{\alpha(\mu-1)} \right)  \notag \\
&=\sum_{\beta=1}^q (-1)^{\beta-1} i_{\frac{\partial \varphi^{\mu-1}}{\partial t}}\left( w_i^{\beta(\mu-1)} \right) w_i^{1(\mu-1)} \wedge \cdots \wedge \widehat{w_i^{\beta(\mu-1)}} \wedge \cdots \wedge w_i^{q(\mu-1)} \wedge \partial w_i^{\alpha(\mu-1)} \notag \\
&+ w_i^{1(\mu-1)} \wedge \cdots \wedge w_i^{q(\mu-1)} \wedge \left( \mathcal{L}_{\frac{\partial \varphi^{\mu-1}}{\partial t} } \left( w_i^{\alpha(\mu-1)} \right) + \partial i_{\frac{\partial \varphi^{\mu-1}}{\partial t}}\left( w_i^{\alpha(\mu-1)} \right)  \right) \notag
\end{align}
Then from $(\ref{ue11}),(\ref{ue12})$ and $(\ref{ue13})$, we obtain
\begin{align*}
\bar{\partial} M_{i|\mu-1}^\alpha \equiv_{\mu-1} (-1)^q\left(  \sum_{\beta=1}^q w_{0i}^1\wedge \cdots \wedge \overbrace{\psi_{i|\mu-1}^\beta}^{\beta-\textnormal{th}} \wedge \cdots w_{0i}^q \wedge dw_{0i}^\alpha +  w_{0i}^1\wedge \cdots \wedge w_{0i}^q \wedge \partial \psi_{i|\mu-1}^\alpha  +  w_{0i}^1\wedge \cdots \wedge w_{0i}^q \wedge \frac{\partial A_{i|\mu}^\alpha}{\partial t}   \right)
\end{align*}

We prove $(\ref{u41})$. In fact, we note that
\begin{align}\label{ue14}
F_{i|\mu}^{\alpha}- \sum_{\beta=1}^q h_{0ij}^{\alpha\beta} \det\left( h_{0ij}^{\gamma\eta} \right) F_{j|\mu}^\beta & = E_{i|\mu}^\alpha + M_{i|\mu-1}^\alpha \wedge dt - \sum_{\beta=1}^q h_{ij}^{\alpha\beta} \det\left( h_{ij}^{\gamma\eta} \right) \left( E_{j|\mu}^\beta + M_{j|\mu-1}^\beta \wedge dt    \right)\\
                        &= E_{i|\mu}^\alpha - \sum_{\beta=1}^q h_{0ij}^{\alpha\beta} \det\left( h_{0ij}^{\gamma\eta} \right)  E_{j|\mu}^\beta +\left( M_{i|\mu-1}^\alpha - \sum_{\beta=1}^q h_{0ij}^{\alpha\beta} \det\left( h_{0ij}^{\gamma\eta} \right) M_{j|\mu-1}^\beta     \right) \wedge dt \notag
\end{align}
and
{\small{\begin{align*}
&F_{i|\mu}^{\alpha}- \sum_{\beta=1}^q h_{0ij}^{\alpha\beta} \det\left( h_{0ij}^{\gamma\eta} \right) F_{j|\mu}^\beta  \\
&\equiv_\mu \left(w_i^{1(\mu-1)}+ v_i^{1(\mu-2)}dt    \right)\wedge \cdots \wedge \left( w_i^{q(\mu-1)} + v_i^{q(\mu-2)}dt \right) \wedge \partial_{M\times B} \left( w_i^{\alpha(\mu-1)}+ v_i^{\alpha(\mu-2)}dt  \right) \\
& - \sum_{\beta=1}^q h_{ij}^{\alpha\beta(\mu-1)}\det\left( h_{ij}^{\gamma\eta(\mu-1)} \right)\left( w_j^{1(\mu-1)} + v_j^{1(\mu-1)}dt \right) \wedge \cdots \wedge \left( w_j^{q(\mu-1)} + v_j^{q(\mu-1)} dt \right) \wedge \partial_{M\times B}\left( w_j^{\beta(\mu-1)} + v_j^{\beta(\mu-1)}dt \right) \\
&\equiv_\mu \left( C_{ij|\mu}^1+ \sum_{\gamma_1=1}^q h_{ij}^{1\gamma_1(\mu-1)} w_j^{\gamma_1(\mu-1)} +\left( N_{ij|\mu-1}^1 + \sum_{\eta_1=1}^q h_{ij}^{1\eta_1(\mu-1)} v_j^{\eta_1 (\mu-2)} \right)dt     \right) \wedge \cdots \\
&\cdots  \wedge \left( C_{ij|\mu}^q+ \sum_{\gamma_q=1}^q h_{ij}^{q\gamma_q(\mu-1)} w_j^{\gamma_q(\mu-1)} +\left( N_{ij|\mu-1}^q + \sum_{\eta_q=1}^q h_{ij}^{q\eta_q(\mu-1)} v_j^{\eta_q (\mu-2)} \right)dt     \right)\\
&\wedge \partial_{M\times B} \left( C_{ij|\mu}^\alpha+ \sum_{\beta=1}^q h_{ij}^{\alpha\beta(\mu-1)} w_j^{\beta(\mu-1)} +\left( N_{ij|\mu-1}^\alpha + \sum_{\delta=1}^q h_{ij}^{\alpha\delta(\mu-1)} v_j^{\delta (\mu-2)} \right)dt     \right)\\
&- \sum_{\beta=1}^q h_{ij}^{\alpha\beta(\mu-1)}\det\left( h_{ij}^{\gamma\eta(\mu-1)} \right)\left( w_j^{1(\mu-1)} + v_j^{1(\mu-1)}dt \right) \wedge \cdots \wedge \left( w_j^{q(\mu-1)} + v_j^{q(\mu-1)} dt \right) \wedge \partial_{M\times B} \left( w_j^{\beta(\mu-1)} + v_j^{\beta(\mu-1)}dt \right)\\
&\equiv_\mu \sum_{\beta=1}^q w_{0i}^1 \wedge \cdots \wedge \overbrace{C_{ij|\mu}^\beta}^{\beta-\textnormal{th}} \wedge \cdots \wedge w_{0i}^q     \wedge dw_{0i}^\alpha + w_{0i}^1\wedge \cdots \wedge w_{0i}^q \wedge \partial C_{ij|\mu}^\alpha \\
&+ \sum_{\beta=1}^q w_{0i}^1\wedge \cdots \wedge \overbrace{ \left( N_{ij|\mu-1}^\beta\right)dt}^{\beta-\textnormal{th}} \wedge \cdots \wedge w_{0i}^q \wedge dw_{0i}^\alpha + w_{0i}^1\wedge \cdots \wedge w_{0i}^q \wedge \partial N_{ij|\mu-1}^\alpha \wedge dt  - w_{0i}^1 \wedge \cdots \wedge w_{0i}^q \wedge \frac{\partial C_{ij|\mu}^\alpha}{\partial t} \wedge dt\\
&\equiv_\mu \sum_{\beta=1}^q w_{0i}^1 \wedge \cdots \wedge \overbrace{C_{ij|\mu}^\beta}^{\beta-\textnormal{th}} \wedge \cdots \wedge w_{0i}^q    \wedge d w_{0i}^\alpha + w_{0i}^1\wedge \cdots \wedge w_{0i}^q \wedge \partial C_{ij|\mu}^\alpha \\
&+ \left( \sum_{\beta=1}^q (-1)^{q-\beta} N_{ij|\mu-1}^\beta w_{0i}^1\wedge \cdots \wedge  \widehat{w_{0i}^\beta} \wedge \cdots \wedge w_{0i}^q \wedge dw_{0i}^\alpha + w_{0i}^1\wedge \cdots \wedge w_{0i}^q \wedge \partial N_{ij|\mu-1}^\alpha - w_{0i}^1 \wedge \cdots \wedge w_{0i}^q \wedge \frac{\partial C_{ij|\mu}^\alpha}{\partial t} \right) \wedge dt
\end{align*}}}
Hence from $(\ref{ue14})$, we have $(\ref{u41})_{\mu-1}$. It remains to show $(\ref{u43})$. In fact, we note that
{\Small{\begin{align*}
&E_{i|\mu}^\alpha + M_{i|\mu-1}^\alpha \wedge dt \\
 &\equiv_\mu \left(w_i^{1(\mu-1)}+ v_i^{1(\mu-2)}dt   \right) \wedge \cdots \wedge \left( w_i^{q(\mu-1)} + v_i^{q(\mu-2)}dt  \right) \wedge \partial_{M\times B} \left( w_i^{\alpha(\mu-1)}+ v_i^{\alpha(\mu-1)}dt  \right)\\ 
    &\equiv_\mu \left(w_i^{1(\mu-1)}+ v_i^{1(\mu-2)}dt   \right) \wedge \cdots \wedge \left( w_i^{q(\mu-1)} + v_i^{q(\mu-2)}dt  \right) \wedge \left( \sum_{\gamma=1}^q\left( a_{i_x}^{\alpha\gamma(\mu-1)} + b_{i_x}^{\alpha\gamma(\mu-2)}dt \right) \wedge \left(w_i^{\gamma(\mu-1)} + v_i^{\gamma(\mu-2)}dt \right)+ L_{i_x|\mu}^\alpha + P_{i_x|\mu-1}^\alpha \wedge dt   \right) \\
    &\equiv_\mu w_{0i}^1\wedge \cdots \wedge w_{0i}^q \wedge L_{i_x|\mu}^\alpha + w_{0i}^1 \wedge \cdots \wedge w_{0i}^q \wedge P_{i_x|\mu-1}^\alpha \wedge dt
\end{align*}}} 
This completes Lemma \ref{ue15}.

\end{proof}

Our purpose is to construct $\varphi^\mu=\varphi^{\mu-1}+ \varphi_\mu, h_{ij}^{\alpha\beta\mu}=h_{ij}^{\alpha\beta(\mu-1)}+ h_{ij|\mu}^{\alpha\beta}, w_i^{\alpha\mu}=w_i^{\alpha(\mu-1)}+ w_{i|\mu}^\alpha, v_i^{\alpha(\mu-1)}=v_i^{\alpha(\mu-2)}+ v_{i|\mu-1}^\alpha , a_{i_x}^{\alpha\beta\mu}=a_{i_x}^{\alpha\beta(\mu-1)}+ a_{i_x|\mu}^{\alpha\beta}$, and $b_{i_x}^{\alpha\beta(\mu-1)}= b_{i_x}^{\alpha\beta(\mu-2)}+ b_{i_x|\mu-1}^{\alpha\beta}$ satisfying $(\ref{u14})_\mu - (\ref{u19})_\mu$.

\begin{lemma}\label{ue16}
$(\ref{u14})_\mu- (\ref{u19})_\mu$ are equivalent to the following equalities, respectively
\begin{align}
\bar{\partial} \varphi_\mu&= - \xi_\mu \label{u45}
\end{align}
\begin{align}
A_{i|\mu}^\alpha&= \bar{\partial} w_{i|\mu}^\alpha -\mathcal{L}_{\varphi_\mu}\left(w_{0i}^\alpha\right) \label{u46}
\end{align}
\begin{align}
\psi_{i|\mu-1}^\alpha &= \bar{\partial} v_{i|\mu-1}^\alpha - i_{\frac{\partial \varphi_\mu}{\partial t}} \left(w_{0i}^\alpha \right) \label{u47} 
\end{align}
\begin{align}
\sum_{\eta=1}^q B_{ij|\mu}^{\alpha\eta} h_{0ij}^{\eta \beta} &= \bar{\partial} h_{ij|\mu}^{\alpha\beta} -\left[ \varphi_\mu, h_{0ij}^{\alpha\beta}  \right] \label{u48} 
\end{align}
\begin{align}
C_{ij|\mu}^\alpha &= - w_{i|\mu}^\alpha + \sum_{\beta=1}^q h_{0ij}^{\alpha\beta} w_{j|\mu}^\beta + \sum_{\beta=1}^q h_{ij|\mu}^{\alpha\beta} w_{0j}^\beta \label{u50}
\end{align}
\begin{align}
N_{ij|\mu-1}^\alpha &= - v_{i|\mu-1}^\alpha + \sum_{\beta=1}^q h_{0ij}^{\alpha\beta} v_{j|\mu-1}^\beta \label{u51}
\end{align}
\begin{align}
- E_{i|\mu}^\alpha &=  \sum_{\beta=1}^q w_{0i}^1 \wedge \cdots \wedge \overbrace{ w_{i|\mu}^\beta }^{\beta-\textnormal{th}} \wedge \cdots \wedge w_{0i}^q  \wedge dw_{0i}^\alpha + w_{0i}^1 \wedge \cdots \wedge w_{0i}^q \wedge dw_{i|\mu}^\alpha \label{u52}
\end{align}
\begin{align}
- M_{i|\mu-1}^\alpha &=  \sum_{\beta=1}^q (-1)^{q-\beta} v_{i|\mu-1}^\beta w_{0i}^1 \wedge \cdots \wedge  \widehat{w_{0i}^\beta} \wedge \cdots \wedge w_{0i}^q \wedge d w_{0i}^\alpha + w_{0i}^1 \wedge \cdots \wedge w_{0i}^q \wedge \partial v_{i|\mu-1}^\alpha - w_{0i}^1 \wedge \cdots \wedge w_{0i}^1 \wedge \frac{\partial w_{i|\mu}^\alpha }{\partial t}    \label{u53}
\end{align}
\begin{align}
L_{i_x|\mu}^\alpha & = - \partial w_{i|\mu}^\alpha + \sum_{\beta=1}^q a_{0i_x}^{\alpha\beta} \wedge w_{i|\mu}^\beta + \sum_{\beta=1}^q a_{i_x|\mu}^{\alpha\beta} \wedge w_{0i}^\beta        \label{u54}\\
P_{i_x|\mu-1}^\alpha &= \frac{\partial w_{i|\mu}^\alpha}{\partial t}  - \partial v_{i|\mu-1}^\alpha + \sum_{\beta=1}^q a_{0i_x}^{\alpha\beta} v_{i|\mu-1}^\beta - \sum_{\beta=1}^q b_{i_x|\mu-1}^{\alpha\beta} w_{0i}^\beta             \label{u55}
\end{align}
\end{lemma}

\begin{proof}
We have shown $(\ref{u45}),(\ref{u46}),(\ref{u48}),(\ref{u50}),(\ref{u52})$ and $(\ref{u54})$ in Part I. We prove $(\ref{u47})$. In fact,
\begin{align*}
&\bar{\partial}\left( v_i^{\alpha(\mu-2)} + v_{i|\mu-1}^\alpha \right) -\left[ \varphi^{\mu-1}+ \varphi_\mu,  v_i^{\alpha(\mu-2)} + v_{i|\mu-1}^\alpha \right] - i_{\frac{\partial (\varphi^{\mu-1}+ \varphi_\mu ) }{\partial t} } \left(  w_i^{\alpha(\mu-1)} + w_{i|\mu}^\alpha \right)     \equiv_{\mu-1} 0 \\
&\iff  - \psi_{i|\mu-1}^\alpha + \bar{\partial} v_{i|\mu-1}^\alpha - i_{\frac{\partial \varphi_\mu}{\partial t}} \left(w_{0i}^\alpha\right) \equiv_{\mu-1} 0
\end{align*}
We prove $(\ref{u51})$. In fact,
\begin{align*}
&w_i^{\alpha(\mu-1)}+ w_{i|\mu}^{\alpha}+ \left( v_i^{\alpha(\mu-2)}+ v_{i|\mu-1}^{\alpha} \right)dt \equiv_\mu \sum_{\beta=1}^q \left( h_{ij}^{\alpha\beta(\mu-1)}+ h_{ij|\mu}^{\alpha\beta}\right) \left( w_j^{\beta(\mu-1)} + w_{j|\mu}^\beta +\left(v_j^{\beta(\mu-2)}+ v_{j|\mu-1}^\beta \right) dt   \right)\\
&\iff C_{ij|\mu}^\alpha= -w_{i|\mu}^\alpha + \sum_{\beta=1}^q h_{0ij}^{\alpha\beta} w_{j|\mu}^\beta + \sum_{\beta=1}^q h_{ij|\mu}^{\alpha\beta} w_{0j}^\beta,\,\,\,\,\,\,\,\,\,\, N_{ij|\mu-1}^\alpha= - v_{i|\mu-1}^\alpha + \sum_{\beta=1}^q h_{0ij}^{\alpha\beta} v_{j|\mu-1}^\beta
\end{align*}

We prove $(\ref{u53})$. In fact,
{\Small{\begin{align*}
0&\equiv_\mu \left( w_i^{1(\mu-1)} + w_{i|\mu}^1+ \left( v_i^{1(\mu-2)} + v_{i|\mu-1}^1 \right) dt  \right) \wedge \cdots \wedge \left( w_i^{q(\mu-1)} + w_{i|\mu}^q + \left( v_i^{q(\mu-2)} + v_{i|\mu-1}^q\right)dt  \right) \\
&\,\,\,\,\,\,\,\,\, \wedge \partial_{M\times B} \left( w_i^{\alpha(\mu-1)} + w_{i|\mu}^\alpha + \left( v_i^{\alpha(\mu-2)} + v_{i|\mu-1}^\alpha \right) dt \right)\\
&\equiv_\mu \left( w_i^{1(\mu-1)} + w_{i|\mu}^1+ \left( v_i^{1(\mu-2)} + v_{i|\mu-1}^1 \right) dt  \right) \wedge \cdots \wedge \left( w_i^{q(\mu-1)} + w_{i|\mu}^q + \left( v_i^{q(\mu-2)} + v_{i|\mu-1}^q\right)dt  \right)\wedge \partial_{M\times B} \left( w_i^{\alpha(\mu-1)} + v_i^{\alpha(\mu-2)} dt \right)\\
&+ \left( w_i^{1(\mu-1)} + w_{i|\mu}^1+ \left( v_i^{1(\mu-2)} + v_{i|\mu-1}^1 \right) dt  \right) \wedge \cdots \wedge \left( w_i^{q(\mu-1)} + w_{i|\mu}^q + \left( v_i^{q(\mu-2)} + v_{i|\mu-1}^q\right)dt  \right)\wedge \left( \partial w_{i|\mu}^\alpha - \frac{\partial w_{i|\mu}^\alpha}{\partial t}\wedge dt + \partial v_{i|\mu-1}^\alpha \wedge dt \right)\\
 &\equiv_\mu E_{i|\mu}^\alpha + M_{i|\mu-1}^\alpha \wedge dt + \sum_{\beta=1}^q w_{0i}^1 \wedge \cdots \wedge \overbrace{w_{i|\mu}^\beta}^{\beta-\textnormal{th}} \wedge \cdots \wedge w_i^q \wedge dw_{0i}^\alpha + w_{0i}^1 \wedge \cdots \wedge w_{0i}^q \wedge \partial w_{i|\mu}^\alpha  \\
 &+\sum_{\beta=1}^q w_{0i}^1 \wedge \cdots \wedge \overbrace{ v_{i|\mu-1}^\beta dt }^{\beta-\textnormal{th}} \wedge \cdots \wedge w_{0i}^q \wedge dw_{0i}^\alpha + w_{0i}^1 \wedge \cdots \wedge w_{0i}^q \wedge \left( -\frac{\partial w_{i|\mu}^\alpha}{\partial t} +  \partial v_{i|\mu -1 }^\alpha   \right),
\end{align*}}}
so that we have
\begin{align*}
- M_{i|\mu-1}^\alpha = \sum_{\beta=1}^q (-1)^{q-\beta} v_{i|\mu-1}^\beta w_{0i}^1 \wedge \cdots \wedge  \widehat{w_i^\beta} \wedge \cdots \wedge w_{0i}^q \wedge d w_{0i}^\alpha + w_{0i}^1 \wedge \cdots \wedge w_{0i}^q \wedge \partial v_{i|\mu-1}^\alpha - w_{0i}^1 \wedge \cdots \wedge w_{0i}^1 \wedge \frac{\partial w_{i|\mu}^\alpha }{\partial t}
\end{align*}

We prove $(\ref{u55})$. In fact,
{\Small{\begin{align*}
&\partial_{M\times B} \left( w_i^{\alpha(\mu-1)}+ w_{i|\mu}^{\alpha} + \left( v_i^{\alpha(\mu-2)} + v_{i|\mu-1}^\alpha \right)dt     \right) \\
&\equiv_\mu \sum_{\beta=1}^q \left(a_{i_x}^{\alpha\beta(\mu-1)}+ a_{i_x|\mu}^{\alpha\beta} +\left(b_{i_x}^{\alpha\beta(\mu-2)}+ b_{i_x|\mu-1}^{\alpha\beta} \right) dt \right) \wedge \left( w_i^{\beta(\mu-1)} + w_{i|\mu}^\beta + \left( v_i^{\beta(\mu-2)} + v_{i|\mu-1}^\beta \right) dt \right) \\
&\iff  \partial_{M\times B} \left(w_i^{\alpha(\mu-1)}+ v_i^{\alpha(\mu-2)}dt \right)+ \partial w_{i|\mu}^{\alpha} - \frac{\partial w_{i|\mu}^\alpha}{\partial t} \wedge dt + \partial v_{i|\mu-1}^\alpha \wedge dt\\
&\equiv_\mu \sum_{\beta=1}^q \left( a_{i_x}^{\alpha\beta(\mu-1)}+ b_{i_x}^{\alpha\beta(\mu-2)}dt\right) \wedge \left( w_i^{\beta(\mu-1)}+ v_i^{\beta(\mu-2)}dt  \right) + \sum_{\beta=1}^q  a_{0 i_x}^{\alpha\beta}\wedge w_{i|\mu}^\beta  +  \sum_{\beta=1}^q a_{i_x|\mu}^{\alpha\beta} \wedge w_{0i}^\beta + \sum_{\beta=1}^q a_{0i_x}^{\alpha\beta} v_{i|\mu-1}^\beta \wedge dt - \sum_{\beta=1}^q b_{i_x|\mu-1}^{\alpha\beta} w_{0i}^\beta \wedge dt \\
&\iff L_{i_x|\mu}^\alpha + P_{i_x|\mu-1}^\alpha \wedge dt = - \partial w_{i|\mu}^\alpha + \sum_{\beta=1}^q a_{0i_x}^{\alpha\beta}\wedge w_{i|\mu}^\beta + \sum_{\beta=1}^q a_{i_x|\mu}^{\alpha\beta} \wedge w_{0i}^\beta + \left(\frac{\partial w_{i|\mu}^\alpha}{\partial t}  - \partial v_{i|\mu-1}^\alpha + \sum_{\beta=1}^q a_{0 i_x}^{\alpha\beta} v_{i|\mu-1}^\beta - \sum_{\beta=1}^q b_{i_x|\mu-1}^{\alpha\beta} w_{0i}^\beta      \right)  \wedge dt,
\end{align*}}}
so that we have $(\ref{u55})$. This completes Lemma \ref{ue16}.

\end{proof}

\begin{lemma}\label{ue23}
Under the hypothesis $\mathbb{H}^2 \left(M, \tilde{\mathcal{N}}_{\mathcal{F}_0}^{*\bullet} \right)=0$, we can find $\varphi_\mu , h_{ij|\mu}^{\alpha\beta}, w_{i|\mu}^\alpha, v_{i|\mu-1}^\alpha, a_{i_x|\mu}^{\alpha\beta}$, and $ b_{i_x|\mu-1}^{\alpha\beta}$ satisfying $(\ref{u45})- (\ref{u55})$.
\end{lemma}

\begin{proof}
We define $D_{ijk|\mu}\in \Gamma\left( U_{ijk}, \mathcal{A}^{0,0}\left( \mathscr{H}om_{\mathcal{O}_M}\left( \mathcal{N}_{\mathcal{F}_0}^*, \mathcal{N}_{\mathcal{F}_0}^* \right) \right) \right)$ by
\begin{align*}
D_{ijk|\mu}:\Gamma\left( U_{ijk},  \mathcal{N}_{\mathcal{F}_0}^* \right) &\to \Gamma\left( U_{ijk}, \mathcal{A}^{0,0}\left( \mathcal{N}_{\mathcal{F}_0}^* \right) \right)\\
     w_{0i}^\alpha &\mapsto \sum_{\xi=1}^q D_{ijk|\mu}^{\alpha \xi} w_{0i}^\xi
\end{align*}
and linearly extends to $\Gamma\left( U_{ijk}, \mathcal{N}_{\mathcal{F}_0}^* \right)$. Then by $(\ref{u30})$, as in Part I, we can find $\left\{D_{ij|\mu} \right\}\in C^1\left( \mathcal{U}, \mathcal{A}^{0,0}\left(\mathscr{H}om_{\mathcal{O}_M}\left( \mathcal{N}_{\mathcal{F}_0}^*, \mathcal{N}_{\mathcal{F}_0}^* \right)   \right) \right)$
with $D_{ij|\mu}\in \Gamma\left( U_{ij}, \mathcal{A}^{0,0}\left( \mathscr{H}om_{\mathcal{O}_M}\left( \mathcal{N}_{\mathcal{F}_0}^*, \mathcal{N}_{\mathcal{F}_0}^* \right) \right) \right)$ defined by
\begin{align*}
D_{ij|\mu}:\Gamma\left( U_{ij} , \mathcal{N}_{\mathcal{F}_0}^* \right) &\to \Gamma\left( U_{ij} , \mathcal{A}^{0,0}\left( \mathcal{N}_{\mathcal{F}_0}^* \right) \right) \\
 w_{0i}^\alpha &\mapsto \sum_{\xi=1}^p D_{ij|\mu}^{\alpha \xi} w_{0i}^\xi
\end{align*}
such that
\begin{align}\label{ne11}
 \sum_{\beta,\eta=1}^q h_{0ij}^{\alpha\beta} D_{jk|\mu}^{\beta \eta} h_{0ji}^{\eta\xi}-D_{ik|\mu}^{\alpha \xi}+ D_{ij|\mu}^{\alpha \xi}=D_{ijk|\mu}^{\alpha\xi}
\end{align}
We define $B_{ij|\mu}\in \Gamma\left(U_{ij}, \mathcal{A}^{0,1}\left(\mathscr{H}om_{\mathcal{O}_M}\left( \mathcal{N}_{\mathcal{F}_0}^* , \mathcal{N}_{\mathcal{F}_0}^* \right) \right) \right)$by
\begin{align*}
B_{ij|\mu}:\Gamma\left( U_{ij}, \mathcal{N}_{\mathcal{F}_0}^* \right) &\to \Gamma\left(U_{ij}, \mathcal{A}^{0,1}\left( \mathcal{N}_{\mathcal{F}_0}^* \right) \right) \\
 w_{0i}^\alpha &\mapsto \sum_{\xi=1}^q B_{ij|\mu}^{\alpha \xi} w_{0i}^\xi
\end{align*}
Then from $(\ref{u33})$ and $(\ref{ne11})$, we can find $\left\{ B_{i|\mu} \right\}\in C^0\left(\mathcal{U}, \mathcal{A}^{0,1}\left( \mathscr{H}om_{\mathcal{O}_M}\left( \mathcal{N}_{\mathcal{F}_0}^*, \mathcal{N}_{\mathcal{F}_0}^* \right) \right) \right)$ with $B_{i|\mu}\in \Gamma\left(U_i, \mathcal{A}^{0,1}\left( \mathscr{H}om_{\mathcal{O}_M}\left(\mathcal{N}_{\mathcal{F}_0}^*, \mathcal{N}_{\mathcal{F}_0}^* \right)\right) \right)$ defined by
\begin{align*}
B_{i|\mu} :\Gamma\left( U_i, \mathcal{N}_{\mathcal{F}_0}^* \right) &\to \Gamma\left( U_i, \mathcal{A}^{0,1}\left( \mathcal{N}_{\mathcal{F}_0}^* \right) \right)\\
 w_{0i}^\alpha &\mapsto \sum_{\xi=1}^q B_{i|\mu}^{\alpha \xi} w_{0i}^\xi
\end{align*}
such that
\begin{align}
\sum_{\eta,\gamma=1}^q h_{0ij}^{\alpha \eta} B_{j|\mu}^{\eta \gamma}h_{0ji}^{\gamma \xi} - B_{i|\mu}^{\alpha \xi}=B_{ij|\mu}^{\alpha \xi}-\bar{\partial} D_{ij|\mu}^{\alpha \xi} \label{n50}
\end{align}

We define $C_{ij|\mu}\in \Gamma\left( U_{ij}, \mathcal{A}^{0,0}\left( \mathscr{H}om_{\mathcal{O}_M}\left( \mathcal{N}_{\mathcal{F}_0}^*, \Omega_M^1 \right)\right) \right)$ by
\begin{align*}
C_{ij|\mu}: \Gamma\left( U_{ij}, \mathcal{N}_{\mathcal{F}_0}^* \right) &\to \Gamma\left( U_{ij} , \mathcal{N}_{\mathcal{F}_0}^* \right) \\
   w_{0i}^\alpha  &\mapsto  C_{ij|\mu}^\alpha
\end{align*}
Then from $(\ref{u34})$ and $(\ref{ne11})$, we have
\begin{align*}
\sum_{\gamma=1}^q h_{0ij}^{\alpha\gamma}\left( C_{jk|\mu}^\gamma - \sum_{\beta=1}^q D_{jk|\mu}^{\gamma \beta} w_{0j}^\beta \right)-\left( C_{ik|\mu}^\alpha -\sum_{\beta=1}^q D_{ik|\mu}^{\alpha\beta} w_{0i}^\beta \right)+\left( C_{ij|\mu}^\alpha -\sum_{\beta=1}^q D_{ij|\mu}^{\alpha\beta} w_{0i}^\beta \right)=0
\end{align*}
This implies that we can find $\left\{ C_{i|\mu} \right\}\in C^0\left( \mathcal{U}, \mathcal{A}^{0,0}\left( \mathscr{H}om_{\mathcal{O}_M}\left( \mathcal{N}_{\mathcal{F}_0}^*, \Omega_M^1 \right) \right) \right)$ with $C_{i|\mu}\in \Gamma\left(U_i, \mathcal{A}^{0,0}\left( \mathscr{H}om_{\mathcal{O}_M} \left( \mathcal{N}_{\mathcal{F}_0}^*, \Omega_M^1   \right) \right)  \right)$ defined by
\begin{align*}
C_{i|\mu} : \Gamma\left( U_i, \mathcal{N}_{\mathcal{F}_0}^* \right) &\to \Gamma\left( U_i, \mathcal{A}^{0,0}\left( \mathcal{N}_{\mathcal{F}_0}^* \right) \right) \\
 w_{0i}^\alpha &\mapsto \sum_{\xi=1}^q C_{i|\mu}^{\alpha \xi} w_{0i}^\xi
\end{align*}
such that
\begin{align}
\sum_{\beta=1}^q h_{0ij}^{\alpha\beta} C_{j|\mu}^\beta - C_{i|\mu}^\alpha= C_{ij|\mu}^\alpha - \sum_{\beta=1}^q  D_{ij|\mu}^{\alpha\beta} w_{0i}^\beta \label{n52}
\end{align}
Then from $(\ref{u38})$ and $(\ref{n52})$, as in Part I, we can show
\begin{align*}
&E_{i|\mu}^\alpha +\sum_{\beta=1}^q w_{0i}^1\wedge \cdots \wedge C_{i|\mu}^\beta \wedge w_{0i}^q \wedge dw_{0i}^\alpha + w_{0i}^1 \wedge \cdots \wedge w_{0i}^q \wedge \partial C_{i|\mu}^\alpha\\
&=\sum_{\beta=1}^q h_{0ij}^{\alpha\beta}\det\left(h_{0ij}\right)\left( E_{j|\mu}^\beta + \sum_{\gamma=1}^q w_{0j}^1 \wedge \cdots \wedge C_{j|\mu}^\gamma\wedge \cdots \wedge w_{0j}^q \wedge dw_{0j}^\beta + w_{0j}^1 \wedge \cdots \wedge w_{0j}^q \wedge \partial C_{j|\mu}^\beta \right)
\end{align*}
This implies that from $(\ref{u42})$
\begin{align*}
I_\mu:=\left\{ I_{i|\mu}\right\} \in A^{0,0}\left(M-S, \mathscr{H}om_{\mathcal{O}_M}\left( \mathcal{N}_{\mathcal{F}_0}^*, \tilde{\mathcal{S}}^2 \right) \right)
\end{align*}
where
\begin{align} \label{ne15}
I_{i|\mu}\left(w_{0i}^\alpha \right)= I_{i|\mu}^\alpha:= E_{i|\mu}^\alpha  +\sum_{\beta=1}^q w_{0i}^1\wedge \cdots \wedge C_{i|\mu}^\beta \wedge w_{0i}^q \wedge dw_{0i}^\alpha + w_{0i}^1 \wedge \cdots \wedge w_{0i}^q \wedge \partial C_{i|\mu}^\alpha\,\,\,\,\,\,\textnormal{on}\,\,\,U_i-S
\end{align}
From $(\ref{u36})$ and $(\ref{n52})$ and $(\ref{n50})$, we can show as in Part I
\begin{align}\label{ne13}
\sum_{\beta=1}^q h_{0ij}^{\alpha\beta}\left( A_{j|\mu}^\beta- \bar{\partial} C_{j|\mu}^\beta +\sum_{\gamma=1}^q B_{j|\mu}^{\beta\gamma} w_{0j}^\gamma \right)-\left( A_{i|\mu}^\alpha - \bar{\partial} C_{i|\mu}^\alpha + \sum_{\gamma=1}^q B_{i|\mu}^{\alpha\gamma} w_{0i}^\gamma \right)=0
\end{align}
We define $\phi_{i|\mu}\in \Gamma\left(U_i, \mathcal{A}^{0,1}\left(\mathscr{H}om_{\mathcal{O}_M}\left( \mathcal{N}_{\mathcal{F}_0}, \Omega_M^1 \right) \right)  \right)$ by
\begin{align}
\phi_{i|\mu}: \Gamma\left( U_i, \mathcal{N}_{\mathcal{F}_0}^* \right) &\to \Gamma\left( U_i, \mathcal{A}^{0,1}\left( \Omega_M^1 \right) \right) \notag \\
  w_{0i}^\alpha &\mapsto \phi_{i|\mu}^\alpha:=  A_{i|\mu}^\alpha- \bar{\partial} C_{i|\mu}^\alpha + \sum_{\gamma=1}^q B_{i|\mu}^{\alpha \gamma} w_{0i}^\gamma \label{ne14}
\end{align}

Then $(\ref{ne13})$ implies that 
\begin{align*}
\phi_\mu:=\left\{ \phi_{i|\mu} \right\} \in A^{0,1}(X, \mathscr{H}om_{\mathcal{O}_X}\left( \mathcal{N}_{\mathcal{F}_0}^*, \Omega_X^1 \right))
\end{align*}
We take $\bar{\phi}_\mu$ to be the image of $\phi_\mu$ in $\frac{A^{0,1}\left( M, \mathscr{H}om_{\mathcal{O}_M}\left( \mathcal{N}_{\mathcal{F}_0}^*, \Omega_M^1 \right) \right)}{A^{0,1}\left( M, \mathscr{H}om_{\mathcal{O}_M}\left( \mathcal{N}_{\mathcal{F}_0}^*, \mathcal{N}_{\mathcal{F}_0}^* \right) \right)}$.

On the other hand, we define $N_{ij|\mu-1}\in \Gamma\left( U_{ij}, \mathcal{A}^{0,0}\left( \mathscr{H}om_{\mathcal{O}_M} \left( \mathcal{N}_{\mathcal{F}_0}^*, \mathcal{O}_M \right) \right)\right)$
\begin{align*}
N_{ij|\mu-1} : \Gamma\left( U_{ij}, \mathcal{N}_{\mathcal{F}_0}^* \right) &\to \Gamma\left( U_{ij}, \mathcal{A}^{0,0} \right) \\
 w_{0i}^\alpha &\mapsto N_{ij|\mu-1}^\alpha
\end{align*}
From  $(\ref{u35})$, we can find $\{ N_{i|\mu-1} \} \in C^0 \left( \mathcal{U},  \mathscr{A}^{0,0}\left( \mathscr{H}om_{\mathcal{O}_M} \left( \mathcal{N}_{\mathcal{F}_0}^*, \mathcal{O}_M \right) \right) \right)$ with $N_{i|\mu}\in \Gamma\left( U_i, \mathscr{H}om_{\mathcal{O}_M}\left(   \mathcal{N}_{\mathcal{F}_0}^*, \mathcal{O}_M  \right) \right)$ defined by
\begin{align*}
N_{i|\mu}:\Gamma \left( U_i, \mathcal{N}_{\mathcal{F}_0}^* \right) &\to \Gamma\left( U_i, \mathcal{A}^{0,0} \right)\\
  w_{0i}^\alpha &\mapsto N_{i|\mu}^\alpha
\end{align*}
 such that
\begin{align} \label{ue17}
\sum_{\beta=1}^q h_{0ij}^{\alpha\beta} N_{j|\mu-1}^\beta - N_{i|\mu-1}^\alpha = N_{ij|\mu-1}^\alpha
\end{align}
Then from $(\ref{ue17})$ and $(\ref{u37})$ we have
\begin{align}\label{ue18}
\sum_{\beta=1}^q h_{0ij}^{\alpha\beta} \left( \psi_{j|\mu-1}^\beta- \bar{\partial} N_{j|\mu-1}^\beta \right) - \left( \psi_{i|\mu-1}^\alpha - \bar{\partial} N_{i|\mu-1}^\alpha  \right) = 0
\end{align}
We define $\sigma_{i|\mu-1}\in \Gamma\left( U_i,  \mathcal{A}^{0,1}\left( \mathscr{H}om_{\mathcal{O}_M} \left( \mathcal{N}_{\mathcal{F}_0}^*, \mathcal{A}^{0,1} \right) \right) \right)$ by
\begin{align*}
\sigma_{i|\mu-1}:\Gamma\left( U_i, \mathcal{N}_{\mathcal{F}_0}^*  \right) & \to \Gamma\left( U_i, \mathcal{A}^{0,1} \right) \\
 w_{0i}^\alpha & \mapsto \sigma_{i|\mu-1}^\alpha:= \psi_{i|\mu-1}^\alpha - \bar{\partial} N_{i|\mu-1}^\alpha
\end{align*}
Then $(\ref{ue18})$ implies that 
\begin{align*}
\sigma_{\mu-1}:= \left\{ \sigma_{i|\mu-1} \right\} \in A^{0,1}\left( \mathscr{H}om_{\mathcal{O}_M} \left(  \mathcal{N}_{\mathcal{F}_0}^* , \mathcal{A}^{0,1}   \right) \right)
\end{align*}
We note that from $(\ref{u41})$ and $(\ref{ue17})$, we have
{\small{\begin{align*}
&M_{i|\mu-1}^\alpha - \sum_{\beta=1}^q h_{0ij}^{\alpha\beta} \det \left(h_{0ij} \right) M_{j|\mu-1}^\beta\\
&= \sum_{\beta=1}^q (-1)^{q-\beta} N_{ij|\mu-1}^\beta w_{0i}^1 \wedge \cdots \wedge \widehat{w_{0i}^\beta} \wedge \cdots \wedge w_{0i}^q \wedge d w_{0i}^\alpha + w_{0i}^1\wedge \cdots \wedge w_{0i}^q \wedge \partial N_{ij|\mu-1}^\alpha - w_{0i}^1 \wedge \cdots \wedge w_{0i}^q \wedge \frac{\partial C_{ij|\mu}^\alpha }{\partial t}\\
&= \sum_{\beta=1}^q (-1)^{q-\beta} \left(\sum_{\eta=1}^q h_{0ij}^{\beta\eta}  N_{j|\mu-1}^\eta - N_{i|\mu-1}^\beta   \right) w_{0i}^1 \wedge \cdots \wedge \widehat{ w_{0i}^\beta } \wedge \cdots \wedge w_{0i}^q \wedge d w_{0i}^\alpha + w_{0i}^1 \wedge \cdots \wedge w_{0i}^q \wedge \partial \left( \sum_{\beta=1}^q h_{0ij}^{\alpha\beta} N_{j|\mu-1}^\beta - N_{i|\mu-1}^\alpha  \right)\\
&- w_{0i}^1 \wedge \cdots \wedge w_{0i}^q \wedge  \left( \sum_{\beta=1}^q  h_{0ij}^{\alpha\beta} \frac{\partial C_{j|\mu}^\beta}{\partial t}- \frac{\partial C_{i|\mu}^\alpha}{\partial t}+ \sum_{\beta=1}^q \frac{\partial D_{ij|\mu}^{\alpha\beta}}{\partial t} w_{0i}^\beta     \right) \\
&= - \sum_{\beta=1}^q (-1)^{q-\beta} N_{i|\mu-1}^\beta w_{0i}^1 \wedge \cdots \wedge \widehat{w_{0i}^\beta} \wedge \cdots \wedge w_{0i}^q \wedge d w_{0i}^\alpha - w_{0i}^1 \wedge \cdots \wedge w_{0i}^q \wedge \partial N_{i|\mu-1}^\alpha + w_{0i}^1 \wedge \cdots \wedge w_{0i}^q \wedge \frac{\partial C_{i|\mu}^\alpha}{\partial t}\\
&+ \sum_{\beta=1}^q (-1)^{q-\beta}\left(\sum_{\eta=1}^q h_{0ij}^{\beta\eta} N_{j|\mu-1}^\eta \right) \left( \sum_{\eta_1=1}^q h_{0ij}^{1\eta_1} w_{0j}^{\eta_1}  \right) \wedge \cdots \wedge \widehat{\left(\sum_{\eta=1}^q h_{0ij}^{\beta\eta} w_{0j}^\eta \right) }  \wedge \cdots \wedge \left( \sum_{\eta_q=1}^q h_{0ij}^{q \eta_q} w_{0j}^{\eta_q}  \right) \wedge d\left( \sum_{\gamma=1}^q h_{0ij}^{\alpha\gamma} w_{0j}^\gamma \right) \\
&+ w_{0i}^1 \wedge \cdots \wedge w_{0i}^q \wedge \partial \left( \sum_{\beta=1}^q h_{0ij}^{\alpha\beta} N_{j|\mu-1}^\beta \right) - \sum_{\beta=1}^q h_{0ij}^{\alpha\beta} \det\left(h_{0ij}\right) w_{0j}^1 \wedge \cdots \wedge w_{0j}^q \wedge \frac{\partial C_{j|\mu}^\beta}{\partial t}\\
&= - \sum_{\beta=1}^q (-1)^{q-\beta} N_{i|\mu-1}^\beta w_{0i}^1 \wedge \cdots \wedge \widehat{w_{0i}^\beta} \wedge \cdots \wedge w_{0i}^q \wedge d w_{0i}^\alpha - w_{0i}^1 \wedge \cdots \wedge w_{0i}^q \wedge \partial N_{i|\mu-1}^\alpha + w_{0i}^1 \wedge \cdots \wedge w_{0i}^q \wedge \frac{\partial C_{i|\mu}^\alpha}{\partial t}\\
&+ \sum_{\beta=1}^q (-1)^{q-\beta} \det\left(h_{0ij}\right) N_{j|\mu-1}^\beta w_{0j}^1 \wedge \cdots \wedge \widehat{w_{0j}^\beta}\wedge \cdots \wedge w_{0j}^q \wedge d\left( \sum_{\gamma=1}^q h_{0ij}^{\alpha\gamma} w_{0j}^\gamma \right) \\
&+ \det\left(h_{0ij} \right) w_{0j}^1 \wedge \cdots \wedge w_{0j}^q \wedge \partial \left( \sum_{\beta=1}^q h_{0ij}^{\alpha\beta} N_{j|\mu-1}^\beta \right) - \sum_{\beta=1}^q h_{0ij}^{\alpha\beta} \det\left(h_{0ij}\right) w_{0j}^1 \wedge \cdots \wedge w_{0j}^q \wedge \frac{\partial C_{j|\mu}^\beta}{\partial t}\\
&= - \sum_{\beta=1}^q (-1)^{q-\beta} N_{i|\mu-1}^\beta w_{0i}^1 \wedge \cdots \wedge \widehat{w_{0i}^\beta} \wedge \cdots \wedge w_{0i}^q \wedge dw_{0i}^\alpha - w_{0i}^1 \wedge \cdots \wedge w_{0i}^q \wedge \partial N_{i|\mu-1}^\alpha + w_{0i}^1 \wedge \cdots \wedge w_{0i}^q \wedge \frac{\partial C_{i|\mu}^\alpha}{\partial t}\\
&+\sum_{\gamma=1}^q h_{0ij}^{\alpha\gamma}\left( \sum_{\beta=1}^q (-1)^{q-\beta} \det\left(h_{0ij}\right) N_{j|\mu-1}^\beta w_{0j}^1 \wedge \cdots \wedge \widehat{w_{0j}^\beta}\wedge \cdots \wedge w_{0j}^q \wedge d  w_{0j}^\gamma \right) \\
&- \sum_{\beta=1}^q (-1)^{q-\beta} \det\left(h_{0ij}\right) N_{j|\mu-1}^\beta w_{0j}^1 \wedge \cdots \wedge \widehat{w_{0j}^\beta} \wedge \cdots \wedge w_{0j}^q \wedge w_{0j}^\beta \wedge d h_{0ij}^{\alpha\beta} \\
&+\sum_{\beta=1}^q h_{0ij}^{\alpha\beta} \det\left(h_{0ij}\right) w_{0j}^1 \wedge \cdots \wedge w_{0j}^q \wedge \partial N_{j|\mu-1}^\beta + \sum_{\beta=1}^q  \det\left(h_{0ij}\right)N_{j|\mu-1}^\beta w_{0j}^1 \wedge \cdots \wedge w_{0j}^q \wedge dh_{0ij}^{\alpha\beta}\\
& - \sum_{\beta=1}^q h_{0ij}^{\alpha\beta} \det\left(h_{0ij}\right) w_{0j}^1 \wedge \cdots \wedge w_{0j}^q \wedge \frac{\partial C_{j|\mu}^\beta}{\partial t}\\
&=  - \sum_{\beta=1}^q (-1)^{q-\beta} N_{i|\mu-1}^\beta w_{0i}^1 \wedge \cdots \wedge \widehat{w_{0i}^\beta} \wedge \cdots \wedge w_{0i}^q \wedge d w_{0i}^\alpha - w_{0i}^1 \wedge \cdots \wedge w_{0i}^q \wedge \partial N_{i|\mu-1}^\alpha + w_{0i}^1 \wedge \cdots \wedge w_{0i}^q \wedge \frac{\partial C_{i|\mu}^\alpha}{\partial t}\\
&+\sum_{\beta=1}^q h_{0ij}^{\alpha\beta}\det\left(h_{0ij}\right)\left(\sum_{\gamma=1}^q (-1)^{q-\gamma} N_{j|\mu-1}^{\gamma} w_{0j}^1 \wedge \cdots \wedge \widehat{w_{0j}^\gamma} \wedge \cdots \wedge w_{0j}^q \wedge d w_{0j}^\beta + w_{0j}^1 \wedge \cdots \wedge w_{0j}^q \wedge \partial N_{j|\mu-1}^\beta- w_{0j}^1 \wedge \cdots \wedge w_{0j}^q\wedge \frac{\partial C_{j|\mu}^\beta}{\partial t }\right)
\end{align*}}}
Then we have
{\Small{\begin{align}\label{ue19}
&M_{i|\mu-1}^\alpha   + \sum_{\beta=1}^q (-1)^{q-\beta} N_{i|\mu-1}^\beta w_{0i}^1 \wedge \cdots \wedge \widehat{w_{0i}^\beta} \wedge \cdots \wedge w_{0i}^q \wedge dw_{0i}^\alpha + w_i^1 \wedge \cdots \wedge w_i^q \wedge \partial N_{i|\mu-1}^\alpha - w_i^1 \wedge \cdots \wedge w_i^q \wedge \frac{\partial C_{i|\mu}^\alpha}{\partial t}\\
&=\sum_{\beta=1}^q h_{0ij}^{\alpha\beta}\det\left(h_{0ij}\right)\left(\sum_{\gamma=1}^q M_{j|\mu-1}^\beta + (-1)^{q-\gamma} N_{j|\mu-1}^{\gamma} w_{0j}^1 \wedge \cdots \wedge \widehat{w_{0j}^\gamma} \wedge \cdots \wedge w_{0j}^q \wedge dw_{0j}^\beta + w_{0j}^1 \wedge \cdots \wedge w_{0j}^q \wedge \partial N_{j|\mu-1}^\beta- w_{0j}^1 \wedge \cdots \wedge w_{0j}^q\wedge \frac{\partial C_{j|\mu}^\beta}{\partial t }\right) \notag
\end{align}}}
We define $\pi_{i|\mu-1}'\in \Gamma\left( U_i, \mathcal{A}^{0,0}\left( \mathscr{H}om_{\mathcal{O}_M}\left( \mathcal{N}_{\mathcal{F}_0}^*, \bigwedge^{q+1} \Omega_M^1\otimes \mathcal{L}_0  \right) \right) \right)$ by
\begin{align*}
\pi_{i|\mu-1}' :&\Gamma\left( U_i, \mathcal{N}_{\mathcal{F}_0}^* \right) \to \Gamma\left( U_i, \bigwedge^{q+1} \Omega_M^1 \otimes \mathcal{L}_0 \right) \\
                w_{0i}^\alpha &\mapsto M_{i|\mu-1}^\alpha   + \sum_{\beta=1}^q (-1)^{q-\beta} N_{i|\mu-1}^\beta w_{0i}^1 \wedge \cdots \wedge \widehat{w_{0i}^\beta} \wedge \cdots \wedge w_{0i}^q \wedge dw_{0i}^\alpha + w_i^1 \wedge \cdots \wedge w_i^q \wedge \partial N_{i|\mu-1}^\alpha - w_i^1 \wedge \cdots \wedge w_i^q \wedge \frac{\partial C_{i|\mu}^\alpha}{\partial t}
\end{align*}
Then $(\ref{ue19})$ implies that $\pi_\mu'=\left\{\pi_{i|\mu}' \right\}\in A^{0,0}\left( M, \mathscr{H}om_{\mathcal{O}_M}\left( \mathcal{N}_{\mathcal{F}_0}^*, \bigwedge^{q+1} \Omega_M^1 \otimes \mathcal{L}_0 \right) \right)$. From $(\ref{u43})$, by restricting $\pi_{\mu}'$ to $M-S$, we obtain
\begin{align*}
\pi_\mu:=\left\{\pi_{i|\mu} \right\} \in A^{0,0}\left(M-S, \mathscr{H}om_{\mathcal{O}_M}\left( \mathcal{N}_{\mathcal{F}_0}^*, \tilde{\mathcal{S}}^1 \right) \right)
\end{align*}
where $\pi_{i|\mu}=\pi_{i|\mu}'|_{U_i-S} :\Gamma\left( U_i-S, \mathcal{N}_{\mathcal{F}_0}^* \right)\to \Gamma\left( U_i-S, \tilde{\mathcal{S}}^1 \right)$. Then we claim that 
{\small{\begin{align}\label{us1}
\left(\frac{\partial I_\mu}{\partial t}, \pi_{\mu-1} \right)\oplus \left( \overline{\frac{\partial \phi_\mu}{\partial t}}, \sigma_{\mu-1} \right) \oplus  - \frac{\partial \xi_\mu}{\partial t} \in &A^{0,0}\left(M - S , \mathscr{H}om_{\mathcal{O}_M}\left( \mathcal{N}_{\mathcal{F}_0}^*, \tilde{\mathcal{S}}^2 \right) \right) \bigoplus A^{0,0}\left(M - S , \mathscr{H}om_{\mathcal{O}_M}\left( \mathcal{N}_{\mathcal{F}_0}^*, \tilde{\mathcal{S}}^1 \right) \right) \\
&\bigoplus \frac{A^{0,1}\left( M, \mathscr{H}om_{\mathcal{O}_M}\left( \mathcal{N}_{\mathcal{F}_0}^*, \Omega_M^1 \right) \right)}{A^{0,1}\left( M, \mathscr{H}om_{\mathcal{O}_M}\left( \mathcal{N}_{\mathcal{F}_0}^*, \mathcal{N}_{\mathcal{F}_0}^* \right) \right)} \bigoplus A^{0,1}\left( M, \mathscr{H}om_{\mathcal{O}_M}\left( \mathcal{N}_{\mathcal{F}_0}^*, \mathcal{O}_M \right) \right) \bigoplus A^{0,2} \left(M, \Theta_M \right) \notag
\end{align}}}
defines a $2$-cocycle in the following Dolbeault type bicomplex associated to $\tilde{\mathcal{N}}_{\mathcal{F}_0}^{*\bullet}$ (see Appendix \ref{up1}):
{\Tiny{\begin{center}
$\begin{CD}
\cdots \\
@A\hat{G}_3AA \\
A^{0,0}\left(M-S, \left(\mathcal{N}_{\mathcal{F}_0}^* \right)^* \otimes \tilde{\mathcal{S}}^3\right)\bigoplus A^{0,0}\left( M-S, \left(\mathcal{N}_{\mathcal{F}_0}^*\right)^*\otimes \tilde{\mathcal{S}}^2  \right) @>-(-1)^q\bar{\partial}>> \cdots \\
@ A \hat{G}_2 AA @A \hat{G}_2 AA \\
A^{0,0}\left(M-S, \left(\mathcal{N}_{\mathcal{F}_0}^*\right)^* \otimes \tilde{\mathcal{S}}^2\right) \bigoplus A^{0,0}\left( M-S, \left(\mathcal{N}_{\mathcal{F}_0}^* \right)^*\otimes \tilde{\mathcal{S}}^1 \right) @>(-1)^q \bar{\partial} >> A^{0,1}\left(M-S, \left(\mathcal{N}_{\mathcal{F}_0}^*\right)^*\otimes \tilde{\mathcal{S}}^2 \right) \bigoplus A^{0,1}\left( M-S, \left( \mathcal{N}_{\mathcal{F}_0}^*\right)^*\otimes \tilde{\mathcal{S}}^1 \right) @>\bar{\partial}>> \cdots\\
@A\hat{G}_1  AA @A\hat{G}_1 AA @A\hat{G}_1AA \\
\frac{A^{0,0}\left(M, \left(\mathcal{N}_{\mathcal{F}_0}^*\right)^* \otimes \Omega_M^1 \right)}{A^{0,0}\left(M, \left( \mathcal{N}_{\mathcal{F}_0}^* \right)^*\otimes \mathcal{N}_{\mathcal{F}_0}^* \right) }\bigoplus A^{0,0}\left(M, \left( \mathcal{N}_{\mathcal{F}_0}^* \right)^* \right) @>-\bar{\partial}>> \frac{A^{0,1}\left(M, \left( \mathcal{N}_{\mathcal{F}_0}^* \right)^*\otimes \Omega_M^1 \right)}{A^{0,1}\left( M, \left(\mathcal{N}_{\mathcal{F}_0}^* \right)^*\otimes \mathcal{N}_{\mathcal{F}_0}^* \right) }\bigoplus A^{0,1}\left(M, \left(\mathcal{N}_{\mathcal{F}_0}^* \right)^* \right) @>\bar{\partial}>> \cdots\\
@A\hat{G}_0AA @A \hat{G}_0 AA @A\hat{G}_0AA \\
A^{0,0}\left(M, \Theta_M\right) @>\bar{\partial}>> A^{0,1}\left( M, \Theta_M\right) @>-\bar{\partial}>> A^{0,2}\left(M, \Theta_M\right) 
\end{CD}$
\end{center}}}
From $(\ref{u29})$, we have
\begin{align*}
\bar{\partial}\left(-\frac{\partial \xi_\mu}{\partial t} \right)= -\frac{\partial}{\partial t}\left( \bar{\partial} \xi_\mu \right)=0
\end{align*}
We show that $\bar{\partial}\left( \overline{\frac{\partial\phi_\mu}{\partial t}} , \sigma_{\mu-1} \right)+ \hat{G}_0\left( -\frac{\partial \xi_\mu}{\partial t} \right)=0$, equivalently, $\bar{\partial}\overline{\frac{\partial \phi_\mu}{\partial t}}+ \hat{E}_0\left( -\frac{\partial \xi_\mu}{\partial t}\right)=0$ and $\bar{\partial}\left(\sigma_{\mu-1} \right) + \hat{G}_0'\left( -\frac{\partial \xi_\mu}{\partial t}\right)=0$. In fact, we have $\bar{\partial}\bar{\phi}_\mu+ \hat{E}_0\left( -\xi_\mu\right)=0$ from Part I. By taking the derivative with respect to $t$, we have $\bar{\partial}\overline{\frac{\partial \phi_\mu}{\partial t}}+ \hat{E}_0\left( -\frac{\partial \xi_\mu}{\partial t}\right)=0$. On the other hand, we note that from $(\ref{u32})$
\begin{align*}
\bar{\partial} \left(\psi_{i|\mu-1}^\alpha- \bar{\partial} N_{i|\mu-1}^\alpha\right) + i_{-\frac{\partial \xi_\mu}{\partial t}}\left( w_{0i}^\alpha \right)= \bar{\partial} \psi_{i|\mu}^\alpha - i_{\frac{\partial \xi_\mu}{\partial t}}\left( w_{0i}^\alpha \right)=0
\end{align*}
We show that $(-1)^q \bar{\partial} \left( \frac{\partial I_\mu}{\partial t},  \pi_{\mu-1} \right) + \hat{G}_1\left( \overline{\frac{\partial \phi_\mu }{\partial t} } , \sigma_{\mu-1} \right)=0$, equivalently, $(-1)^q\bar{\partial} \left(\frac{\partial I_\mu}{\partial t}\right)+ \hat{E}_1\left(\overline{\frac{\partial \phi_\mu}{\partial t}} \right)=0$ and $(-1)^q\bar{\partial} \pi_{\mu-1}^\alpha + \hat{E}_0\left( \sigma_{i|\mu-1}\right)\left( w_{0i}^\alpha \right) - w_{0i}^1 \wedge \cdots \wedge w_{0i}^q \wedge \frac{\partial \phi_{i|\mu}^\alpha}{\partial t}=0$. In fact, we have $(-1)^q\bar{\partial} I_\mu+ \hat{E}_1\left(\bar{\phi}_\mu \right)=0$ from Part I. By taking the derivative with respect $t$, we have $(-1)^q\bar{\partial} \left(\frac{\partial I_\mu}{\partial t}\right)+ \hat{E}_1\left(\overline{\frac{\partial \phi_\mu}{\partial t}} \right)=0$. On the other hand, we note that from $(\ref{u39})$
{\small{\begin{align*}
&\bar{\partial} \left(M_{i|\mu-1}^\alpha + \sum_{\beta=1}^q(-1)^{q-\beta} N_{i|\mu-1}^\beta w_{0i}^1 \wedge \cdots \wedge \widehat{w_{0i}^\beta}\wedge \cdots  \wedge  w_{0i}^q \wedge d w_{0i}^\alpha + w_{0i}^1 \wedge \cdots \wedge w_{0i}^q\wedge \partial N_{i|\mu-1}^\alpha- w_{0i}^1 \wedge \cdots \wedge w_{0i}^q \wedge \frac{\partial C_{i|\mu}^\alpha}{\partial t} \right)\\
&= (-1)^q \sum_{\beta=1}^q w_{0i}^1 \wedge \cdots \wedge \psi_{i|\mu-1}^\beta \wedge \cdots \wedge w_{0i}^q \wedge d w_{0i}^\alpha + (-1)^q w_{0i}^1 \wedge \cdots \wedge w_{0i}^q \wedge \partial \psi_{i|\mu-1}^\alpha  + (-1)^q w_{0i}^1 \wedge \cdots \wedge w_{0i}^q \wedge \frac{\partial A_{i|\mu}^\alpha}{\partial t}\\
&+\sum_{\beta=1}^q (-1)^{q-\beta} \bar{\partial} N_{i|\mu-1}^\beta \wedge w_{0i}^1 \wedge \cdots \wedge \widehat{ w_{0i}^\beta } \wedge \cdots \wedge w_{0i}^q \wedge d w_{0i}^\alpha -(-1)^q w_{0i}^1 \wedge \cdots \wedge w_{0i}^q \wedge \partial \left( \bar{\partial} N_{i|\mu-1}^\alpha \right) - (-1)^q w_{0i}^1 \wedge \cdots \wedge w_{0i}^q \wedge \frac{\partial}{\partial t}\left( \bar{\partial} C_{i|\mu}^\alpha\right)\\
&= (-1)^q \sum_{\beta=1}^q w_{0i}^1 \wedge \cdots \wedge \psi_{i|\mu-1}^\beta \wedge \cdots \wedge w_{0i}^q \wedge dw_{0i}^\alpha + (-1)^q w_{0i}^1 \wedge \cdots \wedge w_{0i}^q \wedge \partial \psi_{i|\mu-1}^\alpha  + (-1)^q w_{0i}^1 \wedge \cdots \wedge w_{0i}^q \wedge \frac{\partial A_{i|\mu}^\alpha}{\partial t}\\
&-\sum_{\beta=1}^q (-1)^{q} w_{0i}^1 \wedge \cdots \wedge \bar{\partial} N_{i|\mu-1}^\beta\wedge \cdots \wedge w_{0i}^q \wedge dw_{0i}^\alpha -(-1)^q w_{0i}^1 \wedge \cdots \wedge w_{0i}^q \wedge \partial \left( \bar{\partial} N_{i|\mu-1}^\alpha \right) - (-1)^q w_i^1 \wedge \cdots \wedge w_i^q \wedge \frac{\partial}{\partial t}\left( \bar{\partial} C_{i|\mu}^\alpha\right)\\
&= (-1)^q \sum_{\beta=1}^q  w_{0i}^1 \wedge \cdots \wedge \left(\psi_{i|\mu-1}^\beta - \bar{\partial} N_{i|\mu-1}^\beta \right) \wedge \cdots \wedge w_{0i}^q \wedge d w_{0i}^\alpha + (-1)^q  w_{0i}^1  \wedge \cdots \wedge w_{0i}^q \wedge \partial \left( \psi_{i|\mu-1}^\alpha - \bar{\partial}  N_{i|\mu-1}^\beta \right) \\
&+ (-1)^q w_i^1 \wedge \cdots \wedge w_i^q \wedge \left(\frac{\partial A_{i|\mu}^\alpha}{\partial t} - \frac{\partial}{\partial t}\left( \bar{\partial} C_{i|\mu}^\alpha\right)\right)
\end{align*}}}
Then we have
{\tiny{\begin{align*}
&(-1)^{q}\bar{\partial} \pi_{\mu-1}^\alpha\\
&= \sum_{\beta=1}^q  w_{0i}^1 \wedge \cdots \wedge \left(\psi_{i|\mu-1}^\beta - \bar{\partial} N_{i|\mu-1}^\beta \right) \wedge \cdots \wedge w_{0i}^q \wedge d w_{0i}^\alpha +   w_{0i}^1  \wedge \cdots \wedge w_{0i}^q \wedge \partial \left( \psi_{i|\mu-1}^\alpha - \bar{\partial}  N_{i|\mu-1}^\beta \right) + w_{0i}^1 \wedge \cdots \wedge w_{0i}^q \wedge  \frac{\partial}{\partial t}\left(  A_{i|\mu}^\alpha - \bar{\partial} C_{i|\mu}^\alpha  + \sum_{\gamma=1}^q B_{i|\mu}^{\alpha\gamma} w_{0i}^\gamma \right)
\end{align*} }}

It remains to show $\hat{G}_2\left( \frac{\partial I_\mu}{ \partial t}, \pi_{\mu-1} \right)=0$, equivalently, $\hat{E}_2\left(\frac{ \partial I_\mu}{ \partial t} \right) = 0 $ and $\hat{E}_1\left(\pi_{\mu-1} \right)  + \frac{\partial I_\mu}{\partial t}=0$. In fact, we have $\hat{E}_2\left( I_\mu \right)=0$ from Part I. By taking the derivative with respect to $t$, we have $\hat{E}_2\left( \frac{\partial I_\mu}{\partial t} \right)=0$. On the other hand, we show that $\hat{E}_1\left(\pi_{\mu-1} \right)  + \frac{\partial I_\mu}{\partial t}=0$. It is sufficient to show that $\hat{E}_1\left( M_{i|\mu}\right) + \frac{\partial E_{i|\mu}}{\partial t}=0$. We note that we have on $U_{i_x}$
{\Small{\begin{align*}
&\partial_{M\times B} \left( E_{i|\mu}^\alpha + M_{i|\mu-1}^\alpha \wedge dt \right)= \partial \left(E_{i|\mu}^\alpha \right) + (-1)^q \frac{\partial E_{i|\mu}}{\partial t} \wedge dt + \partial  \left( M_{i|\mu-1}^\alpha \right)\wedge dt \\
& \equiv_\mu \partial_{M\times B} \left(\left(w_i^{1(\mu-1)}+ v_i^{1(\mu-2)}dt \right) \wedge \cdots \wedge \left( w_i^{q(\mu-1)}+ v_i^{q(\mu-2)}dt \right) \wedge \partial_{M\times B}  \left( w_i^{\alpha(\mu-1)}+ v_i^{\alpha(\mu-2)} dt\right)       \right)\\
&\equiv_\mu \sum_{\beta=1}^q (-1)^{\beta-1} \left( w_i^{1(\mu-1)} + v_i^{1(\mu-2)} dt   \right) \wedge \cdots \wedge \partial_{M\times B}\left( w_i^{\beta(\mu-1)}+ v_i^{\beta(\mu-2)}dt \right) \wedge \cdots \wedge \left(  w_i^{q(\mu-1)}+ v_i^{q(\mu-2)}dt \right) \wedge \partial_{M\times B} \left( w_i^{\alpha(\mu-1)}+ v_i^{\alpha(\mu-2)}dt \right)\\
&\equiv_\mu \sum_{\beta=1}^q (-1)^{\beta-1} \left( w_i^{1(\mu-1)} + v_i^{1(\mu-2)} dt   \right) \wedge \cdots \wedge \partial_{M\times B} \left( w_i^{\beta(\mu-1)}+ v_i^{\beta(\mu-2)}dt \right) \wedge \cdots \wedge \left(  w_i^{q(\mu-1)}+ v_i^{q(\mu-2)}dt \right) \\
& \wedge \left( \sum_{\gamma=1}^q \left(a_{i_x}^{\alpha\gamma (\mu-1)}+ b_{i_x}^{\alpha\gamma(\mu-2)}dt \right) \wedge \left( w_i^{\gamma(\mu-1)} + v_i^{\gamma(\mu-2)}dt  \right)   + L_{i_x|\mu}^\alpha + P_{i_x|\mu-1}^\alpha \wedge dt    \right)\\ 
&\equiv_\mu \sum_{\beta=1}^q \left( a_{i_x}^{\alpha\beta(\mu-1)}+ b_{i_x}^{\alpha\beta(\mu-2)}dt\right) \wedge \left( w_i^{1(\mu-1)} + v_i^{1(\mu-2)}dt \right)  \wedge \cdots \wedge \left( w_i^{q(\mu-1)}+ v_i^{q(\mu-2)} dt \right) \wedge \partial_{M\times B} \left( w_i^{\beta(\mu-1)}+  v_i^{\beta(\mu-2)}dt \right) \\
&+ d \left( w_{0i}^1 \wedge \cdots w_{0i}^q\right) \wedge L_{i_x|\mu}^\alpha + d \left( w_{0i}^1 \wedge \cdots \wedge w_{0i}^q \right) \wedge P_{i_x|\mu-1}^\alpha \wedge dt \\
&\equiv_\mu \sum_{\beta=1}^q \left(a_{i_x}^{\alpha\beta(\mu-1)}+  b_{i_x}^{\alpha\beta(\mu-2)}dt \right) \wedge \left( E_{i|\mu}^\beta + M_{i|\mu-1}^\beta \wedge dt    \right) +d \left(w_{0i}^1 \wedge \cdots \wedge w_{0i}^q \right) \wedge L_{i_x|\mu}^\alpha + d\left(w_{0i}^1\wedge \cdots \wedge w_{0i}^q \right) \wedge P_{i_x|\mu-1}^\alpha \wedge dt \\
&\equiv_\mu \sum_{\beta=1}^q a_{0i_x}^{\alpha\beta} \wedge E_{i|\mu}^\beta + \sum_{\beta=1}^q a_{0i_x}^{\alpha\beta} \wedge M_{i|\mu-1}^\beta \wedge dt + d \left(w_{0i}^1  \wedge \cdots \wedge w_{0i}^q \right) \wedge L_{i_x|\mu}^\alpha + d \left(w_{0i}^1\wedge \cdots \wedge w_{0i}^q \right)\wedge P_{i_x|\mu-1}^\alpha \wedge dt
\end{align*}}}

Then we have
\begin{align}
&\partial E_{i|\mu}^\alpha - \sum_{\beta=1}^q a_{0i_x}^{\alpha\beta} \wedge E_{j|\mu}^\beta - d \left(w_{0i}^1\wedge \cdots \wedge w_{0i}^q \right) \wedge L_{i_x|\mu}^\alpha =0 \notag \\
(-1)^q \frac{\partial E_{i|\mu}}{\partial t} + &\partial M_{i|\mu-1}^\alpha - \sum_{\beta=1}^q a_{0i_x}^{\alpha\beta} \wedge M_{j|\mu-1}^\beta - d \left(w_{0i}^1\wedge \cdots \wedge w_{0i}^q \right)\wedge P_{i_x|\mu-1}^\alpha=0 \label{ue20}
\end{align}
From $(\ref{u43})$, $M_{i|\mu-1}^\alpha= w_{0i}^1\wedge \cdots \wedge w_{0i}^q \wedge P_{i_x|\mu-1}^\alpha$ on $U_{i_x}$. Then from $(\ref{ue20})$, we have
\begin{align*}
&\hat{E}_1\left( M_{i|\mu} \right)\left( w_{0i}^\alpha \right) + \frac{\partial E_{i|\mu}^\alpha}{\partial t}= w_{0i}^1 \wedge \cdots \wedge w_{0i}^q\wedge \left(\partial P_{i_x|\mu}^\alpha  - \sum_{\beta=1}^q a_{0i_x}^{\alpha\beta}\wedge  P_{i_x|\mu}^\beta\right) + \frac{\partial E_{i|\mu}^\alpha}{\partial t}\\
&= (-1)^q \partial \left(w_{0i}^1\wedge \cdots \wedge w_{0i}^q \wedge P_{i_x|\mu}^\alpha   \right) - (-1)^q d\left( w_{0i}^1 \wedge \cdots \wedge w_{0i}^q \right) \wedge P_{i_x|\mu}^\alpha - (-1)^q \sum_{\beta=1}^q a_{0i_x}^{\alpha\beta} w_{0i}^1 \wedge \cdots \wedge w_{0i}^q \wedge P_{i_x|\mu}^\beta + \frac{\partial E_{i|\mu}^\alpha}{\partial t}\\
&= (-1)^q \partial M_{i|\mu-1}^\alpha - (-1)^q d\left( w_{0i}^1 \wedge \cdots \wedge w_{0i}^q \right) \wedge P_{i_x}^\alpha - (-1)^q \sum_{\beta=1}^q a_{0i_x}^{\alpha\beta} \wedge M_{i|\mu-1}^\beta + \frac{\partial E_{i|\mu}^\alpha}{\partial t} = 0
\end{align*}
This proves $\left( \frac{\partial I_\mu}{\partial t}, \pi_{\mu-1}, \overline{\frac{\partial \phi_\mu}{\partial t}}, \sigma_{\mu-1},  - \frac{\partial \xi_\mu}{\partial t} \right)$ defines a $2$-cocycle in the Dolbeault type bicomplex associated to $\tilde{\mathcal{N}}_{\mathcal{F}_0}^{*\bullet}$. Then by the hypothesis $\mathbb{H}^2 \left(M, \mathcal{N}_{\mathcal{F}_0}^{*\bullet} \right)=0$, we can find
\begin{align*}
 \left(\overline{\tilde{w}_{\mu-1}'}, v_{\mu-1}' , \tilde{\varphi}_{\mu-1}' \right) \in \frac{A^{0,0}\left( M, \mathscr{H}om_{\mathcal{O}_M}\left( \mathcal{N}_{\mathcal{F}_0}^*, \Omega_M^1 \right) \right)}{A^{0,0}\left( M, \mathscr{H}om_{\mathcal{O}_M}\left( \mathcal{N}_{\mathcal{F}_0}^*, \mathcal{N}_{\mathcal{F}_0}^* \right) \right)}\bigoplus A^{0,0}\left( M, \mathscr{H}om_{\mathcal{O}_M} \left( \mathcal{N}_{\mathcal{F}_0}^*,  \mathcal{O}_M \right) \right) \bigoplus A^{0,1}\left( M, \Theta_M \right)
\end{align*}
such that
\begin{align}
-\bar{\partial} \tilde{\varphi}_{\mu-1}' &= -\frac{\partial \xi_\mu}{\partial t} \label{u62}
\end{align}
\begin{align}
                                  \overline{ -\bar{\partial}\tilde{w}_{\mu-1}' } + \hat{E}_0(\varphi_{\mu-1}')&= \overline{\left\{ w_{0i}^\alpha \mapsto \frac{\partial A_{i|\mu}^\alpha}{\partial t} - \bar{\partial}\left( \frac{\partial C_{i|\mu}^\alpha}{\partial t} \right)  + \sum_{\gamma=1}^q \frac{\partial B_{i|\mu}^{\alpha\gamma}}{\partial t} w_{0i}^\gamma     \right\} } \label{u63}
                                   \end{align}
                                   \begin{align}
\hat{E}_1 \left( \tilde{w}_{\mu-1}' \right)&=\left\{w_{0i}^\alpha \mapsto \frac{\partial E_{i|\mu}^\alpha}{\partial t} + \sum_{\beta=1}^q w_{0i}^1\wedge \cdots \wedge \overbrace{ \frac{\partial C_{i|\mu}^\beta}{\partial t}}^{\beta-\textnormal{th}}\wedge \cdots \wedge w_{0i}^q \wedge d w_{0i}^\alpha + w_{0i}^1\wedge \cdots \wedge w_{0i}^q \wedge \partial \left( \frac{\partial C_{i|\mu}^\alpha}{\partial t} \right)   \right\} \label{ue64}
\end{align}
\begin{align}
 -\bar{\partial} v_{\mu-1}'+ \hat{G}_0' \left( \tilde{\varphi}_{\mu-1}' \right)&=\left\{ w_i^\alpha \mapsto  \psi_{i|\mu-1}^\alpha - \bar{\partial} N_{i|\mu-1}^\alpha   \right\}  \label{u64} 
 \end{align}
{\Small{\begin{align}
&-\tilde{w}_{\mu-1}'+ \hat{E}_0\left(v_{\mu-1}' \right) \label{u65}\\
 &=\left\{ w_{0i}^\alpha \mapsto M_{i|\mu-1}^\alpha + \sum_{\beta=1}^q (-1)^{q-\beta} N_{i|\mu-1}^\beta w_{0i}^1 \wedge \cdots \widehat{w_{0i}^\beta} \wedge \cdots \wedge w_{0i}^q \wedge d w_{0i}^\alpha + w_{0i}^1 \wedge \cdots \wedge w_{0i}^q \wedge \partial N_{i|\mu-1}^\alpha - w_{0i}^1 \wedge \cdots \wedge w_{0i}^q \wedge \frac{\partial C_{i|\mu}^\alpha}{\partial t}    \right\} \notag
\end{align}}}
where $\overline{\tilde{w}_{\mu-1}'}$ is the image of the global section $\tilde{w}_{\mu-1}'=\left\{\tilde{w}_{i|\mu-1}' \right\}\in A^{0,0}\left( M, \mathscr{H}om_{\mathcal{O}_M}\left(  \mathcal{N}_{\mathcal{F}_0}^*, \Omega_M^1   \right) \right)$ in $\frac{A^{0,0}\left( M, \mathscr{H}om_{\mathcal{O}_M}\left( \mathcal{N}_{\mathcal{F}_0}^* ,  \Omega_M^1 \right) \right)}{A^{0,0}\left( M,  \mathscr{H}om_{\mathcal{O}_M}\left(  \mathcal{N}_{\mathcal{F}_0}^* ,  \mathcal{N}_{\mathcal{F}_0}^* \right) \right)}$ such that $\tilde{w}_{i|\mu-1}' : \Gamma\left( U_i, \mathcal{N}_{\mathcal{F}_0}^*\right) \to \Gamma\left( U_i, \mathcal{A}^{0,0}\left( \Omega_M^1 \right)\right), w_{0i}^\alpha \mapsto \tilde{w}_{i|\mu-1}'^\alpha,\alpha=1,...,q$. In particular, we have
\begin{align}\label{ne17}
\tilde{w}_{i|\mu-1}'^\alpha= \sum_{\beta=1}^q h_{0ij}^{\alpha\beta} \tilde{w}_{j|\mu}'^\beta
\end{align}
We note that 
\begin{align*}
\xi_\mu&=\frac{1}{\mu}t\frac{\partial \xi_\mu}{\partial t}\\
A_{i|\mu}^\alpha &= \frac{1}{\mu} t \frac{\partial A_{i|\mu}^\alpha}{\partial t} \\
C_{i|\mu}^\alpha &= \frac{1}{\mu}t \frac{ \partial C_{i|\mu}^\alpha}{\partial t}\\
B_{i|\mu}^{\alpha\gamma}&= \frac{1}{\mu}t \frac{\partial B_{i|\mu}^{\alpha\beta}}{\partial t}\\
E_{i|\mu}^\alpha &=\frac{1}{\mu} t\frac{\partial E_{i|\mu}^\alpha}{\partial t}
\end{align*}

We note that if we set
\begin{align}
\varphi_\mu' &:= \frac{1}{\mu}t \tilde{\varphi}'_{\mu-1} \Longrightarrow \frac{\partial \varphi_\mu'}{\partial t}= \tilde{\varphi}_{\mu-1}' \label{u67}\\
w_\mu' &:= \frac{1}{\mu} t\tilde{w}_{\mu-1}' \Longrightarrow \frac{\partial w_\mu'}{\partial t}= \tilde{w}'_{\mu-1} \label{u68}
\end{align}
Then from $(\ref{u62}),(\ref{u63})$ and $(\ref{ue64})$, we have
\begin{align}
- \bar{\partial} \varphi_\mu' &=- \xi_\mu \label{u60}
\end{align}
\begin{align}
- \bar{\partial} w_\mu' + \hat{E}_0 \left(\varphi_\mu' \right)&= \overline{\left\{ w_{0i}^\alpha \mapsto A_{i|\mu}^{\alpha}- \bar{\partial} C_{i|\mu}^\alpha + \sum_{\gamma=1}^q B_{i|\mu}^{\alpha\gamma} w_{0i}^\gamma  \right\}}
\end{align}
\begin{align}
E_1'(\tilde{w}_\mu')&= \left\{ w_{0i}^\alpha \mapsto  E_{i|\mu}^\alpha+ \sum_{\beta=1}^q w_{0i}^1 \wedge \cdots \wedge C_{i|\mu}^\beta \wedge \cdots \wedge w_{0i}^q \wedge d w_{0i}^\alpha + w_{0i}^1 \wedge \cdots \wedge w_{0i}^q \wedge \partial C_{i|\mu}^\alpha  \right\} \label{u61}
\end{align}
We need to find $\varphi_\mu, h_{ij|\mu}^{\alpha\beta}, w_{i|\mu}^\alpha, a_{i_x|\mu}^{\alpha\beta}$, and $v_{i|\mu-1}^\alpha, b_{i_x|\mu-1}^{\alpha\beta}$ satisfying $(\ref{u45})-(\ref{u55})$.  From $(\ref{u60})-(\ref{u61})$, we already have $\varphi_\mu, h_{ij|\mu}^{\alpha\beta}, w_{i|\mu}^\alpha$ and $a_{i_x|\mu}^{\alpha\beta}$  satisfying $(\ref{u45})-(\ref{u55})$ from Part I. We keep the notations in Part I. We recall from Part I that we take
\begin{align}
\varphi_\mu &:=-\varphi_\mu'  \label{uut2}\\
w_{i|\mu}^\alpha&:= C_{i|\mu}^\alpha - w_{i|\mu}'^\alpha + \sum_{\beta=1}^q P_{i|\mu}^{\alpha\beta} w_i^\beta \label{ue22} \\
h_{ij|\mu}^{\alpha\beta}&:= \sum_{\gamma=1}^q D_{ij|\mu}^{\alpha\gamma} h_{0ij}^{\gamma\beta} + \sum_{\gamma=1}^q P_{i|\mu}^{\alpha\gamma} h_{0ij}^{\gamma\beta} - \sum_{\gamma=1}^q h_{0ij}^{\alpha\gamma} P_{j|\mu}^{\gamma\beta} \label{uut3}
\end{align} 
For the definition of $P_{i|\mu}^{\alpha\beta}$, see Part I. Don't confuse $P_{i|\mu}^{\alpha\beta}$ with $P_{i_x|\mu-1}^\alpha$. It remains to find $v_{i|\mu-1}^\alpha$ and $b_{i_x|\mu-1}^{\alpha\beta}$ satisfying $(\ref{u47}),(\ref{u51}),(\ref{u53})$ and $(\ref{u55})$. We find $v_{i|\mu-1}^\alpha$ satisfying $(\ref{u47})$. From $(\ref{u64})$ and $(\ref{u67})$  we have
\begin{align*}
-\bar{\partial} v_{i|\mu-1}'^\alpha+ i_{\tilde{\varphi}_{\mu-1}'} \left(w_{0i}^\alpha \right) = \psi_{i|\mu-1}^\alpha - \bar{\partial} N_{i|\mu-1}^\alpha & \iff \bar{\partial}(N_{i|\mu-1}^\alpha- v_{i|\mu-1}'^\alpha )+ i_{\frac{\partial \varphi_\mu'}{\partial t}} (w_i^\alpha)= \psi_{i|\mu-1}^\alpha  \\
 &\iff \bar{\partial} ( N_{i|\mu-1}^\alpha - v_{i|\mu-1}'^\alpha ) - i_{\frac{\partial \varphi_\mu}{\partial t}}(w_i^\alpha)= \psi_{i|\mu-1}^\alpha
\end{align*}
Then if we take 
\begin{align}\label{ue24}
v_{i|\mu-1}^\alpha:= N_{i|\mu-1}^\alpha - v_{i|\mu-1}'^\alpha,
\end{align}
 then $(\ref{u47})_{\mu-1}$ is satisfied. Next we check $(\ref{u51})_{\mu-1}$. We recall that we have $\sum_{\beta=1}^q h_{0ij}^{\alpha\beta} N_{j|\mu-1}^\beta- N_{i|\mu-1}^\alpha= N_{ij|\mu-1}^\alpha$ from $(\ref{ue17})$ and by our construction $v_{i|\mu-1}'^\alpha= \sum_{\beta=1}^q  h_{0ij}^{\alpha\beta} v_{j|\mu-1}'^\beta$. Then we have from $(\ref{ue24})$
\begin{align*}
N_{ij|\mu-1}^\alpha= - v_{i|\mu-1}^\alpha + \sum_{\beta=1}^q h_{0ij}^{\alpha\beta} v_{j|\mu-1}^\beta &\iff \sum_{\beta=1}^q h_{0ij}^{\alpha\beta} N_{j|\mu-1}^\beta- N_{i|\mu-1}^\alpha = - N_{i|\mu-1}^\alpha + v_{i|\mu-1}'^\alpha + \sum_{\beta=1}^q h_{0ij}^{\alpha\beta} \left( N_{j|\mu-1}^\beta - v_{j|\mu-1}'^\beta \right) \\
 &\iff v_{i|\mu-1}'^\alpha= \sum_{\beta=1}^q  h_{0ij}^{\alpha\beta} v_{j|\mu-1}'^\beta
\end{align*}

We check $(\ref{u53})_{\mu-1}$. From $(\ref{u65})$ and $(\ref{u68})$
{\small{\begin{align} \label{ue21}
&- w_{0i}^1\wedge \cdots \wedge w_{0i}^q\wedge  \tilde{w}_{i|\mu-1}'^\alpha  + \sum_{\beta=1}^q (-1)^{q-\beta} v_{i|\mu-1}'^\beta w_{0i}^1 \wedge \cdots \wedge \widehat{w_{0i}^\beta } \wedge \cdots \wedge  w_{0i}^q \wedge d w_{0i}^\alpha + w_{0i}^1 \wedge \cdots \wedge w_{0i}^q\wedge \partial v_{i|\mu-1}'^\beta\\
&= M_{i|\mu-1}^\alpha + \sum_{\beta=1}^q (-1)^{q-\beta} N_{i|\mu-1}^\beta w_{0i}^1 \wedge \cdots \wedge \widehat{w_{0i}^\beta} \wedge \cdots \wedge w_{0i}^q \wedge d w_{0i}^\alpha + w_{0i}^1\wedge \cdots \wedge w_{0i}^q \wedge \partial N_{i|\mu-1}^\alpha - w_{0i}^1 \wedge \cdots \wedge w_{0i}^q \wedge \frac{\partial C_{i|\mu}^\alpha}{\partial t} \notag
\end{align}}}
Since $\frac{\partial w_{i|\mu}^\alpha}{\partial t}= \frac{\partial C_{i|\mu}^\alpha}{\partial t}- \tilde{w}_{i|\mu-1}'^\alpha  + \sum_{\beta=1}^q \frac{\partial \tilde{P}_{i|\mu}^{\alpha\beta} }{\partial t} w_{0i}^\beta$ from $(\ref{ue22})$ and $(\ref{u68})$, we have from $(\ref{ue21})$ and $(\ref{ue24})$
\begin{align}\label{ue26}
 M_{i|\mu-1}^\alpha + \sum_{\beta=1}^q (-1)^{q-\beta} v_{i|\mu-1}^\beta w_{0i}^1 \wedge \cdots \wedge \widehat{w_{0i}^\beta} \wedge \cdots \wedge w_{0i}^q \wedge d w_{0i}^\alpha + w_{0i}^1 \wedge \cdots \wedge w_{0i}^q \wedge \partial v_{i|\mu-1}^\alpha - w_{0i}^1 \wedge \cdots \wedge w_{0i}^q \wedge \frac{\partial w_{i|\mu}^\alpha}{\partial t}=0
\end{align}
By definition of $\hat{G}_2$, $(\ref{ue26})$ holds on $U_i-S$. We show that $(\ref{ue26})$ actually holds on $U_i$. In fact, first we note that 
\begin{align}\label{ue28}
0 = i_{\frac{\partial \varphi_\mu}{\partial t}}\left(  w_{0i}^1 \wedge \cdots \wedge w_{0i}^q \wedge d w_{0i}^\alpha  \right) = \sum_{\beta=1}^q w_{0i}^1 \wedge \cdots \wedge  i_{\frac{\partial \varphi_\mu}{\partial t}} \left( w_{0i}^\alpha \right) \wedge \cdots \wedge w_{0i}^q \wedge dw_{0i}^\alpha + w_{0i}^1\wedge \cdots \wedge w_{0i}^q \wedge \left(i_{\frac{\partial \varphi_\mu}{\partial t}}\partial w_{0i}^\alpha \right)
\end{align}

from $(\ref{u39})$ and $(\ref{u47})$ and $(\ref{ue28})$
{\small{\begin{align*}
&\bar{\partial}\left(  M_{i|\mu-1}^\alpha + \sum_{\beta=1}^q (-1)^{q-\beta} v_{i|\mu-1}^\beta w_{0i}^1 \wedge \cdots \wedge \widehat{w_{0i}^\beta} \wedge \cdots \wedge w_{0i}^q \wedge d w_{0i}^\alpha + w_{0i}^1 \wedge \cdots \wedge w_{0i}^q \wedge \partial v_{i|\mu-1}^\alpha - w_{0i}^1 \wedge \cdots \wedge w_{0i}^q \wedge \frac{\partial w_{i|\mu}^\alpha}{\partial t} \right)\\
&=(-1)^q \sum_{\beta=1}^q w_{0i}^1\wedge \cdots \wedge \psi_{i|\mu-1}^\beta \wedge \cdots w_{0i}^q \wedge dw_{0i}^\alpha +  (-1)^q w_{0i}^1\wedge \cdots \wedge w_{0i}^q \wedge \partial \psi_{i|\mu-1}^\alpha  +  (-1)^q w_{0i}^1\wedge \cdots \wedge w_{0i}^q \wedge \frac{\partial A_{i|\mu}^\alpha}{\partial t} \\
&+ \sum_{\beta=1}^q (-1)^{q-\beta} \left( \psi_{i|\mu-1}^\beta + i_{\frac{\partial \varphi_\mu}{\partial t}}\left( w_{0i}^\beta \right) \right) w_{0i}^1 \wedge \cdots \wedge \widehat{w_{0i}^\beta} \wedge \cdots \wedge w_{0i}^q \wedge d w_{0i}^\alpha - (-1 )^q w_{0i}^1 \wedge \cdots \wedge w_{0i}^q \wedge  \partial \left( \psi_{i|\mu-1}^\alpha + i_{\frac{\partial \varphi_\mu}{\partial t}}\left( w_{0i}^\alpha \right) \right) \\
& - (-1)^q w_{0i}^1 \wedge \cdots \wedge w_{0i}^q \wedge \frac{\partial }{\partial t}\left( A_{i|\mu}^\alpha + \mathcal{L}_{\varphi_\mu}\left( w_{0i}^\alpha \right)\right) \\
&= (-1)^{q-1}\sum_{\beta=1}^q w_{0i}^1 \wedge \cdots \wedge i_{\frac{\partial \varphi_\mu}{\partial t}}\left( w_{0i}^\beta \right) \wedge \cdots \wedge w_{0i}^q \wedge dw_{0i}^\alpha- (-1)^q w_{0i}^1 \wedge \cdots \wedge w_{0i}^q \wedge \left(\partial i_{\frac{\partial \varphi_\mu}{\partial t}} \left(w_{0i}^\alpha \right) + \mathcal{L}_{\frac{\partial \varphi_\mu}{\partial t}} \left( w_{0i}^\alpha \right) \right)\\
&= (-1)^q w_{0i}^1 \wedge \cdots \wedge w_{0i}^q \wedge \left( i_{\frac{\partial \varphi_\mu}{\partial t}}\partial\left(w_{0i}^\alpha \right) - \partial  i_{\frac{\partial \varphi_\mu}{\partial t}}\left(w_{0i}^\alpha \right) - \mathcal{L}_{\frac{\partial \varphi_\mu}{\partial t}}\left( w_{0i}^\alpha \right)   \right) = 0
\end{align*}}}
This implies that $ M_{i|\mu-1}^\alpha + \sum_{\beta=1}^q (-1)^{q-\beta} v_{i|\mu-1}^\beta w_{0i}^1 \wedge \cdots \wedge \widehat{w_{0i}^\beta} \wedge \cdots \wedge w_{0i}^q \wedge d w_{0i}^\alpha + w_{0i}^1 \wedge \cdots \wedge w_{0i}^q \wedge \partial v_{i|\mu-1}^\alpha - w_{0i}^1 \wedge \cdots \wedge w_{0i}^q \wedge \frac{\partial w_{i|\mu}^\alpha}{\partial t}$ is holomorphic on $U_i$. Hence $(\ref{ue26})$ implies that $ M_{i|\mu-1}^\alpha + \sum_{\beta=1}^q (-1)^{q-\beta} v_{i|\mu-1}^\beta w_{0i}^1 \wedge \cdots \wedge \widehat{w_{0i}^\beta} \wedge \cdots \wedge w_{0i}^q \wedge d w_{0i}^\alpha + w_{0i}^1 \wedge \cdots \wedge w_{0i}^q \wedge \partial v_{i|\mu-1}^\alpha - w_{0i}^1 \wedge \cdots \wedge w_{0i}^q \wedge \frac{\partial w_{i|\mu}^\alpha}{\partial t}=0$ on $U_i$, so that $(\ref{u53})_{\mu-1}$ is satisfied. If remains to find $b_{i_x|\mu-1}^{\alpha\beta}$ satisfying $(\ref{u55})_{\mu-1}$. From $(\ref{ue26})$ and $(\ref{u43})$
\begin{align*}
&w_{0i}^1 \wedge \cdots \wedge w_{0i}^q \wedge \left(\frac{\partial w_{i|\mu}^\alpha}{\partial t} - P_{i_x|\mu-1}^\alpha - \partial v_{i|\mu-1}^\alpha + \sum_{\beta=1}^q a_{0 i_x}^{\alpha\beta} v_{i|\mu-1}^\beta   \right)\\
&=  M_{i_x|\mu-1}^\alpha - w_{0i}^1\wedge \cdots \wedge w_{0i}^q \wedge P_{i_x|\mu-1}^\alpha =0
\end{align*}
Since $w_{0i}^1,...,w_{0i}^q$ defines a regular foliation on $U_{i_x} \subset U_i-S$, we can find $b_{i_x|\mu}^{\alpha\beta}$ satisfying $(\ref{u55})_\mu$.

This completes the proof of Lemma \ref{ue23}.

\end{proof}

It remains to determine $\varphi_1, w_{i|1}^\alpha, h_{ij|1}^{\alpha\beta}$ and $v_{i|0}^\alpha, a_{i_x|1}^{\alpha\beta}, b_{i_x|0}^{\alpha\beta}$ satisfying $(\ref{u14})_1,(\ref{u15})_1, (\ref{ucc1})_0,(\ref{u16})_1- (\ref{u19})_1$. Let $\phi\in \mathbb{H}^1\left( M, \tilde{\mathcal{N}}_{\mathcal{F}_0}^{* \bullet} \right)$. Then $\phi$ can be represented by using the Dolbeault type bicomplex associated to $\tilde{\mathcal{N}}_{\mathcal{F}_0}^{*\bullet}$ as
\begin{align*}
\left(  \overline{ \kappa } , \pi    , \rho \right)\in \frac{A^{0,0}\left( M, \mathscr{H}om_{\mathcal{O}_M}\left( \mathcal{N}_{\mathcal{F}_0}^*, \Omega_M^1 \right) \right)}{ A^{0,0}\left( M , \mathscr{H}om_{\mathcal{O}_M}\left( \mathcal{N}_{\mathcal{F}_0}^*, \mathcal{N}_{\mathcal{F}_0}^* \right) \right)} \bigoplus A^{0,0}\left( M, \mathscr{H}om_{\mathcal{O}_M}\left( \mathcal{N}_{\mathcal{F}_0}^*, \mathcal{O}_M  \right) \right) \bigoplus A^{0,1}\left( M, \Theta_M \right)
\end{align*}
where $\kappa\in A^{0,0}\left( M, \mathscr{H}om_{\mathcal{O}_M}\left( \mathcal{N}_{\mathcal{F}_0}^*, \Omega_M^1 \right)\right)$ defined by $\kappa\left( w_{0i}^\alpha \right)=\kappa_i^\alpha$ on $U_i$ for $\kappa_i^\alpha \in \Gamma\left( U_i, \mathcal{A}^{0,0}\left( \Omega_M^1 \right)  \right)$ with $\kappa_{i\lambda}^\alpha = \sum_{\beta=1}^q h_{0ij}^{\alpha\beta} \kappa_j^\beta$ on $U_{ij}$ and $\pi\left( w_{0i}^\alpha \right)= \pi_i^\alpha \in \Gamma\left( U_i, \mathcal{A}^{0,0} \right)$ on $U_i$ with $\pi_i^\alpha= \sum_{\beta=1}^q h_{0ij}^{\alpha\beta} \pi_j^\beta$ on $U_{ij}$ such that there exists $\theta_i\in \Gamma\left( U_i, \mathcal{A}^{0,1}\left( \mathscr{H}om_{\mathcal{O}_M}\left( \mathcal{N}_{\mathcal{F}_0}^*, \mathcal{N}_{\mathcal{F}_0}^* \right) \right) \right)$ with $\theta_i\left( w_{0i}^\alpha \right)=\theta_i^\alpha= \sum_{\beta=1}^q \theta_i^{\alpha\beta} w_{0i}^\beta \in \Gamma\left(U_i, \mathcal{A}^{0,1}\left( \mathcal{N}_{\mathcal{F}_0}^* \right) \right)$ satisfying
\begin{align}
&-\bar{\partial} \rho = 0\\
\sum_{\beta=1}^q w_{0i}^1 \wedge \cdots \wedge \kappa_i^\alpha \wedge \cdots \wedge w_{0i}^q & \wedge  dw_{0i}^\alpha + w_{0i}^1 \wedge \cdots \wedge w_{0i}^q \wedge \partial \kappa_i^\alpha = 0 \label{ucc7} \\
\sum_{\beta=1}^q (-1)^{q-\beta} \pi_i^\beta w_{0i}^1 \wedge \cdots \wedge \widehat{w}_{0i}^\beta \wedge \cdots \wedge w_{0i}^q &\wedge dw_{0i}^\alpha + w_{0i}^1\wedge \cdots \wedge w_{0i}^q \wedge \partial \pi_i^\alpha - w_{0i}^1\wedge \cdots \wedge w_{0i}^q \wedge \kappa_i^\alpha=0 \label{ucc6}\\
-\bar{\partial}\kappa_i^\alpha + &\mathcal{L}_\rho\left( w_{0i}^\alpha \right)= \sum_{\beta=1}^q \theta_i^{\alpha\beta} w_{0i}^\beta \label{ucc3} \\
-\bar{\partial} \pi_i^\alpha &+i_\rho\left(w_{0i}^\alpha \right) = 0
\end{align}
We note that since $\bar{\partial} \theta_{i}^{\alpha\beta}=0$ from $(\ref{ucc3})$, there exists $\mathfrak{K}_i^{\alpha\beta}\in \Gamma\left(U_i, \mathcal{A}^{0,0} \right)$ such that $\bar{\partial} \mathfrak{K}_i^{\alpha\beta} = \theta_i^{\alpha\beta}$. We set
\begin{align}
\varphi_1 &= t\rho \label{us11}\\
w_i^{\alpha 1}&= w_{0i}^\alpha + t\left(\kappa_i^\alpha + \mathfrak{K}_i^\alpha \right) = w_{0i}^\alpha + t\left( \kappa_i^\alpha + \sum_{\beta=1}^q \mathfrak{K}_i^{\alpha\beta} w_{0i}^\beta \right)\\
 h_{ij}^{\alpha \beta 1} &= h_{0ij}^{\alpha\beta} + t\left( \sum_{\gamma=1}^q \mathfrak{K}_i^{\alpha \gamma} h_{0ij}^{\gamma \beta} - \sum_{\gamma=1}^q h_{0ij}^{\alpha\gamma} \mathfrak{K}_j^{\gamma \beta} \right)\\
v_{i|0}^\alpha &= \pi_i^\alpha \label{us12}
\end{align}
Then $(\ref{u14})_1, (\ref{u15})_1, (\ref{ucc1})_0, (\ref{u17})_1, (\ref{u18})_1$ and $(\ref{ucc5})_1$ are satisfied. $(\ref{u15})_1, (\ref{ucc1})_0$ and $(\ref{u17})_1$ implies $(\ref{u16})_1$. It remains to find $a_{i_x}^{\alpha\beta 1}$ and $b_{i_x|0}^{\alpha\beta}$. In fact, from $(\ref{ucc7})$ and $(\ref{ucc6})$, we have on $U_i$
\begin{align*}
\partial \kappa_i^\alpha - \sum_{\beta=1}^q a_{0i_x}^{\alpha\beta} \wedge \kappa_i^\beta = \sum_{\beta=1}^q c_{i_x}^{\alpha\beta} \wedge w_{0i}^\beta \,\,\,\,\,\,\textnormal{for some $c_{i_x}^{\alpha\beta}\in \Gamma\left(  U_i, \mathcal{A}^{0,0}\left( \Omega_M^1 \right) \right)$}
\end{align*}
\begin{align*}
\partial \pi_i^\alpha- \sum_{\beta=1}^q a_{0i_x}^{\alpha\beta} \pi_i^\beta - \kappa_i^\alpha = \sum_{\beta=1}^q l_{i_x}^{\alpha\beta} w_{0i}^\beta \,\,\,\,\,\,\,\textnormal{for some $l_{i_x}^{\alpha\beta}\in \Gamma\left( U_i, \mathcal{A}^{0,0}\right)$}
\end{align*}
Then we have
\begin{align*}
&\partial_{M\times B} \left( w_i^{\alpha 1} + v_{i|0}^\alpha dt \right) = d w_{0i}^\alpha + t \left( \partial  \kappa_{i\lambda}^\alpha  + \sum_{\beta=1}^q \partial \mathfrak{K}_{i\lambda}^{\alpha\beta} \wedge w_{0i}^\beta + \sum_{\beta=1}^q \mathfrak{K}_{i\lambda}^{\alpha\beta} d w_{0i}^\beta \right) +\left( - \kappa_i^\alpha - \sum_{\beta=1}^q \mathfrak{K}_i^{\alpha\beta} w_{0i}^\beta + \partial \pi_i^\alpha \right) \wedge dt\\
&\equiv_1 \sum_{\beta=1}^q \left( a_{0i_x}^{\alpha\beta} +t \left( c_{i_x }^{\alpha\beta}  + \partial \mathfrak{K}_{i}^{\alpha\beta} - \sum_{\eta=1}^q a_{0i_x}^{\alpha \eta} \mathfrak{K}_{i}^{\eta \beta}   + \sum_{\eta=1}^q \mathfrak{K}_{i}^{\alpha\eta} a_{0i_x}^{\eta \beta}  \right) + \left( \mathfrak{K}_i^{\alpha\beta}  - l_{i_x}^{\alpha\beta}  \right) dt \right) \wedge \left( w_{0i}^\beta + t \left( \kappa_{i}^\beta + \sum_{\gamma=1}^q \mathfrak{K}_{i}^{\beta \gamma} w_{0i}^\gamma \right)  + \pi_i^\beta dt\right),
\end{align*}
so that if we take
\begin{align*}
a_{i_x}^{\alpha\beta 1}  & : = a_{0i_x}^{\alpha\beta} +  t \left( b_{i_x\lambda }^{\alpha\beta}  + \partial \mathfrak{K}_{i}^{\alpha\beta} - \sum_{\eta=1}^q a_{0i_x}^{\alpha \eta} \mathfrak{K}_{i}^{\eta \beta}   + \sum_{\eta=1}^q \mathfrak{K}_{i}^{\alpha\eta} a_{0i_x}^{\eta \beta}  \right),\\
b_{i_x|0}^{\alpha\beta} & := \mathfrak{K}_i^{\alpha\beta} - l_{i_x}^{\alpha\beta},
\end{align*}
then $(\ref{u19})_1$ is satisfied. This completes the inductive construction of $\varphi, w_i^\alpha, h_{ij}^{\alpha\beta}$, and $ v_i^\alpha, a_{i_x}^{\alpha\beta}, b_{i_x}^{\alpha\beta}$ satisfying $(\ref{u11})-(\ref{u12})$.

\subsection{Proof of convergence}\

We will prove that $\varphi=\sum_{\mu=1}^\infty \varphi_\mu, w_i^\alpha= \sum_{\mu=0}^\infty w_{i|\mu}^\alpha, h_{ij}^{\alpha\beta}= \sum_{\mu=0}^\infty h_{ij|\mu}^{\alpha\beta}$ and $v_i^\alpha= \sum_{\mu=0}^\infty v_{i|\mu}^\alpha$ converge. Before proceeding to the discussion, we introduce another complex of sheaves which also controls unfoldings of $\left( M, \mathcal{N}_{\mathcal{F}_0}^* \right)$ but removes the quotient in the degree $1$. We shall define the following complex of sheaves
\begin{equation}
\begin{CD}
\cdots \\
@AG_3AA \\
\mathscr{H}om_{\mathcal{O}_M}\left( \mathcal{N}_{\mathcal{F}_0}^*, \tilde{\mathcal{S}}^3 \right) \bigoplus \mathscr{H}om_{\mathcal{O}_M}\left( \mathcal{N}_{\mathcal{F}_0}^*, \tilde{\mathcal{S}}^2  \right) \\
@AG_2AA \\
\mathscr{H}om_{\mathcal{O}_M}\left( \mathcal{N}_{\mathcal{F}_0}^*, \tilde{\mathcal{S}}^2 \right) \bigoplus \mathscr{H}om_{\mathcal{O}_M}\left( \mathcal{N}_{\mathcal{F}_0}^*, \tilde{\mathcal{S}}^1 \right) \\
@A\tilde{G}_1''AA \\
\mathscr{H}om_{\mathcal{O}_M}\left( \mathcal{N}_{\mathcal{F}_0}^* , \Omega_M^1 \right) \bigoplus \mathscr{H}om_{\mathcal{O}_M}\left( \mathcal{N}_{\mathcal{F}_0}^*, \mathcal{O}_M \right) \\
@AG_0''AA \\
\mathcal{E}_{\mathcal{N}_{\mathcal{F}_0}^*}
\end{CD}
\end{equation}
where $\mathcal{E}_{\mathcal{N}_{\mathcal{F}_0}^*}$ is the Atiyah extension associated to $\mathcal{N}_{\mathcal{F}_0}^*$ from Part II. We recall $(\ref{us2})-(\ref{us3})$ in the beginning of the proof of Theorem \ref{u56}. We define $\mathcal{E}_{\mathcal{N}_{\mathcal{F}_0}^*}$ locally on $U_i$ by $\mathcal{E}_{\mathcal{N}_{\mathcal{F}_0}^*}|_{U_i}\cong \mathscr{H}om_{\mathcal{O}_M}\left( \mathcal{N}_{\mathcal{F}_0}^* , \mathcal{N}_{\mathcal{F}_0}^* \right)|_{U_i} \bigoplus \Theta_M|_{U_i}$ such that $(\phi_i,X)\in \mathscr{H}om_{\mathcal{O}_X}\left( \mathcal{N}_{\mathcal{F}_0}^* , \mathcal{N}_{\mathcal{F}_0}^* \right)|_{U_i} \bigoplus \Theta_X |_{U_i}$ and $(\phi_j, X)\in \mathscr{H}om_{\mathcal{O}_X}\left( \mathcal{N}_{\mathcal{F}_0}^* , \mathcal{N}_{\mathcal{F}_0}^* \right)|_{U_j} \bigoplus \Theta_X|_{U_j}$ are identified on $U_{ij}$ if $\phi_i\left(w_{0i}^\alpha\right)=\phi_j\left(w_{0i}^\alpha\right)+\sum_{\beta=1}^q \left[X, h_{0ij}^{\alpha\beta}\right] w_{0j}^\beta$. Then we define $G_0':\mathcal{E}_{\mathcal{N}_{\mathcal{F}_0}^*}\to \mathscr{H}om_{\mathcal{O}_M}\left( \mathcal{N}_{\mathcal{F}_0}^*, \Omega_M^1 \right)$ $\bigoplus \mathscr{H}om_{\mathcal{O}_M}\left( \mathcal{N}_{\mathcal{F}_0}^*, \mathcal{O}_M \right)$ locally on $U_i$ by $\tilde{G}_0''\left(\left(\phi_i, X \right) \right)\left( w_{0i}^\alpha \right) =\left(\left(\mathcal{L}_X - \phi_i \right)\left( w_{0i}^\alpha \right) , i_{X}\left( w_{0i}^\alpha \right) \right)$ and then linearly extends to $\Gamma\left( U_i, \mathcal{N}_{\mathcal{F}_0}^* \right)$. We define $\tilde{G}_1''$ by the composition of the natural quotient map $\mathscr{H}om_{\mathcal{O}_M}\left( \mathcal{N}_{\mathcal{F}_0}^*, \Omega_M^1 \right) $ $\bigoplus \mathscr{H}om_{\mathcal{O}_M}\left( \mathcal{N}_{\mathcal{F}_0}^*, \mathcal{O}_M \right) \to \mathscr{H}om_{\mathcal{O}_M}\left( \mathcal{N}_{\mathcal{F}_0}^*, \frac{\Omega_M^1}{ \mathcal{N}_{\mathcal{F}_0}^* } \right) \bigoplus \mathscr{H}om_{\mathcal{O}_M}\left( \mathcal{N}_{\mathcal{F}_0}^*, \mathcal{O}_M \right)$ with $G_1$. Then we have 
\begin{align}\label{ucc10}
\mathbb{H}^i\left( M, \tilde{\mathcal{E}}_{\mathcal{N}_{\mathcal{F}_0}^*}^\bullet \right)  \cong \mathbb{H}^i\left( M, \tilde{\mathcal{N}}_{\mathcal{F}_0}^{*\bullet} \right),\,\,\,\,\,\,\,\, i\geq 0
\end{align}
(for the detail on the complex of sheaves $\tilde{\mathcal{E}}_{\mathcal{N}_{\mathcal{F}_0}^*}^\bullet$, see subsection \ref{us5}). We denote $\mathcal{A}^{0,p}\left( \mathcal{E}_{\mathcal{N}_{\mathcal{F}_0}^*} \right)$ be the sheaf of germs of $C^\infty(0,p)$-forms with coefficients in $\mathcal{E}_{\mathcal{N}_{\mathcal{F}_0}^*}^\bullet$ and denote by $A^{0,p}\left( M, \mathcal{E}_{\mathcal{N}_{\mathcal{F}_0}^*} \right)$ the global section of $\mathcal{A}^{0,p}\left( \mathcal{E}_{\mathcal{N}_{\mathcal{F}_0}^*} \right)$. Then we have the following Dolbeault type bicomplex associated to $\tilde{\mathcal{E}}_{\mathcal{N}_{\mathcal{F}_0}^*}^\bullet$:
{\Tiny{\begin{center}
$\begin{CD}
\cdots \\
@A\hat{G}_3AA \\
A^{0,0}\left(M-S, \left(\mathcal{N}_{\mathcal{F}_0}^* \right)^* \otimes \tilde{\mathcal{S}}^3\right)\bigoplus A^{0,0}\left( M-S, \left(\mathcal{N}_{\mathcal{F}_0}^*\right)^*\otimes \tilde{\mathcal{S}}^2  \right) @>-(-1)^q\bar{\partial}>> \cdots \\
@ A \hat{G}_2 AA @A \hat{G}_2 AA \\
A^{0,0}\left(M-S, \left(\mathcal{N}_{\mathcal{F}_0}^*\right)^* \otimes \tilde{\mathcal{S}}^2\right) \bigoplus A^{0,0}\left( M-S, \left(\mathcal{N}_{\mathcal{F}_0}^* \right)^*\otimes \tilde{\mathcal{S}}^1 \right) @>(-1)^q \bar{\partial} >> A^{0,1}\left(M-S, \left(\mathcal{N}_{\mathcal{F}_0}^*\right)^*\otimes \tilde{\mathcal{S}}^2 \right) \bigoplus A^{0,1}\left( M-S, \left( \mathcal{N}_{\mathcal{F}_0}^*\right)^*\otimes \tilde{\mathcal{S}}^1 \right) @>-(-1)^q\bar{\partial}>> \cdots\\
@A\hat{G}_1''  AA @A\hat{G}_1'' AA @A\hat{G}_1''AA \\
A^{0,0}\left(M, \left(\mathcal{N}_{\mathcal{F}_0}^*\right)^* \otimes \Omega_M^1 \right) \bigoplus A^{0,0}\left(M, \left( \mathcal{N}_{\mathcal{F}_0}^* \right)^* \right) @>-\bar{\partial}>> A^{0,1}\left(M, \left( \mathcal{N}_{\mathcal{F}_0}^* \right)^*\otimes \Omega_M^1 \right) \bigoplus A^{0,1}\left(M, \left(\mathcal{N}_{\mathcal{F}_0}^* \right)^* \right) @>\bar{\partial}>> \cdots\\
@A\hat{G}_0''AA @A \hat{G}_0'' AA @A\hat{G}_0''AA \\
A^{0,0}\left(M, \mathcal{E}_{\mathcal{N}_{\mathcal{F}_0}^*}\right) @>\bar{\partial}>> A^{0,1}\left( M,  \mathcal{E}_{\mathcal{N}_{\mathcal{F}_0}^*} \right) @>-\bar{\partial}>> A^{0,2}\left(M, \mathcal{E}_{\mathcal{N}_{\mathcal{F}_0}^*} \right) 
\end{CD}$
\end{center}}}
Let us denote the $i$-th cohomology group of the complex associated to the above bicomplex by $\textnormal{G}'^i$. Then $\textnormal{G}'^i\cong \textnormal{G}^i \left(i\geq 0 \right)$ $\textnormal{(\ref{us6})}$ and so $\textnormal{G}'^i\cong \mathbb{H}^i\left( M, \tilde{\mathcal{E}}_{\mathcal{N}_{\mathcal{F}_0}^*}^\bullet  \right)$ for $i=0,1,2$.

We reinterpret $(\ref{us1})$ in terms of the above bicomplex associated to $\tilde{\mathcal{E}}_{\mathcal{N}_{\mathcal{F}_0}^*}^\bullet$. From $(\ref{n50})$, we can show as in Part I
\begin{align}\label{us15}
\left( -\xi_\mu, \left\{ \bar{\partial} B_{i|\mu} \right\}\right) \in A^{0,2}\left( M, \mathcal{E}_{\mathcal{N}_{\mathcal{F}_0}^*} \right) \Longrightarrow \left( - \frac{\partial \xi_\mu}{\partial t}, \left\{ \bar{\partial}\left( \frac{\partial B_{i|\mu}}{\partial t} \right) \right\} \right) \in A^{0,2}\left( M, \mathcal{E}_{\mathcal{N}_{\mathcal{F}_0}^*} \right)
\end{align}

Then we see that from $(\ref{us1})$
{\Small{\begin{align}\label{us20}
\left(\frac{\partial I_\mu}{\partial t}, \pi_{\mu-1} \right)\oplus \left( \frac{\partial \phi_\mu}{\partial t}, \sigma_{\mu-1} \right) \oplus  &\left(- \frac{\partial \xi_\mu}{\partial t} , \left\{ \bar{\partial}\left(\frac{\partial B_{i|\mu}}{\partial t} \right)\right\} \right)\in A^{0,0}\left(M - S , \mathscr{H}om_{\mathcal{O}_M}\left( \mathcal{N}_{\mathcal{F}_0}^*, \tilde{\mathcal{S}}^2 \right) \right) \bigoplus A^{0,0}\left(M - S , \mathscr{H}om_{\mathcal{O}_M}\left( \mathcal{N}_{\mathcal{F}_0}^*, \tilde{\mathcal{S}}^1 \right) \right) \\
&\bigoplus A^{0,1}\left( M, \mathscr{H}om_{\mathcal{O}_M}\left( \mathcal{N}_{\mathcal{F}_0}^*, \Omega_M^1 \right) \right)\bigoplus A^{0,1}\left( M, \mathscr{H}om_{\mathcal{O}_M}\left( \mathcal{N}_{\mathcal{F}_0}^*, \mathcal{O}_M \right) \right) \bigoplus A^{0,2} \left(M, \mathcal{E}_{\mathcal{N}_{\mathcal{F}_0}^*} \right) \notag
\end{align}}}
defines a $2$-cocycle in the above Dolbeault type bicomplex associated to $\tilde{\mathcal{E}}_{\mathcal{N}_{\mathcal{F}_0}^*}^\bullet$.

We define the \"Holder norms $\left|- \right|_{k+\alpha}$ (k: an integer $\geq 2, 0<\alpha<a$)    on the sections of $\mathcal{A}^{0,q}\left( \Theta_M \right)$ and $\mathcal{A}^{0,q}\left(\mathscr{H}om_{\mathcal{O}_M}\left( \mathcal{N}_{\mathcal{F}_0}^* , \Omega_M^1 \right) \right)$ and  $\mathcal{A}^{0,q}\left( M, \mathcal{E}_{\mathcal{N}_{\mathcal{F}_0}^*} \right)$ and $A^{0,q}\left( \mathscr{H}om_{\mathcal{O}_M}\left( \mathcal{N}_{\mathcal{F}_0}^* , \bigwedge^{q+2} \Omega_M^1 \otimes \mathcal{L}_0 \right) \right)$ and introduce a harmonic theory as in Part I. We recall that
\begin{align}\label{us17}
\textnormal{We denote by $\mathfrak{d}$ the adjoint operator of $\bar{\partial}$, and $\widetilde{\square}= \mathfrak{d}\bar{\partial} + \bar{\partial} \mathfrak{d}$ and $G$ the Green's operator on $\Theta_M$.}
\end{align}
\begin{align}
\textnormal{We denote by $\mathfrak{d}''$ the adjoint operator of $\bar{\partial}$, and $\square'=\mathfrak{d}''\bar{\partial} + \bar{\partial} \mathfrak{d}''$ and $G''$ the Green's operator on $\mathscr{H}om_{\mathcal{O}_M}\left( \mathcal{N}_{\mathcal{F}_0}^*, \Omega_M^1 \right)$.}
\end{align}
\begin{align}\label{us16}
\textnormal{We denote by $\tilde{\mathfrak{d}}$ the adjoint operator of $\bar{\partial}$, and $\widetilde{\square}=\tilde{\mathfrak{d}}\bar{\partial} + \bar{\partial} \tilde{\mathfrak{d}}$ and $\tilde{G}$ the Green's operator on $\mathcal{E}_{\mathcal{N}_{\mathcal{F}_0}^*}$}
\end{align}

We also define the \"Holder norm on section of $\mathcal{A}^{0,q}\left(\mathscr{H}om_{\mathcal{O}_M}\left(\mathcal{N}_{\mathcal{F}_0}^*, \mathcal{O}_M \right)\right)$ in a similar way. We introduce a harmonic theory on $\mathscr{H}om_{\mathcal{O}_M}\left( \mathcal{N}_{\mathcal{F}_0}^*, \mathcal{O}_M \right)$. 
\begin{align}
\textnormal{We denote by $\mathfrak{d}^*$ the adjoint operator of $\bar{\partial}$, and $\square^*=\mathfrak{d}^*\bar{\partial} + \bar{\partial} \mathfrak{d}^*$ and $G^*$ the Green's operator on $\mathscr{H}om_{\mathcal{O}_M}\left( \mathcal{N}_{\mathcal{F}_0}^*, \mathcal{O}_M \right)$}
\end{align}

Consider a formal power series $\phi= \phi(t)=\phi(t)=\sum \phi_v t^{v}$ of one variable $t$ with coefficients in $A^{0,q}\left( M, \mathcal{G} \right)$ or in $\Gamma\left( U_i, \mathcal{A}^{0,q}\left(\mathcal{G} \right)  \right)$ or in $\Gamma\left(U_{ij}, \mathcal{A}^{0,q}\left( \mathcal{G} \right) \right)$ where $\mathcal{G}=\Theta_M$ or $\Omega_M^1$ or $\mathcal{O}_X$, and a power series $a(t)=\sum a_{v} t^v , a_v \geq 0$. 
\begin{align}\label{d30}
\textnormal{We indicate by $\left|\phi \right|_{k+\alpha} \ll a(t)$ or $\left| \phi \right|_{k+\alpha}^{U_i} \ll a(t)$ that $\left| \phi_{v} \right|_{k+\alpha} \leq a_{v}$   }
\end{align}
In the following we do not indicate explicitly the domain $U_i$ or $U_{ij}$ if no confusion is possible.

We set 
\begin{align}
A(t)=\frac{b}{16c}\sum_{\mu=1}^\infty \frac{1}{\mu^2} c^\mu t^\mu\,\,\,\,\,\,\,\,\,\textnormal{with}\,\,\,b>0\,\,\,\textnormal{and}\,\,\, c>0
\end{align}
Then we have
\begin{align}\label{utt5}
A(t)^2 \ll \frac{b}{c} A(t) \Longrightarrow 2 A(t)\cdot \frac{\partial A(t)}{\partial t} \ll \frac{b}{c} \frac{\partial A(t)}{\partial t} \Longrightarrow A(t)\cdot \frac{\partial A(t)}{\partial t} \ll \frac{b}{c} \frac{\partial A(t)}{\partial t}
\end{align}

With this preparation, we will show that for a fixed integer $k\geq 2$ and $0< \alpha<1$, the inductive construction of $\varphi, w_i^\alpha , h_{ij}^{\alpha\beta}$ and $v_i^\alpha$ in the previous subsection can be carried out in such a way that
\begin{align*}
\left|\varphi \right|_{k+\alpha} \ll A(t) \\
\left| w_i^{\alpha }- w_{0i}^{\alpha} \right|_{k+ \alpha}  \ll  A(t) \\
\left| h_{ij}^{\alpha\beta }- h_{0ij}^{\alpha\beta} \right|_{k+ 1+ \alpha} \ll A(t) \\
\left| v_i^\alpha \right|_{ k +  \alpha} \ll \frac{\partial A(t)}{\partial t}
\end{align*}
Then it suffices to prove that for $\mu=1,2,3,\cdots$,
\begin{align}
\left| \varphi^{\mu}\right|_{k+\alpha}  \ll A(t) \label{us7} \\
\left| w_i^{\alpha\mu}- w_{0i}^{\alpha} \right|_{k+ \alpha}  \ll  A(t) \label{us9} \\
\left| h_{ij}^{\alpha\beta \mu}- h_{0ij}^{\alpha\beta} \right|_{k+ 1+ \alpha} \ll A(t) \label{us10} \\
\left| v_i^{\alpha(\mu-1)} \right|_{k+\alpha} \ll \frac{\partial A(t)}{\partial t} \label{us8}
\end{align}
for some proper choice of constants $c>b>0$. We prove $(\ref{us7})_\mu, (\ref{us9})_\mu, (\ref{us10})_\mu$ and $(\ref{us8})_{\mu-1}$ by induction on $\mu$. For $\mu=1$, we have $(\ref{us11})-(\ref{us12})$, and the linear term of $A(t)=\frac{b}{16}t$ and the constant term of $\frac{\partial A(t)}{\partial t}= \frac{b}{16}$. Therefore $(\ref{us7})_1, (\ref{us9})_1, (\ref{us10})_1$ and $(\ref{us8})_0$ holds if $b$ is sufficiently small.

Now assume that $(\ref{us7})_{\mu-1}, (\ref{us9})_{\mu-1}, (\ref{us10})_{\mu-1}$ and $(\ref{us8})_{\mu-2}$ are satisfied. Then we will derive $(\ref{us7})_\mu, (\ref{us9})_\mu, (\ref{us10})_\mu$ and $(\ref{us8})_{\mu-1}$. In the following $Y_1,Y_2,Y_3, \cdots$ will denote constants which depend only on $k,\alpha, M, \mathcal{N}_{\mathcal{F}_0}^*$. We note that from $(\ref{us7})_{\mu-1}$ and $(\ref{us9})_{\mu-1}$,  we have
\begin{align}\label{utt6}
\left| \frac{\partial  \varphi^{\mu-1}}{\partial t} \right|_{k+\alpha}  \ll  \frac{\partial A(t) }{\partial t},\,\,\,\,\,\,\,\,\,\,\left| \frac{\partial w_i^{\alpha (\mu-1)}}{\partial t} \right|_{k+\alpha} \ll \frac{\partial A(t)}{\partial t}
\end{align}

From $(\ref{us15})$, we have $\bar{\partial}\left(- \frac{\partial \xi_\mu}{\partial t}, \left\{ \bar{\partial}\left( \frac{\partial B_{i|\mu}}{\partial t} \right) \right\} \right)=0$. Then we take $\left( \varphi_{\mu-1}^\sharp , \left\{ \lambda_{i|\mu-1}^\sharp \right\} \right):=\tilde{\mathfrak{d}}\tilde{G}\left( - \frac{\partial \xi_\mu}{\partial t}, \left\{ \bar{\partial}\left( \frac{\partial B_{i|\mu}}{\partial t} \right) \right\}\right) \in A^{0,1}\left( M, \mathcal{E}_{\mathcal{N}_{\mathcal{F}_0}^*} \right)$ from $(\ref{us16})$, where $\varphi_{\mu-1}^\sharp \in A^{0,1}\left(M, \Theta_M \right)$ and $\lambda_{i|\mu-1}^\sharp \in \Gamma\left( U_i, \mathcal{A}^{0,1}\left( \mathscr{H}om_{\mathcal{O}_M}\left( \mathcal{N}_{\mathcal{F}_0}^*, \mathcal{N}_{\mathcal{F}_0}^* \right) \right) \right)$ defined by $\lambda_{i|\mu-1}^\sharp :\Gamma\left( U_i, \mathcal{N}_{\mathcal{F}_0}^* \right) \to \Gamma\left( U_i, \mathcal{N}_{\mathcal{F}_0}^* \right), w_{0i}^\alpha \mapsto \sum_{\xi=1}^q \lambda_{i|\mu}^{\sharp \alpha \xi} w_{0i}^\xi$. Then we have
\begin{align} \label{us18}
\left| \left( \varphi_{\mu-1}^\sharp, \left\{ \lambda_{i|\mu-1}^\sharp \right\} \right) \right|_{k+\alpha} \ll Y_1 \left| \left( - \frac{\partial \xi_\mu}{\partial t} , \left\{  \bar{\partial} \left( \frac{\partial B_{i|\mu}}{\partial t}\right\} \right) \right)    \right|_{k-1 +\alpha}
\end{align}
as in a similar way with \cite{Kod05} Lemma 5.7 p.276.

From $(\ref{us17})$, we have
\begin{align}
\varphi_{\mu-1}^\sharp= \bold{H} \varphi_{\mu-1}^\sharp + \left( \bar{\partial} \mathfrak{d} + \mathfrak{d} \bar{\partial} \right)G \varphi_{\mu-1}^\sharp
\end{align}
It follows that
\begin{align*}
\bar{\partial} \left( \bold{H} \varphi_{\mu-1}^\sharp + \mathfrak{d} \bar{\partial} G \varphi_{\mu-1}^\sharp  \right) & = - \xi_\mu\\
\mathfrak{d}\left( \bold{H} \varphi_{\mu-1}^\sharp + \mathfrak{d} \bar{\partial} G \varphi_{\mu-1}^\sharp \right) & =0
\end{align*}
On the other hand, we define $q_{ij|\mu-1}^\sharp \in \Gamma\left( U_{ij}, \mathcal{A}^{0,0}\left( \mathscr{H}om_{\mathcal{O}_M}\left( \mathcal{N}_{\mathcal{F}_0}^*, \mathcal{N}_{\mathcal{F}_0}^* \right) \right) \right)$ by
\begin{align*}
q_{ij|\mu-1}^\sharp : \Gamma\left( U_{ij}, \mathcal{N}_{\mathcal{F}_0}^*  \right) &\to \Gamma\left( U_{ij}, \mathcal{A}^{0,0}\left( \mathcal{N}_{\mathcal{F}_0}^* \right) \right) \\
 w_{0i}^\alpha &\mapsto \sum_{\beta=1}^q \left[ \mathfrak{d} G \varphi_{\mu-1}^\sharp , h_{0ij}^{\alpha\beta} \right] w_{0j}^\beta
\end{align*}
and linearly extends to $\Gamma\left( U_{ij}, \mathcal{N}_{\mathcal{F}_0}^* \right)$. Then $\left\{ q_{ij|\mu-1}^\sharp \right\}\in C^1\left( \mathcal{U}, \mathcal{A}^{0,0}\left( \mathscr{H}om_{\mathcal{O}_M}\left( \mathcal{N}_{\mathcal{F}_0}^*, \mathcal{N}_{\mathcal{F}_0}^* \right) \right) \right)$ defines a $1$-cocycle. Let $\left\{ \rho_i(z) \right\}$ be a partition of unity subordinate to the covering $\mathcal{U}=\left\{ U_i \right\}$. By setting $q_{i|\mu-1}^\sharp:=\sum_k \rho_k q_{ik|\mu-1}^\sharp$, we see that $\left\{ q_{i|\mu-1}^\sharp \right\} \in C^0\left(\mathcal{U}, \mathcal{A}^{0,0}\left( \mathscr{H}om_{\mathcal{O}_M}\left(\mathcal{N}_{\mathcal{F}_0}^*, \mathcal{N}_{\mathcal{F}_0}^* \right) \right) \right)$ such that $q_{i|\mu-1}^\sharp- q_{j|\mu-1}^\sharp =  q_{ij|\mu-1}^\sharp$ where
\begin{align*}
q_{i|\mu-1}^\sharp : \Gamma\left( U_{i}, \mathcal{N}_{\mathcal{F}_0}^*  \right) &\to \Gamma\left( U_{i}, \mathcal{A}^{0,0}\left( \mathcal{N}_{\mathcal{F}_0}^* \right) \right) \\
 w_{0i}^\alpha &\mapsto \sum_{\gamma=1}^q q_{i|\mu-1}^{\sharp \alpha \gamma} w_{0i}^\gamma, \,\,\,\,\,\, q_{i|\mu-1}^{\sharp\alpha \gamma} \in \Gamma\left( U_i, \mathcal{A}^{0,0} \right)
\end{align*}
Then we see that
\begin{align}\label{utt45}
\left( \varphi_{\mu-1}'' , \left\{ \lambda_{i|\mu-1}'' \right\} \right) \in A^{0,1}\left( M, \mathcal{E}_{\mathcal{N}_{\mathcal{F}_0}^* } \right),\,\,\,\,\,\,\,\, \mathfrak{d} \varphi_{\mu-1}'' = 0,\,\,\,\,\,\,\,\,\,\bar{\partial}\left( \varphi_{\mu-1}'', \left\{ \lambda_{i|\mu-1}'' \right\} \right) = \left( - \frac{\partial \xi_\mu}{\partial t} , \left\{ \bar{\partial}\left( \frac{\partial B_{i|\mu}}{\partial t} \right) \right\} \right)
\end{align}
where $\varphi_{\mu-1}'':= \bold{H} \varphi_{\mu-1}^\sharp + \mathfrak{d} \bar{\partial} G \varphi_{\mu-1}^\sharp$ and $\lambda_{i|\mu-1}''\in \Gamma\left( U_i ,  \mathcal{A}^{0,1}\left( \mathscr{H}om_{\mathcal{O}_M}\left( \mathcal{N}_{\mathcal{F}_0}^*, \mathcal{N}_{\mathcal{F}_0}^*  \right) \right) \right)$ defined by
\begin{align*}
\lambda_{i|\mu-1}'' : \Gamma\left( U_i, \mathcal{N}_{\mathcal{F}_0}^* \right) &\to \Gamma\left( U_i, \mathcal{A}^{0,1}\left( \mathcal{N}_{\mathcal{F}_0}^* \right) \right) \\
 w_{0i}^\alpha &\mapsto \sum_{\beta=1}^q \left( \lambda_{i|\mu}^{\sharp \alpha\beta} - \bar{\partial} q_{i|\mu}^{\sharp \alpha \beta} \right) w_{0i}^\beta
\end{align*}

Then we have from $(\ref{us18})$
\begin{align}\label{utt22}
\left| \varphi_{\mu-1}'' \right|_{k+\alpha},\,\,\,\,\,\,\,\left| \lambda_{i|\mu-1}'' \right|_{k+\alpha} \ll Y_2\left|\left( - \frac{\partial \xi_\mu}{\partial t} , \left\{ \bar{\partial}\left( \frac{\partial B_{i|\mu}}{\partial t} \right) \right\}    \right) \right|_{k-1+\alpha}
\end{align}

From $(\ref{us20})$, we see that
{\small{\begin{align*}
0 \in &A^{0,2}\left(M, \mathcal{E}_{\mathcal{N}_{\mathcal{F}_0}^*}\right) \\
\phi_{\mu-1}'=&\left\{ w_{0i}^\alpha \mapsto  \frac{\partial A_{i|\mu}^\alpha}{\partial t} - \bar{\partial} \left(  \frac{\partial C_{i|\mu}^\alpha}{\partial t}  \right) + \sum_{\gamma=1}^q \frac{\partial B_{i|\mu}^{\alpha \gamma}}{\partial t} w_{0i}^\gamma +\mathcal{L}_{\varphi_{\mu-1}''}\left(w_{0i}^\alpha\right) - \sum_{\xi=1}^q \lambda_{i|\mu-1}''^{\alpha \xi} w_{0i}^\xi \right\} \in A^{0,1}\left(\mathscr{H}om_{\mathcal{O}_M}\left( \mathcal{N}_{\mathcal{F}_0}^*, \Omega_M^1 \right) \right)\\
\sigma_{\mu-1}'=&\left\{  w_{0i}^\alpha \mapsto \psi_{i|\mu-1}^\alpha - \bar{\partial} N_{i|\mu-1}^\alpha + i_{\varphi_{\mu-1}''} (w_{0i}^\alpha)\right\} \in A^{0,1}\left( M, \mathscr{H}om_{\mathcal{O}_M}\left( \mathcal{N}_\mathcal{F}^*, \mathcal{O}_M \right) \right)\\
\frac{\partial I_\mu}{\partial t}=&\left\{ w_{0i}^\alpha \mapsto  \frac{\partial E_{i|\mu}^\alpha}{\partial t} + \sum_{\beta=1}^q w_{0i}^1\wedge \cdots \wedge \frac{\partial C_{i|\mu}^\beta}{\partial t} \wedge \cdots \wedge w_{0i}^q \wedge dw_{0i}^\alpha + w_{0i}^1\wedge \cdots \wedge w_{0i}^q \wedge  d\left( \frac{\partial C_{i|\mu}^\alpha}{\partial t} \right) \right\} \\
&\in A^{0,0}\left( M, \mathscr{H}om_{\mathcal{O}_M}\left(\mathcal{N}_{\mathcal{F}_0}^*,  \bigwedge^{q+2} \Omega_M^1 \otimes \mathcal{L}_0 \right) \right)\\
\pi_{\mu-1}= &\left\{ w_{0i}^\alpha \mapsto M_{i|\mu-1}^\alpha + \sum_{\beta=1}^q (-1)^{q-\beta} N_{i|\mu-1}^\beta w_{0i}^1\wedge \cdots \wedge  \widehat{w_{0i}^\beta} \wedge \cdots \wedge w_{0i}^q \wedge dw_{0i}^\alpha   + w_{0i}^1\wedge \cdots \wedge w_{0i}^q \wedge \partial N_{i|\mu-1}^\alpha - w_{0i}^1\wedge \cdots \wedge w_{0i}^q \wedge \frac{\partial C_{i|\mu}^\alpha}{\partial t} \right\}\\
& \in A^{0,0}\left( M, \mathscr{H}om_{\mathcal{O}_M}\left(\mathcal{N}_{\mathcal{F}_0}^*,  \bigwedge^q  \Omega_M^1\otimes \mathcal{L}_0  \right) \right)
\end{align*}}}
defines a $2$-cocycle in the above Dolbeault type bicomplex associated to $\tilde{\mathcal{E}}_{\mathcal{N}_{\mathcal{F}_0}^*}^\bullet$. Since $\mathbb{H}^2\left( M, \tilde{\mathcal{N}}_{\mathcal{F}_0}^{*\bullet} \right)=\mathbb{H}^2\left( M, \tilde{\mathcal{E}}_{\mathcal{N}_{\mathcal{F}_0}^*}^\bullet \right)=0$ by the assumption of Theorem \ref{u56}, there exist $\left(\chi_{\mu-1}, \left\{ \delta_i^{\chi_{\mu-1}} \right\} \right)\in A^{0,1}\left( M, \mathcal{E}_{\mathcal{N}_{\mathcal{F}_0}^*} \right)$ and $\Sigma_{\mu-1}\in A^{0,0}\left( M, \mathscr{H}om_{\mathcal{O}_M}\left( \mathcal{N}_{\mathcal{F}_0}^* , \Omega_M^1 \right) \right)$ and $\Pi_{\mu-1} \in A^{0,0}\left( M, \mathscr{H}om_{\mathcal{O}_M} \left( \mathcal{N}_{\mathcal{F}_0}^*, \mathcal{O}_M \right) \right)$ such that $-\bar{\partial}\left( \chi_{\mu-1}, \left\{ \delta_i^{\chi_{\mu-1}} \right\} \right)=0$ and $-\bar{\partial} \left( \Sigma_{\mu-1}, \Pi_{\mu-1} \right)+ \hat{G}_0''\left(\chi_{\mu-1}, \left\{ \delta_i^{\chi_{\mu-1}} \right\} \right)=\left( \phi_{\mu-1}', \sigma_{\mu-1}'  ,\right)$ and $\hat{G}_1''\left( \Sigma_{\mu-1}, \Pi_{\mu-1} \right)=\left( \frac{\partial I_\mu}{\partial t}, \pi_{\mu-1} \right)$. By using the following Lemma, we will choose appropriate elements $\left(\chi_{\mu-1}, \left\{ \delta_i^{\chi_{\mu-1}} \right\} \right)$ and $\Sigma_{\mu-1}$ and $\Pi_{\mu-1}$ in a way that the resulting $\varphi, w_i^\alpha, v_i^\alpha$ and $h_{ij}^{\alpha\beta}$ converge.

We define $G_1'^\sharp:\mathscr{H}om_{\mathcal{O}_M}\left( \mathcal{N}_{\mathcal{F}_0}^*, \Omega_M^1 \right) \bigoplus \mathscr{H}om_{\mathcal{O}_M}\left( \mathcal{N}_{\mathcal{F}_0}^*, \mathcal{O}_M \right)\to \mathscr{H}om_{\mathcal{O}_M}\left( \mathcal{N}_{\mathcal{F}_0}^*, \bigwedge^{q+1} \Omega_M^1\otimes \mathcal{L}_0 \right)$ in the following way: $\phi\in \Gamma\left(U_i, \mathscr{H}om_{\mathcal{O}_M}\left( \mathcal{N}_{\mathcal{F}_0}^*, \Omega_M^1 \right) \right), \psi \in \Gamma\left( U_i, \mathscr{H}om_{\mathcal{O}_M}\left( \mathcal{N}_{\mathcal{F}_0}^*,\mathcal{O}_M \right)  \right)$, we define
\begin{align*}
G_1'^\sharp \left(\phi, \psi \right): \Gamma\left( U_i, \mathcal{N}_{\mathcal{F}_0}^* \right) &\to \Gamma\left( U_i, \bigwedge^{q+1} \Omega_M^1 \otimes \mathcal{L}_0 \right) \\
 w_{0i}^\alpha \mapsto \sum_{\beta=1}^q (-1)^{q-\beta} \psi\left(  w_{0i}^\beta \right) &w_{0i}^1 \wedge \cdots \wedge \widehat{w}_{0i}^\beta \wedge \cdots \wedge w_{0i}^q \wedge dw_{0i}^\alpha + w_{0i}^1 \wedge \cdots \wedge w_{0i}^q \wedge d \psi\left( w_{0i}^\alpha \right) - w_{0i}^1 \wedge \cdots \wedge w_{0i}^q \wedge \phi\left( w_{0i}^\alpha \right)
\end{align*}
This induces 
\begin{align*}
G_1'^\sharp : A^{0,p}\left( M, \mathscr{H}om_{\mathcal{O}_M} \left( \mathcal{N}_{\mathcal{F}_0}^*, \Omega_M^1 \right) \right) \bigoplus A^{0,p}\left( M, \mathscr{H}om_{\mathcal{O}_M} \left( \mathcal{N}_{\mathcal{F}_0}^*, \mathcal{O}_M \right) \right) \to A^{0,p}\left( M, \mathscr{H}om_{\mathcal{O}_M} \left( \mathcal{N}_{\mathcal{F}_0}^*, \bigwedge^{q+1} \Omega_M^1 \otimes \mathcal{L}_0 \right) \right)
\end{align*}
Then we have
\begin{align*}
\hat{G}_1''\left(\phi, \psi \right) = \left( E_1^\sharp (\phi), G_1'^\sharp (\phi, \psi) \right) |_{M-S}
\end{align*}
For the definition of $E_1^\sharp$, see Part I. Let us define $G_1^\sharp$ by
\begin{align*}
G_1^\sharp \left(\phi, \psi \right) =\left( E_1^\sharp\left( \phi \right), G_1'^\sharp\left( \phi, \psi \right) \right)
\end{align*}

Then we have
\begin{lemma}\label{utt1}
Suppose that there exist $\left( \varphi, \left\{ \delta_i \right\} \right)\in A^{0,1}\left( M, \mathcal{E}_{\mathcal{N}_{F_0}^*} \right)$, where $\varphi\in A^{0,1}\left( M, \Theta_M \right)$,  and $\delta_i\in \Gamma\left(U_i, \mathcal{A}^{0,1}\left( \mathscr{H}om_{\mathcal{O}_M}\left( \mathcal{N}_{\mathcal{F}_0}^*, \mathcal{N}_{\mathcal{F}_0}^* \right) \right) \right)$ with $\delta_i\left(w_{0i}^\alpha \right)= \delta_j\left( w_{0i}^\alpha \right) + \sum_{\beta=1}^q \left[ \varphi, h_{0ij}^{\alpha\beta} \right] w_{0j}^\beta$ for $\alpha=1,...,q$, and $V\in A^{0,1}\left( M, \mathscr{H}om_{\mathcal{O}_M}\left( \mathcal{N}_{\mathcal{F}_0}^*, \Omega_M^1 \right) \right), W\in A^{0,1}\left( M, \mathscr{H}om_{\mathcal{O}_M}\left(  \mathcal{N}_{\mathcal{F}_0}^*, \mathcal{O}_M \right) \right)$, $I\in A^{0,0}\left(M, \mathscr{H}om_{\mathcal{O}_M}\left( \mathcal{N}_{\mathcal{F}_0}^*, \bigwedge^{q+2} \Omega_M^1 \otimes \mathcal{L}_0 \right) \right)$ and $ J \in A^{0,0}\left( M, \mathscr{H}om_{\mathcal{O}_M}\left( \mathcal{N}_{\mathcal{F}_0}^*, \bigwedge^{q+1} \Omega_M^1\otimes \mathcal{L}_0 \right) \right)$ such that $I |_{M-S}\in A^{0,0}\left(M-S, \mathscr{H}om_{\mathcal{O}_M}\left( \mathcal{N}_{\mathcal{F}_0}^*, \tilde{\mathcal{S}}^2 \right) \right)$ and $J  |_{M-S}\in A^{0,0}\left(M-S , \mathscr{H}om_{\mathcal{O}_M}\left( \mathcal{N}_{\mathcal{F}_0}^*, \tilde{\mathcal{S}}^1 \right) \right)$, and $\left( 0 , \left( \left( V, W\right)+ \hat{G}'_0 \left(\varphi, \left\{ \delta_i\right\}  \right) \right) , I|_{M-S}, J |_{M-S} \right)$ defines a $2$-cocycle in the Dolbeault type bicomplex associated to $\tilde{\mathcal{E}}_{\mathcal{N}_{\mathcal{F}_0}^*}^\bullet$, and moreover $(-1)^q\bar{\partial}\left( I, \sigma \right) + G_1^\sharp\left( V, W\right)=0$. Then we can find $\left(\chi, \left\{ \delta_i^\chi \right\} \right) \in A^{0,1}\left( M, \mathcal{E}_{\mathcal{N}_{\mathcal{F}_0}^*}\right)$, and $\Sigma\in A^{0,0}\left( M, \mathscr{H}om_{\mathcal{O}_M}\left( M, \mathscr{H}om_{\mathcal{O}_M}\left( \mathcal{N}_{\mathcal{F}_0}^*, \Omega_M^1 \right)\right) \right)$ and $\Pi\in A^{0,0}\left( M, \mathscr{H}om_{\mathcal{O}_M}\left( \mathcal{N}_{\mathcal{F}_0}^*, \mathcal{O}_M \right) \right)$ in such a way that
\begin{align}
&\widetilde{\square}\left( \chi, \left\{ \delta_i^\chi \right\} \right)=0 \label{uut26}\\
-\bar{\partial} \Sigma\left(w_{0i}^\alpha \right) +\mathcal{L}_{\chi}\left( w_{0i}^\alpha\right)- & \delta_i^{\chi} \left( w_{0i}^\alpha \right) = V\left( w_{0i}^\alpha \right) + \mathcal{L}_{\varphi}-  \delta_i \left( w_{0i}^\alpha \right) \\
- \bar{\partial} \Pi\left( w_{0i}^\alpha \right)& + i_{\chi}\left( w_{0i}^\alpha \right)= W \left( w_{0i}^\alpha \right) + i_{\varphi}\left( w_{0i}^\alpha \right)\\
\delta_i^\chi\left( w_{0i}^\alpha \right)& = \delta_j^\chi\left( w_{0i}^\alpha \right) + \sum_{\beta=1}^q \left[ \chi, h_{0ij}^{\alpha\beta} \right] w_{0j}^\beta   \\
&G_1^\sharp\left(  \Sigma, \Pi \right)= \left(  I, \sigma \right) \label{uut27}
\end{align}
\begin{align*}
\left|\left( \chi, \left\{ \delta^\chi \right\} \right) \right|_{k+\alpha} &\ll Y\left( \left|\left( \varphi, \left\{ \delta_i \right\}\right) \right|_{k+\alpha} + \left| V\right|_{k-1+\alpha} +\left| W \right|_{k-1+\alpha} + \left| I \right|_{k-1+\alpha} + \left| J \right|_{k-1+\alpha} \right)  \\
\left| \Sigma \right|_{k+\alpha}&\ll Y \left( \left|\left( \varphi, \left\{ \delta_i \right\}\right) \right|_{k+\alpha} + \left| V\right|_{k-1+\alpha} +\left| W \right|_{k-1+\alpha} + \left| I \right|_{k-1+\alpha} + \left| J \right|_{k-1+\alpha} \right) \\
\left| \Pi \right|_{k+\alpha} & \ll Y\left( \left|\left( \varphi, \left\{ \delta_i \right\}\right) \right|_{k+\alpha} + \left| V\right|_{k-1+\alpha} +\left| W \right|_{k-1+\alpha} + \left| I \right|_{k-1+\alpha} + \left| J \right|_{k-1+\alpha} \right) 
\end{align*}
where $Y$ is a constant which is independent of $\left( \varphi, \left\{\delta_i \right\}\right), V, W , I , J$.

\end{lemma}

\begin{proof}
For any quintuple $\zeta=\left( \left(\varphi, \left\{\delta_i\right\}  \right)  ,V, W, I, J   \right)$ as above, let
\begin{align*}
||\zeta || &= \left| \left(\varphi  , \left\{ \delta_i \right\}\right)\right|_{k+\alpha} + \left|V \right|_{k-1+\alpha} + \left| W \right|_{k-1 +\alpha} + \left| I \right|_{k-1+\alpha}  + \left| J \right|_{k-1+\alpha} \\
\iota(\zeta)&= \inf \left( \left|\left(\chi, \left\{\delta_i^\chi \right\} \right)\right|_{k+\alpha}  + \left|\Sigma  \right|_{k+\alpha}  + \left| \Pi \right|_{k+\alpha} \right)
\end{align*}
where $\inf$ is taken with respect to all solutions $\left( \left(\chi, \left\{\eta_i^\chi \right\}\right), \Sigma , \Pi \right)$ of the equalities $(\ref{uut26})-(\ref{uut27})$. It suffices to prove the existence of $Y$ such that
\begin{align*}
\iota\left( \zeta \right)\leq Y ||\zeta ||\,\,\,\,\,\,\,\,\,\textnormal{for all triple $\zeta$}
\end{align*}
which can be shown in a similar way as in Part I.
\end{proof}

In Lemma \ref{utt1}, we set
{\Small{\begin{align} \label{utt21}
\left( \varphi, \left\{ \delta_i \right\} \right)&:= \left(  \varphi_{\mu-1}'', \left\{\lambda_{i|\mu-1}'' \right\} \right)\\
V&:= \left\{ w_{0i}^\alpha \mapsto \frac{\partial A_{i|\mu}^\alpha}{\partial t} - \bar{\partial}\left(\frac{\partial C_{i|\mu}^\alpha}{\partial t} \right) + \sum_{\gamma=1}^q \frac{\partial B_{i|\mu}^{\alpha \gamma}}{\partial t} \right\} \notag \\
W &:= \left\{ w_{0i}^\alpha \mapsto \psi_{i|\mu-1}^\alpha - \bar{\partial} N_{i|\mu-1}^\alpha  \right\}  \notag \\
I &:= \left\{ w_{0i}^\alpha \mapsto \frac{\partial E_{i|\mu}^\alpha}{\partial t}  + \sum_{\beta=1}^q w_{0i}^1 \wedge \cdots \wedge \frac{\partial C_{i|\mu}^\beta}{\partial t} \wedge \cdots \wedge w_{0i}^q \wedge dw_{0i}^\alpha + w_{0i}^1 \wedge \cdots \wedge w_{0i}^q \wedge  \partial \left( \frac{\partial C_{i|\mu}^\alpha}{\partial t}\right)  \right\} \notag \\
J& =\left\{ w_{0i}^\alpha \mapsto M_{i|\mu-1}^\alpha + \sum_{\beta=1}^q(-1)^{q-\beta} N_{i|\mu-1}^\beta w_{0i}^1 \wedge \cdots \wedge \widehat{w}_{0i}^\beta \wedge \cdots \wedge w_{0i}^q \wedge dw_{0i}^\alpha + w_{0i}^1 \wedge \cdots \wedge w_{0i}^q \wedge \partial N_{i|\mu-1}^\alpha - w_{0i}^1 \wedge \cdots \wedge w_{0i}^q \wedge \frac{\partial C_{i|\mu}^\alpha}{\partial t}  \right\} \notag
\end{align}}}

We will estimate $\left( \varphi_{\mu-1}'', \left\{ \lambda_{i|\mu-1}'' \right\} \right), V, W, I$ and $J$. First we will estimate $(\ref{u20})-(\ref{utt2})$ and $(\ref{utt3})$. By induction hypothesis $(\ref{us7})_{\mu-1}-(\ref{us10})_{\mu-1}$ and $(\ref{us8})_{\mu-2}$, we have as in Part I,
\begin{align} \label{utt23}
\left| \xi_\mu \right|_{k-1+\alpha} \ll Y_3\frac{b}{c} A(t) \Longrightarrow \left|\frac{\partial \xi_\mu}{\partial t} \right|_{k-1+\alpha} \ll Y_3 \frac{b}{c} \frac{\partial A(t)}{\partial t}
\end{align}
and we have as in Part I
\begin{align} \label{utt25}
\left| A_{i|\mu}^\alpha \right|_{k-1+\alpha} \ll Y_4 \frac{b}{c} A(t) \Longrightarrow \left| \frac{\partial A_{i|\mu}^\alpha}{\partial t} \right|_{k-1+\alpha} \ll Y_4 \frac{b}{c} \frac{\partial A(t)}{\partial t}
\end{align}
and we have from $(\ref{utt6})$ and $(\ref{utt5})$
\begin{align}
&-\psi_{i|\mu-1}^\alpha \equiv_{\mu-1} \bar{\partial} v_i^{\alpha(\mu-2)} - \left[ \varphi^{\mu-1}, v_i^{\alpha(\mu-2)}  \right] - i_{\frac{\partial \varphi^{\mu-1}}{\partial t}}\left( w_i^{\alpha(\mu-1)} \right) \equiv_{\mu-1}\left[ \left[ \varphi^{\mu-1}, v_i^{\alpha(\mu-2)} \right] - i_{\frac{\partial \varphi^{\mu-1}}{\partial t}}\left( w_i^{\alpha(\mu-1)} - w_{0i}^\alpha \right)  \right]_{\mu-1} \notag \\
&\Longrightarrow \left| \psi_{i|\mu-1}^\alpha \right|_{k-1+\alpha} \leq  Y_5 \left| \varphi^{\mu-1}  \right|_{k-1+\alpha} \left| v_i^{\alpha(\mu-2)}\right|_{k+\alpha} + Y_6\left| \frac{\partial \varphi^{\mu-1}}{\partial t} \right|_{k+\alpha} \left| w_i^{\alpha( \mu-1)}- w_{0i}^\alpha\right|_{k +\alpha} \ll Y_7 A(t)\cdot \frac{\partial A(t)}{\partial t} \ll Y_7 \frac{b}{c} \frac{\partial A(t)}{\partial t} \label{utt27}
\end{align}
and we have as in Part I
\begin{align} \label{utt17}
\left| B_{ij|\mu}^{\alpha \eta} \right|_{k+\alpha}\ll Y_8 \frac{b}{c} A(t) \Longrightarrow \left| \frac{\partial B_{ij|\mu}^{\alpha \eta}}{\partial t} \right|_{k+\alpha} \ll Y_8 \frac{b}{c} \frac{\partial A(t) }{\partial t}
\end{align}
and we have as in Part I
\begin{align} \label{utt19}
\left| C_{ij|\mu}^\alpha \right|_{k+\alpha} \ll Y_9 \frac{b}{c} A(t) \Longrightarrow \left| \frac{\partial C_{ij|\mu}^\alpha}{\partial t} \right|_{k+\alpha} \ll Y_9 \frac{b}{c} \frac{\partial A(t)}{\partial t}
\end{align}
and we have from $(\ref{utt5})$
\begin{align}
&N_{ij|\mu-1}^\alpha \equiv_{\mu-1} v_i^{\alpha(\mu-2)} - \sum_{\beta=1}^q h_{ij}^{\alpha\beta(\mu-1)} v_j^{\beta(\mu-2)}\equiv_{\mu - 1} \left[ - \sum_{\beta=1}^q \left( h_{ij}^{\alpha\beta(\mu-1)} - h_{0ij}^{\alpha\beta} \right)   v_j^{\beta(\mu-2)}   \right]_{\mu-1} \notag \\
&\Longrightarrow \left| N_{ij|\mu-1}^\alpha \right|_{k+\alpha } \leq \sum_{\beta=1}^q \left| h_{ij}^{\alpha\beta(\mu-1)} - h_{0ij}^{\alpha\beta}  \right|_{k+\alpha} \left| v_j^{\beta(\mu-2)}\right|_{k+\alpha} \ll Y_{10} A(t)\frac{\partial A(t)}{\partial t} \ll Y_{10}\frac{b}{c} \frac{\partial A(t)}{\partial t} \label{utt20}
\end{align}
and we have as in Part I
\begin{align}\label{utt16}
\left| D_{ijk|\mu}^{\alpha\beta} \right|_{k+1+\alpha} \ll Y_{11} \frac{b}{c} A(t) \Longrightarrow  \left| \frac{\partial D_{ijk|\mu}^{\alpha\beta}}{\partial t} \right|_{k+1+\alpha} \ll Y_{11}\frac{b}{c} \frac{\partial A(t)}{\partial t}
\end{align}
and we have as in Part I
\begin{align} \label{utt29}
\left| E_{i|\mu}^\alpha \right|_{k-1+\alpha} \ll Y_{12} \frac{b}{c} A(t) \Longrightarrow \left| \frac{\partial E_{i|\mu}^\alpha}{\partial t} \right|_{k-1+\alpha} \ll Y_{12} \frac{b}{c} \frac{\partial A(t)}{\partial t}
\end{align}
It remains to estimate
\begin{align}\label{utt7}
M_{i|\mu-1}^\alpha \equiv_{\mu-1} & \sum_{\beta=1}^q (-1)^{q-\beta} w_i^{1(\mu-1)} \wedge \cdots \wedge v_i^{\beta(\mu-2)}\wedge \cdots \wedge w_i^{q(\mu-2)} \wedge \partial w_i^{\alpha(\mu-1)} + w_i^{1(\mu-1)}\wedge \cdots \wedge w_i^{q(\mu-1)}\wedge \partial v_i^{\alpha(\mu-2)} \\
   & - w_i^{1(\mu-1)} \wedge \cdots \wedge w_i^{q(\mu-1)} \wedge \frac{\partial w_i^{\alpha(\mu-1)}}{\partial t} \notag
\end{align}
We estimate the first term of $(\ref{utt7})$. We note that
\begin{align*}
&M_1^\alpha := \left[ \sum_{\beta=1}^q (-1)^{q-\beta} w_i^{1(\mu-1)} \wedge \cdots \wedge v_i^{\beta(\mu-2)}\wedge \cdots \wedge w_i^{q(\mu-2)} \wedge \partial w_i^{\alpha(\mu-1)}  \right]_{\mu-1}\\
&\equiv_{\mu-1} \sum_{\beta=1}^q (-1)^{q-\beta} \left(w_i^{1(\mu-1)}- w_{0i}^1+ w_{0i}^1 \right) \wedge \cdots \wedge v_i^{\beta(\mu-2)}\wedge \cdots \wedge \left( w_i^{q(\mu-2)} - w_{0i}^q + w_{0i}^q \right) \wedge \partial \left( w_i^{\alpha(\mu-1)}  - w_{0i}^\alpha + w_{0i}^\alpha \right)\\
&\equiv_{\mu-1} \sum_{\beta=1}^q (-1)^{q-\beta} \left(w_i^{1(\mu-1)} - w_{0i}^1 \right)\wedge \cdots \wedge v_i^{\beta(\mu-2)} \wedge\cdots \wedge \left(w_i^{q(\mu-1)} - w_{0i}^q \right) \wedge \partial \left( w_i^{\alpha(\mu-1)} - w_{0i}^\alpha \right) \\
&+ \cdots \\
&+ \sum_{ 1\leq \beta\ne \gamma \leq q} (-1)^{q-\beta} w_{0i}^1 \wedge \cdots \wedge v_i^{\beta(\mu-2)}\wedge \cdots \wedge \left(w_i^{\gamma(\mu-1)} - w_{0i}^\gamma \right) \wedge \cdots \wedge w_{0i}^q \wedge dw_{0i}^\alpha \\
&+ \sum_{\beta=1}^q (-1)^{q-\beta} w_{0i}^1 \wedge \cdots \wedge v_i^{ \beta(\mu-2)} \wedge \cdots \wedge w_{0i}^q \wedge \partial \left( w_i^{\alpha(\mu-1)} - w_{0i}^\alpha \right) 
\end{align*}
Then we have from $(\ref{utt5})$
\begin{align} \label{utt8}
\left| M_1^\alpha \right|_{k-1+\alpha} \ll Y_{13} \left( \frac{b}{c} \right)^{q-1} A(t) \cdot \frac{\partial A(t)}{\partial t} + Y_{14}  \left( \frac{b}{c} \right)^{q-2} A(t)\cdot \frac{ \partial A(t)}{\partial t} + \cdots + Y_{15} A(t)\cdot \frac{\partial A(t)}{\partial t} \ll Y_{16} A(t) \cdot \frac{ \partial A(t)}{\partial t}
\end{align}
We estimate the second term of $(\ref{utt7})$. We note that
\begin{align*}
&M_2^\alpha:= \left[ w_i^{1(\mu-1)} \wedge \cdots \wedge w_i^{q(\mu-1)} \wedge \partial v_i^{\alpha(\mu-2)} \right]_{\mu-1} \\
&\equiv_{\mu-1} \left(w_i^{1(\mu-1)} - w_{0i}^1 + w_{0i}^1 \right) \wedge \cdots \wedge \left(w_i^{q(\mu-1)} - w_{0i}^q + w_{0i}^q \right) \wedge \partial v_i^{\alpha(\mu-2)}\\
&\equiv_{\mu-1} \left(w_i^{1(\mu-1)}- w_{0i}^1 \right) \wedge \cdots \wedge \left(w_i^{q(\mu-1)} - w_{0i}^q \right) \wedge \partial v_i^{\alpha(\mu-2)}\\
& + \cdots\\
&+ \sum_{\beta=1}^q w_{0i}^1 \wedge \cdots \wedge \left( w_{i}^{\beta(\mu-1)} - w_{0i}^\beta \right) \wedge \cdots \wedge w_{0i}^q \wedge \partial v_i^{\alpha(\mu-2)}  
\end{align*}
Then we have from $(\ref{utt5})$
\begin{align} \label{utt9}
\left| M_2^\alpha \right|_{k-1+\alpha} \ll Y_{17} \left( \frac{b}{c} \right)^{q-1} A(t) \cdot \frac{ \partial A(t)}{\partial t} + \cdots + Y_{18} A(t) \cdot \frac{\partial A(t)}{\partial t} \ll Y_{19} \frac{b}{c} \frac{\partial A(t)}{\partial t}
\end{align}
We estimate the third term of $(\ref{utt7})$. We note that
\begin{align*}
&M_3^\alpha : = \left[w_i^{1(\mu-1)} \wedge \cdots \wedge w_i^{q(\mu-1)} \wedge \frac{\partial w_i^{\alpha(\mu-1)}}{\partial t} \right]_{\mu-1}\\
&\equiv_{\mu-1} \left(w_i^{1(\mu-1)}- w_{0i}^1 + w_{0i}^1\right) \wedge \cdots \wedge \left(w_i^{q(\mu-1)}- w_{0i}^q + w_{0i}^q \right)\wedge \frac{\partial w_i^{\alpha(\mu-1)}}{\partial t}\\
&\equiv_{\mu-1} \left(w_i^{1(\mu-1)}- w_{0i}^1 \right)\wedge \cdots \wedge \left(w_i^{q(\mu-1)} -w_{0i}^q \right) \wedge \frac{\partial  w_i^{\alpha(\mu-1)} }{\partial t}\\
& + \cdots \\
&+ \sum_{\beta=1}^q w_{0i}^1 \wedge \cdots \wedge \left( w_{0i}^{\beta(\mu-1)} - w_{0i}^\beta   \right) \wedge \cdots \wedge w_{0i}^q \wedge \frac{ \partial w_i^{\alpha(\mu-1)} }{\partial t} 
\end{align*}
Then we have from $(\ref{utt6})$ and $(\ref{utt5})$
\begin{align}\label{utt10}
\left| M_3^\alpha \right|_{k-1+\alpha} \ll Y_{20}\left( \frac{b}{c} \right)^{q-1} A(t) \cdot \frac{\partial A(t)}{\partial t} + \cdots + Y_{21} A(t) \cdot \frac{\partial A(t)}{\partial t} \ll Y_{22} \frac{b}{c} \frac{\partial A(t)}{\partial t}
\end{align}
Then from $(\ref{utt8}),(\ref{utt9})$ and $(\ref{utt10})$ we have
\begin{align} \label{utt30}
\left| M_{i|\mu-1}^\alpha \right|_{k-1+\alpha} \ll Y_{23} \frac{b}{c} \frac{\partial A(t)}{\partial t}
\end{align}

We recall that from $(\ref{ne11})$ and $(\ref{n50})$ and $(\ref{n52})$ and $(\ref{ue17})$
\begin{align}
 \sum_{\beta,\eta=1}^q h_{0ij}^{\alpha\beta} D_{jk|\mu}^{\beta \eta} h_{0ji}^{\eta\xi}-D_{ik|\mu}^{\alpha \xi}+ D_{ij|\mu}^{\alpha \xi} &=D_{ijk|\mu}^{\alpha\xi} \label{utt12}\\
 \sum_{\eta,\gamma=1}^q h_{0ij}^{\alpha \eta} B_{j|\mu}^{\eta \gamma}h_{0ji}^{\gamma \xi} - B_{i|\mu}^{\alpha \xi} &=B_{ij|\mu}^{\alpha \xi}-\bar{\partial} D_{ij|\mu}^{\alpha \xi} \label{utt13} \\
 \sum_{\beta=1}^q h_{0ij}^{\alpha\beta} C_{j|\mu}^\beta - C_{i|\mu}^\alpha &= C_{ij|\mu}^\alpha - \sum_{\beta=1}^q  D_{ij|\mu}^{\alpha\beta} w_{0i}^\beta \label{utt14} \\
 \sum_{\beta=1}^q h_{0ij}^{\alpha\beta} N_{j|\mu-1}^\beta - N_{i|\mu-1}^\alpha &= N_{ij|\mu-1}^\alpha \label{utt15}
\end{align}
By using a partition of unity subordinate to the covering $\left\{ U_i \right\}$, we may choose $D_{jk|\mu}$ from $D_{ijk|\mu}$, and choose $B_{j|\mu}$ from $B_{ij|\mu}$ and $\bar{\partial} D_{ij|\mu}$, and choose $C_{i|\mu}$ from $C_{ij|\mu}$ and $D_{ij|\mu}$, and choose $N_{i|\mu-1}$ from $N_{ij|\mu-1}$. Hence we may assume that from $(\ref{utt12})$ and $(\ref{utt16})$
\begin{align}\label{utt18}
\left| \frac{\partial D_{ij|\mu}}{\partial t} \right|_{k+1+\alpha} \ll Y_{23}\frac{b}{c} \frac{\partial A(t)}{\partial t}
\end{align}
and from $(\ref{utt13})$ and $(\ref{utt18})$ and $(\ref{utt17})$
\begin{align} \label{utt24}
\left| \frac{\partial B_{i|\mu}}{\partial t}  \right|_{k+\alpha} \ll Y_{24} \frac{b}{c} \frac{\partial A(t)}{\partial t}
\end{align}
and from $(\ref{utt14})$ and $(\ref{utt18})$ and $(\ref{utt19})$
\begin{align} \label{utt26}
\left| \frac{\partial C_{i|\mu}}{\partial t}  \right|_{k+\alpha} \ll Y_{25} \frac{b}{c} \frac{\partial A(t)}{\partial t}
\end{align}
and from $(\ref{utt15})$ and $(\ref{utt20})$
\begin{align} \label{utt28}
\left| N_{i|\mu-1}\right|_{k+\alpha} \ll Y_{26} \frac{b}{c} \frac{\partial A(t)}{\partial t}
\end{align}

With this preparation, now we estimate $(\ref{utt21})$ and then apply Lemma \ref{utt1}. Let us estimate $\left( \varphi_{\mu-1}'', \left\{ \lambda_{i|\mu-1}'' \right\} \right)$. From $(\ref{utt22})$ and $(\ref{utt23})$ and $(\ref{utt24})$, we have
\begin{align} \label{utt41}
\left| \left( \varphi_{\mu-1}'', \left\{ \lambda_{i|\mu-1}'' \right\} \right)  \right|_{k+\alpha} \ll Y_{27} \frac{b}{c} \frac{\partial A(t)}{\partial t}
\end{align}
We estimate $V= \left\{ w_{0i}^\alpha \mapsto \frac{\partial A_{i|\mu}^\alpha}{\partial t} - \bar{\partial}\left(\frac{\partial C_{i|\mu}^\alpha}{\partial t} \right) + \sum_{\gamma=1}^q \frac{\partial B_{i|\mu}^{\alpha \gamma}}{\partial t} \right\} \in A^{0,1}\left( M, \mathscr{H}om_{\mathcal{O}_M}\left( \mathcal{N}_{\mathcal{F}_0}^*, \Omega_M^1 \right) \right)$. From $(\ref{utt25})$ and $(\ref{utt26})$ and $(\ref{utt24})$
\begin{align*}
\left| V \right|_{k-1+\alpha} \ll Y_{28} \frac{b}{c} \frac{\partial A(t) }{\partial t}
\end{align*}
We estimate $W= \left\{ w_{0i}^\alpha \mapsto \psi_{i|\mu-1}^\alpha - \bar{\partial} N_{i|\mu-1}^\alpha  \right\} \in A^{0,1}\left( M, \mathscr{H}om_{\mathcal{O}_M} \left( \mathcal{N}_{\mathcal{F}_0}^*, \mathcal{O}_M \right) \right) $. From $(\ref{utt27})$ and $(\ref{utt28})$
\begin{align*}
\left| W \right|_{k-1 +\alpha} \ll  Y_{29} \frac{b}{c} \frac{\partial A(t)}{\partial t}
\end{align*}
We estimate 
{\Small{\begin{align*}
I = \left\{ w_{0i}^\alpha \mapsto \frac{\partial E_{i|\mu}^\alpha}{\partial t}  + \sum_{\beta=1}^q w_{0i}^1 \wedge \cdots \wedge \frac{\partial C_{i|\mu}^\beta}{\partial t} \wedge \cdots \wedge w_{0i}^q \wedge dw_{0i}^\alpha + w_{0i}^1 \wedge \cdots \wedge w_{0i}^q \wedge  \partial \left( \frac{\partial C_{i|\mu}^\alpha}{\partial t}\right)  \right\}\in A^{0,0}\left( M, \mathscr{H}om_{\mathcal{O}_M} \left( \mathcal{N}_{\mathcal{F}_0}^*, \bigwedge^{q+2} \Omega_M^1 \otimes \mathcal{L}_0 \right) \right)
\end{align*}}}
From $(\ref{utt29})$ and $(\ref{utt26})$, we have
\begin{align*}
\left| I \right|_{k-1+\alpha} \ll Y_{30} \frac{b}{c} \frac{\partial A(t)}{\partial t}
\end{align*}

We estimate 
{\Small{\begin{align*}
J =&\left\{ w_{0i}^\alpha \mapsto M_{i|\mu-1}^\alpha + \sum_{\beta=1}^q(-1)^{q-\beta} N_{i|\mu-1}^\beta w_{0i}^1 \wedge \cdots \wedge \widehat{w}_{0i}^\beta \wedge \cdots \wedge w_{0i}^q \wedge dw_{0i}^\alpha + w_{0i}^1 \wedge \cdots \wedge w_{0i}^q \wedge \partial N_{i|\mu-1}^\alpha - w_{0i}^1 \wedge \cdots \wedge w_{0i}^q \wedge \frac{\partial C_{i|\mu}^\alpha}{\partial t}  \right\}\\
&\in A^{0,0}\left( M, \mathscr{H}om_{\mathcal{O}_M} \left( \mathcal{N}_{\mathcal{F}_0}^*, \bigwedge^{q+1} \Omega_M^1 \otimes \mathcal{L}_0 \right) \right)
\end{align*}}}
From $(\ref{utt30})$ and $(\ref{utt28})$ and $(\ref{utt26})$ we have
\begin{align*}
\left| J \right|_{k-1+\alpha} \ll Y_{31} \frac{b}{c} \frac{\partial A(t) }{ \partial t}
\end{align*}

Then by Lemma \ref{utt1}, there exists $\left( \tilde{\chi}_{\mu-1}, \left\{ \delta_i^{\tilde{\chi}_{\mu-1}} \right\} \right)\in A^{0,1}\left( M, \mathcal{E}_{\mathcal{N}_{\mathcal{F}_0}^*} \right)$ and $\tilde{\Sigma}_{\mu-1} \in A^{0,0}\left( M, \mathscr{H}om_{\mathcal{O}_M}\left( \mathcal{N}_{\mathcal{F}_0}^*, \Omega_M^1 \right) \right)$ and $\tilde{\Pi}_{\mu-1}\in A^{0,0}\left( M,  \mathscr{H}om_{ \mathcal{O}_M} \left( \mathcal{N}_{\mathcal{F}_0}^*, \mathcal{O}_M    \right) \right)$ such that
\begin{align}
\bar{\partial} \tilde{\chi}_{\mu-1} =0\,\,\,\,\,\,\,\bar{\partial} \delta_i^{\tilde{\chi}_{\mu-1}} =0 &\,\,\,\,\,\left( \tilde{\square}\left( \tilde{\chi}_{\mu-1} , \left\{ \delta_i^{\tilde{\chi}_{\mu-1}} \right\} \right)=0 \right)\\
-\bar{\partial} \tilde{\Sigma}_{\mu-1}\left( w_{0i}^\alpha \right) +\mathcal{L}_{\tilde{\chi}_{\mu-1}}\left( w_{0i}^\alpha \right) - \delta_i^{\tilde{\chi}_{\mu-1}}\left( w_{0i}^\alpha \right)=& \frac{\partial  A_{i|\mu}^\alpha}{\partial t} - \bar{\partial}\left( \frac{\partial C_{i|\mu}^\alpha}{\partial t} \right) + \sum_{\gamma=1}^q \frac{\partial B_{i|\mu}^{\alpha \gamma}}{\partial t} w_{0i}^\gamma + \mathcal{L}_{\varphi_{\mu-1}''}\left( w_{0i}^\alpha \right) - \lambda_{i|\mu-1}''\left( w_{0i}^\alpha \right)  \label{utt31}\\
\sum_{\beta=1}^q w_{0i}^1 \wedge \cdots \wedge \tilde{\Sigma}_{\mu-1}\left( w_{0i}^\beta \right) \wedge \cdots \wedge w_{0i}^q \wedge dw_{0i}^\alpha&+ w_{0i}^1\wedge \cdots \wedge w_{0i}^q \wedge \partial\left(\tilde{\Sigma}_{\mu-1}\left( w_{0i}^\alpha \right) \right)\\
= \frac{\partial E_{i|\mu}^\alpha}{\partial t} + \sum_{\beta=1}^q w_{0i}^1 \wedge \cdots \wedge &\frac{\partial C_{i|\mu}^\beta}{\partial t} \wedge \cdots \wedge w_{0i}^q \wedge dw_{0i}^\alpha + w_{0i}^1 \wedge \cdots \wedge w_{0i}^q \wedge \partial \left( \frac{\partial C_{i|\mu}^\alpha}{\partial t} \right) \notag \\
-\bar{\partial}\tilde{\Pi}_{\mu-1} \left( w_{0i}^\alpha \right) + i_{\tilde{\chi}_{\mu-1}}\left( w_{0i}^\alpha \right) &= \psi_{i|\mu-1}^\alpha - \bar{\partial} N_{i|\mu-1}^\alpha  + i_{\varphi_{\mu-1}''}\left( w_{0i}^\alpha \right) \label{utt35} \\
\sum_{\beta=1}^q (-1)^{q-\beta} \tilde{\Pi}_{\mu-1}\left(w_{0i}^\beta \right) w_{0i}^1 \wedge \cdots \wedge \widehat{w}_{0i}^\beta \wedge \cdots \wedge w_{0i}^q \wedge dw_{0i}^\alpha&+ w_{0i}^1\wedge \cdots \wedge w_{0i}^q \wedge \partial \left(\tilde{\Pi}_{\mu-1} \left( w_{0i}^\alpha \right) \right) - w_{0i}^1 \wedge \cdots \wedge w_{0i}^q \wedge \tilde{\Sigma}_{\mu-1} \left( w_{0i}^\alpha \right)\\
=M_{i|\mu-1}^\alpha + \sum_{\beta=1}^q (-1)^{q-\beta} N_{i|\mu-1}^{q-\beta} w_{0i}^1 \wedge \cdots \wedge \widehat{w}_{0i}^\beta  \wedge \cdots \wedge &w_{0i}^q \wedge dw_{0i}^\alpha + w_{0i}^1 \wedge \cdots \wedge w_{0i}^q \wedge \partial N_{i|\mu-1}^\alpha - w_{0i}^1 \wedge \cdots \wedge w_{0i}^q \wedge \frac{\partial C_{i|\mu}^\alpha}{\partial t} \notag \\
\left| \left( \tilde{\chi}_{\mu-1}, \left\{ \delta_i^{\tilde{\chi}_{\mu-1}} \right\} \right)  \right|_{k+\alpha},\,\,\,\,\,\left| \tilde{\Sigma}_{\mu-1} \right|_{k+\alpha},\,\,\,\,\, & \left| \tilde{\Pi}_{\mu-1} \right|_{k+\alpha}\,\,\,\,\, \ll Y_{32} \frac{b}{c} \frac{\partial A(t)}{\partial t} \label{utt32}
\end{align} 

We note that if we set
\begin{align}
\chi_\mu &:= \frac{1}{\mu}t \tilde{\chi}_{\mu-1} \Longrightarrow  \frac{\partial \chi_\mu}{\partial t}=\tilde{\chi}_{\mu-1} \label{utt46} \\
\delta_i^{\chi_\mu} &:= \frac{1}{\mu} t \delta_i^{\tilde{\chi}_{\mu-1}} \Longrightarrow \frac{\partial \delta_i^{\chi_\mu}}{\partial t} = \delta_i^{\tilde{\chi}_{\mu-1}} \label{utt47}\\
\Sigma_\mu &:= \frac{1}{\mu} t \tilde{\Sigma}_{\mu-1} \Longrightarrow \frac{\partial \Sigma_\mu}{\partial t} = \tilde{\Sigma}_{\mu-1} \label{utt48}
\end{align}
Then from $(\ref{utt32})$ we have
\begin{align} \label{utt53}
\left| \left(\chi_\mu, \left\{ \delta_i^{\chi_\mu} \right\} \right) \right|_{k+\alpha},\,\,\,\,\,\,\,\left| \Sigma_\mu \right|_{k+\alpha}\,\,\,\,\,\ll Y_{32} \frac{b}{c} A(t)
\end{align}

In view of the equality
\begin{align} \label{utt36}
\tilde{\chi}_{\mu-1}= \bold{H} \tilde{\chi}_{\mu-1} + \left(\bar{\partial} \mathfrak{d} + \mathfrak{d}\bar{\partial} \right) G \tilde{\chi}_{\mu-1}
\end{align}
we have from $(\ref{utt31})$
\begin{align} \label{utt34}
-\bar{\partial}\left(\tilde{\Sigma}_{i|\mu-1}^\alpha - \mathcal{L}_{\mathfrak{d}G\tilde{\chi}_{\mu-1}}\left( w_{0i}^\alpha  \right) \right)+ \mathcal{L}_{\bold{H}\tilde{\chi}_{\mu-1}-  \varphi_{\mu-1}''}\left( w_{0i}^\alpha \right) - \left( \delta_i^{\tilde{\chi}_{\mu-1}}  -\lambda_{i|\mu-1}'' \right)\left( w_{0i}^\alpha \right) = \frac{\partial A_{i|\mu}^\alpha}{\partial t} - \bar{\partial}\left( \frac{\partial C_{i|\mu}^\alpha}{\partial t} \right) + \sum_{\gamma=1}^q \frac{\partial B_{i|\mu}^{\alpha \gamma}}{\partial t} w_{0i}^\gamma
\end{align}

On the other hand, we note that $\left\{\tilde{U}_{ij|\mu-1} \right\}$ is a $1$-cocycle with coefficient in $\mathcal{A}^{0,1}\left( \mathscr{H}om_{\mathcal{O}_M}\left( \mathcal{N}_{\mathcal{F}_0}^*, \mathcal{N}_{\mathcal{F}_0}^*   \right) \right)$ where
\begin{align*}
\tilde{U}_{ij|\mu-1}: \Gamma\left( U_{ij}, \mathcal{N}_{\mathcal{F}_0}^* \right)& \to \Gamma\left( U_{ij}, \mathcal{A}^{0,0}\left( \mathcal{N}_{\mathcal{F}_0}^* \right) \right)\\
w_{0i}^\alpha &\mapsto \sum_{\beta=1}^q \left[ \mathfrak{d} G \tilde{\chi}_{\mu-1},  h_{0ij}^{\alpha\beta} \right] w_{0j}^\beta
\end{align*}

Let $\left\{ \rho_i(z)\right\}$ be a partition of unity subordinate  to the covering $\mathcal{U}=\left\{ U_i \right\}$. By setting $\tilde{U}_{i|\mu-1}= \sum_k \rho_k(z)\tilde{U}_{ik|\mu-1}$, we see that $\left\{ \tilde{U}_{i|\mu-1} \right\}\in C^0\left( \mathcal{U}, \mathscr{H}om_{\mathcal{O}_M}\left( \mathcal{N}_{\mathcal{F}_0}^*, \mathcal{N}_{\mathcal{F}_0}^* \right)\right)$ such that $\tilde{U}_{i|\mu-1}- U_{j|\mu-1}= U_{ij|\mu-1}$ where
\begin{align*}
\tilde{U}_{i|\mu-1}:\Gamma\left( U_i, \mathcal{N}_{\mathcal{F}_0}^* \right) &\to \Gamma\left( U_i, \mathcal{A}^{0,0}\left( \mathcal{N}_{\mathcal{F}_0}^* \right) \right) \\
      w_{0i}^\alpha &\mapsto \sum_{\gamma=1}^q \tilde{U}_{i|\mu-1}^{\alpha \gamma} w_{0i}^\gamma,\,\,\,\,\,\,\,\,\tilde{U}_{i|\mu-1}^{\alpha \gamma} \in \Gamma\left( U_i, \mathcal{A}^{0,0} \right)
\end{align*}
and we have
\begin{align}\label{utt33}
\sum_{\beta, \gamma=1}^q h_{0ij}^{\alpha\beta} \tilde{U}_{j|\mu-1}^{\beta \gamma} w_{0j}^\gamma - \sum_{\gamma=1}^q \tilde{U}_{i|\mu-1}^{\alpha \gamma} w_{0i}^\gamma = - \sum_{\beta=1}^q \left[ \mathfrak{d}G \tilde{\chi}_{\mu-1}, h_{0ij}^{\alpha\beta} \right] w_{0j}^\beta
\end{align}
Then we have from $(\ref{utt32})$ 
\begin{align}\label{utt42}
\left| \tilde{U}_{i|\mu-1} \right|_{k+1+\alpha} \ll Y_{33}\left| \mathfrak{d} G \tilde{\chi}_{\mu-1} \right|_{k+1+\alpha} \ll Y_{34} \left| \tilde{\chi}_{\mu-1} \right|_{k+\alpha} \ll Y_{35} \frac{b}{c} \frac{\partial A(t)}{\partial t}
\end{align}

On the other hand, $(\ref{utt33})$ implies that
\begin{align*}
\tilde{R}_i:\Gamma\left( U_i, \mathcal{N}_{\mathcal{F}_0}^* \right) &\to \Gamma\left( U_i, \mathcal{A}^{0,0}\left( \Omega_M^1 \right) \right) \\
 w_{0i}^\alpha &\mapsto \mathcal{L}_{ \mathfrak{d} G \tilde{\chi}_{\mu-1}}\left( w_{0i}^\alpha \right) - \tilde{U}_{i|\mu-1}\left( w_{0i}^\alpha \right)
\end{align*}
defines a global section $\tilde{R}=\left\{\tilde{R}_i \right\} \in A^{0,0}\left( M, \mathscr{H}om_{\mathcal{O}_M}\left( \mathcal{N}_{\mathcal{F}_0}^*, \Omega_M^1 \right)\right)$. Then $\tilde{\Sigma}_{\mu-1} - R\in A^{0,0}\left( M, \mathscr{H}om_{\mathcal{O}_M}\left( \mathcal{N}_{\mathcal{F}_0}^*, \Omega_M^1 \right) \right)$. Then we have from $(\ref{utt34})$

{\small{\begin{align}\label{utt51}
&- \bar{\partial}\left(\tilde{\Sigma}_{i|\mu-1}^\alpha - \mathcal{L}_{\mathfrak{d}G \tilde{\chi}_{\mu-1}}\left(w_{0i}^\alpha \right) + \sum_{\xi=1}^q \tilde{U}_{i|\mu-1}^{\alpha \xi} w_{0i}^\xi \right) + \mathcal{L}_{\bold{H}\tilde{\chi}_{\mu-1} -  \varphi_{\mu-1}''}\left(w_{0i}^\alpha \right) - \left(   \delta_i^{\tilde{\chi}_{\mu-1}} - \lambda_{i|\mu-1}'' -\bar{\partial} \tilde{U}_{i|\mu-1}\right)\left( w_{0i}^\alpha \right) \\
&= \frac{\partial A_{i|\mu}^\alpha}{\partial t} - \bar{\partial}\left( \frac{\partial C_{i|\mu}^\alpha  }{\partial t} \right) + \sum_{\gamma=1}^q \frac{\partial B_{i|\mu}^{\alpha \gamma} }{\partial t} w_{0i}^\gamma  \notag
\end{align}}}
We note that
\begin{align}\label{utt50}
&\sum_{\beta, \xi=1}^q h_{0ij}^{\alpha\beta}\left( \delta_j^{\tilde{\chi}_{\mu-1} \beta \xi} - \lambda_{j|\mu-1}''^{\beta \xi} - \bar{\partial} \tilde{U}_{j|\mu-1}^{\beta \xi}  \right) w_{0j}^\xi - \sum_{\xi=1}^q \left( \delta_i^{\tilde{\chi}_{\mu-1}\alpha\xi}  - \lambda_{i|\mu-1}''^{\alpha \xi} - \bar{\partial} \tilde{U}_{i|\mu-1}^{\alpha \xi}\right) w_{0i}^\xi\\
&= - \sum_{\beta=1}^q \left[ \tilde{\chi}_{\mu-1} - \varphi_{\mu-1}'' - \bar{\partial} \mathfrak{d} G \tilde{\chi}_{\mu-1} ,  h_{0ij}^{\alpha\beta} \right] w_{0j}^\beta = - \sum_{\beta=1}^q \left[ \bold{H} \chi_{\mu-1} - \varphi_{\mu-1}'', h_{0ij}^{\alpha\beta} \right] w_{0j}^\beta  \notag
\end{align}
If we set
\begin{align}\label{utt49}
\varphi_\mu'':= \frac{1}{\mu} t \varphi_{\mu-1}'',\,\,\,\,\,\,\,\,\,\,\,\,\lambda_{i|\mu}'':= \frac{1}{\mu} t \lambda_{i|\mu-1}'',\,\,\,\,\,\,\,\,\,\,\, U_{i|\mu}:= \frac{1}{\mu} t \tilde{U}_{i|\mu-1}
\end{align}
then from $(\ref{utt46}),(\ref{utt47}),(\ref{utt48}),(\ref{utt49})$ and $(\ref{utt51})$ and $(\ref{utt50})$, we have
\begin{align}\label{n40}
- \bar{\partial}\left(\Sigma_{i|\mu}^\alpha - \mathcal{L}_{\mathfrak{d}G \chi_\mu}\left(w_{0i}^\alpha \right) + \sum_{\xi=1}^q U_{i|\mu}^{\alpha \xi} w_{0i}^\xi \right) + \mathcal{L}_{\bold{H}\chi_\mu - \varphi_\mu''}\left(w_{0i}^\alpha \right) - \left( \delta_i^{\chi_\mu}-  \lambda_{i|\mu}''  - \bar{\partial} U_{i|\mu}\right)\left( w_{0i}^\alpha \right) = A_{i|\mu}^\alpha - \bar{\partial} C_{i|\mu}^\alpha + \sum_{\gamma=1}^q B_{i|\mu}^{\alpha \gamma} w_{0i}^\gamma 
\end{align}
and
\begin{align}\label{n41}
&\sum_{\beta, \xi=1}^q h_{0ij}^{\alpha\beta}\left( \delta_j^{\chi_\mu \beta \xi} - \lambda_{j|\mu}''^{\beta \xi} - \bar{\partial} U_{j|\mu}^{\beta \xi}  \right) w_{0j}^\xi - \sum_{\xi=1}^q \left( \delta_i^{\chi_\mu\alpha\xi}  - \lambda_{i|\mu}^{\alpha \xi} - \bar{\partial} U_{i|\mu}^{\alpha \xi}\right) w_{0i}^\xi\\
&= - \sum_{\beta=1}^q \left[ \chi_\mu - \varphi_\mu'' - \bar{\partial} \mathfrak{d} G \chi_\mu,  h_{0ij}^{\alpha\beta} \right] w_{0j}^\beta = - \sum_{\beta=1}^q \left[ \bold{H} \chi_\mu - \varphi_\mu'', h_{0ij}^{\alpha\beta} \right] w_{0j}^\beta  \notag
\end{align}

On the other hand, from $(\ref{utt35})$ and $(\ref{utt36})$
\begin{align}
&-\bar{\partial} \tilde{\Pi}_{i|\mu-1}^\alpha + i_{\tilde{\chi}_{\mu-1}}\left( w_{0i}^\alpha \right) = \psi_{i|\mu-1}^\alpha - \bar{\partial} N_{i|\mu-1}^\alpha + i_{ \varphi_{\mu-1}''}\left( w_{0i}^\alpha \right) \notag \\
&\Longrightarrow -\bar{\partial}\left( \tilde{\Pi}_{i|\mu-1}^\alpha - i_{\mathfrak{d}G\tilde{\chi}_{\mu-1}}\left( w_{0i}^\alpha \right)    \right) + i_{ \bold{H} \tilde{\chi}_{\mu-1} - \varphi_{\mu-1}''}\left( w_{0i}^\alpha \right) = \psi_{i|\mu-1}^\alpha- \bar{\partial} N_{i|\mu-1}^\alpha
\end{align}

We note that from $(\ref{utt49})$ and $(\ref{utt41})$
\begin{align} \label{utt54}
\left| \left( \varphi_{\mu}'', \left\{ \lambda_{i|\mu}'' \right\} \right)  \right|_{k+\alpha} \ll Y_{27} \frac{b}{c} A(t)
\end{align}
and from $(\ref{utt49})$ and $(\ref{utt42})$
\begin{align} \label{uut1}
\left| U_{i|\mu} \right|_{k+1 + \alpha} \ll Y_{35} A(t)
\end{align} 

If we set
\begin{align}
\varphi_\mu' &:= \bold{H}{\chi_{\mu}} -\varphi_{\mu}'' \label{utt37}\\
w_{i|\mu}'^\alpha &:= \Sigma_{i|\mu}^\alpha - \mathcal{L}_{\mathfrak{d} G \chi_\mu }\left( w_{0i}^\alpha \right) + U_{i|\mu}\left( w_{0i}^\alpha \right) \label{utt38} \\
v_{i|\mu-1}'^\alpha &:= \tilde{\Pi}_{i|\mu-1}^\alpha - i_{\mathfrak{d}G\tilde{\chi}_{\mu-1}}\left( w_{0i}^\alpha \right) \label{utt39}\\
V_{i|\mu}^{\alpha \xi}&:=- \delta_i^{\chi_{\mu}\alpha \xi} + \lambda_{i|\mu}''^{\alpha \xi} + \bar{\partial} U_{i|\mu}^{\alpha \xi}, \label{utt40}
\end{align}
then we have from $(\ref{utt37})$ and $(\ref{utt53})$ and $(\ref{utt54})$
\begin{align} \label{uut9}
\left| \varphi_\mu' \right|_{k+\alpha} \ll Y_{36} \frac{b}{c} A(t)
\end{align}
and we have from $(\ref{utt38})$ and $(\ref{utt53})$ and $(\ref{uut1})$
\begin{align} \label{uut10}
\left| w_{i|\mu-1}'^\alpha \right|_{k+\alpha} \ll Y_{37} \frac{b}{c} \frac{\partial A(t)}{\partial t}
\end{align}
and we have from $(\ref{utt39})$ and $(\ref{utt32})$
\begin{align} \label{uut11}
\left| v_{i|\mu-1}'^\alpha  \right|_{k+\alpha} \ll Y_{38} \frac{b}{c} \frac{\partial A(t)}{\partial t}
\end{align}
and we have from $(\ref{utt40})$ and $(\ref{utt53})$ and $(\ref{utt54})$ and $(\ref{uut1})$
\begin{align}
\left| V_{i|\mu-1}^{\alpha\xi} \right|_{k+\alpha} \ll Y_{39} \frac{b}{c} A(t)
\end{align}
From $(\ref{utt45})$ and $(\ref{utt46})$ and $(\ref{utt49})$, we have
\begin{align*}
- \bar{\partial} \varphi_\mu' = - \bar{\partial}\left( \bold{H} \chi_{\mu} - \varphi_{\mu}'' \right) = \bar{\partial} \varphi_{\mu}'' = - \xi_\mu ,\,\,\,\,\,\,\, \mathfrak{d}\varphi_{\mu}' = \mathfrak{d}\left( \bold{H} \chi_{\mu} - \varphi_{\mu}'' \right) = - \mathfrak{d} \varphi_{\mu}'' = 0
\end{align*}
and so from $(\ref{uut2})$, we have $\varphi_\mu:=-\varphi_\mu'$, so that
\begin{align}\label{uut12}
\bar{\partial} \varphi_\mu= - \xi_\mu,\,\,\,\,\,\,\,\,\mathfrak{d}\varphi_\mu = 0
\end{align}
Since $\bar{\partial}\left( V_{i|\mu}^{\alpha\beta} - B_{i|\mu}^{\alpha\beta} \right)=0$, there exists $P_{i|\mu}^{\alpha\beta}= V_{i|\mu}^{\alpha\beta}- B_{i|\mu}^{\alpha\beta}$. Then we have as in Part I
\begin{align}
\left| P_{i|\mu}^{\alpha\beta} \right|_{k+1+\alpha} \ll Y_{40} \frac{b}{c} A(t)
\end{align}

Then we have from $(\ref{uut2})$ and $(\ref{ue22})$ and $(\ref{uut3})$ and $(\ref{ue24})$
\begin{align}
\varphi_\mu&:=  -\varphi_\mu'=-\bold{H}\chi_\mu + \varphi_\mu''  \label{uut4}\\
w_{i|\mu}^\alpha:&= C_{i|\mu}^\alpha - w_{i|\mu}'^\alpha + \sum_{\beta=1}^q P_{i|\mu}^{\alpha\beta} w_{0i}^\beta \label{uut5}\\
h_{ij|\mu}^{\alpha\beta}&:= \sum_{\gamma=1}^q D_{ij|\mu}^{\alpha \gamma} h_{0ij}^{\gamma \beta} + \sum_{\gamma=1}^q P_{i|\mu}^{\alpha \gamma} h_{0ij}^{\gamma \beta} - \sum_{\gamma=1}^q h_{0ij}^{\alpha \gamma} P_{j|\mu}^{\gamma \beta} \label{uut6} \\
v_{i|\mu-1}^\alpha&:= N_{i|\mu-1}^\alpha - v_{i|\mu-1}'^\alpha \label{uut7}
\end{align}
Then we have as in Part I
\begin{align}
\left| \varphi_\mu \right|_{k+\alpha},\,\,\,\,\,\,\,\,\, \left| w_{i|\mu}^\alpha \right|_{k+\alpha},\,\,\,\,\,\,\,\,\,\left| h_{ij|\mu}^{\alpha\beta} \right|_{k+1+\alpha}\,\,\,\,\,\, \ll Y_{41} \frac{b}{c} A(t)
\end{align}
and from $(\ref{uut7})$ and $(\ref{utt28})$ and $(\ref{uut11})$
\begin{align}
\left| v_{i|\mu-1}^\alpha  \right|_{k+\alpha} \ll  Y_{42} \frac{b}{c} \frac{\partial A(t) }{\partial t}
\end{align}
Then we can choose $c>b>0$ satisfying $\left(\ref{us7} \right)_\mu-\left(\ref{us10} \right)_\mu$ and $\left(\ref{us8} \right)_{\mu-1}$. Then $\varphi(t)$ and $w_i^\alpha$ and $v_i^\alpha$ and $h_{ij}^{\alpha\beta}$ converge with respect to \" Holder norm $|-|_{k+\alpha}, |-|_{k+\alpha}, |-|_{k+\alpha},|-|_{k+1+\alpha}$, respectively for some sufficiently small neighborhood of $\Delta_\epsilon \subset B \subset \mathbb{C}$ of $0$. Consequently $\varphi(t)$ is a $C^k$ vector $(0,1)$-form on $M\times \Delta_\epsilon$ and $w_i^\alpha+ v_i^\alpha dt$ is a $C^k$-form of the form $w_i^\alpha+ v_i^\alpha dt= \sum_{\beta=1}^n w_i^{\alpha\beta}\left(z_i,t \right)dz_i^\beta + v_i^\alpha dt$ on $U_i\times \Delta_\epsilon$, and $h_{ij}^{\alpha\beta}$ is a $ C^{k+1}$ on $U_{ij}\times \Delta_\epsilon$. From $(\ref{uut12})$, we have $\mathfrak{d}\varphi_\mu =0$ for $\mu\geq 2$. Hence $\varphi(t)$ is a solution of the quasi-linear partial differential equation of order $2$ 
\begin{align}\label{ns53}
 \frac{\partial^2}{\partial t \partial \bar{t}} \varphi(t) +\square \varphi(t) - \mathfrak{d}\left[ \varphi(t), \varphi(t)\right]=\bar{\partial}\mathfrak{d}\varphi_1(t)
\end{align}
We may assume that $(\ref{ns53})$ is a quasi-linear elliptic partial differential equation on $M\times \Delta_\epsilon$ (see \cite{Kod05} p.281). Therefore its solution $\varphi(t)$ is $C^\infty$ on $M\times \Delta_\epsilon$ (see \cite{DN55} Theorem 5 or \cite{Kod05} Appendix \S 8). Then $\varphi(t)$ determines a complex analytic family $\pi: \mathcal{M}\to \Delta_\epsilon$ (see \cite{Kod05} p.281-283). On the other hand, $w_i^\alpha+v_i^\alpha dt$ and $h_{ij}^{\alpha\beta}$ are holomorphic with respect to the complex structure $\varphi(t)$ by $(\ref{u13})$ and $(\ref{uut15})$ and $(\ref{uut16})$, respectively. $(\ref{uut17}),(\ref{uut18})$ and $(\ref{uut19})$ determines a locally free subsheaf $\tilde{\mathcal{N}}_\mathcal{F}^*$ of $\Omega_{\mathcal{M}}^1$ satisfying the integrability condition. From the definition of linear terms, we infer that the foliated Kodaira-Spence map $\varphi_0: T_0\left(\Delta_\epsilon \right)\cong \mathbb{C}\to \mathbb{H}^1\left( M, \tilde{\mathcal{N}}_{\mathcal{F}_0}^{*\bullet} \right)$ send $\frac{\partial}{\partial t}$ to $\phi$. This completes the proof of Theorem \ref{u56}.
\end{proof}

We finish this section by making the conjecture on the complete form of `Theorem of existence for unfoldings of foliated complex analytic structures defined by locally free subsheaves of cotangent sheaves' as follows:
\begin{conjecture}[Theorem of existence for unfoldings of foliated complex analytic structures defined by locally free subsheaves of cotangent sheaves] \label{conjecture1}
Let $\left( M, \mathcal{N}_{\mathcal{F}_0}^* \right)$ be a compact foliated complex manifold with $\mathcal{N}_{\mathcal{F}_0}^*$ locally free. Suppose that $\mathbb{H}^2\left( M, \tilde{\mathcal{N}}_{\mathcal{F}_0}^{*\bullet} \right)=0$. Then there exists an unfolded analytic family $\left( \mathcal{M}, \tilde{\mathcal{N}}_\mathcal{F}^*, B, \pi \right)$ with $0\in B\subset \mathbb{C}^m$ satisfying the following conditions$:$
\begin{enumerate}
\item $\pi^{-1}(0)=\left( M, \mathcal{N}_{\mathcal{F}_0}^* \right)$
\item The unfolded Kodaira-Spencer map $\varphi_0: T_0 B \to \mathbb{H}^1\left( M,  \tilde{\mathcal{N}}_{\mathcal{F}_0}^{* \bullet} \right)$ is an isomorphism.
\end{enumerate}

\end{conjecture}

\section{Theorem of stability for unfoldings of foliated complex analytic structures for one-parameter families}

Our proof of theorem of existence of unfoldings for one parameter families (Theorem \ref{u56}) provides the proof of theorem of stability for unfoldings of foliated complex analytic structures for one-parameter families (Theorem \ref{nst1}). We define a complex of sheaves $\textnormal{Unf}_{\mathcal{N}_{\mathcal{F}_0}^*}^\bullet$ which we truncate the $0$-th term of the unfolded dual leaf complex
\begin{equation}
\begin{CD}
\cdots \\
@AG_3AA \\
\mathscr{H}om_{\mathcal{O}_M}\left( \mathcal{N}_{\mathcal{F}_0}^*, \tilde{\mathcal{S}}^3 \right) \bigoplus \mathscr{H}om_{\mathcal{O}_M}\left( \mathcal{N}_{\mathcal{F}_0}^*, \tilde{\mathcal{S}}^2  \right) \\
@AG_2AA \\
\mathscr{H}om_{\mathcal{O}_M}\left( \mathcal{N}_{\mathcal{F}_0}^*, \tilde{\mathcal{S}}^2 \right) \bigoplus \mathscr{H}om_{\mathcal{O}_M}\left( \mathcal{N}_{\mathcal{F}_0}^*, \tilde{\mathcal{S}}^1 \right) \\
@AG_1 AA \\
\mathscr{H}om_{\mathcal{O}_M}\left( \mathcal{N}_{\mathcal{F}_0}^* , \frac{\Omega_M^1}{\mathcal{N}_{\mathcal{F}_0}^*} \right) \bigoplus \mathscr{H}om_{\mathcal{O}_M}\left( \mathcal{N}_{\mathcal{F}_0}^*, \mathcal{O}_M \right) 
\end{CD}
\end{equation}
We will denote the $i$-th cohomology group of $\textnormal{Unf}_{\mathcal{N}_{\mathcal{F}_0}^*}^\bullet$ by $\mathbb{H}^i\left( M,\textnormal{Unf}_{\mathcal{N}_{\mathcal{F}_0}^*}^\bullet \right)$.

\begin{theorem}[Theorem of stability for unfoldings of foliated complex analytic structures for one-parameter families] \label{nst1}
Let $\left( M, \mathcal{N}_{\mathcal{F}_0}^* \right)$ be a compact foliated complex manifold with $\mathcal{N}_{\mathcal{F}_0}^*$ locally free. Assume that the first cohomology $\mathbb{H}^1\left( M, \textnormal{Unf}_{\mathcal{N}_{\mathcal{F}_0}^*}^\bullet \right)=0$. Then for any complex analytic family $\pi:\mathcal{M}\to B$ of deformations of $\pi^{-1}(0)=M, 0\in B\subset \mathbb{C}$ in the sense of Kodaira-Spencer, there exists an open neighborhood $N\subset B$ of $0$ and a locally free subsheaf $\tilde{\mathcal{N}}_{\mathcal{F}}^*$ of $\Omega_{\mathcal{M}|_N }^1$ such that $\left(\mathcal{M}|_N, \tilde{\mathcal{N}}_{\mathcal{F}}^*\right)$ defines an unfolded analytic family $\pi|_N: \left(\mathcal{M}|_N , \tilde{\mathcal{N}}_{\mathcal{F}}^* \right) \to N$ of unfoldings of $\left(M, \mathcal{N}_{\mathcal{F}_0}^* \right)=\pi^{-1}(0)$.
\end{theorem}

\begin{proof}
We copy the proof of Theorem \ref{u56}. The difference is that we assume that $\varphi(t)$ (determined by $\mathcal{M}$ over some neighborhood $N'\subset B \subset \mathbb{C}$ of $0$) is given from the beginning and replace $\varphi^\mu$ by $\varphi(t)$ in the proof of Theorem \ref{u56}. Accordingly we replace $\xi_\mu$ and $\varphi_\mu$ by $0$. Then the obstructions for solving $(\ref{u21})-(\ref{uut25})$ are in $\mathbb{H}^1\left( M, \textnormal{Unf}_{\mathcal{N}_{\mathcal{F}_0}^*}^\bullet \right)=0$. Hence we can construct $\tilde{\mathcal{N}}_{\mathcal{F}}^*$ on $\mathcal{M}|_N$ for some neighborhood $N\subset N'$ of $0$. 
\end{proof}

\section{Theorem of completeness for unfoldings of foliated complex analytic structures}

\subsection{Change of parameters}\

Consider an unfolded complex analytic family $\left( \mathcal{M}, \tilde{\mathcal{N}}_\mathcal{F}^*, B, \pi \right)$ with $\tilde{\mathcal{N}}_\mathcal{F}^*$ locally free of unfoldings of $\left( M_t, \mathcal{N}_{\mathcal{F}_t}^* \right)$ which $\pi^{-1}(t),t\in B$ induces, where $B$ is a domain of $\mathbb{C}^m$. Let $D$ be a domain of $\mathbb{C}^{m'}$ and $s:u\to t=s(u), u\in D$, a holomorphic map of $D$ into $B$. Then by changing the parameter from $t$ to $u$, we will construct an unfolded analytic family $\left\{ \left( M_{s(u)}, \tilde{\mathcal{N}}_{\mathcal{F}_{s(u)}}^* \right) | u \in D \right\}$ on the parameter space $D$ in the following.

Let $\mathcal{M}\times_B D:=\left\{\left(p,u\right)\in \mathcal{M}\times D |\pi(p)=s(u) \right\}$. Then we have the following commutative diagram
\begin{center}
$\begin{CD}
\mathcal{M}\times_B D @>p>> \mathcal{M} \\
@V\pi' VV @VV \pi V\\
D @>s >> B 
\end{CD}$
\end{center}
such that $\left( \mathcal{M}\times_B D, D, \pi' \right)$ is a complex analytic family in the sense of Kodaira-Spencer and $\pi'^{-1}(u)= M_{s(u)}$. On the other hand, $\Pi^* \tilde{\mathcal{N}}_\mathcal{F}^*$ is a locally free subsheaf of $\Omega_{\mathcal{M}\times D}^1$ where $\Pi:\mathcal{M}\times D \to \mathcal{M}$ is the natural projection. Since $\mathcal{M}\times_B D$ is a complex submanifold of $\mathcal{M}\times D$, we have an exact sequence  $0\to \mathcal{I}/\mathcal{I}^2 \xrightarrow{d} \Omega_{\mathcal{M}\times D}^1|_{\mathcal{M}\times_B D} \to \Omega_{\mathcal{M}\times_B D}^1\to 0$ where $\mathcal{I}$ is the ideal sheaf of $\mathcal{O}_{\mathcal{M}\times D}$ defining $\mathcal{M}\times_B D$. The ideal sheaf $\mathcal{I}$ is locally defined by the equation $t=s(u)$, i.e. $t_\alpha=s_\alpha\left(u_1,..., u_{m'}\right),\alpha=1,...,m$. Then the image $d\left(\mathcal{I}/\mathcal{I}^2 \right)$ of $\mathcal{I}/\mathcal{I}^2$ in $\Omega_{\mathcal{M}\times D}^1|_{\mathcal{M}\times_B D}$ is locally generated by the form $dt_\alpha= \sum_{\beta=1}^{m'} \frac{\partial s_\alpha}{\partial u_\beta}(u)du_\beta$, so that $d\left( \mathcal{I}/\mathcal{I}^2\right)\bigcap p^* \tilde{\mathcal{N}}_{\mathcal{F}}^* = 0$. Here we note that $\Pi^*\tilde{\mathcal{N}}_\mathcal{F}^*|_{\mathcal{M}\times_B D}=p^*\tilde{\mathcal{N}}_\mathcal{F}^*$. Hence $p^*\tilde{\mathcal{N}}_\mathcal{F}^*$ is isomorphically mapped under $\Omega_{\mathcal{M}\times D}^1|_{\mathcal{M}\times_B D}\to \Omega_{\mathcal{M}\times_B D}^1$ so that $p^*\tilde{\mathcal{N}}_\mathcal{F}^*$ can be considered to be a locally free subsheaf of $\Omega_{\mathcal{M}\times_B D}^1$. Let $S$ be the locus of $\mathcal{M}$ such that $\frac{\Omega_\mathcal{M}^1}{\tilde{\mathcal{N}}_\mathcal{F}^*}$ is not locally free, and let $\tilde{S}=S\times_B D$. Then $p$ induces $p: \left(\mathcal{M}\times_B D \right) - \tilde{S} \to \mathcal{M}- S$. Since the pullback commutes with the differential and wedge products, $d\tilde{\mathcal{N}}_\mathcal{F}^*\subset \tilde{\mathcal{N}}_\mathcal{F}^*\wedge \Omega_\mathcal{M}^1$ on $\mathcal{M}- S$ implies that $d\left(p^* \tilde{\mathcal{N}}_\mathcal{F}^* \right) \subset p^*\tilde{\mathcal{N}}_\mathcal{F}^* \wedge \Omega_{\mathcal{M}\times_B D}^1 $ on $\left(\mathcal{M}\times_B D \right) - \tilde{\mathcal{S}}$. From the exact sequence $0\to \pi'^* \Omega_D^1\to \Omega_{\mathcal{M}\times_B D}^1 \xrightarrow{\alpha'} \Omega_{\mathcal{M}\times_B D/D}^1\to 0$, and $\pi'^* \Omega_D^1\cap p^* \tilde{\mathcal{N}}_\mathcal{F}^*=0$, we see that $p^*\tilde{\mathcal{N}}_\mathcal{F}^* \cong \alpha'\left( p^* \tilde{\mathcal{N}}_\mathcal{F}^* \right)= p^* \mathcal{N}_\mathcal{F}^*$, where $\mathcal{N}_\mathcal{F}^*$ is the isomorphic image of $\tilde{\mathcal{N}}_\mathcal{F}^*$ under $\alpha: \Omega_\mathcal{M}^1\to \Omega_{\mathcal{M}/B}^1$. Moreover for each $u\in D$, the composition $p^*\tilde{\mathcal{N}}_\mathcal{F}^*|_{\pi'^{-1}(u)}\cong p^* \mathcal{N}_{\mathcal{F}}^* |_{\pi^{-1}(u)}= \mathcal{N}_{\mathcal{F}_{s(u)}}^* \to \Omega_{\mathcal{M}\times_B D}^1|_{\pi'^{-1}(u)} \to \Omega_{\mathcal{M}\times_B D/D}^1|_{\pi'^{-1}(u)} = \Omega_{M_{s(u)}}^1$ is injective, so that $p^*\tilde{\mathcal{N}}_\mathcal{F}^*|_{\pi'^{-1}(u)} \to \Omega_{\mathcal{M}\times_B D}^1|_{\pi'^{-1}(u)} $ is injective. Hence $\frac{\Omega_{\mathcal{M}\times_B D}^1}{p^* \tilde{\mathcal{N}}_\mathcal{F}^*}$ is flat over $D$. We note that we have the exact sequence $0\to \pi'^* \Omega_D^1 \to \frac{\Omega_{\mathcal{M}\times_B D}^1}{p^*\tilde{\mathcal{N}}_\mathcal{F}^*}\to \frac{\Omega_{\mathcal{M}\times_B D/D}^1}{p^* \mathcal{N}_\mathcal{F}^*}\to 0$ on $\mathcal{M}\times_B D $ and so on $\left(\mathcal{M}\times_B D \right) - \tilde{S}$. Let $\pi'': \left( \mathcal{M}\times_B D \right)- S \to D$ be the restriction of $\pi'$ to $\left( \mathcal{M}\times_B D  \right)-\tilde{S}$. Then for each $u\in D$, $\frac{\Omega_{\mathcal{M}\times_B D/D}^1}{p^* \mathcal{N}_\mathcal{F}^*}|_{\pi''^{-1}(u)}= \frac{\Omega_{M_{s(u)}}^1}{\mathcal{N}_{\mathcal{F}_{s(u)}}^*}|_{M_{s(u)}- S_{s(u)}}$ is locally free, so that $\frac{\Omega_{\mathcal{M}\times_B D/D}^1}{p^* \mathcal{N}_\mathcal{F}^*}$ is locally free on $\left(\mathcal{M}\times_B D \right) - \tilde{\mathcal{S}}$. Hence $\frac{\Omega_{\mathcal{M}\times_B D}^1}{p^* \tilde{\mathcal{N}}_\mathcal{F}^*}$ is locally free on $\left( \mathcal{M} \times_B D \right) - \tilde{S}$. This implies that $\left( \mathcal{M}\times_B D, p^* \tilde{\mathcal{N}}_\mathcal{F}^*, D, \pi' \right)$ forms an unfolded complex analytic family.

\begin{definition}
The unfolded complex analytic family $\left( \mathcal{M}\times_B D, p^* \tilde{\mathcal{N}}_\mathcal{F}^*, D, \pi' \right)$ is called the unfolded analytic family induced from $\left( \mathcal{M}, \tilde{\mathcal{N}}_\mathcal{F}^*, B, \pi \right)$ by the holomorphic map $s:D\to B$.
\end{definition}

To investigate the relation of the infinitesimal unfolding of $\left(M_{s(u)}, \mathcal{N}_{\mathcal{F}_{s(u)}}^* \right)$ and that of $\left( M_t , \mathcal{N}_{\mathcal{F}_t}^* \right)$, we assume that $0\in B, 0\in D$ and $s(0)=0$. Taking a sufficiently small coordinate polydisk $\Delta$ with $0\in \Delta \subset B$, we represent $\left( \mathcal{M}_\Delta, \tilde{\mathcal{N}}_{\mathcal{F}_\Delta}^* \right)=\pi^{-1}\left(\Delta \right)$ in the form
\begin{align*}
\left(\mathcal{M}_\Delta, \tilde{\mathcal{N}}_{\mathcal{F}_\Delta}^* \right) = \bigcup \left( U_j \times \Delta, \tilde{\mathcal{N}}_{\mathcal{F}_\Delta}^*|_{U_j\times \Delta} \right)
\end{align*}
where $\left(z_j,t \right)\in U_j\times \Delta$ and $\left( z_k, t\right)\in U_k\times \Delta$ are the same point of $\mathcal{M}_\Delta$ if $z_j= f_{jk}(z_k,t)$, and $\Gamma\left( U_j\times \Delta,\tilde{\mathcal{N}}_\mathcal{F}^*\right)$ is generated by $\omega_j^1\left(z_j,t\right),..., \omega_j^q(z_j,t)$, where $\omega_j^\alpha\left(z_j,t\right)= w_j^\alpha(z_j,t) + v_j^\alpha(z_j,t)= \sum_{\beta=1}^n w_j^{\alpha\beta}(z_j,t_j)dz_j^\beta + \sum_{\eta=1}^m v_j^{\alpha\eta}(z_j,t_j)dt_j^\eta$ for some $w_j^{\alpha\beta}(z_j,t_j)$ and $v_j^{\alpha\beta}\left(z_j,t_j\right)\in \Gamma\left( U_j\times \Delta, \mathcal{O}_\mathcal{M}\right),\alpha=1,...,q$ with $\omega_j^\alpha\left(z_j,t_j \right)=\sum_{\beta=1}^q h_{jk}^{\alpha\beta}\left(z_k,t\right) \omega_k^\beta(z_k,t)$ for some $h_{jk}^{\alpha\beta}(z_k,t)\in \Gamma\left( \left(U_j\times \Delta \right)\bigcap\left( U_k \times \Delta \right), \mathcal{O}_\mathcal{M}\right)$ for $\alpha=1,...,q$, and the integrability condition $\omega_j^1\left(z_j,t\right) \wedge \cdots \wedge \omega_j^q\left(z_j,t\right) \wedge d \omega_j^\alpha \left(z_j,t\right)=0$ for $\alpha=1,...,q$. Take a polydisk $\Delta'$ with $0\in \Delta' \subset D$. Then $\left( \mathcal{M}\times_B D|_{\Delta'}, p^* \tilde{\mathcal{N}}_\mathcal{F}^*|_{\Delta'}, \Delta', \pi' \right)$ is represented by the form 
\begin{align*}
\left( \mathcal{M}\times_B D|_{\Delta'}, p^* \tilde{\mathcal{N}}_\mathcal{F}^* |_{\Delta'} \right)=  \bigcup \left( U_j\times \Delta', p^* \tilde{\mathcal{N}}_\mathcal{F}^*|_{U_j\times \Delta'} \right)
\end{align*}
where $(z_j, u)\in U_j\times \Delta'$ and $(z_k, u)\in U_k\times \Delta'$ are the same point of $\mathcal{M}\times_B D|_{\Delta'}$ if $z_j= f_{jk}\left(z_k, s(u)\right)$ and $\Gamma\left( U_j\times \Delta', p^* \tilde{\mathcal{N}}_\mathcal{F}^* \right)$ is generated by $\omega_j^1\left(z_j, s(u)\right),..., w_j^q\left(z_j, s(u)\right)$, where
\begin{align*}
\omega_j^\alpha(z_j, s(u))=w_j^\alpha\left(z_j, s(u)\right)+ v_j^\alpha \left(z_j, s(u)\right)= \sum_{\beta=1}^n w_j^{\alpha\beta}\left(z_j, s(u) \right)  dz_j^\beta + \sum_{\delta=1}^{m'} \sum_{\eta=1}^m v_j^{\alpha\eta}\left(z_j, s(u) \right) \frac{\partial s^\eta (u)}{\partial u_\delta} d u_\delta 
\end{align*}
and $\omega_j^1\left( z_j, s(u) \right) \wedge \cdots \wedge \omega_j^q\left( z_j, s(u) \right) \wedge d \omega_j^\alpha\left( z_j, u(s) \right)=0$ for $\alpha=1,...,q$. Then we can show that

\begin{theorem}
For any tangent vector $\frac{\partial }{\partial u}=c_1\frac{\partial}{\partial u_1} + \cdots + c_{m'}\frac{\partial}{\partial u_{m'}}\in T_u \left(D \right) $, the infinitesimal unfolding of $\left(M_{s(u)}, \mathcal{N}_{\mathcal{F}_{s(u)}}^* \right)$ along $\frac{\partial}{\partial u}$ is given by
\begin{align*}
\frac{\partial \left(M_{s(u)},  \tilde{N}_{\mathcal{F}_{s(u)}}^*\right)}{\partial u}&=\left(\sum_{\gamma=1}^m \frac{\partial t_\gamma}{\partial u}\frac{\partial M_t}{\partial t_\gamma},\left\{w_j^\alpha\left(z_j,s(u)\right) \mapsto -\sum_{\gamma=1}^m \frac{\partial t_\gamma}{\partial u} \frac{\partial w_j^\alpha}{\partial t_\gamma} \right\}, \left\{ w_j^\alpha \left(z_j, s(u)\right)\mapsto  - \sum_{\gamma=1}^m \frac{\partial t_\gamma}{\partial u} i_{\frac{\partial}{\partial t_\gamma}}\left(v_j^\alpha \right)\right\} \right)
           \\
           & =\sum_{\gamma=1}^m \frac{\partial t_\gamma}{\partial u} \frac{\partial \left(M_t, \tilde{\mathcal{N}}_{\mathcal{F}_t}^* \right) }{\partial t_\gamma}
\end{align*}
\end{theorem}

\subsection{Theorem of completeness for unfoldings of foliated complex analytic structures}\

\begin{definition}
Let $\left( \mathcal{M}, \tilde{\mathcal{N}}_\mathcal{F}^*, B, \pi \right)$ with $\tilde{\mathcal{N}}_\mathcal{F}^*$ locally free be an unfolded analytic family of unfoldings of compact foliated complex manifolds, and $t^0\in B$. Then $\left( \mathcal{M} , \tilde{\mathcal{N}}_\mathcal{F}^*, B, \pi \right)$ is called complete at $t^0 \in B$ if for any unfolded complex analytic family $\left( \mathcal{M}', \tilde{\mathcal{N}}_{\mathcal{F}'}, D, \pi' \right)$ such that $D$ is a domain of $\mathbb{C}^{m'}$ containing $0$ and that $\pi'^{-1}(0)= \pi^{-1}(t^0)$ induces $\left(M, \mathcal{N}_{\mathcal{F}_0}^* \right)$, there are a sufficiently small domain $\Delta$ with $0\in \Delta \subset D$, and a holomorphic map $s:u\to t=s(u)$ with $s(0)=t^0$ such that $\left( \mathcal{M}'_{\Delta}, \tilde{\mathcal{N}}_{\mathcal{F}'|_\Delta}^*, \Delta, \pi' \right)$ is the unfolded complex analytic family induced from $\left( \mathcal{M}, \tilde{\mathcal{N}}_\mathcal{F}^*, B, \pi \right)$ by $h$ where $\left( \mathcal{M}'_{\Delta}, \tilde{\mathcal{N}}_{\mathcal{F}'|_\Delta}^*, \Delta, \pi'   \right)$ is the restriction of $\left( \mathcal{M}', \tilde{\mathcal{N}}_{\mathcal{F}'}^*, D, \pi' \right)$ to $\Delta$.
\end{definition}

\begin{theorem}\label{uc1}
Let $\left(\mathcal{M}, \tilde{\mathcal{N}}_\mathcal{F}^*, B, \pi \right)$ be an unfolded analytic family of unfdolgins of a compact foliated complex manifold $\left(M , \mathcal{N}_{\mathcal{F}_0}^* \right)=\omega^{-1}(0)$ such that $\mathcal{N}_{\mathcal{F}_0}^*$ is locally free, and $B$ is a domain of $\mathbb{C}^r$ containing $0$. If the unfolded Kodaira-Spence map $\varphi_0: T_0(B)\to \mathbb{H}^1\left(M, \tilde{\mathcal{N}}_{\mathcal{F}_0}^{ * \bullet } \right)$ is surjective, then the unfolded analytic family $(\mathcal{M}, \tilde{\mathcal{N}}_\mathcal{F}^*, B, \pi)$ is complete at $0 \in B$.
\end{theorem}

\begin{proof}
Let $\left(\mathcal{M}, \tilde{\mathcal{N}}_\mathcal{F}^*, B, \pi' \right)$ be an unfolded analytic family which is represented as in Remark \ref{u3}.  We keep the notations in Remark \ref{u3}. Since the problem is local with respect to $B$, we may assume that $B=\left\{ t\in \mathbb{C}^m | |t| <1\right\}$ is a polydisk, and $\mathcal{M}$ is written in the form
\begin{align}
\mathcal{M}=\bigcup_j \,\,\, \mathcal{U}_j,\,\,\,\,\,\,\,\,\mathcal{U}_j=\left\{\left( z_j, t\right) \in \mathbb{C}^n \times B\,\,\, |\,\,\, |z_j| <1 \right\}
\end{align}

Set $U_i:=M \cap \mathcal{U}_i=\left\{ z_i\in \mathbb{C}^n | |z_i|<1 \right\}$. Then $\left( M, \mathcal{N}_{\mathcal{F}_0}^* \right)=\pi^{-1}(0)$ is described in terms of covering $\mathcal{U}:=\{U_i\}=\left\{ M \cap \mathcal{U}_i \right\}$ in the following way: (1) local coordinates $z_i= \left(z_i^1,..., z_i^n \right)$ on $U_i$ with $z_i= f_{0ij}^\alpha(z_j)= \left( f_{0ij}^1(z_j),..., f_{0ij}^n(z_j) \right)$ where we set $f_{0ij}(z_j):= f_{ij}(z_j,0)$ and (2) $\Gamma\left(U_i, \mathcal{N}_{\mathcal{F}_0}^* \right)$ is generated by $w_{0i}^\alpha(z_i):= \sum_{\beta=1}^n w_i^{\alpha\beta}(z_i,0)dz_i^\beta \in \Gamma\left(U_i, \mathcal{N}_{\mathcal{F}_0}^* \right)$ with $w_{0i}^\alpha= \sum_{\beta=1}^q h_{0ij}^{\alpha\beta}(z_j) w_{0j}^\beta$ where we set $h_{0ij}^{\alpha\beta}(z_j):=h_{ij}^{\alpha\beta}(z_j,0)$, and $w_{0i}^1\wedge \cdots \wedge w_{0i}^q \wedge dw_{0i}^\alpha=0$.

Let $\left(\mathcal{M}', \tilde{\mathcal{N}}_{\mathcal{F}'}^*, D, \pi' \right)$ be an another foliated analytic family such that $\pi'^{-1}(0')= \left(M, \mathcal{N}_{\mathcal{F}_0}^* \right)$. We may assume the following:
\begin{enumerate}
\item $D\subset \left\{ u\in \mathbb{C}^{m'} ||u|<1 \right\}$ is a sufficiently small polydisk in $\mathbb{C}^{m'}$ with a system of coordinates $u=(u_1,..., u_{m'})$ centered at $0'$.
\item $\mathcal{M}'$ is covered by a finite number of coordinate neighborhood $\mathcal{U}_i'=\left\{ \left(\xi_i, u\right) \in \mathbb{C}^n \times D\,\,\,|\,\,\,|\xi_i|<1\right\}$ with a system of coordinate $(\xi_i, u)$ such that $\pi'(\xi_i, u)=u$.
\item $(\xi_i, u)$ coincides with $(\xi_j, u)$ if and only if $\xi_i = f_{ij}'(\xi_i, u)$.
\item $\tilde{\mathcal{N}}_{\mathcal{F}'}^*$ on $\mathcal{U}_i'$ is generated by
\begin{align*}
\omega_i'^\alpha(\xi_i, u):=w_i'^\alpha(\xi_i, u)+ v_i'^\alpha(\xi_i, u)=\sum_{\beta=1}^n w_i'^{\alpha\beta}(\xi_i, u)d\xi_i^\beta + \sum_{\eta=1}^{m'} v_i'^{\alpha\eta} (\xi_i, u)du_\eta,\,\,\,\,\,\,\,\alpha=1,...,q
\end{align*}
\item on $\mathcal{U}_i' \cap \mathcal{U}_j'\ne \emptyset$, we have $\omega_i'^\alpha(\xi_i,u)= \sum_{\beta=1}^q h_{ij}'(\xi_j, u) \omega_j'^\beta(\xi_j, u)$ for $h_{ij}'^{\alpha\beta}(\xi_j, u)\in \Gamma\left( \mathcal{U}_i'\cap \mathcal{U}_j', \mathcal{O}_\mathcal{M} \right)$.
\item $\omega_i'^1(\xi_i,u) \wedge \cdots \wedge \omega_i'^q(\xi_i, u) \wedge d \omega_i'^\alpha(\xi_i, u)=0,\alpha=1,...,q$.
\item $\pi'^{-1}(0)\cap \mathcal{U}_i'=U_i$, and $\xi_i= z_i$ on $U_i$ and $f_{0ij}(z_j)=f_{0ij}'(\xi_j)$ where $f_{0ij}'(\xi_j) : =f_{ij}'(\xi_j,0)$.
\item Setting $w_{0i}'^\alpha(z_i):= w_i'^\alpha(z_i,0)$, i.e.
\begin{align*}
w_{0i}'^\alpha(z_i):= w_i'^\alpha(z_i,0)=\sum_{\beta=1}^n w_i'^{\alpha\beta}(z_i,0)dz_i^\beta,\,\,\,\,\,\alpha=1,...,q,
\end{align*}
we can find $c_{0i}^{\alpha\beta}(z_i)\in \Gamma\left( U_i, \mathcal{O}_M\right),\alpha,\beta=1,...,q$ such that
\begin{align*}
w_{0i}^\alpha(z_i)= \sum_{\beta=1}^q c_{0i}^{\alpha\beta}(z_i) w_{0i}'^\beta (z_i)
\end{align*}
\item we set $h_{0ij}'^{\alpha\beta}(\xi_j):= h_{ij}'^{\alpha\beta}(\xi_j,0)$.
\end{enumerate}
We note that  as in Part I
\begin{align*}
 \sum_{\beta=1}^q c_{0i}^{\alpha\beta} h_{0ij}'^{\beta\gamma} =  \sum_{\beta=1}^q h_{0ij}^{\alpha\beta} c_{0j}^{\beta\gamma}
\end{align*}

In order to prove Theorem \ref{uc1} it suffices to construct holomorphic functions
\begin{align*}
\varphi_i&:\mathcal{U}_i'\to \mathbb{C}^n \\
s &:D \to \mathbb{C}^{m}
\end{align*}
and a matrix function
\begin{center}
$\left[\begin{matrix}
c_{i}^{11}(\xi_i,u) & c_{i}^{12}(\xi_i,u) & \cdots & c_{i}^{1q}(\xi_i, u)\\
c_{i}^{21}(\xi_i, u) & c_{i}^{22}(\xi_i, u) & \cdots & c_{i}^{2q}(\xi_i, u)\\
\cdot & \cdot & \cdots & \cdot\\
\cdot & \cdot & \cdots & \cdot\\
\cdot & \cdot & \cdots & \cdot\\
c_{i}^{q1}(\xi_i, u) & c_{i}^{q2}(\xi_i, u) & \cdots & c_{i}^{qq}(\xi_i, u)\\
\end{matrix}\right]$
\end{center}
where $c_i^{\alpha\beta}:\mathcal{U}_i'\to \mathbb{C}$ such that
\begin{align}
\varphi_i(\xi_i, 0)=\xi_i, \,\,\,\,\,s(0)&=0,\,\,\,\,\, c_i^{\alpha\beta}(\xi_i,0)= c_{0i}^{\alpha\beta}(\xi_i)\\
\varphi_i \left(f_{ij}', u \right) &= f_{ij}\left(\varphi_j, s(u)\right) \label{uc2}
\end{align}
\begin{align}
\sum_{\beta=1}^q c_i^{\alpha\beta}\left(f_{ij}' (\xi_j, u), u\right) h_{ij}'^{\beta\gamma}(\xi_j , u) = \sum_{\beta=1}^q h_{ij}^{\alpha\beta} \left (\varphi_j(\xi_j, u), s(u) \right) c_j^{\beta\gamma}(\xi_j, u)\label{uc3}
\end{align}
By setting
\begin{align*}
\left(\varphi_i, s \right): \mathcal{U}_i' \to \mathbb{C}^n \times\mathbb{C}^m,
\end{align*}
we have
\begin{align*}
&\left(\varphi_i,s\right)^*\left(w_i^\alpha+ v_i^\alpha\right)=\left(\varphi_i,s\right)^*\left(\sum_{\beta=1}^q w_i^{\alpha\beta}(z_i,t)dz_i^\beta+ \sum_{\eta=1}^r v_i^{\alpha\eta}(z_i,t)dt_\eta \right) \\
 &= \sum_{\beta=1}^q w_i^{\alpha\beta}\left(\varphi_i(\xi_i,u),s(u) \right)d \varphi_i^\beta(\xi_i,u) + \sum_{\eta=1}^r v_i^{\alpha\eta}\left(\varphi_i(\xi_i,u),s(u)\right)d s^\eta(u)\\
 &= \sum_{\beta=1}^q \sum_{\gamma=1}^n w_i^{\alpha \beta}\left( \varphi_i(\xi_i, u), s(u)\right)\frac{\partial \varphi_i^\beta}{\partial \xi_i^\gamma}d\xi_i^\gamma + \sum_{\beta=1}^q \sum_{\delta=1}^{m'} w_i^{\alpha\beta} \left( \varphi_i(\xi_i, u), s(u) \right) \frac{\partial \varphi_i^\beta}{\partial u_\delta} du_\delta + \sum_{\eta=1}^r \sum_{\delta=1}^{m'} v_i^{\alpha\eta}\left( \varphi_i(\xi_i, u),s(u)\right)\frac{\partial s^\eta}{\partial u_\delta} du_\delta,
\end{align*}
so that the following is satisfied: for $\alpha=1,...,q$,
\begin{align}
\sum_{\beta=1}^q\sum_{\gamma=1}^n c_i^{\alpha\beta}(\xi_i, u) w_i'^{\beta\gamma}(\xi_i, u)d\xi_i^\gamma &= \sum_{\beta=1}^n \sum_{\gamma=1}^n w_i^{\alpha\beta}\left(\varphi_i(\xi_i, u), s(u)\right)\frac{\partial \varphi_i^\beta}{\partial \xi_i^\gamma} d\xi_i^\gamma\label{uc4}
\end{align}
{\small{\begin{align}
\sum_{\beta=1}^q\sum_{\delta=1}^{m'} c_i^{\alpha\beta}(\xi_i, u) v_i'^{\beta\delta}(\xi_i, u)du_\delta & = \sum_{\delta=1}^{m'}\left(  \sum_{\beta=1}^q w_i^{\alpha\beta}\left( \varphi_i(\xi_i, u), s(u)\right)\frac{\partial \varphi_i^\beta}{\partial u_\delta} + \sum_{\eta=1}^r v_i^{\alpha\eta}\left(\varphi_i(\xi_i, u),s(u)\right)\frac{\partial s^\eta}{\partial u_\delta} \right) du_\delta,\label{uc5}
\end{align}}}

First we prove the existence of formal solution of $(\ref{uc2})-(\ref{uc5})$. We recall Notation \ref{te27}. Then $(\ref{uc2})-(\ref{uc5})$ are equivalent to the following systems of congruences
\begin{align}
\varphi_i^\mu \left(f_{ij}'(\xi_j, u), u \right)& \equiv_\mu f_{ij}\left(\varphi_j^\mu(\xi_j, u), s^\mu(u)\right) \label{uc6}\\
\sum_{\beta=1}^q c_i^{\alpha\beta\mu}\left(f_{ij}'(\xi_j, u),u\right) h_{ij}'^{\beta\gamma}& (\xi_j, u)  \equiv_\mu \sum_{\beta=1}^q  h_{ij}^{\alpha\beta}\left(\varphi_j^\mu(\xi_j, u), s^\mu(u) \right) c_j^{\beta\gamma\mu}(\xi_j, u) \label{uc7} \\
\sum_{\beta=1}^q\sum_{\gamma=1}^n c_i^{\alpha\beta \mu}(\xi_i, u) w_i'^{\beta\gamma}(\xi_i, u)d\xi_i^\gamma &= \sum_{\beta=1}^n \sum_{\gamma=1}^n w_i^{\alpha\beta}\left(\varphi_i^\mu (\xi_i, u), s^\mu (u) \right)\frac{\partial \varphi_i^{ \beta \mu}}{\partial \xi_i^\gamma}  d\xi_i^\gamma \label{uc8}
\end{align}
\begin{align}
\sum_{\beta=1}^q\sum_{\delta=1}^{m'} c_i^{\alpha\beta\mu}(\xi_i, u) v_i'^{\beta\delta}(\xi_i, u)du_\delta & \equiv_\mu \sum_{\delta=1}^{m'}\left(  \sum_{\beta=1}^q w_i^{\alpha\beta}\left(\varphi_i^\mu(\xi_i, u), s^\mu(u) \right)\frac{\partial \varphi_i^{\beta\mu}}{\partial u_\delta} + \sum_{\eta=1}^m v_i^{\alpha\eta}\left( \varphi_i^\mu(\xi_i, u),s^\mu(u) \right)\frac{\partial s^{\eta\mu}}{\partial u_\delta} \right) du_\delta \label{uc9}
\end{align}
for each $\mu=1,2,...$.  We shall construct $\varphi_i^\mu, s^\mu$, and $c_i^{\alpha\beta \mu}$ satisfying $(\ref{uc6})_\mu- (\ref{uc9})_\mu$ by induction on $\mu$. We assume that $\varphi_i^{\mu-1}, s^{\mu-1}$ and $c_i^{\alpha\beta(\mu-1)}$ are already determined. Then we define homogenous polynomials $\Gamma_{ij|\mu}^\alpha, A_{ij|\mu}^{\alpha\beta}, W_{i|\mu}^{\alpha \gamma}$ and $V_{i|\mu-1}^{\alpha \delta}$ of degree $\mu$ by the following congruences:
{\small{\begin{align}
\Gamma_{ij|\mu}^\alpha &\equiv_\mu \varphi_i^{\alpha(\mu-1)}\left(f_{ij}'(\xi_j, u), u \right)- f_{ij}^{\alpha}\left( \varphi_j^{\mu-1}(\xi_j, u), s^{\mu-1}(u)\right) \label{uc10} \\
\sum_{\beta,\sigma=1}^q A_{ij|\mu}^{\alpha\beta} c_{0i}^{\beta\sigma} h_{0ij}'^{\sigma \gamma} &\equiv_\mu \sum_{\beta=1}^q c_i^{\alpha\beta (\mu-1)}\left(f_{ij}'(\xi_j, u), u\right) h_{ij}'^{\beta\gamma}(\xi_j, u) - \sum_{\beta=1}^q h_{ij}^{\alpha\beta} \left(\varphi_j^{\mu-1}(\xi_j, u\right), s^{\mu-1}(u))c_j^{\beta\gamma(\mu-1)}(\xi_j, u)\label{uc11}\\
W_{i|\mu}^{\alpha\gamma}& \equiv_{\mu} \sum_{\beta=1}^q  c_i^{\alpha\beta (\mu-1)}(\xi_i, u) w_i'^{\beta\gamma}(\xi_i, u)- \sum_{\beta=1}^q w_i^{\alpha\beta}\left(\varphi_i^{\mu-1}(\xi_i, u), s^{\mu-1}(u) \right) \frac{\partial \varphi_i^{\beta (\mu-1)}}{\partial \xi_i^\gamma} \label{uc12}
\end{align}}}
{\small{\begin{align}
V_{i|\mu-1}^{\alpha \delta} &\equiv_{\mu-1} \sum_{\beta=1}^q c_i^{\alpha\beta (\mu-1)}(\xi_i, u) v_i'^{\beta\delta}(\xi_i, u)- \sum_{\beta=1}^q w_i^{\alpha\beta}\left(\varphi_i^{\mu-1}(\xi_i, u), s^{\mu-1}(u) \right)\frac{\partial \varphi_i^{\beta (\mu-1)} }{\partial u_\delta} - \sum_{\eta=1}^m v_i^{\alpha\eta} \left(\varphi_i^{\mu-1}(\xi_i , u), s^{\mu-1}(u)\right) \frac{\partial s^{\eta (\mu-1)}}{\partial u_\delta} \label{uc13}
\end{align}}}
We set 
\begin{align*}
\Gamma_{ij|\mu}&= \sum_{\alpha=1}^n \Gamma_{ij|\mu}^{\alpha}\frac{\partial}{\partial z_i^\gamma},\,\,\,\,\,\,\,\,\,\,\,\, W_{i|\mu}^\alpha = \sum_{\gamma=1}^n W_{i|\mu}^{\alpha\gamma} dz_i^\gamma
\end{align*}

\begin{lemma}\label{uc14}
We have the following equalities$:$ for $\alpha=1,...,q$ and $\delta=1,..., m'$
\begin{align}
&\,\,\,\,\,\,\,\,\Gamma_{jk|\mu}- \Gamma_{ik|\mu} + \Gamma_{ij|\mu}=0\label{uc15}\\
W_{i|\mu}^\alpha &=  \sum_{\beta=1}^q A_{ij|\mu}^\beta w_{0i}^\beta + \sum_{\beta=1}^q h_{0ij}^{\alpha\beta} W_{j|\mu}^\beta- \mathcal{L}_{\Gamma_{ij|\mu}}\left(w_{0i}^\alpha \right) \label{uc16}\\
\sum_{\beta=1}^q w_{0i}^1 \wedge \cdots \wedge &\overbrace{\left(-  W_{i|\mu}^\beta \right)}^{\beta-\textnormal{th}} \wedge \cdots \wedge w_{0i}^q \wedge d w_{0i}^\alpha + w_{0i}^1 \wedge \cdots \wedge w_{0i}^q \wedge d \left( - W_{i|\mu}^\alpha \right)=0 \label{uc17}\\
V_{i|\mu-1}^{\alpha \delta}= \sum_{\beta=1}^q h_{0ij}^{\alpha\beta}& V_{j|\mu-1}^{\beta\delta} - i_{\sum_{\beta=1}^n \frac{\partial \Gamma_{ij|\mu}^\beta}{\partial u_\delta} \frac{\partial }{\partial z_i^\beta}} \left( w_{0i}^\alpha \right) = \sum_{\beta=1}^q h_{0ij}^{\alpha \beta} V_{j|\mu-1}^{\beta\delta}- i_{\frac{\partial \Gamma_{ij|\mu}}{\partial u_\delta}}\left(w_{0i}^\alpha \right) \label{uc18}
\end{align}
{\small{\begin{align} \label{uc19}
 \sum_{\beta=1}^q (-1)^{q-\beta} \left(- V_{i|\mu-1}^{\beta \delta}  \right) w_{0i}^1 \wedge \cdots \wedge &\widehat{ w_{0i}^{\beta} } \wedge \cdots \wedge  w_{0i}^q \wedge d w_{0i}^\alpha   + w_{0i}^1 \wedge \cdots \wedge w_{0i}^q \wedge d \left(- V_{i|\mu-1}^{\alpha \delta} \right)    - w_{0i}^1 \wedge \cdots \wedge w_{0i}^q \wedge \frac{\partial \left(- W_{i|\mu}^\alpha \right) }{\partial u_\delta }=0 
\end{align}}}
\begin{align}\label{uc20}
\frac{\partial V_{i|\mu-1}^{\alpha \delta}}{\partial u_\eta}= \frac{\partial  V_{i|\mu-1}^{\alpha \eta}}{\partial u_\delta}\,\,\,\,\,\,\,\,\,\textnormal{for $\delta\ne \eta$}
\end{align}

\end{lemma}

\begin{proof}
$(\ref{uc15})$ follows from \cite{Kod05} p.292-p.294. $(\ref{uc16})-(\ref{uc17})$ follows from Part I. We prove $(\ref{uc18})$. In fact, from $(\ref{uc13})$ and $(\ref{uc11})$ and $(\ref{uc10})$, and $(\ref{uc12})$,
{\small{\begin{align*}
&V_{i|\mu-1}^{\alpha\delta} \equiv_{\mu-1} \sum_{\beta=1}^q c_i^{\alpha\beta (\mu-1)}(\xi_i, u) v_i'^{\beta\delta}(\xi_i, u)- \sum_{\beta=1}^q w_i^{\alpha\beta}\left( \varphi_i^{\mu-1}(\xi_i, u), s^{\mu-1}(u) \right)\frac{\partial \varphi_i^{\beta (\mu-1)} }{\partial u_\delta} - \sum_{\eta=1}^m v_i^{\alpha\eta} \left( \varphi_i^{\mu-1}(\xi_i , u), s^{\mu-1}( u ) \right) \frac{\partial s^{\eta (\mu-1)}}{\partial u_\delta} \\
         &\equiv_{\mu-1} \sum_{\beta=1}^q c_i^{\alpha\beta(\mu-1)}\left( \sum_{\sigma=1}^q h_{ij}'^{\beta\sigma}v_j'^{\sigma\delta} - \sum_{\sigma=1}^n w_i'^{\beta\sigma}\left(f_{ij}', u \right)\frac{\partial f_{ij}'^\sigma}{\partial u_\delta}  \right) - \sum_{\beta=1}^q w_i^{\alpha\beta}\left( \varphi_i^{\mu-1}\left(f_{ij}',u \right), s^{\mu-1}\right)\frac{\partial \varphi_i^{\beta (\mu-1)} }{\partial u_\delta} \left(f_{ij}', u \right) \\
         & - \sum_{\eta=1}^{m} v_i^{\alpha\eta} \left( \varphi_i^{\mu-1} \left( f_{ij}',u \right) , s^{\mu-1} \right) \frac{\partial s^{\eta (\mu-1)}}{\partial u_\delta}\\
         &\equiv_{\mu-1} \sum_{\sigma=1}^q \left( \sum_{\beta, \gamma=1}^q A_{ij|\mu}^{\alpha \beta} h_{0ij}'^{\beta\gamma} c_{0j}^{\gamma \sigma} + \sum_{\beta=1}^q h_{ij}^{\alpha\beta}\left( \varphi_j^{\mu-1}, s^{\mu-1} \right) c_j^{\beta \sigma (\mu-1)}   \right) v_j'^{\sigma \delta} -  \sum_{\beta=1}^q \sum_{\sigma=1}^n c_i^{\alpha\beta(\mu-1)} w_i'^{\beta\sigma}\left(f_{ij}', u\right) \frac{\partial f_{ij}'^\sigma}{\partial u_\delta} \\
         & - \sum_{\beta=1}^q w_i^{\alpha\beta}\left(\Gamma_{ij|\mu}+ f_{ij}\left( \varphi_j^{\mu-1}, s^{\mu-1}\right)  , s^{\mu-1}\right) \frac{\partial \varphi_i^{\beta(\mu-1)}}{\partial u_\delta} - \sum_{\eta=1}^{m} v_i^{\alpha\eta}\left( \Gamma_{ij|\mu}+ f_{ij}\left( \varphi_j^{\mu-1}, s^{\mu-1}\right) , s^{\mu-1} \right) \frac{\partial s^{\eta(\mu-1)}}{\partial u_\delta}\\
         &\equiv_{\mu-1} \sum_{\beta=1}^q h_{ij}^{\alpha\beta} \left(\varphi_j^{\mu-1}, s^{\mu-1}\right) \left( V_{j|\mu-1}^{\beta \delta} +\sum_{\gamma=1}^n w_j^{\beta\gamma}\left( \varphi_j^{\mu-1}, s^{\mu-1} \right)\frac{\partial \varphi_j^{\gamma(\mu-1)}}{\partial u_\delta}   + \sum_{\eta=1}^m v_j^{\beta\eta} \left( \varphi_j^{\mu-1}, s^{\mu-1} \right) \frac{\partial s^{\eta(\mu-1)}}{\partial u_\delta} \right) \\
         &- \sum_{\sigma=1}^n \left(W_{i|\mu}^{\alpha\sigma} + \sum_{\beta=1}^q w_i^{\alpha\beta} \left( \varphi_i^{\mu-1}, s^{\mu-1} \right)\frac{\partial \varphi_i^{\beta(\mu-1)}}{\partial \xi_i^\sigma} \right) \frac{\partial f_{ij}'^\sigma}{\partial u_\delta}\\
         & - \sum_{\beta=1}^q w_i^{\alpha\beta}\left(f_{ij}\left(\varphi_j^{\mu-1}, s^{\mu-1}\right), s^{\mu-1} \right) \frac{\partial \varphi_i^{\beta(\mu-1)}}{\partial u_\delta}\left( f_{ij}', u \right) - \sum_{\eta=1}^{m} v_i^{\alpha\eta}\left(f_{ij} \left(\varphi_j^{\mu-1}, s^{\mu-1}\right), s^{\mu-1} \right)\frac{\partial s^{\eta(\mu-1)}}{\partial u_\delta}\\
         &\equiv_{\mu-1} \sum_{\beta=1}^q h_{0ij}^{\alpha\beta} V_{j|\mu-1}^{\beta\delta} + \sum_{\sigma,\gamma=1}^n w_i^{\alpha \sigma}\left(f_{ij}\left( \varphi_j^{\mu-1}, s^{\mu-1}\right), s^{\mu-1}\right) \frac{\partial f_{ij}^\sigma}{\partial z_j^\gamma} \left( \varphi_j^{\mu-1}, s^{\mu-1} \right) \frac{\partial \varphi_j^{\gamma(\mu-1)}}{\partial u_\delta} \\
          &+ \sum_{\sigma=1}^n  \sum_{\eta=1}^m w_i^{\alpha\sigma}\left(f_{ij}\left( \varphi_j^{\mu-1}, s^{\mu-1}\right), s^{\mu-1}\right) \frac{\partial f_{ij}^\sigma}{\partial t_\eta} \left( \varphi_j^{\mu-1}, s^{\mu-1}\right) \frac{\partial s^{\eta(\mu-1)}}{\partial u_\delta}           + \sum_{\eta=1}^m v_i^{\alpha\eta}\left(f_{ij}\left( \varphi_j^{\mu-1}, s^{\mu-1}\right), s^{\mu-1} \right) \frac{\partial s^{\eta(\mu-1)}}{\partial u_\delta}\\
          & - \sum_{\sigma=1}^n \sum_{\beta=1}^q w_i^{\alpha\beta} \left( \varphi_i^{\mu-1}, s^{\mu-1}\right) \frac{\partial \varphi_i^{\beta( \mu-1)}}{\partial \xi_i^\sigma} \frac{\partial f_{ij}'^\sigma }{\partial u_\delta} \\
                & - \sum_{\beta=1}^q w_i^{\alpha\beta}\left(f_{ij} \left(\varphi_j^{\mu-1}, s^{\mu-1} \right), s^{\mu-1} \right) \frac{\partial \varphi_i^{\beta(\mu-1)}}{\partial u_\delta}\left(f_{ij}', u\right) - \sum_{\eta=1}^{m} v_i^{\alpha\eta}\left(f_{ij}\left(\varphi_j^{\mu-1}, s^{\mu-1} \right), s^{\mu-1} \right)\frac{\partial s^{\eta(\mu-1)}}{\partial u_\delta}\\
                &\equiv_{\mu-1} \sum_{\beta=1}^q h_{0ij}^{\alpha\beta} V_{j|\mu-1}^{\beta\delta} +  \sum_{\sigma=1}^n w_i^{\alpha\sigma}\left(f_{ij} \left(\varphi_j^{\mu-1}, s^{\mu-1} \right), s^{\mu-1} \right)\frac{\partial f_{ij}^\sigma \left( \varphi_j^{\mu-1}, s^{\mu-1}\right)}{\partial u_\delta}   \\
                & - \sum_{\sigma=1}^n \sum_{\beta=1}^q w_i^{\alpha\beta} \left( \varphi_i^{\mu-1}, s^{\mu-1}\right) \frac{\partial \varphi_i^{\beta(\mu-1)}}{\partial \xi_i^\sigma} \frac{\partial f_{ij}'^{\sigma}}{\partial u_\delta}   - \sum_{\beta=1}^q w_i^{\alpha\beta}\left( \varphi_i^{(\mu-1)}\left(f_{ij}', u \right)- \Gamma_{ij|\mu}, s^{\mu-1}\right) \frac{\partial \varphi_i^{\beta(\mu-1)}}{\partial u_\delta}\left(f_{ij}', u \right)\\
                &\equiv_{\mu-1} \sum_{\beta=1}^q h_{0ij}^{\alpha\beta} V_{j|\mu-1}^{\beta\delta} + \sum_{\beta=1}^n w_i^{\alpha\beta} \left( \varphi_i^{(\mu-1)}\left(f_{ij}', u \right) - \Gamma_{ij|\mu} , s^{\mu-1}\right)\frac{\partial\left( \varphi_i^{\beta (\mu-1)}\left(f_{ij}', u \right)- \Gamma_{ij|\mu}^\beta \right)}{\partial u_\delta} - \sum_{\beta=1}^q w_i^{\alpha\beta} \left( \varphi_i^{\mu-1}, s^{\mu-1}\right) \frac{\partial \varphi_i^{\beta(\mu-1)}\left(f_{ij}', u\right)}{\partial u_\delta} \\
                & \equiv_{\mu-1} \sum_{\beta=1}^q h_{0ij}^{\alpha\beta} V_{j|\mu-1}^{\beta\delta} - \sum_{\beta=1}^n w_{0i}^{\alpha\beta} \frac{\partial \Gamma_{ij|\mu}^\beta}{\partial u_\delta}
\end{align*}}}

We prove $(\ref{uc19})$ and $(\ref{uc20})$. In fact, from $\omega_i^1 \wedge \cdots \wedge \omega_i^q \wedge d\omega_i^\alpha=0$, i.e. 
\begin{align}\label{uc21}
&\left( \sum_{\beta_1=1}^n w_i^{1\beta_1}(z_i,t) dz_i^{\beta_1}  + \sum_{\eta_1=1}^m v_i^{1\eta_1} (z_i,t) dt_{\eta_1} \right) \wedge \cdots \wedge \left(\sum_{\beta_q=1}^n w_i^{q\beta_q}(z_i,t) dz_i^{\beta_q} + \sum_{\eta_q=1}^m v_i^{q\eta_q} (z_i,t) d t_{\eta_q}    \right) \\
&\wedge d_{\mathcal{M}}\left(\sum_{\beta=1}^q w_i^{\alpha\beta}(z_i,t) dz_i^\beta + \sum_{\eta=1}^m v_i^{\alpha \eta}(z_i,t) dt_\eta   \right) =0 \notag
\end{align} 
By pulling back $(\ref{uc21})$ under $\left( \varphi_i^{\mu-1}, s^{\mu-1}\right)$, we have
\begin{align*}
&\left( \sum_{\beta_1=1}^n w_i^{1\beta_1}\left(\varphi_i^{\mu-1},s^{\mu-1}\right) d_\mathcal{M} \varphi_i^{\beta_1(\mu-1)}  + \sum_{\eta_1=1}^m v_i^{1\eta_1} (z_i,t) d_\mathcal{M} s_i^{\eta_1(\mu-1)} \right) \wedge \cdots \\
& \wedge \left(\sum_{\beta_q=1}^n w_i^{q\beta_q}\left( \varphi_i^{\mu-1},s^{\mu-1} \right) d_\mathcal{M} \varphi_i^{\beta_q (\mu-1)}+ \sum_{\eta_q=1}^m v_i^{q\eta_q} \left( \varphi_i^{\mu-1}, s^{\mu-1} \right) d_\mathcal{M} s_i^{\eta_q(\mu-1)}    \right) \\
&\wedge d_{\mathcal{M}'}\left(\sum_{\beta=1}^q w_i^{\alpha\beta}\left( \varphi_i^{\mu-1},s^{\mu-1}\right) d_\mathcal{M} \varphi_i^{\beta(\mu-1)} + \sum_{\eta=1}^m v_i^{\alpha \eta} \left( \varphi_i^{\mu-1},s^{\mu-1}\right) d_\mathcal{M} s_i^{\eta(\mu-1)}   \right) =0
\end{align*}
Then we have from $(\ref{uc12})$ and $(\ref{uc13})$,
\begin{align*}
&\left(\sum_{\beta=1}^n c_i^{1\beta(\mu-1)} \left(\sum_{\gamma=1}^n w_i'^{\beta \gamma}d\xi_i^\gamma + \sum_{\delta=1}^{m'} v_i'^{\beta \delta}du_\delta \right)  - \sum_{\gamma=1}^n W_{i|\mu}^{1\gamma} dz_i^\gamma  - \sum_{\delta=1}^{m'} V_{i|\mu-1}^{1\delta} du_\delta    \right) \wedge \cdots  \\
&\wedge \left(\sum_{\beta=1}^n c_i^{q\beta(\mu-1)} \left(\sum_{\gamma=1}^n w_i'^{\beta \gamma}d\xi_i^\gamma + \sum_{\delta=1}^{m'} v_i'^{\beta \delta}du_\delta \right)  - \sum_{\gamma=1}^n W_{i|\mu}^{q\gamma} dz_i^\gamma  - \sum_{\delta=1}^{m'} V_{i|\mu-1}^{q\delta} du_\delta    \right)  \\
& \wedge d_{\mathcal{M}'} \left(\sum_{\beta=1}^n c_i^{\alpha\beta(\mu-1)} \left(\sum_{\gamma=1}^n w_i'^{\beta \gamma}d\xi_i^\gamma + \sum_{\delta=1}^{m'} v_i'^{\beta \delta}du_\delta \right)  - \sum_{\gamma=1}^n W_{i|\mu}^{\alpha \gamma} dz_i^\gamma  - \sum_{\delta=1}^{m'} V_{i|\mu-1}^{\alpha \delta} du_\delta       \right)  = 0 
\end{align*}
In other words, we have
{\small{\begin{align}
&\left(\sum_{\beta=1}^n c_i^{1\beta(\mu-1)} \left(w_i'^\beta + v_i'^\beta \right)  - W_{i|\mu}^{1}   - \sum_{\delta=1}^{m'} V_{i|\mu-1}^{1\delta} du_\delta    \right) \wedge \cdots  \wedge \left(\sum_{\beta=1}^n c_i^{q\beta(\mu-1)} \left(w_i'^\beta + v_i'^\beta \right)  -  W_{i|\mu}^q   - \sum_{\delta=1}^{m'} V_{i|\mu-1}^{q\delta} du_\delta    \right) \label{uc22} \\
& \wedge d_{\mathcal{M}'} \left(\sum_{\beta=1}^n c_i^{\alpha\beta(\mu-1)} \left(w_i'^\beta + v_i'^\beta \right)  - W_{i|\mu}^\alpha  - \sum_{\delta=1}^{m'} V_{i|\mu-1}^{\alpha \delta} du_\delta       \right)  = 0 \notag
\end{align}}}
Since $\left(w_i'^1 + v_i'^1 \right)\wedge \cdots \left(w_i'^q+ v_i'^q \right) \wedge d_{\mathcal{M}'}\left(w_i'^\beta + v_i'^\beta \right)=0,\beta=1,...,q$, $(\ref{uc22})$ implies that in particular
\begin{align*}
&\sum_{\beta=1}^q w_{0i}^1 \wedge \cdots \wedge \overbrace{\left(-  W_{i|\mu}^\beta \right)}^{\beta-\textnormal{th}} \wedge \cdots \wedge w_{0i}^q \wedge d w_{0i}^\alpha + w_{0i}^1 \wedge \cdots \wedge w_{0i}^q \wedge d \left( - W_{i|\mu}^\alpha \right) - \sum_{\delta=1}^{m'} w_{0i}^1 \wedge \cdots \wedge w_{0i}^q \wedge \frac{\partial \left(- W_{i|\mu}^\alpha \right) }{\partial u_\delta } \wedge d u_\delta \\
&+ \sum_{\delta=1}^{m'} \sum_{\beta=1}^q (-1)^{q-\beta} \left(- V_{i|\mu-1}^{\beta \delta}  \right) w_{0i}^1 \wedge \cdots \wedge \widehat{ w_{0i}^{\beta} } \wedge \cdots \wedge  w_{0i}^q \wedge d w_{0i}^\alpha \wedge du_\delta  + \sum_{\delta=1}^{m'} w_{0i}^1 \wedge \cdots \wedge w_{0i}^q \wedge d \left(- V_{i|\mu-1}^{\alpha \delta} \right) \wedge du_\delta \\
& - w_{0i}^1 \wedge \cdots \wedge w_{0i}^q \wedge \left( \sum_{\delta, \eta=1}^{m'} \frac{\partial V_{i|\mu-1}^{\alpha \delta}}{\partial u_\eta} du_\eta \wedge du_\delta     \right) =0 
\end{align*}
Hence for $\delta=1,...,m'$, we have
{\small{\begin{align*}
-w_{0i}^1 \wedge \cdots \wedge w_{0i}^q \wedge \left( -\frac{\partial W_{i|\mu}^\alpha}{\partial u_\delta}   \right) +  \sum_{\beta=1}^q (-1)^{q-\beta} \left(- V_{i|\mu-1}^{\beta \delta}  \right) w_{0i}^1 \wedge \cdots \wedge \widehat{ w_{0i}^{\beta} } \wedge \cdots \wedge  w_{0i}^q \wedge d w_{0i}^\alpha + w_{0i}^1 \wedge \cdots \wedge w_{0i}^q \wedge d \left(- V_{i|\mu-1}^{\alpha \delta} \right)=0
\end{align*}}}
and 
\begin{align*}
 \sum_{\delta, \eta=1}^{m'} \frac{\partial V_{i|\mu-1}^{\alpha \delta}}{\partial u_\eta} du_\eta \wedge du_\delta =0
\end{align*}
Then for $\delta \ne \eta$, we have
\begin{align*}
\frac{\partial V_{i|\mu-1}^{\alpha \delta}}{\partial u_\eta}= \frac{\partial  V_{i|\mu-1}^{\alpha \eta}}{\partial u_\delta}
\end{align*}
This completes the proof of Lemma \ref{uc14}.

\end{proof}

Our purpose is to determine $\varphi^\mu= \varphi^{\mu-1}+ \varphi_{i|\mu}, s^\mu= s^{\mu-1}+ s_\mu$, and $c_i^{\alpha\beta(\mu-1)}+ c_{i|\mu}^{\alpha\beta}$ satisfying $(\ref{uc6})_\mu-(\ref{uc9})_\mu$.

\begin{lemma} \label{uc27}
$(\ref{uc6})_\mu -(\ref{uc9})_\mu$ are equivalent to the following equalities$:$ for $\alpha=1,...,q $ and $\delta=1,...,m'$,
\begin{align}
\Gamma_{ij|\mu}&= \varphi_{j|\mu}-  \varphi_{i | \mu} + \sum_{\lambda=1}^{m} s_\mu^\lambda \rho_{ij\lambda} \label{uc23}\\
\sum_{\beta, \sigma=1}^q A_{ij|\mu}^{\alpha\beta} h_{0ij}^{\beta\sigma} c_{0j}^{\sigma \gamma} &=  -\sum_{\beta=1}^q c_{i|\mu}^{\alpha\beta} h_{0ij}'^{\beta\gamma} + \sum_{\beta=1}^q \left[ \varphi_{j|\mu} , h_{0ij}^{\alpha\beta}  \right]  c_{0j}^{\beta\gamma} + \sum_{\beta=1}^q \sum_{\lambda=1}^m \frac{\partial h_{ij}^{\alpha\beta}}{\partial t_\lambda }|_{t=0} s_\mu^\lambda c_{0j}^{\beta \gamma} - \sum_{\beta=1}^q h_{0ij}^{\beta\gamma} c_{j|\mu}^{\beta\gamma}    \label{uc24}   \\
W_{i|\mu}^{\alpha } &=  - \sum_{\gamma=1}^n \sum_{\beta=1}^q c_{i|\mu}^{\alpha\beta} w_{0i}'^{\beta} + \mathcal{L}_{\varphi_{i|\mu}}\left( w_{0i}^\alpha  \right) - \sum_{\lambda=1}^m s_\mu^\lambda \beta_{i\lambda}^\alpha   \label{uc25} \\
V_{i|\mu-1}^{\alpha \delta } &= i_{\frac{\partial \varphi_{i|\mu}}{\partial u_\delta}}\left( w_{0i}^\alpha \right)   - \sum_{\lambda=1}^m  \gamma_{i\lambda}^\alpha \frac{\partial s_\mu^\lambda}{ \partial u_\delta }   \label{uc26}
\end{align}
where
\begin{align*}
\varphi_{i|\mu}=\sum_{\alpha=1}^n \varphi_{i|\mu}^\alpha \frac{\partial}{\partial z_i^\alpha},\,\,\,\,\,\,\,\rho_{ij\lambda}=\sum_{\alpha=1}^n \frac{\partial f_{ij}^\alpha}{\partial t_\lambda}|_{t=0} \frac{\partial}{\partial z_i^\alpha},\,\,\,\,\,\,\,\,\,\,\,\, \beta_{i\lambda}^\alpha = - \frac{\partial w_i^\alpha (z_i,t) }{\partial t_\lambda}|_{t=0},\,\,\,\,\,\,\,\, \gamma_{i\lambda}^\alpha =-i_{\frac{\partial}{\partial t_\lambda}}\left( v_i^\alpha(z_i,t)   \right)|_{t=0}= - v_{0i}^{\alpha \lambda}
\end{align*}
\end{lemma}

\begin{proof}
$(\ref{uc23})$ follows from \cite{Kod05} p.290. $(\ref{uc24})$ and $(\ref{uc25})$ follows from Part I. We prove $(\ref{uc26})$. In fact, for $\delta=1,...,m'$,
{\small{\begin{align*}
&\sum_{\beta=1}^q \left( c_i^{\alpha\beta (\mu-1)}+ c_{i|\mu}^{\alpha\beta} \right) v_i'^{\beta\delta}- \sum_{\beta=1}^q w_i^{\alpha\beta}\left( \varphi_i^{\mu-1} +  \varphi_{i|\mu}, s^{\mu-1}+ s_\mu \right) \frac{\partial \left( \varphi_i^{\beta(\mu-1)}+ \varphi_{i|\mu}^\beta \right)}{\partial u_\delta} - \sum_{\eta=1}^r v_i^{\alpha\eta}\left(\varphi_i^{\mu-1}+ \varphi_{i|\mu}, s^{\mu-1}+ s_\mu \right) \frac{\partial \left(s^{\eta(\mu-1)}+ s_\mu^\eta \right) }{\partial u_\delta} \\
&\equiv_{\mu-1} V_{i|\mu-1}^{\alpha\delta}  - \sum_{\beta=1}^q w_{0i}^{\alpha\beta} \frac{\partial \varphi_{i|\mu}^\beta}{\partial u_\delta} -\sum_{\lambda=1}^r v_{0i}^{\alpha \lambda } \frac{\partial s_\mu^\lambda }{\partial u_\delta} \equiv_{\mu-1} V_{i|\mu-1}^{\alpha\delta}  - \sum_{\beta=1}^q w_{0i}^{\alpha\beta} \frac{\partial \varphi_{i|\mu}^\beta}{\partial u_\delta} -\sum_{\lambda=1}^r v_{0i}^{\alpha \lambda } \frac{\partial s_\mu^\lambda }{\partial u_\delta} 
\end{align*}}}
This completes the proof of Lemma \ref{uc27}.
\end{proof}

For $\delta=1,..., m'$, we define an element $\overline{\frac{\partial W_{i|\mu}}{\partial u_\delta}}\in \Gamma\left( U_i, \mathscr{H}om_{\mathcal{O}_M} \left( \mathcal{N}_{\mathcal{F}_0}^*, \frac{\Omega_M^1}{\mathcal{N}_{\mathcal{F}_0}^*}  \right) \right)$ by
\begin{align*}
\overline{\frac{\partial W_{i|\mu}}{\partial t_\delta}}: \Gamma\left( U_i, \mathcal{N}_{\mathcal{F}_0}^*  \right) &\to \Gamma\left(  U_i  ,   \frac{\Omega_M^1}{\mathcal{N}_{\mathcal{F}_0}^*}   \right) \\
                                       w_{0i}^\alpha &\mapsto \overline{\frac{\partial W_{i|\mu}^\alpha}{\partial t_\delta} }
\end{align*}
and linearly extends to $\Gamma\left( U_i, \mathcal{N}_{\mathcal{F}_0}^* \right)$. For $\delta=1,...,q'$, we define an element $V_{i|\mu}^\delta\in \Gamma\left( U_i, \mathscr{H}om_{\mathcal{O}_M}\left( \mathcal{N}_{\mathcal{F}_0}^*, \mathcal{O}_M \right)\right)$ by
\begin{align*}
V_{i|\mu}^\delta :\Gamma\left( U_i,  \mathcal{N}_{\mathcal{F}_0}^* \right) &\to \Gamma\left( U_i, \mathcal{O}_M \right)\\
                  w_{0i}^\alpha &\mapsto V_{i|\mu}^{\alpha \delta}
\end{align*}
and linearly extends to $\Gamma\left( U_i, \mathcal{N}_{\mathcal{F}_0}^* \right)$. Then $(\ref{uc15})-(\ref{uc19})$ implies that
\begin{align}\label{ucc12}
\left(  \left\{ -  \overline{\frac{\partial W_{i|\mu}}{\partial t_\delta}}\right\}   , \left\{- V_{i|\mu}^\delta \right\}  \left\{  \frac{\partial \Gamma_{ij|\mu}}{\partial t_\delta} \right\} \right) \in C^0\left( \mathcal{U}, \mathscr{H}om_{\mathcal{O}_M}\left( \mathcal{N}_{\mathcal{F}_0}^* ,  \frac{\Omega_M^1}{\mathcal{N}_{\mathcal{F}_0}^*} \right)   \right) \bigoplus C^0\left( \mathcal{U}, \mathscr{H}om_{\mathcal{O}_M}\left( \mathcal{N}_{\mathcal{F}_0}^*, \mathcal{O}_M \right) \right) \bigoplus C^1 \left( \mathcal{U}, \Theta_M  \right)
\end{align}
defines a $1$-cocycle in the following \v Cech resolution of $\tilde{\mathcal{N}}_{\mathcal{F}_0}^{* \bullet}$:
{\small{\begin{center}
$\begin{CD}
\cdots \\
@AG_2AA \\
C^0\left( \mathcal{U}, \left(\left(\mathcal{N}_{\mathcal{F}_0}^* \right)^*\otimes \tilde{\mathcal{S}}^2 \right) \bigoplus \left( \mathcal{N}_{\mathcal{F}_0}^* \right)^*\otimes \tilde{\mathcal{S}}^1 \right) @>-\delta >> \cdots \\
@AG_1AA @AG_1AA \\
C^0\left(\mathcal{U},\left( \left(\mathcal{N}_{\mathcal{F}_0}^* \right)^* \otimes \frac{\Omega_{M}^1}{ \mathcal{N}_{\mathcal{F}_0}^*} \right) \bigoplus \left( \mathcal{N}_{\mathcal{F}_0}^* \right)^* \right) @>\delta >> C^1\left(\mathcal{U}, \left( \left(\mathcal{N}_{\mathcal{F}_0}^* \right)^* \otimes \frac{\Omega_{M}^1}{ \mathcal{N}_{\mathcal{F}_0}^*} \right) \bigoplus   \left( \mathcal{N}_{\mathcal{F}_0}^* \right)^* \right) @>-\delta>> \cdots \\
@AG_0AA @AG_0AA @AG_0AA \\
C^0\left(\mathcal{U}, \Theta_{M}\right) @>-\delta>> C^1\left(\mathcal{U}, \Theta_{M}\right) @>\delta>> C^2\left(\mathcal{U}, \Theta_{M}\right) @>-\delta>> \cdots
\end{CD}$
\end{center}}}
By hypothesis that the unfolded Kodaira-Spencer map $\varphi_0: T_0 (B)\to \mathbb{H}^1\left(M, \tilde{\mathcal{N}}_\mathcal{F}^{*\bullet}\right)$ is surjective, we can find homogeneous polynomials $\tilde{s}_{\mu-1}^{\delta \lambda}$ such that
\begin{align*}
\varphi_0\left( \sum_{\lambda=1}^r \tilde{s}_{\mu-1}^{\delta\lambda} \frac{\partial}{\partial t_\lambda}  \right)=\left( \left\{ \overline{-\frac{\partial W_{i|\mu} }{\partial u_\delta}  } \right\}, \left\{  - V_{i|\mu-1}^{ \delta}  \right\},  \left\{ \frac{ \Gamma_{ij|\mu}}{\partial u_\delta} \right\}   \right)
\end{align*}
Since we have
\begin{align*}
\varphi_0\left(\frac{\partial}{\partial t_\lambda} \right)=\left( \{\beta_{i\lambda}\} = \left\{ w_{0i}^\alpha \mapsto \overline{- \frac{\partial w_i^\alpha(z_i,t)}{\partial t_\lambda}|_{t=0} } \right\}, \{\gamma_{i\lambda}\}= \left\{w_{0i}^\alpha \mapsto  - v_{0i}^{\alpha \lambda}(z_i) \right\} , \{\rho_{ij\lambda}\}= \left\{\sum_{\alpha=1}^n \frac{\partial f_{ij}^\alpha}{\partial t_\lambda}|_{t=0}\frac{\partial}{\partial z_i^\alpha}  \right\}    \right),
\end{align*}
there exists $\left\{\tilde{\varphi}_{i|\mu-1}^\delta:=\sum_{\alpha=1}^n \tilde{\varphi}_{i|\mu-1}^{\delta\alpha} \frac{\partial}{\partial z_i^\alpha} \right\} \in C^0 \left(\mathcal{U} , \Theta_M \right)$ and $\left\{ \tilde{c}_{i|\mu-1}^\delta \right\}\in C^0\left( \mathcal{U}, \mathscr{H}om_{\mathcal{O}_M} \left( \mathcal{N}_{\mathcal{F}_0}^*, \mathcal{N}_{\mathcal{F}_0}^*    \right) \right)$ where $\tilde{c}_{i|\mu-1}^\delta$ is defined by $\tilde{c}_{i|\mu-1}^\delta:\Gamma\left( U_i, \mathcal{N}_{\mathcal{F}_0}^* \right)\to \Gamma\left(U_i, \mathcal{N}_{\mathcal{F}_0}^* \right), w_{0i}^\alpha \mapsto \sum_{\gamma=1}^q \tilde{c}_{i|\mu-1}^{\delta \alpha \gamma} w_{0i}'^\gamma$ such that
\begin{align}
\tilde{\varphi}_{i|\mu-1}^\delta- \tilde{\varphi}_{j|\mu-1}^\delta&=- \frac{\partial \Gamma_{ij|\mu}}{\partial u_\delta} +  \sum_{\lambda=1}^m \tilde{s}_{\mu-1}^{\delta\lambda} \rho_{ij\lambda} \label{uc30} \\
\mathcal{L}_{\tilde{\varphi}_{i|\mu-1}^\delta}\left(w_{0i}^\alpha\right)&= \frac{\partial W_{i|\mu}^\alpha}{\partial u_\delta} + \sum_{\gamma=1}^q \tilde{c}_{i|\mu-1}^{ \delta \alpha \gamma  } w_{0i}'^\gamma +\sum_{\lambda=1}^m \tilde{s}_{\mu-1}^{\delta \lambda} \beta_{i\lambda}^\alpha \label{uc31} \\
i_{\tilde{\varphi}_{i|\mu-1}^\delta}\left( w_{0i}^\alpha \right) &=V_{i|\mu-1}^{\alpha \delta} + \sum_{\lambda=1}^m \tilde{s}_{\mu-1}^{\delta \lambda} \gamma_{i\lambda}^\alpha \label{uc32}
\end{align}
By setting
\begin{align*}
f_{i|\mu} := \frac{1}{\mu} \sum_{\delta=1}^{m'} u_\delta \tilde{\varphi}_{i|\mu-1}^\delta,\,\,\,\,\,\,\,\,\,\,s_\mu^\lambda : = \frac{1}{\mu} \sum_{\delta=1}^{m'} u_\delta s_{\mu-1}^{\delta \lambda},\,\,\,\,\,\,\,\,\,\,\, c_{i|\mu}^{\alpha \gamma} := \frac{1}{\mu} \sum_{\delta=1}^{m'} u_\delta \tilde{c}_{i|\mu-1}^{\delta \alpha \gamma },
\end{align*}
$(\ref{uc30}),(\ref{uc31})$ and $(\ref{uc32})$ imply that
\begin{align}
\varphi_{i|\mu} - \varphi_{j|\mu} & = - \Gamma_{ij|\mu}+ \sum_{\lambda=1}^m s_\mu^{\lambda} \rho_{ij\lambda} \label{uc33}\\
\mathcal{L}_{\varphi_{i|\mu} } \left( w_{0i}^\alpha \right) &= W_{i|\mu}^\alpha + \sum_{\gamma=1}^q c_{i|\mu}^{\alpha\gamma} w_{0i}'^\gamma  + \sum_{\lambda=1}^m s_\mu^{\lambda} \beta_{i\lambda}^\alpha \label{uc34} \\
i_{\varphi_{i|\mu}}\left(w_{0i}^\alpha \right)&= \frac{1}{\mu}\sum_{\delta=1}^{r'} u_\delta V_{i|\mu-1}^{\alpha \delta} + \sum_{\lambda=1}^m s_\mu^{\lambda} \gamma_{i\lambda}^\alpha \label{uc35}
\end{align}
Then $(\ref{uc33})$ and $(\ref{uc34})$ implies $(\ref{uc23}), (\ref{uc24})$, and $(\ref{uc25})$ from Part I. It remains to show $(\ref{uc26})$. In fact, first we note that from $(\ref{uc20})$, 
\begin{align} \label{uc36}
\frac{\partial}{\partial u_\tau }\left( \frac{1}{\mu}\sum_{\delta=1}^{m'} u_\delta V_{i|\mu-1}^{\alpha \delta}  \right) & = \frac{1}{\mu} V_{i|\mu-1}^{\alpha \tau} + \frac{1}{\mu}\sum_{\delta=1}^{m'} u_\delta \frac{\partial V_{i|\mu-1}^{\alpha \delta}}{\partial u_\tau}= \frac{1}{\mu} V_{i|\mu-1}^{\alpha \tau} + \frac{1}{\mu} \sum_{\delta=1}^{m'} u_\delta \frac{\partial V_{i|\mu-1}^{\alpha \tau}}{\partial u_\delta} \\
&= \frac{1}{\mu} V_{i|\mu-1}^{\alpha \tau} + \frac{1}{\mu}\left( (\mu-1) V_{i|\mu-1}^{\alpha \tau} \right) = V_{i|\mu-1}^{\alpha \tau} \notag
\end{align}
Then $(\ref{uc35})$ and $(\ref{uc36})$ imply $(\ref{uc26})$. This completes the inductive construction of $\varphi_i^\mu, s^\mu$ and $c_i^{\alpha\beta \mu}$.

\begin{remark}\label{ucc13}
We recall from $(\ref{ucc10})$ that we have $\mathbb{H}^1\left( M, \tilde{\mathcal{N}}_{\mathcal{F}_0}^{*\bullet} \right) \cong \mathbb{H}^1\left( M, \tilde{\mathcal{E}}_{\mathcal{N}_{\mathcal{F}_0}^*}^\bullet \right)$. We will reinterpret $(\ref{ucc12})$ in $\mathbb{H}^1\left( M ,  \tilde{\mathcal{N}}_{\mathcal{F}_0}^{* \bullet} \right) $ in terms of $\mathbb{H}^1\left( M,  \tilde{\mathcal{E}}_{\mathcal{N}_{\mathcal{F}_0}^*}^\bullet \right)$ in the following \v Cech resolution of $\tilde{\mathcal{E}}_{ \mathcal{N}_{\mathcal{F}_0}^*}^\bullet:$
{\small{\begin{equation}  \label{ucc15}
\begin{CD}
\cdots \\
@AG_2AA \\
C^0\left( \mathcal{U}, \left( \mathcal{N}_{\mathcal{F}_0}^* \right)^* \otimes \tilde{\mathcal{S}}^1 \right)\bigoplus C^0\left( \mathcal{U}, \left( \mathcal{N}_{\mathcal{F}_0}^* \right)^* \otimes \tilde{\mathcal{S}}^1 \right) @>\delta>> \cdots \\
@AG_1''AA @AG_1''AA \\
C^0\left( \mathcal{U}, \left( \mathcal{N}_{\mathcal{F}_0}^* \right)^*\otimes \Omega_M^1 \right) \bigoplus C^0\left( \mathcal{U}, \left( \mathcal{N}_{\mathcal{F}_0}^* \right)^* \right) @>\delta>> C^1\left( \mathcal{U}, \left( \mathcal{N}_{\mathcal{F}_0}^* \right)^*\otimes \Omega_M^1 \right) \bigoplus C^1\left( \mathcal{U}, \left( \mathcal{N}_{\mathcal{F}_0}^* \right)^* \right) @>-\delta>> \cdots \\
@AG_0''AA @AG_0''AA @AG_0''AA \\
C^0\left(\mathcal{U} , \mathcal{E}_{\mathcal{N}_{\mathcal{F}_0}^*} \right) @>-\delta>> C^1\left( \mathcal{U}, \mathcal{E}_{\mathcal{N}_{\mathcal{F}_0}^*} \right) @>\delta>> C^2\left( \mathcal{U}, \mathcal{E}_{\mathcal{N}_{\mathcal{F}_0}^*} \right) @>- \delta>> \cdots
\end{CD}
\end{equation}}}
We set $A_{ij|\mu}\in \Gamma\left( U_{ij}, \mathscr{H}om_{\mathcal{O}_M} \left( \mathcal{N}_{\mathcal{F}_0}^*, \mathcal{N}_{\mathcal{F}_0}^* \right) \right)$ by
\begin{align}
A_{ij|\mu}: \Gamma\left( U_{ij}, \mathcal{N}_{\mathcal{F}_0}^*\right) &\to \Gamma\left( U_{ij}, \mathcal{N}_{\mathcal{F}_0}^* \right)\\
 w_{0i}^\alpha &\mapsto \sum_{\beta=1}^q A_{ij|\mu}^{\alpha\beta} w_{0i}^\beta \notag
\end{align}
and linearly extends to $\Gamma\left( U_{ij}, \mathcal{N}_{\mathcal{F}_0}^* \right)$. Then we have $\left\{ \left( \Gamma_{ij|\mu}, A_{ij|\mu}   \right)\right\} \in C^1\left( \mathcal{U}, \mathcal{E}_{\mathcal{N}_{\mathcal{F}_0}^*} \right)$ and $\delta\left( \left\{\left( \Gamma_{ij|\mu}, A_{ij|\mu} \right) \right\} \right)=0$ as in \textnormal{Part I}, and so $\delta\left( \left\{\left( \frac{\partial \Gamma_{ij|\mu}}{\partial u_\delta},  \frac{\partial A_{ij|\mu}  }{\partial u_\delta} \right\} \right) \right)=0$. Hence
{\small{\begin{align}
\left( \left\{ - \frac{\partial W_{i|\mu} }{\partial u_\delta} \right\}, \left\{ - V_{i|\mu}^\delta\right\} , \left\{ \left( \frac{\partial \Gamma_{ij|\mu}}{\partial u_\delta}, \frac{\partial A_{ij|\mu}}{\partial u_\delta} \right) \right\}  \right) \in C^0\left( \mathcal{U}, \mathscr{H}om_{\mathcal{O}_M}\left(  \mathcal{N}_{\mathcal{F}_0}^*, \Omega_M^1 \right) \right) \bigoplus C^0\left( \mathcal{U}, \mathscr{H}om_{\mathcal{O}_M} \left( \mathcal{N}_{\mathcal{F}_0}^*, \mathcal{O}_M \right) \right) \bigoplus C^1 \left( \mathcal{U}, \mathcal{E}_{\mathcal{N}_{\mathcal{F}_0}^*} \right)
\end{align}}}
defines a $1$-cocycle in the above \v Cech resolution of $\tilde{\mathcal{E}}_{\mathcal{N}_{\mathcal{F}_0}^*}^\bullet$.

\end{remark}

\subsection{Proof of convergence}\

We prove that we can choose appropriate solutions $\varphi_{i|\mu}, s_\mu$ and $c_{i|\mu}^{\alpha\beta}$ satisfying $(\ref{uc23})-(\ref{uc26})$ in each inductive step so that
\begin{align*}
s(u)&=s_1(u)+ s_2(u)+\cdots + s_\mu(u)+ \cdots\\
\varphi_i\left( \xi_i, u\right)&= \xi_i+ \varphi_{i|1}\left(\xi_i, u\right) +  \varphi_{i|2}\left( \xi_i, u \right) + \cdots + \varphi_{i|\mu}(\xi_i, u) + \cdots \\
c_i^{\alpha\beta}(\xi_i, u)&= c_{0i}^{\alpha\beta}(\xi_i) + c_{i|1}^{\alpha\beta}(\xi_i,u) + \cdots + c_{i|\mu}^{\alpha\beta} (\xi_i, u) + \cdots
\end{align*}
converge absolutely and uniformly for $|u|< \epsilon$ if $\epsilon>0$ is sufficiently small. 

Consider a power series 
\begin{align*}
g(u)=\sum g_{v_1,..., v_{m'}} u_1^{v_1} u_2^{v_2} \cdots u_{m'}^{v_{m'}}
\end{align*}
whose coefficients $g_{v_1,..., v_{m'}}$ are vectors, and a power series
\begin{align*}
a(u)=\sum a_{v_1,..., v_{m'}} u_1^{v_1} u_2^{v_2} \cdots u_{m'}^{v_{m'}}
\end{align*}
with non-negative coefficients $a_{v_1\cdots v_{m'}}\geq 0$. We recall in Part I that
\begin{align}\label{ncc1}
\textnormal{We indicate by writing $g(u)\ll a(u)$ that $\left| g_{v_1,..., v_{m'}} \right| \leq a_{v_1\cdots v_{m'}}$}
\end{align}
We recall that
\begin{align}\label{ncc2}
A(u)= \frac{b}{16c} \sum_{\mu=1}^\infty \frac{1}{\mu^2}c^\mu\left( u_1 + u_2 + \cdots + u_{m'} \right)^\mu \,\,\,\textnormal{with $b,c>0$}\,\,\,\,\, \Longrightarrow A(u)^v \ll \left(\frac{b}{c} \right)^{v-1} A(u),\,\,\,\,\,v=2,3,\cdots
\end{align}

It suffices to prove the estimates $s(u)\ll A(u), \varphi_i(\xi_i,u)-\xi_i\ll A(u)$ and $c_i^{\alpha\beta}(\xi_i, u)- c_{0i}^{\alpha\beta}(\xi_i) \ll A(u)$ for suitable constants $b$ and $c$, equivalently
\begin{align}\label{uuc1}
s^\mu(u)\ll A(u),\,\,\,\,\,\,\,\,\,\varphi_i^\mu(\xi_i, u)- \xi_i\ll A(u),\,\,\,\,\,\,\,\,\,\,\,c_i^{\alpha\beta \mu}(\xi_i, u)- c_{0i}^{\alpha\beta}(\xi_i) \ll A(u)
\end{align}
for $\mu=1,2,3,\cdots$. We will prove $(\ref{uuc1})$ by induction on $\mu=1,2,3,\cdots$. For $\mu=1$, since the linear term of $A(u)$ is $\frac{b}{16}\left( u_1+ \cdots + u_{r'} \right)$, the estimates $\left(\ref{uuc1} \right)_1$ holds if $b$ is sufficiently large. Let $\mu\geq  2$ and assume that the induction $\left(\ref{uuc1}\right)_{\mu-1}$ holds for $\mu-1$, i.e.
\begin{align*}
s^{\mu-1}(u)\ll A(u),\,\,\,\,\,\,\,\,\,\varphi_i^{\mu-1}(\xi_i, u)- \xi_i \ll A(u),\,\,\,\,\,\,\,\,\,\,\,c_i^{\alpha\beta(\mu-1)}(\xi_i, u) - c_{0i}^{\alpha\beta}(\xi_i) \ll A(u)
\end{align*}
We will prove that $(\ref{uuc1})_\mu$ holds. We estimate $(\ref{uc10}),(\ref{uc11}),(\ref{uc12})$ and $(\ref{uc13})$ in the following Lemma.

\begin{lemma} \label{uuc13}
\begin{align}
\Gamma_{ij|\mu} &\ll \left( \frac{K_1}{b}+ \frac{K_2}{c} + \frac{K_3b}{c} \right) A(u)\,\,\,\,\,\,\,\,\,\textnormal{on}\,\,\, U_{ij} \label{uuc2} \\
A_{ij|\mu}^{\alpha\beta}\left(\xi_j, u \right) &\ll  \left( \frac{K_4}{c} + \frac{K_5}{c} + \frac{K_6 b}{c}     \right) A(u)  \,\,\,\,\,\,\,\,\textnormal{on}\,\,\,U_{ij} \label{uuc3}\\
W_{i|\mu}^{\alpha \gamma} &\ll \left( \frac{K_8}{c} + \frac{K_9b}{c}  \right) A(u)          \,\,\,\,\,\,\,\,\textnormal{on} \,\,\,U_i^\delta \label{uuc4}\\
V_{i|\mu-1}^{\alpha \delta} &\ll \left( \frac{K_{15}}{b} + \frac{K_{16}}{c} + \frac{K_{17} b}{c} \right) \frac{\partial A(u)}{\partial u_\delta},\,\,\,\,\,\delta=1,..., m'. \,\,\,\,\,\,\textnormal{on}\,\,\, U_i  \label{uuc5}
\end{align}
where $K_1,K_2,K_3, \cdots $ are constants independent of $\mu$.
\end{lemma}

\begin{proof}

The estimate $(\ref{uuc2})$ follows from \cite{Kod05} p.302. The estimate $(\ref{uuc3})$ and $(\ref{uuc4})$ follows from Part II. It remains to estimate $V_{i|\mu-1}^{\alpha \delta}$ from $(\ref{uc13})$
\begin{align}\label{uuc6}
\left[\sum_{\beta=1}^q c_i^{\alpha\beta (\mu-1)}(\xi_i, u) v_i'^{\beta\delta}(\xi_i, u)- \sum_{\beta=1}^q w_i^{\alpha\beta}\left(\varphi_i^{\mu-1}(\xi_i, u), s^{\mu-1}(u) \right)\frac{\partial \varphi_i^{\beta (\mu-1)} }{\partial u_\delta} - \sum_{\eta=1}^m v_i^{\alpha\eta} \left(\varphi_i^{\mu-1}(\xi_i , u), s^{\mu-1}(u)\right) \frac{\partial s^{\eta (\mu-1)}}{\partial u_\delta}\right]_{\mu-1}
\end{align}

First we estimate the first term of $(\ref{uuc6})$. We note that
\begin{align}\label{uuc7}
&\left[ c_i^{\alpha\beta(\mu-1)}(\xi_i, u) v_i'^{\beta \delta}(\xi_i, u)   \right]_{\mu-1} = \left[ \left( c_i^{\alpha\delta(\mu-1)}(\xi_i,u) - c_{0i}^{\alpha\delta} + c_{0i}^{\alpha \delta}\right) v_i'^{\beta \delta}(\xi_i, u) \right]_{\mu-1}\\
&= \left[ \left( c_i^{\alpha\delta(\mu-1)}(\xi_i,u) - c_{0i}^{\alpha\delta}  \right) v_i'^{\beta \delta}(\xi_i, u)  \right]_{\mu-1}  + \left[ c_{0i}^{\alpha \delta} v_i'^{\beta \delta}\left(  \xi_i, u \right) \right]_{\mu-1} \notag
\end{align}
Since $v_i'^{\beta \delta}(\xi_i, u)$ is a holomorphic function, we may assume that
\begin{align}\label{uuc8}
v_i'^{\beta \delta}(\xi_i, u) \ll \frac{\partial A_0(u)}{\partial u_\delta} \ll \frac{b_0}{b} \frac{ \partial A(u)}{\partial u_\delta},\,\,\,\,\,\,\,\,\,\,\,\frac{\partial A_0(u)}{\partial u_\delta}= \frac{b_0}{16 c_0} \sum_{v=1}^\infty \frac{ c_0^v\left(u_1 + \cdots + u_{m'} \right)^{v-1}}{v}
\end{align}
holds for $\xi_i\in U_i$ with $b_0>0$ and $c_0>0$ and $b>b_0$ and $c>c_0$. Then from $(\ref{uuc7})$ and $(\ref{uuc8})$ and the induction hypothesis, we have
\begin{align}\label{uuc12}
\left[\sum_{\beta=1}^q c_i^{\alpha\beta (\mu-1)}(\xi_i, u) v_i'^{\beta\delta}(\xi_i, u) \right]_{\mu-1} \ll K_{10} \frac{b_0}{b}A(u)\frac{\partial A(u)}{\partial u_\delta} + K_{11}\frac{b_0}{b} \frac{\partial A(u)}{\partial u_\delta} \ll \left( K_{10}\frac{1}{2} \frac{b_0}{c} + K_{11} \frac{b_0}{b} \right) \frac{\partial A(u)}{\partial u_\delta}
\end{align}

We estimate the second term of $(\ref{uuc6})$. We expand $w_i^{\alpha\beta}\left(\xi_i+ y, t \right)$ into power series in $y_1,...,y_n, t_1,...,t_m$. Then we may assume that
\begin{align}
w_i^{\alpha\beta}(\xi_i + y, t)- w_{0i}^{\alpha\beta} \ll \sum_{v=1}^\infty a_1^v\left( y_1+ \cdots + y_n + t_1 + \cdots + t_m \right)^v,\,\,\, a_1 >0
\end{align}
If we set $y=\varphi_i^{\mu-1}(\xi_i, u)- \xi_i$ and $t= s^{\mu-1}(u)$, then since $y \ll A(u)$ and $ t \ll A(u)$ by the induction hypothesis, we obtain
\begin{align}
w_i^{\alpha\beta}\left(\varphi_i^{\mu-1}(\xi_i, u), s^{\mu-1}(u) \right) - w_{0i}^\beta \ll \sum_{v=1}^\infty a_a^v (n+m)^v A(u)^v \ll \sum_{v=1}^{\infty} a_a^v (m+n)^v \left(\frac{b}{c} \right)^{v-1}A(u)
\end{align}
Then we have
\begin{align}
&\left[ \sum_{\beta=1}^q w_i^{\alpha\beta}\left(\varphi_i^{\mu-1}(\xi_i, u), s^{\mu-1}(u) \right)\frac{\partial \varphi_i^{\beta (\mu-1)} }{\partial u_\delta} \right]_{\mu-1} =\left[ \sum_{\beta=1}^q \left( w_i^{\alpha\beta}\left(\varphi_i^{\mu-1}(\xi_i, u), s^{\mu-1}(u) \right)- w_{0i}^{\alpha\beta} + w_{0i}^{\alpha\beta} \right)\frac{\partial \varphi_i^{\beta (\mu-1)} }{\partial u_\delta} \right]_{\mu-1}\\
&= \left[ \sum_{\beta=1}^q \left( w_i^{\alpha\beta}\left(\varphi_i^{\mu-1}(\xi_i, u), s^{\mu-1}(u) \right)- w_{0i}^{\alpha\beta} \right) \frac{\partial \left(\varphi_i^{\beta (\mu-1)} - \xi_i^\beta \right) }{\partial u_\delta} \right]_{\mu-1} \ll K_{12}\sum_{v=1}^{\infty} a_1^v(m+n)^v\left(\frac{b}{c} \right)^{v-1}A(u) \frac{\partial A(u)}{\partial u_\delta}      \notag \\
&\ll \frac{1}{2}K_{12} \frac{a_1b(m+n)}{c}  \sum_{v=0}^\infty \left(\frac{a_1 b(m+n)}{c} \right)^v \frac{\partial A(u)}{\partial u_\delta}   \notag
\end{align}
If we take a constant $c$ such that $\frac{ba_2(m+n)}{c}<\frac{1}{2}$, we obtain
\begin{align}\label{uuc10}
\left[ \sum_{\beta=1}^q w_i^{\alpha\beta}\left(\varphi_i^{\mu-1}(\xi_i, u), s^{\mu-1}(u) \right)\frac{\partial \varphi_i^{\beta (\mu-1)} }{\partial u_\delta} \right]_{\mu-1} \ll  K_{12} \frac{a_1 b(m+n)}{c} \frac{\partial A(u)}{\partial u_\delta} \ll K_{13}\frac{b}{c} \frac{\partial A(u)}{\partial u_\delta}
\end{align}
We estimate the third term of $(\ref{uuc6})$. In a similar way with $(\ref{uuc10})$
\begin{align} \label{uuc11}
\left[ \sum_{\eta=1}^m v_i^{\alpha\eta} \left(\varphi_i^{\mu-1}(\xi_i , u), s^{\mu-1}(u)\right) \frac{\partial s^{\eta (\mu-1)}}{\partial u_\delta} \right]_{\mu-1} \ll K_{14} \frac{b}{c} \frac{\partial A(u)}{\partial u_\delta}
\end{align}
Then from $(\ref{uuc12})$ and $(\ref{uuc10})$ and $(\ref{uuc11})$, we have
\begin{align}
V_{i|\mu-1}^{\alpha \delta} \ll \left( \frac{K_{15}}{b} + \frac{K_{16}}{c} + \frac{K_{17}b}{c} \right) \frac{\partial A(u)}{\partial u_\delta}
\end{align}
This completes the proof of Lemma \ref{uuc13}.
\end{proof}

From Lemma \ref{uuc13} we obtain
\begin{align}
\frac{\partial \Gamma_{ij|\mu}}{\partial u_\delta} &\ll \left( \frac{K_1}{b}+ \frac{K_2}{c} + \frac{K_3b}{c} \right) \frac{\partial A(u)}{u_\delta} \,\,\,\,\,\,\,\,\,\textnormal{on}\,\,\, U_{ij} \label{uuc14} \\
\frac{\partial A_{ij|\mu}^{\alpha\beta}}{\partial u_\delta}\left(\xi_j, u \right) &\ll  \left( \frac{K_4}{b} + \frac{K_5}{c} + \frac{K_6 b}{c}     \right) \frac{ \partial A(u) }{u_\delta}  \,\,\,\,\,\,\,\,\textnormal{on}\,\,\,U_{ij} \label{uuc15}\\
\frac{\partial W_{i|\mu}^{\alpha \gamma}}{\partial u_\delta} &\ll \left( \frac{K_8}{c} + \frac{K_9b}{c}  \right) \frac{\partial A(u)   }{u_\delta}        \,\,\,\,\,\,\,\,\textnormal{on} \,\,\,U_i^\delta \label{uuc16}\\
V_{i|\mu-1}^{\alpha \delta} &\ll \left( \frac{K_{15}}{b} + \frac{K_{16}}{c} + \frac{K_{17} b}{c} \right) \frac{\partial A(u)}{\partial u_\delta}, \,\,\,\,\,\,\textnormal{on}\,\,\, U_i  \label{uuc17}
\end{align}
for $\delta=1,...,m'$.

We recall Remark \ref{uuc13}. For any quadruple $\sigma=\left( -{W} , - V, \left( {\Gamma}, {A} \right) \right)=\left( \left\{- {W}_i \right\} ,\left\{ - V_i \right\} , \left\{ \left( {\Gamma}_{ij},  {A}_{ij} \right) \right\} \right)\in C^0\left( \mathcal{U}, \mathscr{H}om_{\mathcal{O}_M}\left( \mathcal{N}_{\mathcal{F}_0}^*, \Omega_M^1 \right) \right)\bigoplus C^0\left( \mathcal{U}, \mathscr{H}om_{\mathcal{O}_M}\left( \mathcal{N}_{\mathcal{F}_0}^*, \mathcal{O}_M \right) \right) \bigoplus C^1\left( \mathcal{U}, \mathcal{E}_{\mathcal{N}_{\mathcal{F}_0}^*} \right)$ which is a $1$-cocycle in the \v Cech resolution $(\ref{ucc15})$ of $\tilde{\mathcal{E}}_{\mathcal{N}_{\mathcal{F}_0}^*}^\bullet$, where where $W_i\in \Gamma\left( U_i, \mathscr{H}om_{\mathcal{O}_M}\left( \mathcal{N}_{\mathcal{F}_0}^*, \Omega_M^1 \right) \right)$ defined by
\begin{align} \label{ssc13}
W_i : \Gamma\left( U_i, \mathcal{N}_{\mathcal{F}_0}^*  \right) &\to \Gamma\left( U_i, \Omega_M^1  \right)\\
 w_{0i}^\alpha &\mapsto W_i^\alpha:=\sum_{\gamma=1}^n W_i^{\alpha \gamma} dz_i^\gamma \notag
\end{align}
and $V_i\in \Gamma\left(U_i, \mathscr{H}om_{\mathcal{O}_M} \left( \mathcal{N}_{\mathcal{F}_0}^*, \mathcal{O}_M \right) \right)$ defined by
\begin{align}
V_i: \Gamma\left( U_i, \mathcal{N}_{\mathcal{F}_0}^* \right) &\to \Gamma\left( U_i, \mathcal{O}_M  \right)\\
w_{0i}^\alpha &\mapsto V_i^\alpha \notag
\end{align}
and $\Gamma_{ij}\in \Gamma\left( U_{ij}, \Theta_M \right)$ and $A_{ij}\in \Gamma\left( U_{ij}, \mathscr{H}om_{\mathcal{O}_M}\left( \mathcal{N}_{\mathcal{F}_0}^*, \mathcal{N}_{\mathcal{F}_0}^*   \right) \right)$ defined by
\begin{align} 
A_{ij}: \Gamma\left( U_{ij}, \mathcal{N}_{\mathcal{F}_0}^* \right) &\to \Gamma\left( U_{ij} , \mathcal{N}_{\mathcal{F}_0}^* \right) , \\
 w_{0i}^\alpha &\mapsto A_{ij}^\alpha:= \sum_{\beta=1}^q A_{ij}^{\alpha\beta} w_{0i}^\beta \notag
\end{align}

we define the norm $\left| \left| \sigma \right| \right|$
\begin{align*}
\left| \left| \sigma \right| \right| =\left| \left| \Gamma \right| \right| + \left|\left| A \right| \right| + \left|\left| W\right|\right| + \left| \left| V \right| \right|
\end{align*}
where
\begin{align*}
&\left|\left| \Gamma \right|\right|=\max_{i,j} \sup_{\xi_i\in U_i\cap U_j} \left| \Gamma_{ij} (\xi_i)\right| , \,\,\,\,\,\,\,\,\,\,\,\,\,\,\,\,\,\,\,\,\left|\left| A \right|\right|=\max_{i,j}\max_{\alpha,\beta}\sup_{\xi_i  \in U_i\cap U_j}\left|A_{ij}^{\alpha\beta} \right| \\
&\left|\left| W \right| \right| = \max_{i }\max_{\alpha,\gamma}\sup_{\xi_i\in U_i^\delta} \left| W_i^{\alpha\gamma} (\xi_i) \right|,\,\,\,\,\,\,\,\,\left|\left| V \right|\right|= \max_i \max_\alpha \sup_{\xi_i\in U_i} \left| V_i^\alpha \right|
\end{align*}

\begin{lemma}\label{uuc18}
For any quadruple $\sigma=\left( -{W}, -V , \left( {\Gamma} , {A} \right) \right)=\left( \left\{ - {W}_i \right\}, \left\{ - V_i \right\}, \left\{ \left( \Gamma_{ij}, A_{ij} \right)  \right\} \right)$ which is a $1$-cocycle in the \v Cech resolution of $\tilde{\mathcal{E}}_{\mathcal{N}_{\mathcal{F}_0}^*}^\bullet$, we can find $\varphi_i(\xi), s^\lambda$ and $c_i^{\alpha \gamma}(\xi_i)$ satisfying
\begin{align}
\varphi_i- \varphi_j & = - \Gamma_{ij} + \sum_{\lambda=1}^m s^\lambda \rho_{ij\lambda} \\
\mathcal{L}_{\varphi_i}\left( w_{0i}^\alpha \right) - \sum_{\beta=1}^q & c_{i|\mu}^{\alpha\beta} w_{0i}'^\beta = W_i^\alpha  - \sum_{\lambda=1}^m s_\mu^\lambda \beta_{i\lambda}^\alpha \\
i_{\varphi_i}\left( w_{0i}^\alpha \right) &= V_i^\alpha + \sum_{\lambda=1}^m s^\lambda \gamma_{i\lambda}^\alpha \\
\left|\varphi_i\left(\xi_i \right) \right| \leq K \left|\left| \sigma \right|\right|,\,\,\,\,\,\,\,\,&\left|s \right| \leq K \left|\left| \sigma \right|\right|,\,\,\,\,\,\,\,\,\,\,\,\left| c_i^{\alpha \gamma}(\xi_i) \right| \leq K \left|\left| \sigma \right|\right|  
\end{align}
where $K$ is a constant independent of $\sigma$.
\end{lemma}

\begin{proof}
We can prove Lemma \ref{uuc18} in a similar way as in Part I.

\end{proof}

By $(\ref{uuc14})-(\ref{uuc17})$ and Lemma \ref{uuc18} we can choose solutions $\tilde{\varphi}_{i|\mu-1}^\delta(\xi_i, u)$ and $\tilde{s}_{\mu-1}^{\delta \lambda}$ and $\tilde{c}_{i|\mu-1}^{\delta \alpha\beta}$ for each $\delta=1,...,m'$ of equations $(\ref{uc30})$ and $(\ref{uc31})$ and $(\ref{uc32})$
such that
\begin{align*}
\tilde{\varphi}_{i|\mu-1}^\delta \ll K K^* \frac{\partial A(u)}{\partial u_\delta} ,\,\,\,\,\,\,\,\,\,\, \tilde{s}_{\mu-1}^\delta \ll K K^*  \frac{\partial A(u)}{\partial u_\delta},\,\,\,\,\,\,\,\,\,\, \tilde{c}_{i|\mu-1}^{\delta \alpha \beta}\ll K K^* \frac{\partial A(u)}{\partial u_\delta}
\end{align*}
Then we have
\begin{align*}
\varphi_{i|\mu}:=\frac{1}{\mu}\sum_{\delta=1}^{m'} u_\delta \tilde{\varphi}_{i|\mu-1}^\delta \ll KK^* \sum_{\delta=1}^{m'} \frac{1}{\mu} u_\delta \left[\frac{\partial A(u)}{\partial u_\delta}\right]_{\mu-1} =KK^*[A(u)]_\mu \ll KK^* A(u)\\
s_{i|\mu}=\frac{1}{\mu}\sum_{\delta=1}^{m'} u_\delta \tilde{s}_{i|\mu-1}^\delta \ll KK^* \sum_{\delta=1}^{m'} \frac{1}{\mu} u_\delta \left[\frac{\partial A(u)}{\partial u_\delta}\right]_{\mu-1} =KK^*[A(u)]_\mu \ll KK^* A(u)\\
c_{i|\mu}^{\alpha \gamma}=\frac{1}{\mu}\sum_{\delta=1}^{m'} u_\delta \tilde{a}_{i|\mu-1}^{\delta\alpha \gamma} \ll K_6K^* \sum_{\delta=1}^{m'} \frac{1}{\mu} u_\delta \left[\frac{\partial A(u)}{\partial u_\delta}\right]_{\mu-1} =KK^*[A(u)]_\mu \ll KK^* A(u)
\end{align*}
where $K^*=\frac{K_1 + K_4 + K_{15}}{b} + \frac{K_2 + K_5 + K_8 + k_{16}}{c} + \frac{(K_3 + K_6 + K_9 + K_{17})b}{c}$. We choose $b$ and $c$ in a way that $KK^*<1$. Then we have $\varphi_{i|\mu}(\xi_i, u)\ll A(u), s_\mu(u)\ll A(u)$ and $c_{i|\mu}^{\alpha\beta}\ll A(u)$. Hence the induction $(\ref{ucc1})_\mu$ holds for $\mu$. This completes the proof of Theorem \ref{uc1}.

\end{proof}

\appendix

\section{Unfolded dual leaf complex controlling unfoldings of compact foliated complex manifolds defined by locally free sheaves of cotangent sheaves}\label{AppendixA}

Let $\left(M, \mathcal{N}_{\mathcal{F}_0}^*\right)$ be a compact foliated complex manifold with $\mathcal{N}_{\mathcal{F}_0}^*$ locally free of rank $q$. We note that we have a section $\omega_0\in \Gamma\left( M, \bigwedge^q\Omega_M^1\otimes \mathcal{L}_0 \right)$ where $\mathcal{L}_0= \bigwedge^q \left( \mathcal{N}_{\mathcal{F}_0}^* \right)^*$, which is locally decomposable and satisfy integrability condition. Let $\mathcal{U}=\{U_i\}$ be an open covering of $M$ by coordinate neighborhoods such that $\Gamma(U_i, \mathcal{N}_{\mathcal{F}_0}^*)$  is generated by $w_i^1,..., w_i^q\in \Gamma(U_i, \mathcal{N}_{\mathcal{F}_0}^*)$ with the relation $w_i^\alpha= \sum_{\beta=1}^q h_{ij}^{\alpha\beta}w_j^\beta$ on $U_i\cap U_j$ for $h_{ij}^{\alpha\beta}\in \Gamma(U_{ij}, \mathcal{O}_M)$ and satisfy $w_i^1\wedge \cdots \wedge w_i^q \wedge dw_i^\alpha=0,\alpha=1,...,q$. We set $w_i:= w_i^1\wedge \cdots \wedge w_i^q \in \Gamma\left( U_i, \bigwedge^q \Omega_M^1 \right)$. Then $w_i=\det\left( h_{ij}^{\alpha\beta} \right) w_j$ on $U_i\cap U_j$. Then $w_0$ is defined by $\{w_i\}$. We recall from Part I that we defined the dual leaf complex associated to $\mathcal{N}_{\mathcal{F}_0}^\bullet$:

\begin{align*}
\mathcal{N}_{\mathcal{F}_0}^{*\bullet}: \Theta_M \xrightarrow{E_0} \mathscr{H}om_{\mathcal{O}_M}\left( \mathcal{N}_{\mathcal{F}_0}^*, \frac{\Omega_M^1}{\mathcal{N}_{\mathcal{F}_0}^*} \right) \xrightarrow{E_1} \mathscr{H}om_{\mathcal{O}_M}\left( \mathcal{N}_{\mathcal{F}_0}^*, \tilde{\mathcal{S}}^2 \right) \xrightarrow{E_2} \mathscr{H}om_{\mathcal{O}_M}\left( \mathcal{N}_{\mathcal{F}_0}^*, \tilde{\mathcal{S}}^3\right) \xrightarrow{E_3} \to \cdots
\end{align*}

We shall define the unfolded dual leaf complex for $\left( M, \mathcal{N}_{\mathcal{F}_0}^* \right)$ with $\mathcal{N}_{\mathcal{F}_0}^*$ locally free:
{\small{\begin{align}\label{Ap1}
\tilde{\mathcal{N}}_{\mathcal{F}_0}^{*\bullet}: \Theta_M  \xrightarrow{G_0} \mathscr{H}om_{\mathcal{O}_M}\left(\mathcal{N}_{\mathcal{F}_0}^*, \frac{\Omega_M^1}{\mathcal{N}_{\mathcal{F}_0}^*} \right) \bigoplus & \mathscr{H}om_{\mathcal{O}_M}\left(\mathcal{N}_{\mathcal{F}_0}^*, \mathcal{O}_M \right) \xrightarrow{G_1} \mathscr{H}om_{\mathcal{O}_M}\left(\mathcal{N}_{\mathcal{F}_0}, \tilde{\mathcal{S}}^2 \right) \bigoplus \mathscr{H}om_{\mathcal{O}_M}\left( \mathcal{N}_{\mathcal{F}_0}, \tilde{\mathcal{S}}^1 \right) \\
&\xrightarrow{G_2} \mathscr{H}om_{\mathcal{O}_M}\left( \mathcal{N}_{\mathcal{F}_0}^*, \tilde{\mathcal{S}}^3 \right) \bigoplus \mathscr{H}om_{\mathcal{O}_M}\left( \mathcal{N}_{\mathcal{F}_0}^*, \tilde{\mathcal{S}}^2 \right) \xrightarrow{G_3} \cdots \notag
\end{align}}}
For $X\in \Theta_M$, $G_0(X)$ is defined by $\left(E_0(X), G_0'(X) \right)$ where $G_0'(X)$ is locally defined on $U_i$ by
\begin{align*}
G_0'(X)\left(w_i^\alpha\right)= i_X\left(w_i^\alpha\right),\,\,\,\alpha=1,...,q
\end{align*}
and linearly extends to $\Gamma\left(U_i, \mathcal{N}_{\mathcal{F}_0}^* \right)$. This is well-defined since $G_0'(X)(w_i^\alpha)- G_0'(X)\left(\sum_{\beta=1}^q h_{ij}^{\alpha\beta} w_j^\beta \right)=0$, and $E_1$ is defined in the following way: for $\phi\in \mathscr{H}om_{\mathcal{O}_M}\left(\mathcal{N}_{\mathcal{F}_0}^*, \frac{\Omega_M^1}{\mathcal{N}_{\mathcal{F}_0}^*} \right)$ and $\psi\in \mathscr{H}om_{\mathcal{O}_M}\left(\mathcal{N}_{\mathcal{F}_0}^*, \mathcal{O}_M \right)$, we define $G_1(\phi, \psi)=\left(E_1(\phi), G_1'(\phi, \psi)\right)$. Here $G_1'(\phi, \psi)$ is locally defined on $U_i$ by
\begin{align*} 
&G_1'(\phi, \psi)\left(w_i^\alpha\right)\\
&= - w_i^1\wedge \cdots \wedge w_i^q \wedge \widetilde{\phi\left(w_i^\alpha\right)}+ \sum_{\beta=1}^q (-1)^{q-\beta} \psi\left(w_i^\beta \right) w_i^1\wedge \cdots \wedge \widehat{w}_i^\beta \wedge \cdots \wedge w_i^q\wedge dw_i^\alpha + w_i^1\wedge \cdots \wedge w_i^q \wedge d\psi\left( w_i^\alpha \right) \\
&\in \Gamma\left(U_i, \bigwedge^{q+1} \Omega_M^1 \right)\cong \Gamma\left(U_i, \bigwedge^{q+1} \Omega_M^1 \otimes \mathcal{L}_0 \right)
\end{align*}
for $\alpha=1,...,q$, and linearly extends to $\Gamma\left(U_i, \mathcal{N}_{\mathcal{F}_0}^* \right)$. Here $\widetilde{\phi(w_i^\alpha)}$ is a lifting of $\phi(w_i^\alpha)$ under $\Omega_M^1\to \frac{\Omega_M^1}{\mathcal{N}_{\mathcal{F}_0}^*}$.\footnote{Strictly speaking $\Omega_M^1\to \frac{\Omega_M^1}{\mathcal{N}_\mathcal{F}^*}$ is not surjective on $U_i$. More precisely for $x\in U_i$, choose a neighborhood $U_x$ in $U_i$ such that $\phi(w_i^\alpha)$ on $U_x$ has a lifting $\widetilde{\phi_x(w_i^\alpha)}$. Then we compute $w_i^1\wedge \cdots \wedge w_i^q \wedge \widetilde{\phi_x(w_i^\alpha)}$ locally on $U_x$ for each $x\in U_i$. They are equal on $U_x\cap U_y\ne \emptyset$ since $\widetilde{\phi_x(w_i^\alpha)}- \widetilde{\phi_y(w_i^\alpha)}\in \mathcal{N}_{\mathcal{F}_0}^*$. Hence we can glue together $w_i^1\wedge \cdots w_i^q \wedge \widetilde{\phi_x(w_i^\alpha)}$ to define $w_i^1\wedge \cdots \wedge w_i^q \wedge \widetilde{\phi(w_i^\alpha)}$.} $G_1'(\phi, \psi)(w_i^\alpha)$ is independent of a lifting of $\phi(w_i^\alpha)$ since for $A\in \mathcal{N}_{\mathcal{F}_0}^*$, $w_i^1\wedge \cdots \wedge w_i^q\wedge A=0$. We show that $G_1'(\phi, \psi)\left(w_i^\alpha \right)$ actually lies in $\Gamma\left(U_i, \tilde{\mathcal{S}}^1 \right)$. For each $x\in U_i-S$ where $S=\textnormal{Sing}\left(\mathcal{F}_0\right)$, choose a neighborhood $U_x\subset U_i-S$ such that $dw_i^\alpha= \sum_{\beta=1}^q a_{i_x}^{\alpha\beta} \wedge w_i^\beta$ for some $a_{i_x}^{\alpha\beta}\in \Gamma\left(U_x, \Omega_M^1 \right)$. Then we see that on $U_x$
\begin{align*}
&\sum_{\beta=1}^q  (-1)^{q-\beta} \psi\left(w_i^\beta \right) w_i^1\wedge \cdots \wedge \widehat{w}_i^\beta \wedge \cdots \wedge w_i^q \wedge dw_i^\alpha + w_i^1\wedge \cdots \wedge w_i^q \wedge d \psi\left(w_i^\alpha \right)\\
&=\sum_{\beta=1}^q (-1)^{q-\beta} \psi\left(w_i^\beta \right) w_i^1 \wedge \cdots \hat{w}_i^\beta \wedge \cdots \wedge w_i^q\wedge \left( \sum_{\gamma=1}^q a_{i_x}^{\alpha \gamma} \wedge w_i^\gamma \right) + w_i^1\wedge \cdots \wedge w_i^q \wedge d \psi\left(w_i^\alpha \right)\\
&= w_i^1\wedge \cdots \wedge w_i^q \wedge \left( d\psi(w_i^\alpha) - \sum_{\beta=1}^q \psi\left(w_i^\beta \right)a_{i_x}^{\alpha\beta} \right)\in \Gamma\left(U_x, \tilde{\mathcal{S}}^1 \right)
\end{align*}
Next we show that $G_1' \left(\phi, \psi \right)$ is well-defined: in other words, $G_1'\left(\phi, \psi\right)$ on $U_i$ and $G_1'\left(\phi, \psi\right)$ on $U_j$ defines the same homomorphism on $U_i\cap U_j$. We recall that $w_i^\alpha= \sum_{\beta=1}^q h_{ij}^{\alpha\beta} w_j^\beta$. By abuse of notation, we write $\phi(w_i^\alpha)$ instead of $\widetilde{\phi(w_i^\alpha)}$. We set $H_{ij}=\det\left( h_{ij}^{\alpha\beta} \right)$.
{\small{\begin{align*}
&G_1'\left(\phi, \psi \right)\left(w_i^\alpha\right)\\
&= - w_i^1\wedge \cdots \wedge w_i^q \wedge \phi\left(w_i^\alpha\right) + \sum_{\beta=1}^q (-1)^{q-\beta} \psi\left(w_i^\beta \right) w_i^1\wedge \cdots \wedge \widehat{w}_i^\beta \wedge \cdots \wedge w_i^q\wedge dw_i^\alpha + w_i^1\wedge \cdots \wedge w_i^q \wedge d\psi\left( w_i^\alpha \right)\\
&= - H_{ij} w_j^1 \wedge \cdots \wedge w_j^q \wedge \phi\left( \sum_{\beta=1}^q h_{ij}^{\alpha\beta} w_j^\beta \right) + H_{ij} \sum_{\beta=1}^q (-1)^{q-\beta} \psi\left(w_j^\beta \right) w_j^1 \wedge \cdots \wedge \widehat{w}_j^\beta \wedge \cdots \wedge w_j^q \wedge d\left( \sum_{\gamma=1}^q h_{ij}^{\alpha\gamma} w_j^\gamma \right)\\
&+ H_{ij} w_j^1 \wedge \cdots \wedge w_j^q \wedge d\psi\left(\sum_{\beta=1}^q h_{ij}^{\alpha\beta} w_j^\beta \right)\\
&= - H_{ij} w_j^1 \wedge \cdots \wedge w_j^q \wedge \phi\left( \sum_{\beta=1}^q h_{ij}^{\alpha\beta} w_j^\beta \right) + H_{ij} \sum_{\beta=1}^q (-1)^{q-\beta} \psi\left(w_j^\beta \right) w_j^1 \wedge \cdots \wedge \widehat{w}_j^\beta \wedge \cdots \wedge w_j^q \wedge \left( \sum_{\gamma=1}^q dh_{ij}^{\alpha\gamma} \wedge w_j^\gamma + \sum_{\gamma=1}^q h_{ij}^{\alpha \gamma} dw_j^\gamma \right)\\
&+ H_{ij} w_j^1 \wedge \cdots \wedge w_j^q \wedge \left( \sum_{\beta=1}^q \psi \left( w_j^\beta \right) dh_{ij}^{\alpha\beta} + \sum_{\beta=1}^q h_{ij}^{\alpha\beta} d\psi\left( w_j^\beta \right) \right)\\
&= H_{ij}\sum_{\beta=1}^q h_{ij}^{\alpha\beta}\left( - w_j^1 \wedge \cdots \wedge w_j^q \wedge \phi(w_j^\beta) + \sum_{\gamma=1}^q (-1)^{q-\gamma} \psi\left(w_j^\gamma \right) w_j^1 \wedge \cdots \wedge \widehat{w}_j^\gamma\wedge \cdots \wedge w_j^q\wedge dw_j^\beta + w_j^1 \wedge \cdots \wedge w_j^q \wedge d \psi\left( w_j^\beta \right) \right)\\
&= H_{ij}\sum_{\beta=1}^q h_{ij}^{\alpha\beta} G_1'\left(\phi, \psi\right)\left(w_j^\beta\right)
\end{align*}}}
Hence $G_1'$ is well-defined.

For $r\geq 2$, we define {\Small{$G_r:\mathscr{H}om_{\mathcal{O}_M}\left( \mathcal{N}_{\mathcal{F}_0}^* , \tilde{\mathcal{S}}^r \right) \bigoplus \mathscr{H}om_{\mathcal{O}_M}\left( \mathcal{N}_{\mathcal{F}_0}^*, \tilde{\mathcal{S}}^{r-1} \right) \to \mathscr{H}om_{\mathcal{O}_M}\left( \mathcal{N}_{\mathcal{F}_0}^*, \tilde{\mathcal{S}}^{r+1} \right)\bigoplus \mathscr{H}om_{\mathcal{O}_M}\left( \mathcal{N}_{\mathcal{F}_0}^*, \tilde{\mathcal{S}}^r \right)$}} in the following way: for $\phi\in \mathscr{H}om_{\mathcal{O}_M}\left(\mathcal{N}_{\mathcal{F}_0}^*, \tilde{\mathcal{S}}^r \right)$ and $\psi\in \mathscr{H}om_{\mathcal{O}_M}\left( \mathcal{N}_{\mathcal{F}_0}^*, \tilde{\mathcal{S}}^{r-1}\right)$, we define $G_r\left( \phi, \psi \right)= \left( E_r(\phi) , (-1)^r \phi + E_{r-1}(\psi) \right)$.

We show that $G_1\circ G_0=0$. In fact, for $X\in \Theta_X$,
\begin{align*}
G_1(G_0(X))= G_1\left( E_0(X), G_0'(X) \right) = \left(E_1(E_0(X)), G_1'\left(E_0(X), G_0'(X) \right)    \right) = \left( 0, G_1'\left(E_0(X), G_0'(X) \right)  \right)
\end{align*}
We note that since $0=i_X\left(w_i^1 \wedge \cdots \wedge w_i^q \wedge dw_i^\alpha \right)=\sum_{\beta=1}^q (-1)^{\beta-1} i_X\left(w_i^\beta \right)w_i^1 \wedge \cdots \wedge \widehat{w}_i^\beta \wedge \cdots \wedge w_i^q \wedge dw_i^\alpha + (-1)^q w_i^1\wedge \cdots \wedge w_i^q \wedge i_X\left(dw_i^\alpha\right)$, we have
\begin{align*}
&G_1'\left(E_0(X), G_0'(X) \right) (w_i^\alpha)\\
&= - w_i^1 \wedge \cdots \wedge w_i^q \wedge \mathcal{L}_X\left(w_i^\alpha\right) + \sum_{\beta=1}^q (-1)^{q-\beta} i_X\left(w_i^\beta\right) w_i^1 \wedge \cdots \wedge \widehat{w}_i^\beta \wedge \cdots \wedge w_i^q \wedge dw_i^\alpha + w_i^1 \wedge \cdots \wedge w_i^q \wedge di_X\left(w_i^\alpha\right)\\
&= w_i^1 \wedge \cdots \wedge w_i^q \wedge \left(- \mathcal{L}_X \left(w_i^\alpha \right) + i_X \left( dw_i^\alpha \right) + di_X\left(w_i^\alpha \right) \right) = 0
\end{align*}
We show that for $r\geq1$, $G_{r+1} \circ G_r=0$. In fact, $G_{r+1}\left( G_r(\phi, \psi) \right)= G_{r+1}\left(E_r(\phi), (-1)^r \phi + E_{r-1}(\psi) \right) = \left( E_{r+1}E_r(\phi),  (-1)^{r+1} E_r(\phi) +  E_r \left( (-1)^r\phi\right)+ E_r E_{r-1}(\phi) \right) =0 $. Hence $\textnormal{(\ref{Ap1})}$ defines a complex of sheaves.

\section{Dolbeault type bicomplex associated to the unfolded dual leaf complex $\tilde{\mathcal{N}}_{\mathcal{F}_0}^{*\bullet}$} \label{up1}

Let $\left(  M,  \mathcal{N}_{\mathcal{F}_0}^* \right)$ be a compact foliated complex manifold with $\mathcal{N}_{\mathcal{F}_0}^*$ locally free. We will construct a Dolbeault type bicomplex assocaited to the dual leaf complex $\mathcal{N}_{\mathcal{F}_0}^{*\bullet}$ from $\textnormal{(\ref{Ap1})}$. Since $\tilde{\mathcal{S}}^r,r\geq 1$ is not necessarily locally free, we do not have the Dolbeault resolution of $\tilde{\mathcal{S}}^r$ on $M$. However, since $\tilde{\mathcal{S}}^r$ is locally free on $M-S$, we have the Dolbeault resolution of $\tilde{\mathcal{S}}^r$ on $M-S$. As in a similar way to construct the Dolbeault type bicomplex associated to the dual leaf complex $\mathcal{N}_{\mathcal{F}_0}^{*\bullet}$ in Part I, we have the following bicomplex 
{\Tiny{\begin{center}
$\begin{CD}
\cdots \\
@A\hat{G}_3AA \\
A^{0,0}\left(M-S, \left(\mathcal{N}_{\mathcal{F}_0}^* \right)^* \otimes \tilde{\mathcal{S}}^3\right)\bigoplus A^{0,0}\left( M-S, \left(\mathcal{N}_{\mathcal{F}_0}^*\right)^*\otimes \tilde{\mathcal{S}}^2  \right) @>-(-1)^q\bar{\partial}>> \cdots \\
@ A \hat{G}_2 AA @A \hat{G}_2 AA \\
A^{0,0}\left(M-S, \left(\mathcal{N}_{\mathcal{F}_0}^*\right)^* \otimes \tilde{\mathcal{S}}^2\right) \bigoplus A^{0,0}\left( M-S, \left(\mathcal{N}_{\mathcal{F}_0}^* \right)^*\otimes \tilde{\mathcal{S}}^1 \right) @>(-1)^q \bar{\partial} >> A^{0,1}\left(M-S, \left(\mathcal{N}_{\mathcal{F}_0}^*\right)^*\otimes \tilde{\mathcal{S}}^2 \right) \bigoplus A^{0,1}\left( M-S, \left( \mathcal{N}_{\mathcal{F}_0}^*\right)^*\otimes \tilde{\mathcal{S}}^1 \right) @>\bar{\partial}>> \cdots\\
@A\hat{G}_1  AA @A\hat{G}_1 AA @A\hat{G}_1AA \\
\frac{A^{0,0}\left(M, \left(\mathcal{N}_{\mathcal{F}_0}^*\right)^* \otimes \Omega_M^1 \right)}{A^{0,0}\left(M, \left( \mathcal{N}_{\mathcal{F}_0}^* \right)^*\otimes \mathcal{N}_{\mathcal{F}_0}^* \right) }\bigoplus A^{0,0}\left(M, \left( \mathcal{N}_{\mathcal{F}_0}^* \right)^* \right) @>-\bar{\partial}>> \frac{A^{0,1}\left(M, \left( \mathcal{N}_{\mathcal{F}_0}^* \right)^*\otimes \Omega_M^1 \right)}{A^{0,1}\left( M, \left(\mathcal{N}_{\mathcal{F}_0}^* \right)^*\otimes \mathcal{N}_{\mathcal{F}_0}^* \right) }\bigoplus A^{0,1}\left(M, \left(\mathcal{N}_{\mathcal{F}_0}^* \right)^* \right) @>\bar{\partial}>> \cdots\\
@A\hat{G}_0AA @A \hat{G}_0 AA @A\hat{G}_0AA \\
A^{0,0}\left(M, \Theta_M\right) @>\bar{\partial}>> A^{0,1}\left( M, \Theta_M\right) @>-\bar{\partial}>> A^{0,2}\left(M, \Theta_M\right) 
\end{CD}$
\end{center}}}
where
{\Small{\begin{align*}
\hat{G}_1:\frac{A^{0,p}\left(M, \mathscr{H}om_{\mathcal{O}_M}\left( \mathcal{N}_{\mathcal{F}_0}^*, \Omega_M^1 \right) \right)}{A^{0,p}\left( M, \mathscr{H}om_{\mathcal{O}_M}\left( \mathcal{N}_{\mathcal{F}_0}^*, \mathcal{N}_{\mathcal{F}_0}^* \right) \right)}\bigoplus A^{0,p}\left( M, \left( \mathcal{N}_{\mathcal{F}_0}^*\right)^*  \right) \to A^{0,p}\left(M-S,\left( \mathcal{N}_{\mathcal{F}_0}^* \right)^* \otimes\tilde{\mathcal{S}}^2\right)\bigoplus \mathcal{A}^{0,p}\left( M-S,\left( \mathcal{N}_{\mathcal{F}_0}^* \right)^* \otimes \tilde{\mathcal{S}}^1 \right)
\end{align*}}}
 is defined in the following way: first we have the composition 
{\Tiny{\begin{align*}
\mathscr{H}om_{\mathcal{O}_M}\left( \mathcal{N}_{\mathcal{F}_0}^*, \frac{\Omega_M^1}{\mathcal{N}_{\mathcal{F}_0}^*} \right)\bigoplus \mathscr{H}om_{\mathcal{O}_M}\left( \mathcal{N}_{\mathcal{F}_0}^*, \mathcal{O}_M \right)  \xrightarrow{G_1}  \left( \mathcal{N}_{\mathcal{F}_0}^*\right)^* \otimes \tilde{\mathcal{S}}^2  \bigoplus \left( \mathcal{N}_{\mathcal{F}_0}^* \right)^* \otimes \tilde{\mathcal{S}}^1 \hookrightarrow \mathscr{H}om_{\mathcal{O}_M}\left(\mathcal{N}_{\mathcal{F}_0}^*, \bigwedge^{q+2} \Omega_M^1\otimes \mathcal{L}_0 \right)\bigoplus \mathscr{H}om_{\mathcal{O}_M}\left( \mathcal{N}_{\mathcal{F}_0}^*, \bigwedge^{ q+1} \Omega_M^1 \otimes \mathcal{L}_0 \right)
\end{align*}}}
where $\mathcal{L}_0:= \bigwedge^q \left( \mathcal{N}_{\mathcal{F}_0}^* \right)^*$, which induces
{\Tiny{\begin{align*}
G_1^\sharp: \frac{A^{0,p}\left( M,  \mathscr{H}om_{\mathcal{O}_M}\left( \mathcal{N}_{\mathcal{F}_0}^*, \Omega_M^1 \right) \right)}{A^{0,p}\left( M , \mathscr{H}om_{\mathcal{O}_M}\left( \mathcal{N}_{\mathcal{F}_0}^*, \mathcal{N}_{\mathcal{F}_0}^* \right) \right)} \bigoplus  A^{0,p}\left( M, \mathscr{H}om_{\mathcal{O}_M} \left( \mathcal{N}_{\mathcal{F}_0}^*, \mathcal{O}_M \right) \right) \to A^{0,p}\left( M ,  \mathscr{H}om_{\mathcal{O}_M}\left( \mathcal{N}_{\mathcal{F}_0}^*, \bigwedge^{q+2} \Omega_M^1\otimes \mathcal{L}_0 \right)\right) \bigoplus \mathscr{H}om_{\mathcal{O}_M}\left( \mathcal{N}_{\mathcal{F}_0}^*,  \bigwedge^{q+1} \Omega_M^1 \otimes \mathcal{L}_0 \right)
\end{align*}}}
We note that when we restrict the image of $G_1^\sharp$ on $M-S$, the restriction lies in $A^{0,p}\left( M - S , \mathscr{H}om_{\mathcal{O}_M} \left( \mathcal{N}_{\mathcal{F}_0}^*, \tilde{\mathcal{S}}^2 \right) \right)$ $ \bigoplus  A^{0,0}\left( M-S , \mathscr{H}om_{\mathcal{O}_M}\left( \mathcal{N}_{\mathcal{F}_0}^*, \tilde{\mathcal{S}}^1 \right) \right)$. Then we define $\hat{G}_1$ to be the composition of $G_1^\sharp$ and the restriction on $M-S$.

We will denote the $i$-th cohomology group of the complex associated to the above bicomplex by $\textnormal{G}^i$. Then we will show in Part IV-1 that for $i=0,1,2$,
\begin{align}\label{us6}
\textnormal{G}^i \cong \mathbb{H}^i\left( M, \tilde{\mathcal{N}}_{\mathcal{F}_0}^{*\bullet} \right)
\end{align}
by explicitly constructing isomorphisms between $\textnormal{G}^i$ and $i$-th cohomology group from the \v Cech resolution of $\tilde{\mathcal{N}}_{\mathcal{F}_0}^{*\bullet}$ for $i= 0,1,2$.

We will describe $\hat{G}_0$ and $\hat{G}_1$ and $\hat{G}_2$ explicitly which we use in the main body of the paper. Let $\mathcal{U}=\left\{U_i\right\}$ be a Stein open covering of $M$ by coordinate neighborhoods with local coordinates $\left( z_i^1,..., z_i^n \right)$ on $U_i$ such that $\Gamma\left( U_j, \mathcal{N}_{\mathcal{F}_0}^* \right)$ is generated by $w_j^1,..., w_j^q\in \Gamma\left( U_j, \Omega_M^1 \right)$ with the relation $w_j^\alpha= \sum_{\beta=1}^q h_{jk}^{\alpha\beta} w_j^\beta$ for some $h_{jk}^{\alpha\beta}\in \Gamma\left( U_j \cap U_k, \mathcal{O}_M \right)$ and $w_j^1 \wedge \cdots \wedge w_j^q \wedge dw_i^\alpha=0$ for $\alpha= 1,...,q$.

\subsubsection{\textnormal{Description of $\hat{G}_0$}}\

We describe $\hat{G}_0: A^{0,p}\left( M, \Theta_M \right) \to \frac{ A^{0,p}\left( M, \mathscr{H}om_{\mathcal{O}_M}\left( \mathcal{N}_{\mathcal{F}_0}^*, \Omega_M^1 \right) \right)}{A^{0,p}\left( M, \mathscr{H}om_{\mathcal{O}_M}\left( \mathcal{N}_{\mathcal{F}_0}^*, \mathcal{N}_{\mathcal{F}_0}^* \right) \right)}\bigoplus A^{0,p}\left( M, \mathscr{H}om_{\mathcal{O}_M}\left( \mathcal{N}_{\mathcal{F}_0}^*, \mathcal{O}_M \right) \right)$. For $\phi\in A^{0,p}\left( M, \Theta_M \right)$, we define $\hat{G}_0\left( \phi \right)=\left(\hat{E}_0\left( \phi\right), \hat{G}_0'\left( \phi \right) \right)$, where $\hat{E}_0$ is defined as in Part I, and $\hat{G}_0'(\phi)$ is defined in the following way: locally on $U_i$, we define
\begin{align*}
\hat{G}_0'(\phi)\left( w_i^\alpha \right) = i_{\phi}\left( w_i^\alpha \right)
\end{align*}

\subsubsection{ \textnormal{Description of $\hat{G}_1$}}\

We describe 
{\tiny{\begin{align*}
\hat{G}_1:\frac{ A^{0,p}\left( M, \mathscr{H}om_{\mathcal{O}_M}\left( \mathcal{N}_{\mathcal{F}_0}^*, \Omega_M^1 \right) \right)}{A^{0,p}\left( M, \mathscr{H}om_{\mathcal{O}_M}\left( \mathcal{N}_{\mathcal{F}_0}^*, \mathcal{N}_{\mathcal{F}_0}^* \right) \right)}\bigoplus A^{0,p}\left( M, \mathscr{H}om_{\mathcal{O}_M}\left( \mathcal{N}_{\mathcal{F}_0}^*, \mathcal{O}_M \right) \right) \to A^{0,p}\left( M-S, \mathscr{H}om_{\mathcal{O}_M}\left( \mathcal{N}_{\mathcal{F}_0}^*, \tilde{\mathcal{S}}^2 \right) \right) \bigoplus A^{0,p}\left(  M-S, \mathscr{H}om_{\mathcal{O}_M}\left( \mathcal{N}_{\mathcal{F}_0}^*, \tilde{\mathcal{S}}^1 \right) \right)
\end{align*}}}
It is sufficient to describe 
{\small{\begin{align*}
\hat{G}_1^\sharp &:\frac{ A^{0,p}\left( M, \mathscr{H}om_{\mathcal{O}_M}\left( \mathcal{N}_{\mathcal{F}_0}^*, \Omega_M^1 \right) \right)}{A^{0,p}\left( M, \mathscr{H}om_{\mathcal{O}_M}\left( \mathcal{N}_{\mathcal{F}_0}^*, \mathcal{N}_{\mathcal{F}_0}^* \right) \right)}\bigoplus A^{0,p}\left( M, \mathscr{H}om_{\mathcal{O}_M}\left( \mathcal{N}_{\mathcal{F}_0}^*, \mathcal{O}_M \right) \right) \\
&\to A^{0,p}\left( M , \mathscr{H}om_{\mathcal{O}_M}\left( \mathcal{N}_{\mathcal{F}_0}^*, \bigwedge^{q+2} \Omega_M^1 \otimes \mathcal{L}_0 \right) \right) \bigoplus A^{0,p}\left(  M, \mathscr{H}om_{\mathcal{O}_M}\left( \mathcal{N}_{\mathcal{F}_0}^* , \bigwedge^{q+1} \Omega_M^1 \otimes \mathcal{L}_0 \right) \right)
\end{align*}}}
For $\left(\overline{\phi}, \psi \right) \in \frac{ A^{0,p}\left( M, \mathscr{H}om_{\mathcal{O}_M}\left( \mathcal{N}_{\mathcal{F}_0}^*, \Omega_M^1 \right) \right)}{A^{0,p}\left( M, \mathscr{H}om_{\mathcal{O}_M}\left( \mathcal{N}_{\mathcal{F}_0}^*, \mathcal{N}_{\mathcal{F}_0}^* \right) \right)}\bigoplus A^{0,p}\left( M, \mathscr{H}om_{\mathcal{O}_M}\left( \mathcal{N}_{\mathcal{F}_0}^*, \mathcal{O}_M \right) \right)$, we let 
\begin{align*}
\phi\left( w_i^\alpha \right) &= \sum_{\gamma_1< \cdots <\gamma_p} d\bar{z}_i^{\gamma_1} \wedge \cdots \wedge d \bar{z}_i^{\gamma_p}\wedge A_{\gamma_1,..., \gamma_p}^\alpha,\,\,\,\,\,\,\,\,\,\,\,\,A_{\gamma_1,..., \gamma_p}^\alpha\in \Gamma\left( U_i, \mathcal{A}^{0,0}\left(  \Omega_M^1\right) \right) \\
\psi\left( w_i^\alpha \right) &= \sum_{\gamma_1<\cdots <\gamma_p} d \bar{z}_i^{\gamma_1}\wedge \cdots \wedge d \bar{z}_i^{\gamma_p} \cdot B_{\gamma_1,..., \gamma_p}^\alpha,\,\,\,\,\,\,\,\,\,\,\ B_{\gamma_1,..., \gamma_p}^\alpha \in \Gamma\left( U_i, \mathcal{A}^{0,0} \right)
\end{align*}
Then we have
\begin{align*}
\hat{G}_1^\sharp \left( \overline{\phi}, \psi \right) \left( w_i^\alpha \right)=\left( \hat{E}_1^\sharp\left( \overline{\phi}  , \psi \right)\left( w_i^\alpha \right), \hat{G}_1'^\sharp \left( \phi\right)\left( w_i^\alpha \right) \right)
\end{align*}
where  $\hat{E}_1^\sharp\left( \overline{\phi} \right)$ is defined as in Part I, and 
{\tiny{\begin{align*}
&\hat{G}_1^\sharp\left( \overline{\phi}, \psi\right)\left( w_i^\alpha \right) \\
&= (-1)^{pq} \sum_{\gamma_1 < \cdots < \gamma_p} d\bar{z}_i^{\gamma_1} \wedge \cdots \wedge d\bar{z}_i^{\gamma_p} \left(  \sum_{\beta=1}^q (-1)^{q-\beta} B_{\gamma_1,..., \gamma_p}^\beta w_i^1\wedge \cdots  \wedge  \widehat{w}_i^\beta \wedge \cdots \wedge w_i^q \wedge d w_i^\alpha + w_i^1 \wedge \cdots \wedge w_i^q \wedge \partial \left( B_{\gamma_1,..., \gamma_p }^\alpha\right) - w_i^1 \wedge \cdots \wedge w_i^q \wedge A_{\gamma_1,..., \gamma_p}^\alpha    \right)
\end{align*}}}
In particular, for $p=1$, we have
{\Small{\begin{align*}
\hat{G}_1'\left( \overline{\phi}, \psi \right)\left( w_i^\alpha \right) = (-1)^q \left( \sum_{\beta=1}^q (-1)^{q-1} w_i^1 \wedge \cdots \wedge \psi\left( w_i^\beta \right)  \wedge \cdots \wedge w_i^q \wedge d w_i^\alpha - (-1)^q w_i^1 \wedge \cdots \wedge w_i^q \wedge \partial\left( \psi\left( w_i^\alpha \right) \right) - (-1)^q w_i^1 \wedge \cdots \wedge w_i^q \wedge \phi\left( w_i^\alpha  \right)\right)
\end{align*}}}

\subsubsection{\textnormal{Description of $\hat{G}_2$}}\

We describe
{\small{\begin{align*}
\hat{G}_2&: A^{0,p}\left( M-S, \mathscr{H}om_{\mathcal{O}_M} \left( \mathcal{N}_{\mathcal{F}_0}^*, \tilde{\mathcal{S}}^2   \right) \right) \bigoplus A^{0,p}\left( M-S, \mathscr{H}om_{\mathcal{O}_M}\left( \mathcal{N}_{\mathcal{F}_0}^*, \tilde{\mathcal{S}}^1   \right) \right)  \\
&\to A^{0,p}\left(M-S, \mathscr{H}om_{\mathcal{O}_M}\left( \mathcal{N}_{\mathcal{F}_0}^*, \tilde{\mathcal{S}}^3 \right) \right) \bigoplus A^{0,p}\left( M-S, \mathscr{H}om_{\mathcal{O}_M} \left( \mathcal{N}_{\mathcal{F}_0}^*, \tilde{\mathcal{S}}^2  \right) \right)
\end{align*}}}
For $x\in U_x \subset U_i - S$, we locally compute on $U_x$. We note that we have $dw_i^\alpha - \sum_{\beta=1}^q a_{i_x}^{\alpha\beta} \wedge w_i^\beta$ for some $a_{i_x}^{\alpha\beta} \in \Gamma\left( U_x , \Omega_M^1 \right)$. For $\left(C,D\right)\in A^{0,p}\left(M-S, \mathscr{H}om_{\mathcal{O}_M}\left( \mathcal{N}_{\mathcal{F}_0}^*, \tilde{\mathcal{S}}^2 \right) \right) \bigoplus A^{0,p}\left( M-S, \mathscr{H}om_{\mathcal{O}_M}\left( \mathcal{N}_{\mathcal{F}_0}^*, \tilde{\mathcal{S}}^1 \right)\right)$, we can write $C\left( w_i^\alpha \right)= w_{i}^1 \wedge \cdots \wedge w_i^q \wedge C_{i_x}^\alpha $ on $U_x$, where $C_{i_x}^\alpha=  \sum_{\gamma_1< \cdots < \gamma_r} d\bar{z}_i^{\gamma_1}\wedge \cdots \wedge d \bar{z}_i^{\gamma_p} \wedge C_{\gamma_1,..., \gamma_p}^\alpha$ with $C_{\gamma_1,..,\gamma_p}^\alpha\in \Gamma\left( U_x, \mathcal{A}^{0,0}\left(\bigwedge^2 \Omega_M^1\right) \right)$ and we can write $D\left( w_i^\alpha \right)= w_i^1 \wedge \cdots \wedge w_i^q \wedge D_{i_x}^\alpha$ on $U_x$, where $D_{i_x}^{\alpha}=\sum_{\gamma_1<\cdots <\gamma_p} d \bar{z}_i^{\gamma_1} \wedge \cdots \wedge d \bar{z}_i^{\gamma_p}\wedge D_{\gamma_1,..., \gamma_p}^\alpha $ with $D_{\gamma_1,...,\gamma_p}^\alpha \in \Gamma\left( U_x, \mathcal{A}^{0,0}\left( \Omega_M^1 \right) \right)$. Then we define on $U_x$
\begin{align*}
\hat{G}_2'\left(C, D \right) \left( w_i^\alpha \right) &= w_i^1 \wedge \cdots \wedge w_i^q \wedge \left( \sum_{\gamma_1<\cdots <\gamma_p} d\bar{z}_i^{\gamma_1}\wedge \cdots \wedge d \bar{z}_i^{\gamma_p}\wedge \left( \partial\left(  D_{\gamma_1,..., \gamma_p}^\alpha  \right) - \sum_{\beta=1}^q a_{i_x}^{\alpha\beta} \wedge D_{\gamma_1,..., \gamma_p}^\beta   +  C_{\gamma_1,..., \gamma_p}^\alpha  \right) \right) \\
&= (-1)^p w_i^1 \wedge \cdots \wedge w_i^q \wedge \left( \partial\left( D_{i_x}^\alpha \right)  - \sum_{\beta=1}^q a_{i_x}^{\alpha\beta} \wedge D_{i_x}^\beta  + (-1)^p C_{i_x}^\alpha \right)
\end{align*}

\section{Generalities on unfoldings over local Artin rings}\

 Our presentation is based on \cite{Ser06}. We will always denote by $k=\mathbb{C}$ the field of complex numbers. All schemes will be assumed to be defined over $k$, locally-noetherian and separated. We denote by $\bold{Art}$ the category of local artinian $k$-algebras with residue field $k$.

\begin{definition}[small extension]
We say that for $\left( \tilde{A}, \tilde{\mathfrak{m}} \right), \left(A, \mathfrak{m} \right)\in \bold{Art}$, an exact sequence of the form $0\to (t)\to \tilde{A} \to A \to 0$ is a small extension if $t\in \tilde{\mathfrak{m}}$ is annihilated by $\tilde{\mathfrak{m}}$ $($i.e. $t\cdot \tilde{\mathfrak{m}}=0)$, so that $(t)$ is an one dimensional $k$-vector space.
\end{definition}

\subsection{Generalities on infinitesimal unfoldings over $A\in \bold{Art}$} \label{b2}\

We extend the presentation \cite{Ser06} p.20-p.21 in the context of infinitesimal unfoldings. Let $(X, \mathcal{N}_{\mathcal{F}_0}^*)$ be a foliated nonsingular variety, so that we have an exact sequence
\begin{align*}
0\to \mathcal{N}_{\mathcal{F}_0}^*\to \Omega_X^1\to \Omega_X^1/\mathcal{N}_{\mathcal{F}_0}^*\to 0
\end{align*}
where $\mathcal{N}_{\mathcal{F}_0}^*$ is a coherent subsheaf of $\Omega_X^1$ such that $\Omega_X^1/\mathcal{N}_{\mathcal{F}_0}^*$ is torsion-free, and $d\left(\mathcal{N}_{\mathcal{F}_0, x}^* \right)\subset \mathcal{N}_{\mathcal{F}_0,x}^* \bigwedge \Omega_X^1$ for $x \in X- \textnormal{Sing}\left( \mathcal{F}_0 \right)$. A cartesian diagram of morphism of algebraic schemes
\begin{center}
$\eta:\begin{CD}
\left(X, \mathcal{N}_{\mathcal{F}_0}^* \right)@>>> \left(\mathcal{X}, \tilde{\mathcal{N}}_\mathcal{F}^* \right)\\
@VVV @VV\pi V\\
\textnormal{Spec}(k)@>>> \textnormal{Spec}(A)
\end{CD}$
\end{center}
where $\pi$ is flat, smooth and surjective, and $A\in \bold{Art}$, is called an infinitesimal unfolding of $\left(X, \mathcal{N}_{\mathcal{F}_0}^* \right)$ over $\textnormal{Spec}(A)$ if there exists an exact sequence 
\begin{align}\label{b3}
0\to  \tilde{\mathcal{N}}_\mathcal{F}^* \to \Omega_{\mathcal{X}/k}^1 \to \Omega_{\mathcal{X}/k}^1 / \tilde{\mathcal{N}}_\mathcal{F}^* \to 0
\end{align}
such that  $d_{\mathcal{X}} \left( \tilde{\mathcal{N}}_{\mathcal{F},x}^*\right) \subset \tilde{\mathcal{N}}_{\mathcal{F},x}^* \bigwedge \Omega_{\mathcal{X}/k,x}^1 $ for $x\in X - \textnormal{Sing}\left( \mathcal{F}_0\right)$. Moreover, from the exact sequence $0\to \pi^*\Omega_{\textnormal{Spec}(A)/k}^1\to \Omega_{\mathcal{X}/k}^1\xrightarrow{\alpha} \Omega_{\mathcal{X}/A}^1\to 0$, let $\mathcal{N}_{\mathcal{F}}^*$ be the image of $\tilde{\mathcal{N}}_{\mathcal{F}}^*$ under $\alpha$. Then $\tilde{\mathcal{N}}_{\mathcal{F}}^*  \cong \alpha\left( \tilde{\mathcal{N}}_{\mathcal{F}}^* \right):= \mathcal{N}_{\mathcal{F}}^*$ $($i.e $\alpha$ is isomorphic on $\tilde{\mathcal{N}}_{\mathcal{F}}^*$$)$ and $\left(\mathcal{X}, \mathcal{N}_{\mathcal{F}}^* \right)$ is an infinitesimal foliated deformation of $\left(X, \mathcal{N}_{\mathcal{F}_0}^* \right)$ over $A$ in the sense of Part II. In other words, we have an induced exact sequence from $\textnormal{(\ref{b3})}$ 
\begin{align}\label{b1}
0\to \mathcal{N}_{\mathcal{F}}^*\to \Omega_{\mathcal{X}/A}^1\to \Omega_{\mathcal{X}/A}^1/\mathcal{N}_{\mathcal{F}}^*\to 0
\end{align}
such that $\Omega_{\mathcal{X}/A}^1/ \mathcal{N}_{\mathcal{F}}^*$ is flat over $A$. An infinitesimal unfolding over $A=k[\epsilon]$ is called first-order unfolding. Given another infinitesimal unfolding of $\left( X, \mathcal{N}_{\mathcal{F}_0}^*\right)$ over $A\in \bold{Art}$
\begin{center}
$\xi: \begin{CD}
\left( X, \mathcal{N}_{\mathcal{F}_0}^* \right) @>>>  \left( \mathcal{Y}, \tilde{\mathcal{N}}_\mathcal{G}^*\right)\\
@VVV @VVV\\
\textnormal{Spec}(k) @>>> \textnormal{Spec}(A)
\end{CD}$
\end{center}
an isomorphism between infinitesimal unfoldings $\eta$ and $\xi$ is an unfolded $A$-isomorphism $\phi:\left(\mathcal{X}, \tilde{\mathcal{N}}_{\mathcal{F}}^* \right)\to \left(  \mathcal{Y}, \tilde{\mathcal{N}}_{\mathcal{G}}^* \right)$ with $\tilde{\mathcal{N}}_\mathcal{F}^*=\phi^* \tilde{\mathcal{N}}_{\mathcal{G}}^*$ as a subsheaf of $\Omega_{\mathcal{X}/k}^1$ and $\mathcal{N}_\mathcal{F}^* =\phi^* \mathcal{N}_{\mathcal{G}}^*$ inducing the identity on $\left( X, \mathcal{N}_{\mathcal{F}_0}^*\right)$, or equivalently the commutative diagrams
{\small{\begin{center}
$\begin{CD}
0 @>>> \tilde{\mathcal{N}}_\mathcal{F}^* @>>> \Omega_{\mathcal{X}/k}^1 \\
@. @V\cong VV @| \\
0@>>> \phi^* \tilde{\mathcal{N}}_\mathcal{G}^* @>>> \phi^* \Omega_{\mathcal{Y}/k}^1=\Omega_{\mathcal{X}/k}^1
\end{CD}$\,\,\,\,\,\,\,\,\,\,
$\begin{CD}
0 @>>>  \mathcal{N}_\mathcal{F}^* @>>> \Omega_{\mathcal{X}/A}^1 \\
@. @V\cong VV @| \\
0@>>> \phi^* \mathcal{N}_\mathcal{G}^* @>>> \phi^* \Omega_{\mathcal{Y}/A}^1=\Omega_{\mathcal{X}/A}^1
\end{CD}$\,\,\,\,\,\,\,
$\begin{CD}
\tilde{\mathcal{N}}_\mathcal{F}^* @>\alpha>\cong>\mathcal{N}_\mathcal{F}^* \\
@V\cong VV @VV\cong V\\
\phi^*\tilde{\mathcal{N}}_\mathcal{G}^* @>\alpha'>\cong> \phi^*\mathcal{N}_\mathcal{G}^*
\end{CD}$
\end{center}}}

\begin{remark}
We note that $\tilde{\mathcal{N}}_\mathcal{F}^*\cong \mathcal{N}_\mathcal{F}^*$ is flat over $A\in \bold{Art}$. We note that we did not assume that $\frac{\Omega_{\mathcal{X}/k}^1}{\tilde{\mathcal{N}}_{\mathcal{F}}^*}$ is flat over $A$. We have the exact sequence $0\to \pi^* \Omega_{\textnormal{Spec}(A)/k}^1 \to \frac{\Omega_{\mathcal{X}/k}^1}{\tilde{\mathcal{N}}_{\mathcal{F}}^*} \to \frac{\Omega_{\mathcal{X}/A}^1}{\mathcal{N}_\mathcal{F}^*}\to 0 $. If $\pi^* \Omega_{\textnormal{Spec}(A)/k}^1$ is not flat over $A$, then $\frac{\Omega_{\mathcal{X}/k}^1}{\tilde{\mathcal{N}}_\mathcal{F}^*}$ is not flat over $A$.
\end{remark}

\begin{definition}[Trivial unfolding]
Let $\left( X, \mathcal{N}_{\mathcal{F}_0}^* \right)$ be a nonsingular foliated variety over $k$.  We define the trivial unfoldig induced by $\left( X, \mathcal{N}_{\mathcal{F}_0}^* \right)$ and $A\in \bold{Art}$ to be the following infinitesimal unfolding of $\left( X, \mathcal{N}_{\mathcal{F}_0}^* \right)$ over $A$,
\begin{center}
$\begin{CD}
\left( X, \mathcal{N}_{\mathcal{F}_0}^* \right) @>>> \left(  X\times_k \textnormal{Spec}(A)    , \pi_1^* \mathcal{N}_{\mathcal{F}_0}^* \right)\\
@VVV @VV\pi_2 V\\
\textnormal{Spec}(k) @>>> \textnormal{Spec}(A)
\end{CD}$
\end{center}
where $\pi_1:X\times \textnormal{Spec}(A)\to X$ and $\pi_2:X\times \textnormal{Spec}(A) \to \textnormal{Spec}(A)$ are the projections. A flat unfolding of $\left( X, \mathcal{N}_{\mathcal{F}_0}^* \right)$ over $A$ is called trivial if it is isomorphic to the trivial unfolding over $A$.
\end{definition}

\begin{definition}[Rigid unfoldings]
Let $\left( X, \mathcal{N}_{\mathcal{F}_0}^* \right)$ be a foliated nonsingular variety. $\left(X, \mathcal{N}_{\mathcal{F}_0}^*\right)$ is called rigid in unfoldings if every infinitesimal unfolding of $\left( X, \mathcal{N}_{\mathcal{F}_0}^*\right)$ over $A$ is trivial for every $A \in \bold{Art}$.

\end{definition}

\section{Unfoldings of $\left( X, \mathcal{N}_{\mathcal{F}_0}^* \right)$ with $\mathcal{N}_{\mathcal{F}_0}^*$ locally free such that the underlying flat foliated deformations of $\left(X, \mathcal{N}_{\mathcal{F}}^* \right)$ are trivial}

Let $\left( X, \mathcal{N}_{\mathcal{F}_0}^* \right)$ be a foliated nonsingular variety such that $\mathcal{N}_{\mathcal{F}_0}^*$ is a locally free subsheaf of $\Omega_X^1$, so that we have an exact sequence $0\to \mathcal{N}_{\mathcal{F}_0}^*\to \Omega_X^1 \to \frac{\Omega_X^1}{\mathcal{N}_{\mathcal{F}_0}^*}\to 0$ where $\frac{\Omega_X^1}{\mathcal{N}_{\mathcal{F}_0}^*}$ is a torsion-free $\mathcal{O}_X$-module. We start from the most simple unfoldings. In this section we study flat unfoldings of $\left( X, \mathcal{N}_{\mathcal{F}_0}^*\right)$ such that the underlying flat foliated deformations of $\left(X, \mathcal{N}_{\mathcal{F}_0}^* \right)$ are trivial.

\begin{definition}
Let $\left( X, \mathcal{N}_{\mathcal{F}_0}^* \right)$ be a foliated nonsingular variety with $\mathcal{N}_{\mathcal{F}_0}^*$ locally free as above. An infinitesimal unfoldings $\left(\mathcal{X}=X\times \textnormal{Spec}(A), \tilde{\mathcal{N}}_{\mathcal{F}}^* \right)$ over $A\in \bold{Art}$ of $\left(X, \mathcal{N}_{\mathcal{F}_0} \right)$ such that the underlying foliated deformation of $\left(X, \mathcal{N}_{\mathcal{F}_0}^* \right)$ is trivial is a locally free subsheaf $\tilde{\mathcal{N}}_\mathcal{F}^*$ of $\Omega_{X\times \textnormal{Spec}(A)}^1$ such that $d\left(\tilde{\mathcal{N}}_{\mathcal{F},x}^* \right) \subset \tilde{\mathcal{N}}_{\mathcal{F},x}^* \bigwedge \Omega_{\mathcal{X}/k,x}^1$ for $x\in X- \textnormal{Sing}\left(\mathcal{F}_0 \right)$, and $\tilde{\mathcal{N}}_\mathcal{F}^* \cong \alpha\left( \tilde{\mathcal{N}}_{\mathcal{F}}^* \right)=\mathcal{N}_\mathcal{F}^* = \mathcal{N}_{\mathcal{F}_0}^*\otimes A$ under $\alpha:\Omega_{\mathcal{X}/k}^1 \to \Omega_{\mathcal{X}/A}^1$ and $\left(\mathcal{X},\mathcal{N}_\mathcal{F}^* \right)$ is the trivial foliated deformation and, $\tilde{\mathcal{N}}_{\mathcal{F}}^*$ induces $\mathcal{N}_{\mathcal{F}_0}^*$ over $\textnormal{Spec}(k)$. Then we have a functor of Artin rings
\begin{align*}
\textnormal{Unf}_{\left\{X, \mathcal{N}_{\mathcal{F}_0}^* \right\}_{\textnormal{fixed}}}: \bold{Art} \to &\left( \textnormal{Sets}\right) \\
 A \to  &\left\{ \textnormal{infinitesimal unfoldings of $\left(X, \mathcal{N}_{\mathcal{F}_0}^* \right)$ over $A$ such that}\right.\\
 &\left.  \textnormal{the underlying foliated deformation of $\left( X, \mathcal{N}_{\mathcal{F}_0}^* \right)$ is trivial}   \right\}
\end{align*}
 
\end{definition}
Then we will show that
{\small{\begin{align*}
\mathscr{H}om_{\mathcal{O}_X}\left( \mathcal{N}_{\mathcal{F}_0}^* , \mathcal{O}_X \right)^\bullet :\mathscr{H}om_{\mathcal{O}_X}\left( \mathcal{N}_{\mathcal{F}_0}^*, \mathcal{O}_X \right) \xrightarrow{G_0'} \mathscr{H}om_{\mathcal{O}_X}\left( \mathcal{N}_{\mathcal{F}_0}^*, \tilde{\mathcal{S}}^1 \right) \xrightarrow{G_1'} \mathscr{H}om_{\mathcal{O}_X}\left( \mathcal{N}_{\mathcal{F}_0}^* , \tilde{\mathcal{S}}^2 \right) \xrightarrow{E_2} \mathscr{H}om_{\mathcal{O}_X}\left( \mathcal{N}_{\mathcal{F}_0}^*, \tilde{\mathcal{S}}^3  \right) \xrightarrow{E_3} \cdots
\end{align*}}}
controls infinitesimal unfoldings of $\left( X , \mathcal{N}_{\mathcal{F}_0}^* \right)$ such that the underlying foliated deformations of $\left( X , \mathcal{N}_{\mathcal{F}_0}^* \right)$ are trivial. We will denote the $i$-th cohomology group of $\mathscr{H}om_{\mathcal{O}_X}\left( \mathcal{N}_{\mathcal{F}_0}^*, \mathcal{O}_X  \right)^\bullet$ by $\mathbb{H}^i\left( X, \mathscr{H}om_{\mathcal{O}_X}\left( \mathcal{N}_{\mathcal{F}_0}^*, \mathcal{O}_X \right)^\bullet \right)$. 

\begin{proposition}\label{uu5}
Let $\left( X, \mathcal{N}_{\mathcal{F}_0}^* \right)$ be a foliated nonsingular variety over $k$ with $\mathcal{N}_{\mathcal{F}_0}^*$ locally free. Then
\begin{enumerate}
\item there is a $1-1$ correspondence 
\begin{align*}
\textnormal{Unf}_{ \left\{ X, \mathcal{N}_{\mathcal{F}_0}^* \right\}_{\textnormal{fixed}}} \left(k[\epsilon] \right) \cong \mathbb{H}^0\left( X, \mathscr{H}om_{\mathcal{O}_X}\left(\mathcal{N}_{\mathcal{F}_0}^*, \mathcal{O}_X \right)^\bullet \right)
\end{align*}
\item Let $\eta=\left( \mathcal{X}, \mathcal{N}_\mathcal{F}^* \right)$ be an infinitesimal foliated deformations over $A:=k[x]/\left(x^m \right)$ of $\left( X, \mathcal{N}_{\mathcal{F}_0}^* \right)$ such that the underlying foliated deformation is trivial.  To small extension $0\to (t)=(x^{m}) \to \tilde{A}:=k[x]/(x^{m+1}) \to A:=k[x]/(x^m) \to 0$, we can associate an element $o_\eta(e)\in \mathbb{H}^1\left( X, \mathscr{H}om_{\mathcal{O}_X}\left( \mathcal{N}_{\mathcal{F}_0}^*,  \mathcal{O}_X \right)^\bullet \right)$ called the obstruction lifting of $\eta$ to $\tilde{A}$, which is $0$ if and only if a lifting of $\eta$ to $\tilde{A}$ exists.
\end{enumerate}
\end{proposition}

\begin{proof}
Given a first-order unfolding $\left(X\times \textnormal{Spec}\left( k[\epsilon]\right), \tilde{\mathcal{N}}_\mathcal{F}^*\right)$ of $\left( X, \mathcal{N}_{\mathcal{F}_0}^* \right)$ over $k[\epsilon]$ under the trivial foliated deformation of $\left(X, \mathcal{N}_{\mathcal{F}_0}^* \right)$, we choose an affine open covering $\mathcal{U}=\{U_i\}$ of $X$ such that $\Gamma\left( U_i, \mathcal{N}_{\mathcal{F}_0}^* \right)$ is generated by $w_{0i}^1,..., w_{0i}^q\in \Gamma\left( U_i, \Omega_X^1 \right)$ with $w_{0i}^\alpha= \sum_{\beta=1}^q h_{0ij}^{\alpha\beta} w_{0j}^\beta$ for some $h_{0ij}^{\alpha\beta}\in \Gamma\left( U_{ij}, \mathcal{O}_X \right)$ and $w_{0i}^1 \wedge \cdots \wedge w_{0i}^q \wedge d w_{0i}^\alpha =0, \alpha=1,...,q$, and $\Gamma\left( U_i, \tilde{\mathcal{N}}_\mathcal{F}^* \right)$ is generated by $W_i^\alpha=w_{0i}^\alpha  + v_i^\alpha d\epsilon$ for some $v_i^\alpha \in \Gamma\left( U_i, \mathcal{O}_X \right),\alpha=1,...,q$ with $ w_{0i}^\alpha + v_i^\alpha d\epsilon = \sum_{\beta=1}^q h_{0ij}^{\alpha \beta} \left( w_{0j}^\beta +  v_j^\beta d\epsilon \right)$ and $\left(w_{0i}^1+ v_i^1 d\epsilon \right)\wedge \cdots \left( w_{0i}^q + v_i^q d\epsilon \right)\wedge d \left(w_{0i}^\alpha + v_i^\alpha d\epsilon \right)=0,\alpha=1,...q$. Then we have
\begin{align} \label{uu1}
v_i^\alpha = \sum_{\beta=1}^q h_{0ij}^{\alpha\beta} v_j^\beta
\end{align}
and we have
\begin{align*}
0&=W_i^1\wedge \cdots \wedge W_i^q\wedge dW_i^\alpha= \left(w_{0i}^1+ v_i^1d\epsilon\right)\wedge \cdots \wedge \left(w_{0i}^q+ v_i^q d\epsilon \right)\wedge \left(dw_{0i}^\alpha+ dv_i^\alpha \wedge d\epsilon \right)\\
  &= w_{0i}^1\wedge \cdots \wedge w_{0i}^q\wedge dv_i^\alpha \wedge d\epsilon + \sum_{\gamma=1}^q w_{0i}^1\wedge \cdots \wedge \overbrace{ v_i^\gamma d\epsilon}^{\gamma-\textnormal{th}} \wedge \cdots \wedge w_{0i}^q\wedge dw_{0i}^\alpha\\
  &=w_{0i}^1\wedge \cdots \wedge w_{0i}^q\wedge dv_i^\alpha \wedge d\epsilon+\sum_{\gamma=1}^q(-1)^{q-\gamma} v_i^\gamma w_{0i}^1\wedge \cdots \wedge \hat{w}_{0i}^\gamma \wedge \cdots \wedge w_i^q\wedge dw_{0i}^\alpha \wedge d\epsilon 
\end{align*}
so that we have 
\begin{align} \label{uu2}
w_{0i}^1\wedge \cdots \wedge w_{0i}^q\wedge dv_i^\alpha+ \sum_{\gamma=1}^q (-1)^{q-\gamma}v_i^\gamma w_{0i}^1\wedge \cdots \wedge \hat{w}_{0i}^\gamma \wedge \cdots \wedge w_{0i}^q\wedge d w_{0i}^\alpha=0
\end{align}
We define $v_i\in \Gamma\left( U_i, \mathscr{H}om_{\mathcal{O}_X}\left( \mathcal{N}_{\mathcal{F}_0}^*, \mathcal{O}_X \right) \right)$ in the following way:
\begin{align}\label{uu12}
v_i: \Gamma\left( U_i, \mathcal{N}_{\mathcal{F}_0}^* \right) &\to \Gamma\left( U_i, \mathcal{O}_X \right) \\
 w_{0i}^\alpha &\mapsto v_i^\alpha \notag
\end{align}
and linearly extends to $\Gamma\left( U_i, \mathcal{N}_{\mathcal{F}_0}^* \right)$. Then $\textnormal{(\ref{uu1})}$ implies that $v=\left\{ v_i  \right\}\in H^0\left( X, \mathscr{H}om_{\mathcal{O}_X}\left( \mathcal{N}_{\mathcal{F}_0}^*, \mathcal{O}_X \right)\right)$ and $\textnormal{(\ref{uu2})}$ implies that $G_0'(v)=0$, so that $v\in \mathbb{H}^0\left( X, \mathscr{H}om_{\mathcal{O}_X}\left( \mathcal{N}_{\mathcal{F}_0}^*, \mathcal{O}_X \right)^\bullet \right)$

We define an inverse map. Given a $v\in \mathbb{H}^0\left(X, \mathscr{H}om_{\mathcal{O}_X}\left( \mathcal{N}_{\mathcal{F}_0}^* , \mathcal{O}_X \right)^\bullet \right)$, we have $\textnormal{(\ref{uu1})}$ and $\textnormal{(\ref{uu2})}$. By reversing the above process, we obtain a first-order unfolding of $\left(X, \mathcal{N}_{\mathcal{F}_0}^* \right)$ whose underlying foliated deformation is trivial.

Next we identify obstructions. Given a small extension $0\to (t)=(x^m)\to \tilde{A}=k[x]/(x^{m+1})\to A:=k[x]/(x^m)\to 0$, and an infinitesimal unfolding $\eta=\left( X\times \textnormal{Spec}(A) , \tilde{\mathcal{N}}_{\mathcal{F}}^* \right)$ with the underlying foliated deformation trivial over $A$, we choose an affine open covering $\mathcal{U}=\{U_i\}$ of $X$ such that $\Gamma\left( U_i, \tilde{\mathcal{N}}_{\mathcal{F}}^* \right)$ is generated by  $W_i^1,...,W_i^q\in \left( \Gamma\left(U_i, \Omega_X^1 \right)\otimes A \right) \bigoplus \left(\Gamma\left( U_i, \mathcal{O}_X\right)\otimes \Omega_{A/k}^1\right)$ such that $W_i^\alpha=w_{0i}^\alpha+v_i^\alpha$ where $v_i^\alpha\in \Gamma\left(U_i, \mathcal{O}_X\right)\otimes \Omega_{A/k}^1$ and we have $W_i^\alpha=\sum_{\beta=1}^q h_{0ij}^{\alpha\beta} W_j^\beta$. In order to see if there exists a lifting of $\left( \mathcal{X}, \tilde{\mathcal{N}}_{\mathcal{F}}^* \right)$ to $\tilde{A}$, we choose an arbitrary lifting $\tilde{v}_i^\alpha \in \Gamma\left( U_i, \mathcal{O}_X \right)\otimes \Omega_{\tilde{A}/k}^1$ of $v_i^\alpha$. We note that we have a surjection $\Omega_{\tilde{A}/k}^1\otimes_{\tilde{A}} A\to \Omega_{A/k}^1$. We let $\tilde{W}_i^\alpha:= w_{0i}^\alpha + \tilde{v}_i^\alpha \in \Gamma\left( U_i,  \Omega_{U_i\times \textnormal{Spec}\left( \tilde{A} \right)/k}^1 \right)$. Then we have
\begin{align*}
&\tilde{W}_i^\alpha-\sum_{\beta=1}^q h_{0ij}^{\alpha\beta}\tilde{W}_i^\beta= w_{0i}^\alpha+\tilde{v}_i^\alpha -\sum_{\beta=1}^q h_{0ij}^{\alpha\beta}\left( w_{0j}^\beta+\tilde{v}_j^\beta \right)=G_{ij}^\alpha dt \\
& \iff \tilde{v}_i^\alpha-\sum_{\beta=1}^q h_{0ij}^{\alpha\beta} \tilde{v}_j^\beta=G_{ij}^\alpha dt \,\,\,\,\,\,\textnormal{for some}\,\,\, G_{ij}^\alpha \in \Gamma\left( U_{ij}, \mathcal{O}_X \right)
\end{align*}
We define $G_{ij} \in \Gamma\left( U_{ij}, \mathscr{H}om_{\mathcal{O}_X}\left( \mathcal{N}_{\mathcal{F}_0}^*, \mathcal{O}_X \right) \right)$ by
\begin{align}\label{uu17}
G_{ij}:\Gamma\left( U_{ij}, \mathcal{N}_{\mathcal{F}_0}^* \right) &\to \Gamma\left( U_{ij}, \mathcal{O}_X \right)\\
 w_{0i}^\alpha &\mapsto G_{ij}^\alpha \notag
\end{align}
and linearly extends to $\Gamma\left( U_{ij}, \mathcal{N}_{\mathcal{F}_0}^* \right)$. We also note that we have
\begin{align*}
\left(w_{0i}^1+\tilde{v}_i^1\right)\wedge \cdots \wedge \left(w_{0i}^q+\tilde{v}_i^q \right)\wedge d \left(w_{0i}^\alpha+\tilde{v}_i^\alpha \right)= t B_i^\alpha+ F_i^\alpha \wedge dt,\,\,\,\,\,\,\,\, \textnormal{for some}\,\,\, B_i^\alpha \in \Gamma\left( U_i, \bigwedge^{q+2} \Omega_X^1 \right),\,\,\,\,\,\,\, F_i^\alpha \in \Gamma\left( U_i, \bigwedge^{q+1} \Omega_X^1 \right)
\end{align*}
In particular, we have $w_{0i}^1\wedge \cdots \wedge w_{0i}^q \wedge dw_{0i}^\alpha = t B_i^\alpha$, so that $B_i^\alpha=0$. We show that $F_i^\alpha \in \Gamma\left(U_i, \tilde{\mathcal{S}}^1 \right),\alpha=1,...,p $. In fact, for $x\in U_i-S$ with $S=\textnormal{Sing}\left(  \mathcal{F}_0 \right)$, we choose a neighborhood $U_x\subset U_i-S$ of $x$ such that $dw_{0i}^\alpha - \sum_{\beta=1}^q a_{0i_x}^{\alpha\beta} \wedge w_{0i}^\beta$ for some $a_{0i_{i_x}}^{\alpha\beta}\in \Gamma\left( U_i, \Omega_X^ 1\right)$. Then we can write
\begin{align*}
d W_i^\alpha - \sum_{\beta=1}^q \left( a_{i_x}^{\alpha\beta} + b_{i_x}^{\alpha\beta} \right) \wedge W_i^\beta=0 \iff d\left( w_{0i}^\alpha + v_i^\alpha \right) - \sum_{\beta=1}^q \left( a_{i_x}^{\alpha\beta} + b_{i_x}^{\alpha\beta} \right) \wedge \left( w_{0i}^\beta + v_i^\beta \right) = 0
\end{align*}
for some $a_{i_x}^{\alpha\beta}\in \Gamma\left( U_i, \Omega_X^1  \right)\otimes A$ which induces $a_{0i_x}^{\alpha\beta}$ by modulo $\mathfrak{m}_A$, and $b_{i_x}^{\alpha\beta}\in \Gamma\left( U_i,\mathcal{O}_X \right)\otimes \Omega_{A/k}^1$. Let $\tilde{a}_{i_x}^{\alpha\beta}\in \Gamma\left(U_x, \Omega_X^1 \right)\otimes \tilde{A}$ be an arbitrary lifting of $a_{i_x}^{\alpha\beta}$, $\tilde{b}_{i_x}^{\alpha\beta} \in \Gamma\left( U_i, \mathcal{O}_X \right)\otimes \Omega_{\tilde{A}/k}^1$ be an arbitrary lifting of $b_{i_x}^{\alpha\beta}$. Then we can write
\begin{align}\label{uu3}
 d \left( w_{0i}^\alpha + \tilde{v}_i^\alpha \right)- \sum_{\beta=1}^q \left( \tilde{a}_{i_x}^{\alpha\beta} + \tilde{b}_{i_x}^{\alpha\beta} \right) \wedge \left( w_{0i}^\beta + \tilde{v}_i^\beta \right)= t C_{i_x}^\alpha + F_{i_x}^\alpha \wedge dt \,\,\,\,\,\,\,\,\textnormal{for some}\,\,\, C_{i_x}^{\alpha}\in \Gamma\left(U_x, \bigwedge^2 \Omega_X^1 \right),\,\,\,\,\,F_{i_x}^\alpha \in \Gamma\left( U_x, \Omega_X^1 \right)
\end{align}
In particular, we have $dw_{0i}^\alpha - \sum_{\beta=1}^q \tilde{a}_{i_x}^{\alpha\beta}\wedge w_{0i}^\beta = t C_{i_x}^\alpha$. This implies that $w_{0i}^1\wedge \cdots \wedge w_{0i}^q\wedge C_{i_x}^{\alpha}=0$. Then we have
\begin{align*}
F_i^\alpha \wedge dt& =\left( w_{0i}^1 + \tilde{v}_i^1\right) \wedge \cdots \wedge \left( w_{0i}^q + \tilde{v}_i^q \right) \wedge d\left( w_{0i}^\alpha + \tilde{v}_i^\alpha \right)\\
&= \left( w_{0i}^1 + \tilde{v}_i^1\right) \wedge \cdots \wedge \left( w_{0i}^q + \tilde{v}_i^q \right) \wedge \left( t C_{i_x}^\alpha + F_{i_x}^\alpha \wedge dt \right)= w_{0i}^1 \wedge \cdots \wedge w_{0i}^q \wedge F_{i_x}^\alpha \wedge dt \\
\Longrightarrow & F_i^\alpha = w_{0i}^1 \wedge \cdots \wedge w_{0i}^q \wedge F_{i_x}^\alpha
\end{align*}
This implies that $F_i^\alpha \in \Gamma\left( U_i, \tilde{\mathcal{S}}^1 \right)$. We define $F_i\in \Gamma\left( U_i, \mathscr{H}om_{\mathcal{O}_X}\left( \mathcal{N}_{\mathcal{F}_0}^*, \tilde{\mathcal{S}}^1 \right) \right)$ by
\begin{align}
F_i: \Gamma\left( U_i, \mathcal{N}_{\mathcal{F}_0}^* \right) &\to \Gamma\left( U_i, \tilde{\mathcal{S}}^1 \right) \label{uu18} \\
 w_{0i}^\alpha &\mapsto F_i^\alpha \notag
\end{align}
Then we claim that $\left( \left\{ F_i\right\}, \left\{ - G_{ij} \right\} \right)\in C^0\left( \mathcal{U} , \mathscr{H}om_{\mathcal{O}_X}\left( \mathcal{N}_{\mathcal{F}_0}^*, \tilde{\mathcal{S}}^1 \right)  \right)\bigoplus C^1\left( \mathcal{U}, \mathscr{H}om_{\mathcal{O}_X}\left(  \mathcal{N}_{\mathcal{F}_0}^*, \mathcal{O}_X \right) \right)$ defines a $1$-cocycle in the following \v Cech resolution of $\mathscr{H}om_{\mathcal{O}_X}\left( \mathcal{N}_{\mathcal{F}_0}^*, \mathcal{O}_X \right)^\bullet$:
{\small{\begin{center}
$\begin{CD}
\cdots \\
@AAA \\
C^0\left( \mathcal{U}, \mathscr{H}om_{\mathcal{O}_X}\left( \mathcal{N}_{\mathcal{F}_0}^*, \tilde{\mathcal{S}}^2 \right) \right) @>\delta>> \cdots\\
@AAA @AAA \\
C^0\left( \mathcal{U}, \mathscr{H}om_{\mathcal{O}_X}\left( \mathcal{N}_{\mathcal{F}_0}^*, \tilde{\mathcal{S}}^1 \right) \right) @>-\delta>> C^1\left( \mathcal{U}, \mathscr{H}om_{\mathcal{O}_X}\left( \mathcal{N}_{\mathcal{F}_0}^*,  \tilde{\mathcal{S}}^1 \right) \right) @>\delta>> \cdots \\
@AAA @AAA @AAA \\
C^0\left( \mathcal{U}, \mathscr{H}om_{\mathcal{O}_X}\left( \mathcal{N}_{\mathcal{F}_0}^*, \mathcal{O}_X \right) \right) @>\delta >> C^1\left( \mathcal{U}, \mathscr{H}om_{\mathcal{O}_X}\left( \mathcal{N}_{\mathcal{F}_0}^*, \mathcal{O}_X \right) \right) @>-\delta >> C^2\left( \mathcal{U}, \mathscr{H}om_{\mathcal{O}_X}\left( \mathcal{N}_{\mathcal{F}_0}^*, \mathcal{O}_X \right) \right) @>\delta>> \cdots
\end{CD}$
\end{center}}}
We note that by taking $d$ on $\textnormal{(\ref{uu3})}$, we have
\begin{align}\label{uu4}
- \sum_{\beta=1}^q \left( d \tilde{a}_{i_x}^{\alpha\beta} + d\tilde{b}_{i_x}^{\alpha\beta} \right) \wedge \left( w_{0i}^\beta + \tilde{v}_i^\beta \right) + \sum_{\beta=1}^q \left( \tilde{a}_{i_x}^{\alpha\beta} + \tilde{b}_{i_x}^{\alpha\beta} \right) \wedge d\left( w_{0i}^\beta + \tilde{v}_i^\beta \right)= dt \wedge C_{i_x}^\alpha +t dC_{i_x}^\alpha + d F_{i_x}^\alpha \wedge dt
\end{align}
Then by multiplying $\left( w_{0i}^1 + \tilde{v}_i^1 \right) \wedge \cdots \wedge \left( w_{0i}^q + \tilde{v}_i^q \right)$ on $\textnormal{(\ref{uu4})}$, we have
\begin{align*}
\sum_{\beta=1}^q \left( w_{0i}^1 + \tilde{v}_i^1 \right)\wedge \cdots \wedge \left( w_{0i}^q + \tilde{v}_i^q \right)\wedge \left( \tilde{a}_{i_x}^{\alpha\beta} + \tilde{b}_{i_x}^{\alpha\beta} \right)\wedge \left( t C_{i_x}^\beta+ F_{i_x}^\beta \wedge dt  \right)=t w_{0i}^1\wedge \cdots \wedge w_{0i}^q \wedge d C_{i_x}^\alpha+ w_{0i}^1\wedge \cdots \wedge w_{0i}^q \wedge d F_{i_x}^\alpha \wedge dt,
\end{align*}
so that in particular by considering the coefficient of $dt$, we have 
\small{\begin{align*}
w_{0i}^1 \wedge \cdots \wedge w_{0i}^q \left( d F_{i_x}^\alpha - \sum_{\beta=1}^q a_{0i_x}^{\alpha\beta}\cdot F_{i_x}^\beta    \right)=0
\end{align*}
This implies that $G_1'\left(\left\{ F_i \right\} \right)=0$. By setting $h_{0ij}= \det\left( h_{0ij}^{\alpha\beta} \right)$, we have
{\begin{align*}
&F_i^\alpha \wedge dt= \left(w_{0i}^1+\tilde{v}_i^1\right)\wedge \cdots \wedge \left(w_{0i}^q+\tilde{v}_i^q \right) \wedge d \left(w_{0i}^\alpha + \tilde{v}_i^\alpha\right)\\
                                  &=\left(G_{ij}^1dt+\sum_{\beta_1=1}^q h_{0ij}^{1\beta}\left(w_{0j}^{\beta_1}+\tilde{v}_j^{\beta_1}\right) \right) \wedge \cdots \wedge \left( G_{ij}^qdt+\sum_{\beta_q=1}^q h_{0ij}^{q\beta_q}\left( w_{0j}^{\beta_q}+\tilde{v}_j^{\beta_q}\right) \right) \wedge d\left( G_{ij}^\alpha dt+\sum_{\beta=1}^q h_{0ij}^{\alpha\beta} \left(w_{0j}^\beta+ \tilde{v}_j^\beta\right) \right)\\
                                  &=\sum_{\gamma=1}^q w_{0i}^1\wedge \cdots \wedge G_{ij}^\gamma dt\wedge \cdots \wedge w_{0i}^q\wedge dw_{0i}^\alpha +w_{0i}^1\wedge \cdots \wedge w_{0i}^q \wedge d G_{ij}^\alpha \wedge dt      + h_{0ij} \left(w_{0j}^1+\tilde{v}_j^1\right)\wedge \cdots \wedge \left(w_{0j}^q+\tilde{v}_j^q \right) \wedge d\left(\sum_{\beta=1}^q h_{0ij}^{\alpha\beta}\left(w_{0j}^\beta+\tilde{v}_j^\beta\right) \right)\\
                                  &= \sum_{\gamma=1}^q w_{0i}^1\wedge \cdots \wedge G_{ij}^\gamma dt\wedge \cdots \wedge w_{0i}^q\wedge dw_{0i}^\alpha +w_{0i}^1\wedge \cdots \wedge w_{0i}^q \wedge d G_{ij}^\alpha \wedge dt   + h_{0ij} \sum_{\beta=1}^q h_{0ij}^{\alpha\beta} F_j^\beta \wedge dt
\end{align*}}}
This implies that
\begin{align*}
F_i^\alpha - h_{0ij}\sum_{\beta=1}^q h_{0ij}^{\alpha\beta} F_j^\beta - \sum_{\gamma=1}^q (-1)^{q-\gamma} G_{ij}^\gamma w_{0i}^1\wedge \cdots \wedge \hat{w}_{0i}^\gamma \wedge \cdots \wedge w_{0i}^q \wedge dw_{0i}^q - w_{0i}^1\wedge \cdots \wedge w_{0i}^q\wedge dG_{ij}^\alpha  =0,
\end{align*}
so that $-\delta\left( \left\{ F_i \right\} \right) + G_0'\left(\left\{- G_{ij} \right\} \right)=0$. We note that
\begin{align*}
G_{ij}^\alpha dt= \tilde{v}_i^\alpha-\sum_{\beta=1}^q h_{0ij}^{\alpha\beta} \tilde{v}_j^\beta=\tilde{v}_i^\alpha-\sum_{\beta=1}^q h_{0ij}^{\alpha\beta}  \left( G_{jk}^{\beta}dt+\sum_{\gamma=1}^q h_{0jk}^{\beta\gamma} \tilde{v}_k^\gamma \right)= G_{ik}^\alpha dt -\sum_{\beta=1}^q h_{0ij}^{\alpha\beta} G_{jk}^\beta dt
\end{align*}
This implies that $-\delta\left( \left\{ - G_{ij} \right\}\right)=0$. This proves the claim. We show that for another lifting $\tilde{v}_i'^\alpha$ of $v_i^\alpha,\alpha=1,...,q$, the associated $1$-cocycle $\beta:=\left( \left\{ F_i'\right\}, \left\{ - G_{ij}' \right\} \right)$ is cohomologous to $\alpha:=\left(\left\{ F_i\right\}, \left\{ - G_{ij} \right\} \right)$. Then $\tilde{v}_i'^\alpha = \tilde{v}_i^\alpha + r_i^\alpha dt$ for some $r_i^\alpha\in \Gamma\left( U_i, \mathcal{O}_X \right)$. We define $r_i\in \Gamma\left( U_i, \mathscr{H}om_{\mathcal{O}_X} \left( \mathcal{N}_{\mathcal{F}_0}^*, \mathcal{O}_X \right) \right)$ by
\begin{align*}
r_i:\Gamma\left( U_i, \mathcal{N}_{\mathcal{F}_0}^* \right) &\to \Gamma\left( U_i, \mathcal{O}_X \right) \\
 w_{0i}^\alpha &\mapsto r_i^\alpha
\end{align*}
Then we see that
\begin{align*}
\left( F_i'^\alpha - F_i^\alpha \right) \wedge dt& = \left(w_{0i}^1+\tilde{v}_i^1 + r_i^1 dt\right)\wedge \cdots \wedge \left(w_{0i}^q+\tilde{v}_i^q+ r_i^q dt\right) \wedge d\left(w_{0i}^\alpha + \tilde{v}_i^\alpha + r_i^\alpha dt\right) - \left(w_{0i}^1+\tilde{v}_i^1\right)\wedge \cdots \wedge \left(w_{0i}^q+\tilde{v}_i^q\right) \wedge d\left(w_{0i}^\alpha + \tilde{v}_i^\alpha\right)\\
&= \sum_{\beta=1}^q (-1)^{q-\beta} r_i^\beta w_{0i}^1 \wedge \cdots \wedge \widehat{w}_{0i}^\beta \wedge \cdots \wedge w_{0i}^q \wedge dw_{0i}^\alpha \wedge dt + w_{0i}^1 \wedge \cdots \wedge w_{0i}^q \wedge d r_i^\alpha \wedge dt
\end{align*}
This implies that $G_0'\left(  \left\{ r_i \right\} \right)= \left\{F_i'- F_i \right\}$. On the other hand, we have
\begin{align*}
\left( -G_{ij}'^\alpha+ G_{ij}^\alpha \right)dt = - \tilde{v}_i^\alpha  - r_i^\alpha dt + \sum_{\beta=1}^q h_{0ij}^{\alpha\beta}\left( \tilde{v}_j^\beta + r_j^\beta dt \right) + \tilde{v}_i^\alpha - \sum_{\beta=1}^q h_{0ij}^{\alpha\beta} \tilde{v}_j^\beta = \left( \sum_{\beta=1}^q h_{0ij}^{\alpha\beta} r_j^\beta - r_i^\alpha \right)dt
\end{align*}
This implies that $\delta\left( \left\{ r_i \right\} \right)=\left\{ - G_{ij}' - \left(- G_{ij} \right) \right\}$. This proves that $\alpha$ is cohomologous to $\beta$. So given a small extension $e: 0 \to (t)=(x^m) \to \tilde{A}= k[x]/(x^{m+1})\to A= k[x]/(x^m) \to 0$ and an infinitesimal unfolding $\eta$ of $\left( X, \mathcal{N}_{\mathcal{F}_0}^* \right)$ with the underlying foliated deformation trivial, we can associate an element $o_\eta(e):=$ the cohomology class of $\alpha \in \mathbb{H}^1\left( X, \mathscr{H}om_{\mathcal{O}_X}\left( \mathcal{N}_{\mathcal{F}_0}^*, \mathcal{O}_X \right)^\bullet \right)$. We note that $o_\eta(e)=0$ if and only if there exists $\left\{ M_i \right\}\in C^0\left( \mathcal{U}, \mathscr{H}om_{\mathcal{O}_X}\left( \mathcal{N}_{\mathcal{F}_0}^*, \mathcal{O}_X \right) \right)$ where $M_i:\Gamma\left( U_i, \mathcal{N}_{\mathcal{F}_0}^*, \mathcal{O}_X \right), w_{0i}^\alpha\mapsto M_i^\alpha$ such that $M_j- M_i = - G_{ij}$ and $G_0'\left( M_i \right)= F_i$, so that
\begin{align*}
&\left( \sum_{\beta=1}^q h_{0ij}^{\alpha\beta} M_j^\beta - M_i^\alpha \right) dt = - \tilde{v}_i^\alpha + \sum_{\beta=1}^q h_{0ij}^{\alpha\beta} \tilde{v}_j^\beta \Longrightarrow w_{0i}^\alpha + \tilde{v}_i^\alpha - M_i^\alpha dt = \sum_{\beta=1}^q h_{0ij}^{\alpha\beta}\left( w_{0j}^\beta + \tilde{v}_j^\beta - M_j^\beta dt \right)\\
&\left(\sum_{\beta=1}^q (-1)^{q-\beta} M_i^\beta w_{0i}^1 \wedge \cdots \wedge \widehat{w}_{0i}^\beta \wedge \cdots \wedge w_{0i}^q \wedge d w_{0i}^\alpha + w_{0i}^1 \wedge \cdots \wedge w_{0i}^q \wedge d M_i^\alpha \right) \wedge dt= \left( w_{0i}^1 + \tilde{v}_i^1 \right) \wedge \cdots \wedge \left( w_{0i}^q + \tilde{v}_i^q \right) \wedge d\left( w_{0i}^\alpha + \tilde{v}_i^\alpha \right)\\
&\Longrightarrow \left( w_{0i}^1 + \tilde{v}_i^\alpha - M_i^1 dt \right)\wedge \cdots \wedge \left( w_{0i}^q + \tilde{v}_i^q - M_i^q dt \right) \wedge d\left( w_{0i}^\alpha + \tilde{v}_i^\alpha - M_i^\alpha dt \right)=0
\end{align*}
Hence a lifting $\tilde{\eta}$ of $\eta$ to $\tilde{A}$ exists. This completes the proof of Proposition \ref{uu5}.

\end{proof}

\section{Unfoldings of $\left(X, \mathcal{N}_{\mathcal{F}_0}^* \right)$ with $\mathcal{N}_{\mathcal{F}_0}^*$ locally free such that the underlying flat deformations of $\left( X, \mathcal{N}_{\mathcal{F}_0}^* \right)$ are trivial}

Let $\left( X, \mathcal{N}_{\mathcal{F}_0}^* \right)$ be a foliated nonsingular variety with $\mathcal{N}_{\mathcal{F}_0}^*$ locally free. In the previous section we study infinitesimal unfoldings of $\left( X, \mathcal{N}_{\mathcal{F}_0}^* \right)$ such that the underlying foliated deformations of $\left( X, \mathcal{N}_{\mathcal{F}_0}^* \right)$ are trivial. In this section, in addition we study the second simple case, i.e. infinitesimal unfoldings of $\left(X, \mathcal{N}_{\mathcal{F}_0}^* \right)$ in a way that the underlying deformation of a pair of $X$ and $\mathcal{N}_{\mathcal{F}_0}^*$ is trivial.

\begin{definition}
Let $\left( X, \mathcal{N}_{\mathcal{F}_0}^* \right)$ be a foliated nonsingular variety with $\mathcal{N}_{\mathcal{F}_0}^*$ locally free as above. An infinitesimal unfolding $\left(\mathcal{X}=X\times \textnormal{Spec}(A), \tilde{\mathcal{N}}_{\mathcal{F}}^* \right)$ over $A\in \bold{Art}$ of $\left(X, \mathcal{N}_{\mathcal{F}_0} \right)$ such that the underlying deformation of $\left(X, \mathcal{N}_{\mathcal{F}_0}^* \right)$ is trivial, is a locally free subsheaf $\tilde{\mathcal{N}}_\mathcal{F}^*$ of $\Omega_{X\times \textnormal{Spec}(A)}^1$ such that $d\left(\tilde{\mathcal{N}}_{\mathcal{F},x}^* \right) \subset \tilde{\mathcal{N}}_{\mathcal{F},x}^* \bigwedge \Omega_{\mathcal{X}/k,x}^1$ for $x\in X- \textnormal{Sing}\left(\mathcal{F}_0 \right)$ and $\tilde{\mathcal{N}}_\mathcal{F}^* \cong \alpha\left( \tilde{\mathcal{N}}_{\mathcal{F}}^* \right):= \mathcal{N}_\mathcal{F}^* \cong \mathcal{N}_{\mathcal{F}_0}^*\otimes A$ under $\alpha:\Omega_{\mathcal{X}/k}^1\to \Omega_{\mathcal{X}/A}^1$ and $\left( \mathcal{X}, \mathcal{N}_\mathcal{F}^* \right)$ defines a flat foliated deformation over $A$ with $\mathcal{N}_{\mathcal{F}_0}^*$ and $X$ fixed, and $\tilde{\mathcal{N}}_{\mathcal{F}}^*$ induces $\mathcal{N}_{\mathcal{F}_0}^*$ over $\textnormal{Spec}(k)$. Equivalently an infinitesimal unfolding over $A$ of $\left(X, \mathcal{N}_{\mathcal{F}_0}^*\right)$  with $\mathcal{N}_{\mathcal{F}_0}^*$ and $X$ fixed is an exact sequence on $\mathcal{X}= X\times_k \textnormal{Spec}(A)$
\begin{align}
0\to \mathcal{N}_{\mathcal{F}_0}^* \otimes A \xrightarrow{l} \Omega_{\mathcal{X}/k}^1 \to \frac{\Omega_{\mathcal{X}/k}^1}{ l\left(\mathcal{N}_{\mathcal{F}_0}^* \otimes A\right)}\to 0
\end{align}
which induces $0\to \mathcal{N}_{\mathcal{F}_0}^* \to \Omega_X^1 \to \frac{\Omega_X^1}{\mathcal{N}_{\mathcal{F}_0}^*}\to 0$ over $\textnormal{Spec}(k)$ such that $d\left(l\left(\mathcal{N}_{\mathcal{F}_0}^*\otimes A \right)_x \right)\subset l\left( \mathcal{N}_{\mathcal{F}_0}^*\otimes A\right)_x \bigwedge \Omega_{\mathcal{X}/k,x}^1$ for $x\in X-\textnormal{Sing}\left( \mathcal{F}_0 \right)$, and $\Omega_{\mathcal{X}/k}^1/l\left(\mathcal{N}_{\mathcal{F}_0}^* \otimes A \right)$ is flat over $\textnormal{Spec}(A)$. Then we have a functor of Artin rings
\begin{align*}
\textnormal{Unf}_{\left(X, \mathcal{N}_{\mathcal{F}_0}^*\right)_{\textnormal{fixed}}}:\bold{Art}&\to \textnormal{Sets}\\
                                                                              A &\mapsto        \{\textnormal{infinitesimal unfoldings of $\left(X, \mathcal{N}_{\mathcal{F}_0}^*\right)$ over $A$  with $\mathcal{N}_{\mathcal{F}_0}^*$ and $X$ fixed} \}
\end{align*}

\end{definition}
Then we will show that
{\small{\begin{align*}
\textnormal{Unf}_{\left( X, \mathcal{N}_{\mathcal{F}_0}^* \right)_{\textnormal{fixed}}}^\bullet: \mathscr{H}om_{\mathcal{O}_X}\left( \mathcal{N}_{\mathcal{F}_0}^*, \Omega_X^1 \right) \bigoplus \mathscr{H}om_{\mathcal{O}_X}\left( \mathcal{N}_{\mathcal{F}_0}^*, \mathcal{O}_X \right) \xrightarrow{G_1''} \mathscr{H}om_{\mathcal{O}_X}\left( \mathcal{N}_{\mathcal{F}_0}^* \tilde{\mathcal{S}}^2 \right) \bigoplus \mathscr{H}om_{\mathcal{O}_X}\left( \mathcal{N}_{\mathcal{F}_0}^*, \tilde{\mathcal{S}}^1 \right) \\
 \xrightarrow{G_2}  \mathscr{H}om_{\mathcal{O}_X}\left( \mathcal{N}_{\mathcal{F}_0}^* \tilde{\mathcal{S}}^3 \right) \bigoplus \mathscr{H}om_{\mathcal{O}_X}\left( \mathcal{N}_{\mathcal{F}_0}^*, \tilde{\mathcal{S}}^2 \right)  \xrightarrow{G_3} \cdots
\end{align*}}}
controls infinitesimal unfoldings of $\left( X , \mathcal{N}_{\mathcal{F}_0}^*  \right)$ with $X$ and $\mathcal{N}_{\mathcal{F}_0}^*$ fixed. Here $G_1''$ is defined by the composition of the natural map $\mathscr{H}om_{\mathcal{O}_X}\left( \mathcal{N}_{\mathcal{F}_0}^*, \Omega_X^1  \right) \bigoplus \mathscr{H}om_{\mathcal{O}_X}\left( \mathcal{N}_{\mathcal{F}_0}^*, \mathcal{O}_X \right) \to \mathscr{H}om_{\mathcal{O}_X}\left( \mathcal{N}_{\mathcal{F}_0}^*, \frac{\Omega_X^1}{ \mathcal{N}_{\mathcal{F}_0}^*} \right) \bigoplus \mathscr{H}om_{\mathcal{O}_X}\left( \mathcal{N}_{\mathcal{F}_0}^*, \mathcal{O}_X \right)$ induced by the quotient with $G_1$. We will denote the $i$-th cohomology group of $\textnormal{Unf}_{\left( X, \mathcal{N}_{\mathcal{F}_0}^* \right)_{\textnormal{fixed}}}^\bullet$ by $\mathbb{H}^i\left( X, \textnormal{Unf}_{\left( X, \mathcal{N}_{\mathcal{F}_0}^* \right)_{ \textnormal{fixed}}}^\bullet \right)$.

\begin{proposition}\label{uu11}
Let $\left( X, \mathcal{N}_{\mathcal{F}_0}^* \right)$ be a foliated nonsingular variety over $k$ with $\mathcal{N}_{\mathcal{F}_0}^*$ locally free. Then
\begin{enumerate}
\item there is a $1-1$ correspondence 
\begin{align*}
\textnormal{Unf}_{\left( X, \mathcal{N}_{\mathcal{F}_0}^*\right)_{\textnormal{fixed}}} \left(k[\epsilon] \right) \cong \mathbb{H}^0\left( X, \textnormal{Unf}_{\left( X, \mathcal{N}_{\mathcal{F}_0}^* \right)_{\textnormal{fixed}}}^\bullet \right)
\end{align*}
\item Let $\eta=\left( \mathcal{X}, \mathcal{N}_\mathcal{F}^* \right)$ be an infinitesimal unfoldings over $A:=k[x]/\left( x^m\right)$ of $\left( X, \mathcal{N}_{\mathcal{F}_0}^* \right)$ with $X$ and $\mathcal{N}_{\mathcal{F}_0}^*$ fixed. To small extension $0\to (t)=(x^m) \to \tilde{A} : = k[x]/ \left( x^{m+1} \right) \to  A : = k[x]/\left( x^m \right) \to 0$, we can associate an element $o_\eta(e)\in \mathbb{H}^1\left( X, \textnormal{Unf}_{\left( X , \mathcal{N}_{\mathcal{F}_0}^* \right)_{\textnormal{fixed}}}^\bullet \right)$ called the obstruction lifting of $\eta$ to $\tilde{A}$, which is $0$ if and only if a lifting of $\eta$ to $\tilde{A}$ exists.
\end{enumerate}
\end{proposition}

\begin{proof}
Given a first-order unfolding $\left(\mathcal{X}, \tilde{\mathcal{N}}_\mathcal{F}^* \right)$ such that the underlying flat deformation of $\left(X, \mathcal{N}_{\mathcal{F}_0}^* \right)$ is trivial, we choose an affine open covering $\mathcal{U}= \left\{ U_i \right\}$ of $X$ such that $\Gamma \left( U_i, \mathcal{N}_{\mathcal{F}_0}^* \right)$ is generated by $w_{0i}^1,..., w_{0i}^q\in \Gamma\left( U_i, \Omega_X^1 \right)$ with $w_{0i}^\alpha = \sum_{\beta=1}^q h_{0ij}^{\alpha\beta} w_{0j}^\beta$ for some $h_{0ij}^{\alpha\beta} \in \Gamma\left( U_{ij}, \mathcal{O}_X \right),\alpha=1,...,q$ and $w_{0i}^1\wedge \cdots \wedge w_{0i}^q \wedge dw_{0i}^\alpha=0,\alpha=1,...,q$, and $\Gamma\left( U_i, \tilde{\mathcal{N}}_{\mathcal{F}}^*\right)$ is generated by $W_i^\alpha= w_{0i}^\alpha + \epsilon w_i'^\alpha + v_i^\alpha d\epsilon$ for some $w_i'^\alpha \in \Gamma\left( U_i, \Omega_X^1 \right)$ and $v_i^\alpha\in \Gamma\left(  U_i, \mathcal{O}_X \right)$ with $w_{0i}^\alpha + \epsilon w_i'^\alpha + v_i^\alpha d\epsilon =\sum_{\beta=1}^q h_{0ij}^{\alpha\beta}\left( w_{0j}^\beta + \epsilon w_j'^\beta + v_j^\beta d\epsilon \right)$, and $\left(w_{0i}^1 + \epsilon w_i'^1 + v_i^1 d\epsilon \right) \wedge \cdots \wedge \left( w_{0i}^q + \epsilon w_i'^q + v_i^q d \epsilon \right) \wedge d\left( w_{0i}^\alpha + \epsilon w_i'^\alpha + v_i^\alpha d \epsilon \right)=0,\alpha=1,...q$. Then we have
\begin{align}
&w_i'^\alpha = \sum_{\beta=1}^q h_{0ij}^{\alpha\beta} w_j'^\beta \label{uu13}\\
&v_i^\alpha = \sum_{\beta=1}^q h_{0ij}^{\alpha\beta} v_j^\beta \label{uu14}
\end{align}
and we have
\begin{align*}
0&=\left(w_{0i}^1 + \epsilon w_i'^1 + v_i^1 d\epsilon \right) \wedge \cdots \wedge \left( w_{0i}^q + \epsilon w_i'^q + v_i^q d \epsilon \right) \wedge d\left( w_{0i}^\alpha + \epsilon w_i'^\alpha + v_i^\alpha d \epsilon \right)\\
&= \left(w_{0i}^1 + \epsilon w_i'^1 + v_i^1 d\epsilon \right) \wedge \cdots \wedge \left( w_{0i}^q + \epsilon w_i'^q + v_i^q d \epsilon \right) \wedge \left( dw_{0i}^\alpha + d\epsilon \wedge w_i'^\alpha  + \epsilon d w_i'^\alpha + dv_i^\alpha \wedge d \epsilon \right) \\
&= \epsilon \left(  \sum_{\beta=1}^q w_{0i}^1 \wedge \cdots \wedge w_i'^\beta \wedge \cdots \wedge w_{0i}^q \wedge d w_{0i}^\alpha + w_{0i}^1 \wedge \cdots \wedge w_{0i}^q \wedge d w_i'^\alpha \right) \\
&+ \left( \sum_{\beta=1}^q (-1)^{q-\beta} v_i^\beta w_{0i}^1 \wedge \cdots \wedge \widehat{w}_{0i}^\beta \wedge \cdots \wedge w_{0i}^q \wedge d w_{0i}^\alpha + w_{0i}^1 \wedge \cdots \wedge w_{0i}^q \wedge d v_i^\alpha - w_{0i}^1\wedge \cdots \wedge w_{0i}^q \wedge w_i'^\alpha \right) \wedge d\epsilon
\end{align*}
This implies that
\begin{align}
&\sum_{\beta=1}^q w_{0i}^1 \wedge \cdots \wedge w_i'^\beta \wedge \cdots \wedge w_{0i}^q \wedge d w_{0i}^\alpha + w_{0i}^1 \wedge \cdots \wedge w_{0i}^q \wedge d w_i'^\alpha=0 \label{uu15}\\
 \sum_{\beta=1}^q (-1)^{q-\beta} & v_i^\beta  w_{0i}^1 \wedge \cdots \wedge  \widehat{w}_{0i}^\beta \wedge \cdots \wedge w_{0i}^q \wedge d w_{0i}^\alpha + w_{0i}^1 \wedge \cdots \wedge w_{0i}^q \wedge d v_i^\alpha - w_{0i}^1\wedge \cdots \wedge w_{0i}^q \wedge w_i'^\alpha \label{uu16}
\end{align}

We define $w_i'\in \Gamma\left( U_i, \mathscr{H}om_{\mathcal{O}_X}\left( \mathcal{N}_{\mathcal{F}_0}^*, \Omega_X^1 \right) \right)$ in the following way:
\begin{align}
w_i':\Gamma\left( U_i, \mathcal{N}_{\mathcal{F}_0}^* \right) &\to \Gamma\left( U_i, \Omega_X^1 \right)\\
 w_{0i}^\alpha &\mapsto w_i'^\alpha \notag
\end{align}
and linearly extends to $\Gamma\left(U_i, \mathcal{N}_{\mathcal{F}_0}^* \right)$. We define $v_i\in \Gamma\left( U_i, \mathscr{H}om_{\mathcal{O}_X}\left( \mathcal{N}_{\mathcal{F}_0}^*, \mathcal{O}_X \right)  \right)$ as in $\textnormal{(\ref{uu12})}$. Then $\textnormal{(\ref{uu13})},\textnormal{(\ref{uu14})}, \textnormal{(\ref{uu15})}$ and $\textnormal{(\ref{uu16})}$ imply $\left(w'=\left\{ w_i'\right\}, v= \left\{ v_i \right\} \right)\in \mathbb{H}^0 \left(X, \textnormal{Unf}_{\left( X, \mathcal{N}_{\mathcal{F}_0}^* \right)_{\textnormal{fixed}}}^\bullet \right)$.

We define an inverse map. Given a $\left(w', v \right) \in \mathbb{H}^0 \left(X, \textnormal{Unf}_{\left( X, \mathcal{N}_{\mathcal{F}_0}^* \right)_{\textnormal{fixed}}}^\bullet \right)$, we have $\textnormal{(\ref{uu13})},\textnormal{(\ref{uu14})}, \textnormal{(\ref{uu15})}$ and $\textnormal{(\ref{uu16})}$. By reversing the above process, we obtain a first-order unfolding of $\left(X, \mathcal{N}_{\mathcal{F}_0}^* \right)$ whose underlying flat deformation of $\left(X, \mathcal{N}_{\mathcal{F}_0}^* \right)$ is trivial.

Next we identify obstructions. Given a small extension $0\to (t)= (x^m) \to \tilde{A}= k[x]/(x^{m+1}) \to A= k[x]/(x^m) \to 0$, and an infinitesimal unfolding $\eta=\left( X\times \textnormal{Spec}(A), \tilde{\mathcal{N}}_\mathcal{F}^* \right)$ with the the underlying flat deformation of $\left( X, \mathcal{N}_{\mathcal{F}_0}^* \right)$ trivial, we choose an affine open covering $\mathcal{U}=\left\{ U_i \right\}$ of $X$ such that $\Gamma\left( U_i, \tilde{\mathcal{N}}_\mathcal{F}^* \right)$ is generated by $W_i^\alpha= w_i^\alpha + v_i^\alpha \in \left( \Gamma\left( U_i, \Omega_X^1 \right) \otimes A \right) \bigoplus \left( \Gamma\left( U_i, \mathcal{O}_X \right)\otimes \Omega_{A/k}^1 \right),\alpha=1,...,q$ with the relation $W_i^\alpha = \sum_{\beta=1}^q h_{0ij}^{\alpha\beta} W_j^\beta$. In order to see if there exists a lifting of $\left( X\times \textnormal{Spec}(A), \tilde{\mathcal{N}}_\mathcal{F}^* \right)$ to $\tilde{A}$, we choose an arbitrary lifting $\tilde{w}_i^\alpha \in \Gamma\left( U_i, \Omega_X^1 \right)\otimes \tilde{A}$, and an arbitrary lifting $\tilde{v}_i^\alpha\in \Gamma\left( U_i, \mathcal{O}_X \right)\otimes \Omega_{\tilde{A}/k}^1$ of $v_i^\alpha$. We let $\tilde{W}_i^\alpha:= \tilde{w}_i^\alpha + \tilde{v}_i^\alpha \in \Gamma\left( U_i, \Omega_{U_i\times \textnormal{Spec}\left( \tilde{A} \right)/k}^1 \right)$. Then we have
\begin{align*}
&\tilde{W}_i^\alpha - \sum_{\beta=1}^q h_{0ij}^{\alpha\beta} \tilde{W}_j^\beta= \tilde{w}_i^\alpha + \tilde{v}_i^\alpha - \sum_{\beta=1}^q h_{0ij}^{\alpha\beta}\left( \tilde{w}_j^\beta + \tilde{v}_j^\beta \right) = t R_{ij}^\alpha + G_{ij}^\alpha dt \\
&\iff \tilde{w}_i^\alpha - \sum_{\beta=1}^q h_{0ij}^{\alpha\beta} \tilde{w}_j^\beta = t R_{ij}^\alpha,\,\,\,\,\,\,\,\,\, \tilde{v}_i^\alpha - \sum_{\beta=1}^q h_{0ij}^{\alpha\beta} \tilde{v}_j^\beta = G_{ij}^\alpha dt \,\,\,\,\,\textnormal{for some}\,\,\, R_{ij}^\alpha \in \Gamma\left(U_i, \Omega_X^1 \right),\,\,\,G_{ij}^\alpha \in \Gamma\left(  U_i, \mathcal{O}_X \right)
\end{align*}
We define $R_{ij}\in \Gamma\left( U_{ij}, \mathscr{H}om_{\mathcal{O}_X}\left( \mathcal{N}_{\mathcal{F}_0}^*, \Omega_X^1 \right) \right)$ by 
\begin{align*}
R_{ij}: \Gamma\left( U_{ij}, \mathcal{N}_{\mathcal{F}_0}^* \right) &\to \Gamma\left( U_{ij}, \Omega_X^1 \right) \\
 w_{0i}^\alpha &\mapsto R_{ij}^\alpha
\end{align*}
and linearly extends to $\Gamma\left( U_{ij}, \mathcal{N}_{\mathcal{F}_0}^* \right)$. We define $G_{ij}\in \Gamma\left( U_{ij}, \mathscr{H}om_{\mathcal{O}_X}\left( \mathcal{N}_{\mathcal{F}_0}^*, \mathcal{O}_X  \right) \right)$ as in $\textnormal{(\ref{uu17})}$. On the other hand, we note that
\begin{align*}
\left(\tilde{w}_i^1+\tilde{v}_i^1\right)\wedge \cdots \wedge \left(\tilde{w}_i^q+\tilde{v}_i^q \right)\wedge d \left(\tilde{w}_i^\alpha+\tilde{v}_i^\alpha \right)= t B_i^\alpha+ F_i^\alpha \wedge dt,\,\,\,\,\,\,\,\, \textnormal{for some}\,\,\, B_i^\alpha \in \Gamma\left( U_i, \bigwedge^{q+2} \Omega_X^1 \right),\,\,\,\,\,\,\, F_i^\alpha \in \Gamma\left( U_i, \bigwedge^{q+1} \Omega_X^1 \right)
\end{align*}
We show that $B_i^\alpha \in \Gamma\left(U_i, \tilde{\mathcal{S}}^2 \right)$ and $F_i^\alpha \in \Gamma\left( U_i, \tilde{\mathcal{S}}^1 \right),\alpha=1,...,p$. In fact, for $x\in U_i-S$, we choose a neighborhood $U_x\subset U_i-S$ of $x$ such that $dw_{0i}^\alpha- \sum_{\beta=1}^q a_{0i_x}^{\alpha\beta} \wedge w_{0i}^\beta$ for some $a_{0_{i_x}}^{\alpha\beta}\in \Gamma\left( U_i, \Omega_X^1 \right)$. Then we can write $d\left( w_i^\alpha + v_i^\alpha \right)- \sum_{\beta=1}^q \left( a_{i_x}^{\alpha\beta} + b_{i_x}^{\alpha\beta} \right) \wedge \left( w_i^\beta + v_i^\beta \right)=0$ for some $a_{i_x}^{\alpha\beta}\in \Gamma\left( U_i, \Omega_X^1\right)\otimes A$ which induces $a_{0i_x}^{\alpha\beta}$ by modulo $\mathfrak{m}_A$, and $b_{i_x}^{\alpha\beta}\in \Gamma\left( U_i, \mathcal{O}_X \right)\otimes \Omega_{A/k}^1$. Let $\tilde{a}_{i_x}^{\alpha\beta}\in \Gamma\left(U_x, \Omega_X^1 \right)\otimes \tilde{A}$ be an arbitrary lifting of $a_{i_x}^{\alpha\beta}$, $\tilde{b}_{i_x}^{\alpha\beta} \in \Gamma\left( U_i, \mathcal{O}_X \right)\otimes \Omega_{\tilde{A}/k}^1$ be an arbitrary lifting of $b_{i_x}^{\alpha\beta}$. Then we can write
\begin{align}\label{uu19}
 d \left( \tilde{w}_i^\alpha + \tilde{v}_i^\alpha \right)- \sum_{\beta=1}^q \left( \tilde{a}_{i_x}^{\alpha\beta} + \tilde{b}_{i_x}^{\alpha\beta} \right) \wedge \left( \tilde{w}_i^\beta + \tilde{v}_i^\beta \right)= t C_{i_x}^\alpha + F_{i_x}^\alpha \wedge dt \,\,\,\,\,\,\,\,\textnormal{for some}\,\,\, C_{i_x}^{\alpha}\in \Gamma\left(U_x, \bigwedge^2 \Omega_X^1 \right),\,\,\,\,\,F_{i_x}^\alpha \in \Gamma\left( U_x, \Omega_X^1 \right)
\end{align}
Then we have
\begin{align*}
tB_i^\alpha + F_i^\alpha \wedge dt&= \left(\tilde{w}_i^1+\tilde{v}_i^1\right)\wedge \cdots \wedge \left(\tilde{w}_i^q+\tilde{v}_i^q \right)\wedge d \left(\tilde{w}_i^\alpha+\tilde{v}_i^\alpha \right)\\
&= \left(\tilde{w}_i^1+\tilde{v}_i^1\right)\wedge \cdots \wedge \left(\tilde{w}_i^q+\tilde{v}_i^q \right) \wedge \left( t C_{i_x}^\alpha + F_{i_x}^\alpha \wedge dt \right)= t w_{0i}^1\wedge \cdots \wedge w_{0i}^q \wedge C_{i_x}^\alpha + w_{0i}^1 \wedge \cdots \wedge w_{0i}^q\wedge F_{i_x}^\alpha \wedge dt\\
\Longrightarrow & B_i^\alpha = w_{0i}^1\wedge \cdots \wedge w_{0i}^q \wedge C_{i_x}^\alpha,\,\,\,\,\,\,\,\,F_i^\alpha = w_{0i}^1 \wedge \cdots \wedge w_{0i}^q \wedge F_{i_x}^\alpha
\end{align*}
This implies that $B_i^\alpha \in \Gamma\left(U_i, \tilde{\mathcal{S}}^2 \right)$ and $F_i^\alpha \in \Gamma\left( U_i, \tilde{\mathcal{S}}^1 \right)$. We define $B_i \in \Gamma\left( U_i, \mathscr{H}om_{\mathcal{O}_X}\left( \mathcal{N}_{\mathcal{F}_0}^*, \tilde{\mathcal{S}}^2 \right) \right)$ by
\begin{align*}
B_i: \Gamma\left( U_i, \mathcal{N}_{\mathcal{F}_0}^* \right) &\to \Gamma\left( U_i, \tilde{\mathcal{S}}^2 \right) \\
 w_{0i}^\alpha &\mapsto B_i^\alpha
\end{align*}
and we define $F_i \in \Gamma\left( U_i, \mathscr{H}om_{\mathcal{O}_X}\left(  \mathcal{N}_{\mathcal{F}_0}^*, \tilde{\mathcal{S}}^1 \right)\right)$ as in $\textnormal{(\ref{uu18})}$. Then we claim that $\left( \left\{ B_i \right\}, \left\{ F_i \right\}, \left\{ - R_{ij} \right\}, \left\{ - G_{ij} \right\} \right)\in C^0\left( \mathcal{U}, \mathscr{H}om_{\mathcal{O}_X}\left( \mathcal{N}_{\mathcal{F}_0}^*, \tilde{\mathcal{S}}^2 \right) \right) \bigoplus C^0\left( \mathcal{U}, \mathscr{H}om_{\mathcal{O}_X}\left( \mathcal{N}_{\mathcal{F}_0}^*  , \tilde{\mathcal{S}}^1 \right) \right) \bigoplus C^1\left( \mathcal{U}, \mathscr{H}om_{\mathcal{O}_X}\left( \mathcal{N}_{\mathcal{F}_0}^*, \Omega_X^1 \right) \right) \bigoplus C^1\left( \mathcal{U}, \mathscr{H}om_{\mathcal{O}_X}\left( \mathcal{N}_{\mathcal{F}_0}^*, \mathcal{O}_X \right) \right)$ defines a $1$-cocycle in the following \v Cech resolution of $\textnormal{Unf}_{\left( X, \mathcal{N}_{\mathcal{F}_0}^* \right)_{\textnormal{fixed}}}^\bullet$:
{\Small{\begin{center}
$\begin{CD}
\cdots \\
@AAA\\
C^0\left( \mathcal{U}, \left( \mathcal{N}_{\mathcal{F}_0}^* \right)^*\otimes \tilde{\mathcal{S}}^2 \right) \bigoplus C^0\left( \mathcal{U}, \left(  \mathcal{N}_{\mathcal{F}_0}^* \right)^* \otimes \tilde{\mathcal{S}}^1 \right) @>-\delta>> \cdots \\
@AAA @AAA \\
C^0\left( \mathcal{U}, \left(\mathcal{N}_{\mathcal{F}_0}^* \right)^* \otimes \Omega_X^1 \right) \bigoplus C^0\left( \mathcal{U}, \left(\mathcal{N}_{\mathcal{F}_0}^* \right)^* \right) @>\delta>> C^1\left( \mathcal{U}, \left(\mathcal{N}_{\mathcal{F}_0}^* \right)^* \otimes \Omega_X^1 \right) \bigoplus C^0\left( \mathcal{U}, \left(\mathcal{N}_{\mathcal{F}_0}^* \right)^* \right) @>-\delta>> \cdots
\end{CD}$
\end{center}}}
We note that
\begin{align}
&tR_{ij}^\alpha+ G_{ij}^\alpha dt - \left(tR_{ik}^\alpha+ G_{ik}^\alpha dt \right)+\sum_{\beta=1}^q h_{0ij}^{\alpha\beta} \left(tR_{jk}^\beta + G_{jk}^\beta dt \right) \\
&= \left(\tilde{w}_i^\alpha + \tilde{v}_i^\alpha \right) -\sum_{\beta=1}^q h_{0ij}^{\alpha\beta} \left( \tilde{w}_j^\beta + \tilde{v}_j^\beta  \right) - \left(\tilde{w}_i^\alpha + \tilde{v}_i^\alpha \right)+\sum_{\beta=1}^q h_{0ik}^{\alpha\beta} \left( \tilde{w}_k^\beta + \tilde{v}_i^\beta \right) +\sum_{\beta=1}^q h_{0ij}^{\alpha\beta} \left( \tilde{w}_j^\beta + \tilde{v}_j^\beta \right) -\sum_{\beta,\gamma=1}^q h_{0ij}^{\alpha\beta} h_{0jk}^{\beta\gamma} \left( \tilde{w}_k^\gamma + \tilde{v}_k^\gamma \right) = 0 \notag 
\end{align}
This implies that $-\delta\left( \left\{ - R_{ij} \right\} ,\left\{ - G_{ij} \right\} \right) = 0$. We note that by taking $d$ on $\textnormal{(\ref{uu19})}$, we have
\begin{align}\label{uu20}
- \sum_{\beta=1}^q \left( d \tilde{a}_{i_x}^{\alpha\beta} + d\tilde{b}_{i_x}^{\alpha\beta} \right) \wedge \left( \tilde{w}_{i}^\beta + \tilde{v}_i^\beta \right) + \sum_{\beta=1}^q \left( \tilde{a}_{i_x}^{\alpha\beta} + \tilde{b}_{i_x}^{\alpha\beta} \right) \wedge d\left( \tilde{w}_{i}^\beta + \tilde{v}_i^\beta \right)= dt \wedge C_{i_x}^\alpha +t dC_{i_x}^\alpha + d F_{i_x}^\alpha \wedge dt
\end{align}
Then by multiplying $\left( \tilde{w}^1 + \tilde{v}_i^1 \right) \wedge \cdots \wedge \left( \tilde{w}^q + \tilde{v}_i^q \right)$ on $\textnormal{(\ref{uu20})}$, we have
\begin{align*}
&\sum_{\beta=1}^q \left( \tilde{w}_{i}^1 + \tilde{v}_i^1 \right)\wedge \cdots \wedge \left( \tilde{w}_{i}^q + \tilde{v}_i^q \right)\wedge \left( \tilde{a}_{i_x}^{\alpha\beta} + \tilde{b}_{i_x}^{\alpha\beta} \right)\wedge \left( t C_{i_x}^\beta+ F_{i_x}^\beta \wedge dt  \right)\\
&= w_{0i}^1 \wedge \cdots \wedge w_{0i}^q \wedge C_{i_x}^\alpha \wedge dt + t w_{0i}^1\wedge \cdots \wedge w_{0i}^q \wedge d C_{i_x}^\alpha+ w_{0i}^1\wedge \cdots \wedge w_{0i}^q \wedge d F_{i_x}^\alpha \wedge dt,\\
&\iff t w_{0i}^1 \wedge \cdots \wedge w_{0i}^q \wedge \left( \sum_{\beta=1}^q a_{0i_x}^{\alpha\beta} \wedge C_{i_x}^\beta \right) + w_{0i}^1 \wedge \cdots \wedge w_{0i}^q \wedge \left( \sum_{\beta=1}^q a_{0i_x}^{\alpha\beta} \wedge F_{i_x}^\beta  \right) \wedge dt\\
&\,\,\,\,\,\,\,\,\,\,\,\,\,= w_{0i}^1 \wedge \cdots \wedge w_{0i}^q \wedge C_{i_x}^\alpha \wedge dt + t w_{0i}^1\wedge \cdots \wedge w_{0i}^q \wedge d C_{i_x}^\alpha+ w_{0i}^1\wedge \cdots \wedge w_{0i}^q \wedge d F_{i_x}^\alpha \wedge dt
\end{align*}
Then we have
\begin{align*}
w_{0i}^1 \wedge \cdots \wedge w_{0i}^q \wedge \left(  d C_{i_x}^\alpha - \sum_{\beta=1}^q a_{0i_x}^{\alpha\beta} \wedge C_{i_x}^\beta \right) = 0,\,\,\,\,\,\,\,\, w_{0i}^1 \wedge\cdots \wedge w_{0i}^q \wedge \left( d F_{i_x}^\alpha - \sum_{\beta=1}^q a_{0i_x}^{\alpha\beta} \wedge F_{i_x}^\beta + C_{i_x} \right) = 0
\end{align*}
This implies that $G_2 \left( \left\{ B_i \right\}  , \left\{ F_i \right\} \right)=0$. By setting $h_{0ij}= \det\left( h_{0ij}^{\alpha\beta} \right)$, we have
{\tiny{\begin{align*}
&t B_i^\alpha+ F_i^\alpha \wedge dt=\left(\tilde{w}_i^1+ \tilde{v}_i^1\right)\wedge \cdots \wedge \left(\tilde{w}_i^q+ \tilde{v}_i^q\right) \wedge d\left( \tilde{w}_i^\alpha+ \tilde{v}_i^\alpha\right)\\
     &=\left(tR_{ij}^1+G_{ij}^1 dt +\sum_{\beta_1=1}^q h_{0ij}^{\alpha\beta_1}\left( \tilde{w}_j^{\beta_1}+ \tilde{v}_j^{\beta_1}\right) \right) \wedge \cdots \wedge \left( t R_{ij}^q+ G_{ij}^q dt +\sum_{\beta_q=1}^q h_{0ij}^{\alpha\beta_q} \left(\tilde{w}_j^{\beta_q} + \tilde{v}_j^{\beta_q} \right)\right) \wedge d \left( tR_{ij}^\alpha+ G_{ij}^\alpha dt +\sum_{\beta=1}^q h_{0ij}^{\alpha\beta}\left( \tilde{w}_j^\beta  + \tilde{v}_j^\beta\right) \right)\\
          &=\left(tR_{ij}^1+G_{ij}^1 dt +\sum_{\beta_1=1}^q h_{0ij}^{\alpha\beta_1}\left( \tilde{w}_j^{\beta_1}+ \tilde{v}_j^{\beta_1}\right) \right) \wedge \cdots \wedge \left( t R_{ij}^q+ G_{ij}^q dt +\sum_{\beta_q=1}^q h_{0ij}^{\alpha\beta_q} \left(\tilde{w}_j^{\beta_q} + \tilde{v}_j^{\beta_q} \right)\right) \wedge  \left( dt\wedge R_{ij}^\alpha+tdR_{ij}^\alpha+ dG_{ij}^\alpha \wedge dt +d\left(\sum_{\beta=1}^q h_{0ij}^{\alpha\beta}\left( \tilde{w}_j^\beta  + \tilde{v}_j^\beta\right)\right) \right)\\
     &=t \sum_{\gamma=1}^q w_{0i}^1 \wedge \cdots \wedge R_{ij}^\gamma \wedge \cdots \wedge w_{0i}^q\wedge d w_{0i}^\alpha +t w_{0i}^1\wedge \cdots \wedge w_{0i}^q \wedge dR_{ij}^\alpha\\
     &+ \left( \sum_{\gamma=1}^q(-1)^{q-\gamma} G_{ij}^\gamma w_{0i}^1\wedge \cdots \wedge \hat{w}_{0i}^\gamma \wedge \cdots \wedge w_{0i}^q\wedge  d w_{0i}^\alpha + w_{0i}^1 \wedge \cdots \wedge w_{0i}^q \wedge dG_{ij}^\alpha  - w_{0i}^1 \wedge \cdots \wedge w_{0i}^q \wedge R_{ij}^\alpha\right) \wedge dt\\
     & + h_{0ij} \sum_{\beta=1}^q h_{0ij}^{\alpha\beta} \left(\tilde{w}_j^1+\tilde{v}_j^1\right) \wedge \cdots \wedge \left( \tilde{w}_j^q +\tilde{v}_j^q \right) \wedge d \left( \tilde{w}_j^\beta+\tilde{v}_j^\beta \right)\\
     &= t \sum_{\gamma=1}^q w_{0i}^1 \wedge \cdots \wedge R_{ij}^\gamma \wedge \cdots \wedge w_{0i}^q\wedge d w_{0i}^\alpha +t w_{0i}^1\wedge \cdots \wedge w_{0i}^q \wedge dR_{ij}^\alpha\\
     &+ \left( \sum_{\gamma=1}^q(-1)^{q-\gamma} G_{ij}^\gamma w_{0i}^1\wedge \cdots \wedge \hat{w}_{0i}^\gamma \wedge \cdots \wedge w_{0i}^q\wedge  d w_{0i}^\alpha +w_{0i}^1 \wedge \cdots \wedge w_{0i}^q \wedge dG_{ij}^\alpha  - w_{0i}^1 \wedge \cdots \wedge w_{0i}^q \wedge R_{ij}^\alpha\right) \wedge dt  + h_{0ij}\sum_{\beta=1}^q h_{0ij}^{\alpha\beta} \left( tB_j^\beta+ F_j^\beta \wedge dt \right)
\end{align*}}}
This implies that
\begin{align*}
B_i^\alpha-  h_{0ij}\sum_{\beta=1}^q h_{0ij}^{\alpha\beta} B_j^\beta&=\sum_{\gamma=1}^q w_{0i}^1 \wedge \cdots \wedge R_{ij}^\gamma \wedge \cdots \wedge w_{0i}^q \wedge dw_{0i}^\alpha+ w_{0i}^1 \wedge \cdots \wedge w_{0i}^q\wedge dR_{ij}^\gamma\\
F_i^\alpha- h_{0ij}\sum_{\beta=1}^q h_{0ij}^{\alpha\beta} F_j^\beta &= \sum_{\gamma=1}^q (-1)^{q-\gamma} G_{ij}^\gamma w_{0i}^1\wedge \cdots \wedge \hat{w}_{0i}^\gamma \wedge \cdots \wedge w_{0i}^q \wedge dw_{0i}^\alpha + w_{0i}^1\wedge \cdots \wedge w_{0i}^q\wedge dG_{ij}^\alpha-w_{0i}^1\wedge \cdots \wedge w_{0i}^q \wedge R_{ij}^\alpha
\end{align*}
This implies that $-\delta\left( \left\{ B_i\right\}, \left\{ F_i \right\} \right) + G_1\left( \left\{ - R_{ij}\right\}, \left\{ - G_{ij} \right\} \right)=0$. This proves the claim. We show that for another lifting $\tilde{w}_i'^\alpha$ of $w_i^\alpha$ and another lifting of $\tilde{v}_i'^\alpha$ of $v_i^\alpha,\alpha=1,...,q$, the associated $1$-cocycle $\beta:=\left( \left\{ B_i'\right\}, \left\{ F_i' \right\} , \left\{ - R_{ij}' \right\} , \left\{ - G_{ij}' \right\} \right)$ is cohomologous to $\alpha:=\left( \left\{ B_i \right\}, \left\{ F_i \right\}, \left\{ - R_{ij} \right\} , \left\{ - G_{ij} \right\} \right)$. Then we have $\tilde{w}_i'^\alpha= \tilde{w}_i^\alpha+ t V_i^\alpha$ for some $V_i^\alpha\in \Gamma\left(U_i, \Omega_X^1 \right)$ and $\tilde{v}_i'^\alpha= \tilde{v}_i^\alpha + r_i^\alpha dt$ for some $r_i^\alpha \in \Gamma\left( U_i, \mathcal{O}_X \right)$. We define $V_i\in \Gamma\left( U_i, \mathscr{H}om_{\mathcal{O}_X}\left( \mathcal{N}_{\mathcal{F}_0}^*, \Omega_X^1 \right)\right)$ by $V_i:\Gamma\left( U_i, \mathcal{N}_{\mathcal{F}_0}^* \right) \to \Gamma\left( U_i, \Omega_X^1 \right), w_{0i}^\alpha \mapsto V_i^\alpha$, and we define $r_i\in \Gamma\left( U_i , \mathscr{H}om_{\mathcal{O}_X}\left( \mathcal{N}_{\mathcal{F}_0}^*, \mathcal{O}_X \right) \right)$ by $r_i:\Gamma\left( U_i, \mathcal{N}_{\mathcal{F}_0}^* \right) \to \Gamma\left( U_i, \mathcal{O}_X \right), w_{0i}^\alpha \mapsto r_i^\alpha$. Then we see that
{\small{\begin{align}\label{uu21}
&t\left(B_i'^\alpha - B_i^\alpha \right) + \left( F_i'^\alpha - F_i^\alpha \right) \wedge dt \\
&= \left( \tilde{w}_i^1 + \tilde{v}_i^1 + tV_i^1 + r_i^1 dt \right) \wedge \cdots \wedge \left( \tilde{w}_i^q + \tilde{v}_i^q +t V_i^q + r_i^q dt \right) \wedge d\left( \tilde{w}_i^\alpha + \tilde{v}_i^\alpha + t V_i^\alpha + r_i^\alpha dt \right) - \left( \tilde{w}_i^1 + \tilde{v}_i^1 \right) \wedge \cdots \wedge \left( \tilde{w}_i^q + \tilde{v}_i^q \right) \wedge d \left( \tilde{w}_i^\alpha + \tilde{v}_i^\alpha \right)     \notag \\
&=t \sum_{\beta=1}^q w_{0i}^1 \wedge \cdots \wedge V_i^\beta \wedge \cdots \wedge w_{0i}^q \wedge dw_{0i}^\alpha +t w_{0i}^1 \wedge \cdots \wedge w_{0i}^q \wedge d V_i^\alpha \notag \\
& + \left(\sum_{\beta=1}^q (-1)^{q-\beta} r_i^\beta w_{0i}^1 \wedge \cdots \wedge \widehat{w}_{0i}^\beta \wedge \cdots \wedge w_{0i}^q \wedge d w_{0i}^\alpha + w_{0i}^1 \wedge \cdots \wedge w_{0i}^q \wedge d r_i^\alpha - w_{0i}^1 \wedge \cdots \wedge w_{0i}^q \wedge V_i^\alpha \right) \wedge dt    \notag
\end{align}}}
This implies that $G_1'\left( \left\{ V_i \right\}, \left\{ r_i \right\} \right)= \left( \left\{ B_i'-B_i\right\}, \left\{ F_i' - F_i \right\} \right)$. On the other hand, 
\begin{align*}
&t\left(- R_{ij}'^\alpha + R_{ij}^\alpha \right)+ \left( - G_{ij}'^\alpha + G_{ij}^\alpha \right) dt \\
&=  - \tilde{w}_i^\alpha  - t V_i^\alpha - \tilde{v}_i^\alpha - r_i^\alpha dt + \sum_{\beta=1}^q h_{0ij}^{\alpha\beta}\left( \tilde{w}_j^\beta + t V_j^\beta + \tilde{v}_j^\beta +  r_j^\beta dt   \right) + \tilde{w}_i^\alpha + \tilde{v}_i^\alpha - \sum_{\beta=1}^q h_{0ij}^{\alpha\beta} \left( \tilde{w}_j^\beta + \tilde{v}_j^\beta \right)\\
&= t \left( \sum_{\beta=1}^q h_{0ij}^{\alpha\beta} V_j^\beta - V_i^\alpha  \right) + \left(\sum_{\beta=1}^q h_{0ij}^{\alpha\beta} r_j^\beta - r_i^\alpha \right) dt
\end{align*}
This implies that $\delta\left( \left\{ V_i \right\}\right)=\left\{ - R_{ij}' - \left( - R_{ij}\right) \right\}$ and $\delta\left(\left\{ r_i \right\} \right)= \left\{ -G_{ij}' - \left( - G_{ij} \right) \right\}$. This proves $\alpha$ is cohomologous to $\beta$. So given a small extension $e:0\to (t)=(x^m)\to \tilde{A}=k[x]/(x^{m+1})\to A=k[x]/(x^m)\to 0$ and an infinitesimal unfolding $\eta$ over $A$ of $\left(X, \mathcal{N}_{\mathcal{F}_0}^* \right)$ with $X$ fixed, we can associate an element $o_\eta(e):=$ the cohomology class of $\alpha\in \mathbb{H}^1\left(X, \textnormal{Unf}_{ \left( X, \mathcal{N}_{\mathcal{F}_0^*}\right)_{\textnormal{fixed}}}^\bullet \right)$. We note that $o_\eta(e)=0$ if and only if there exist $\left\{ N_i \right\} \in C^0\left( \mathcal{U}, \mathscr{H}om_{\mathcal{O}_X}\left( \mathcal{N}_{\mathcal{F}_0}^*, \Omega_X^1 \right) \right)$ where $N_i :\Gamma\left( U_i, \mathcal{N}_{\mathcal{F}_0}^* \right)\to \Gamma\left( U_i, \Omega_X^1 \right), w_{0i}^\alpha \mapsto N_i^\alpha $ and $ \left\{ M_i \right\} \in C^0\left( \mathcal{U}, \mathscr{H}om_{\mathcal{O}_X}\left( \mathcal{N}_{\mathcal{F}_0}^*, \mathcal{O}_X \right) \right)$ where $M_i:\Gamma\left( U_i, \mathcal{N}_{\mathcal{F}_0}^* \right)\to \Gamma\left( U_i, \mathcal{O}_X \right),w_{0i}^\alpha \mapsto M_i^\alpha$ such that $N_j - N_i=- R_{ij}$ and $M_j- N_i = - G_{ij}$, and $G_1'\left( N_i, M_i \right)= \left( B_i, F_i \right)$, so that
\begin{align*}
&t\left( \sum_{\beta=1}^q h_{0ij}^{\alpha\beta} N_j^\beta - N_i^\alpha \right)+ \left( \sum_{\beta=1}^q h_{0ij}^{\alpha\beta} M_j^\beta - M_i^\alpha \right) \wedge dt = - \tilde{w}_i^\alpha - \tilde{v}_i^\alpha + \sum_{\beta=1}^q h_{0ij}^{\alpha\beta}\left( \tilde{w}_j^\beta + \tilde{v}_j^\beta \right)\\
&\Longrightarrow  \tilde{w}_i^\alpha + \tilde{v}_i^\alpha - t N_i^\alpha - M_i^\alpha \wedge dt = \sum_{\beta=1}^q h_{0ij}^{\alpha\beta}\left( \tilde{w}_j^\beta + \tilde{v}_j^\beta - t N_j^\beta - M_j^\beta \wedge dt  \right)
\end{align*}
and we have
\begin{align}
& t \sum_{\beta=1}^q w_{0i}^1 \wedge \cdots \wedge N_i^\beta \wedge \cdots \wedge w_{0i}^q \wedge d w_{0i}^\alpha + t w_{0i}^1 \wedge \cdots \wedge w_{0i}^q \wedge d N_i^\alpha \notag \\
& + \left( \sum_{\beta=1}^q (-1)^{q-\beta} M_i^\beta w_{0i}^1 \wedge \cdots \wedge \widehat{w}_{0i}^\beta \wedge \cdots \wedge w_{0i}^q \wedge d w_{0i}^\alpha + w_{0i}^1 \wedge \cdots \wedge w_{0i}^q \wedge d M_i^\alpha - w_{0i}^1 \wedge \cdots \wedge w_{0i}^q \wedge N_i^\alpha \right) \wedge dt \notag \\
&= t B_i^\alpha + F_i^\alpha \wedge dt = \left( \tilde{w}_i^1 + \tilde{v}_i^1 \right) \wedge \cdots \wedge \left( \tilde{w}_i^q + \tilde{v}_i^q \right) \wedge d \left(  \tilde{w}_i^\alpha + \tilde{v}_i^\alpha \right) \notag \\
&\Longrightarrow  \left( \tilde{w}_i^1 + \tilde{v}_i^1 - t N_i^1 - M_i^1dt \right) \wedge \cdots \wedge \left( \tilde{w}_i^q + \tilde{v}_i^q - t N_i^q - M_i^q dt \right) \wedge d \left(  \tilde{w}_i^\alpha + \tilde{v}_i^\alpha  - t N_i^\alpha - M_i^\alpha dt \right) = 0. \label{uu22}
\end{align}
Hence a lifting of $\eta$ to $\tilde{A}$ exists. This completes the proof of Proposition \ref{uu11}.

\end{proof}

\subsection{Relationships between $\textnormal{Unf}_{\left\{ X, \mathcal{N}_{\mathcal{F}_0}^* \right\}_{\textnormal{fixed}}}, \textnormal{Unf}_{\left( X, \mathcal{N}_{\mathcal{F}_0}^* \right)_{\textnormal{fixed}}}$ and $\textnormal{FDef}_{\left( X, \mathcal{N}_{\mathcal{F}_0}^* \right)_{\textnormal{fixed}}}$}\

Let $\left( X, \mathcal{N}_{\mathcal{F}_0}^* \right)$ be a foliated nonsingular variety with $\mathcal{N}_{\mathcal{F}_0}^*$ locally free. Then we have the canonical morphism of deformation functors
\begin{align*}
\textnormal{Unf}_{\left\{X, \mathcal{N}_{\mathcal{F}_0}^*\right\}_{\textnormal{fixed}}} \to \textnormal{Unf}_{\left(X, \mathcal{N}_{\mathcal{F}_0}^* \right)_{\textnormal{fixed}}}
\end{align*}
and the morphism of deformation functors from unfoldings to the underlying foliated deformations (Part II)
\begin{align*}
\textnormal{Unf}_{\left(X, \mathcal{N}_{\mathcal{F}_0}^* \right)_{\textnormal{fixed}}}  \to \textnormal{FDef}_{\left( X, \mathcal{N}_{\mathcal{F}_0}^* \right)}
\end{align*}
Then the relationship of those three functors $\textnormal{Unf}_{\left\{X, \mathcal{N}_{\mathcal{F}_0} \right\}_{\textnormal{fixed}}} , \textnormal{Unf}_{\left( X , \mathcal{N}_{\mathcal{F}_0}^* \right)_{\textnormal{fixed}}}$ and $\textnormal{FDef}_{\left(X, \mathcal{N}_{\mathcal{F}_0}^* \right)_{\textnormal{fixed}}}$ is expressed in the following exact sequence of complexes of sheaves
{\Small{\begin{center}
$\begin{CD}
@. \cdots @. \cdots @. \cdots \\
@. @AAA @AAA @AAA \\
0 @>>> \mathscr{H}om_{\mathcal{O}_X}\left( \mathcal{N}_{\mathcal{F}_0}^*, \tilde{\mathcal{S}}^2 \right) @>>> \mathscr{H}om_{\mathcal{O}_X}\left( \mathcal{N}_{\mathcal{F}_0}^*, \tilde{\mathcal{S}}^3  \right) \bigoplus \mathscr{H}om_{\mathcal{O}_X}\left( \mathcal{N}_{\mathcal{F}_0}^*, \tilde{\mathcal{S}}^2 \right) @>>> \mathscr{H}om_{\mathcal{O}_X}\left( \mathcal{N}_{\mathcal{F}_0}^*, \tilde{\mathcal{S}}^3 \right) @>>> 0 \\
@. @AAA @AAA @AAA \\
0 @>>> \mathscr{H}om_{\mathcal{O}_X}\left( \mathcal{N}_{\mathcal{F}_0}^*, \tilde{\mathcal{S}}^1 \right) @>>> \mathscr{H}om_{\mathcal{O}_X}\left( \mathcal{N}_{\mathcal{F}_0}^*, \tilde{\mathcal{S}}^2  \right) \bigoplus \mathscr{H}om_{\mathcal{O}_X}\left( \mathcal{N}_{\mathcal{F}_0}^*, \tilde{\mathcal{S}}^1 \right) @>>> \mathscr{H}om_{\mathcal{O}_X}\left( \mathcal{N}_{\mathcal{F}_0}^*, \tilde{\mathcal{S}}^2 \right) @>>> 0 \\
@. @AAA @AAA @AAA \\
0 @>>> \mathscr{H}om_{\mathcal{O}_X}\left( \mathcal{N}_{\mathcal{F}_0}^*, \mathcal{O}_X \right) @>>> \mathscr{H}om_{\mathcal{O}_X}\left( \mathcal{N}_{\mathcal{F}_0}^*, \Omega_X^1 \right) \bigoplus \mathscr{H}om_{\mathcal{O}_X}\left( \mathcal{N}_{\mathcal{F}_0}^*, \mathcal{O}_X  \right) @>>> \mathscr{H}om_{\mathcal{O}_X}\left( \mathcal{N}_{\mathcal{F}_0}^*, \Omega_X^1 \right) @>>> 0 
\end{CD}$
\end{center}}}
Then from the long exact sequence of cohomology groups, we have
\begin{align*}
 \mathbb{H}^0\left(X, \mathscr{H}om_{\mathcal{O}_X}\left(\mathcal{N}_{\mathcal{F}_0}^*, \mathcal{O}_X \right)^\bullet \right) \to \mathbb{H}^0\left( X, \textnormal{Unf}_{\left( X, \mathcal{N}_{\mathcal{F}_0}^* \right)_{\textnormal{fixed}}}^\bullet \right) \to H^0 \left(X, \mathscr{H}om_{\mathcal{O}_X}\left(\mathcal{N}_{\mathcal{F}_0}^*, \Omega_X^1 \right)^\bullet \right)
\end{align*}
which gives a morphism of first-order deformations 
\begin{align*}
\textnormal{Unf}_{\left\{ X, \mathcal{N}_{\mathcal{F}_0}^* \right\}_{\textnormal{fixed}}}\left(k[\epsilon] \right)\to \textnormal{Unf}_{ \left(X, \mathcal{N}_{\mathcal{F}_0}^* \right)_{\textnormal{fixed}}} \left( k[\epsilon] \right) \to \textnormal{FDef}_{\left(X, \mathcal{N}_{\mathcal{F}_0}^*\right)_{\textnormal{fixed}}} \left( k[\epsilon ] \right)
\end{align*}

and the long exact sequence also gives an exact sequence 
\begin{align*}
 \mathbb{H}^1\left(X, \mathscr{H}om_{\mathcal{O}_X}\left(\mathcal{N}_{\mathcal{F}_0}^*, \mathcal{O}_X \right)^\bullet \right) \to \mathbb{H}^1 \left( X, \textnormal{Unf}_{\left( X, \mathcal{N}_{\mathcal{F}_0}^* \right)_{\textnormal{fixed}}}^\bullet \right) \to H^1 \left(X, \mathscr{H}om_{\mathcal{O}_X}\left(\mathcal{N}_{\mathcal{F}_0}^*, \Omega_X^1 \right)^\bullet \right)
\end{align*}
which gives obstruction maps of morphisms of deformation functors $\textnormal{Unf}_{\left\{ X, \mathcal{N}_{\mathcal{F}_0}^* \right\}_{\textnormal{fixed}}} \to \textnormal{Unf}_{\left( X, \mathcal{N}_{\mathcal{F}_0}^* \right)_{\textnormal{fixed}}} \to \textnormal{FDef}_{\left( X, \mathcal{N}_{\mathcal{F}_0}^* \right)_{\textnormal{fixed}}}$ for small extensions of the form $0\to \left(x^m\right) \to k[x]/\left(x^{m+1}\right)\to k[x]/\left(x^m\right) \to 0$ in $\bold{Art}$.

\section{Unfolidngs of $\left(X, \mathcal{N}_{\mathcal{F}_0}^* \right)$ with $\mathcal{N}_{\mathcal{F}_0}^*$ locally free under the trivial deformation of $X$}

Let $\left(X, \mathcal{N}_{\mathcal{F}_0}^* \right)$ be a foliated nonsingular variety over $k$. In the previous section we study flat unfoldings of $\left(X, \mathcal{N}_{\mathcal{F}_0}^* \right)$ such that the underling flat deformations of $\left(X,\mathcal{N}_{\mathcal{F}_0}^* \right)$ is trivial. In this section, in addition, we vary $\mathcal{N}_{\mathcal{F}_0}^*$ with $X$ fixed.

\begin{definition}
Let $\left( X, \mathcal{N}_{\mathcal{F}_0}^* \right)$ be a foliated nonsingular variety over $k$ with $\mathcal{N}_{\mathcal{F}_0}^*$ locally free, so that we have an exact sequence $0\to \mathcal{N}_{\mathcal{F}_0}^* \to \Omega_X^1 \to \frac{\Omega_X^1}{\mathcal{N}_{\mathcal{F}_0}^*}\to 0$.  An infinitesimal unfoldings of $\left(X, \mathcal{N}_{\mathcal{F}_0}^* \right)$ under
 the trivial deformation of $X$ over $A\in \bold{Art}$ is a  locally free subsheaf  $\tilde{\mathcal{N}}_{\mathcal{F}}^*$ of $\Omega_{\mathcal{X}/k}^1=\Omega_{X\times \textnormal{Spec}(A)/k}^1$ such that $ d \left( \tilde{\mathcal{N}}_{\mathcal{F},x}^* \right) \subset \tilde{\mathcal{N}}_{\mathcal{F},x}^*\bigwedge \Omega_{\mathcal{X}/k,x}^1 $ for $x\in X- \textnormal{Sing}\left( \mathcal{F}_0 \right)$ and $\tilde{\mathcal{N}}_{\mathcal{F}}^*\cong \alpha\left(\tilde{\mathcal{N}}_{\mathcal{F}}^* \right) \cong \mathcal{N}_\mathcal{F}^*$ under $\alpha:\Omega_{\mathcal{X}/k}^1\to \Omega_{\mathcal{X}/A}^1$, and $\left( \mathcal{X}, \mathcal{N}_{\mathcal{F}}^* \right)$ defines a flat foliated deformation of $\left(X, \mathcal{N}_{\mathcal{F}_0}^* \right)$ over $A$ with $X$ fixed, i.e. $\Omega_{\mathcal{X}/A}^1/ \mathcal{N}_{\mathcal{F}}^*$  is flat over $\textnormal{Spec}(A)$, and the exact sequence on $\mathcal{X}:=X\times_{\textnormal{Spec}(k)} \textnormal{Spec}(A)$ 
\begin{align}
0\to \mathcal{N}_{\mathcal{F}}^* \to \Omega_{\mathcal{X}/k}^1  \to \frac{\Omega_{\mathcal{X}/k}^1}{\mathcal{N}_{\mathcal{F}}^*}\to 0
\end{align}
induces $0\to \mathcal{N}_{\mathcal{F}_0}^* \to \Omega_X^1 \to \frac{\Omega_X^1}{\mathcal{N}_{\mathcal{F}_0}^*}\to 0$ over $\textnormal{Spec}(k)$.
Then we have a functor of Artin rings
\begin{align*}
\textnormal{Unf}_{\mathcal{N}_{\mathcal{F}_0}^*}: \bold{Art} &\to \textnormal{(Sets)} \\
                             A &\to \left\{ \textnormal{Infinitesimal unfoldings of $\left( X, \mathcal{N}_{\mathcal{F}_0}^*\right)$ with $X$ fixed}   \right\}
\end{align*}
\end{definition}

We will show that the following complex of sheaves truncated from the unfolded dual leaf complex
{\small{\begin{align*}
\textnormal{Unf}_{\mathcal{N}_{\mathcal{F}_0}^*}^\bullet :   \mathscr{H}om_{\mathcal{O}_X}\left( \mathcal{N}_{\mathcal{F}_0}^* ,   \frac{\Omega_X^1}{\mathcal{N}_{\mathcal{F}_0}^*}  \right) \bigoplus \mathscr{H}om_{\mathcal{O}_X}\left(  \mathcal{N}_{\mathcal{F}_0}^*, \mathcal{O}_X \right) \xrightarrow{G_1} \mathscr{H}om_{\mathcal{O}_X}\left( \mathcal{N}_{\mathcal{F}_0}^* , \tilde{\mathcal{S}}^2 \right) \bigoplus \mathscr{H}om_{\mathcal{O}_X}\left( \mathcal{N}_{\mathcal{F}_0}^*, \tilde{\mathcal{S}}^1 \right) \\
\xrightarrow{G_2}  \mathscr{H}om_{\mathcal{O}_X}\left( \mathcal{N}_{\mathcal{F}_0}^* , \tilde{\mathcal{S}}^2 \right) \bigoplus \mathscr{H}om_{\mathcal{O}_X}\left( \mathcal{N}_{\mathcal{F}_0}^*, \tilde{\mathcal{S}}^1 \right)  \xrightarrow{G_3} \cdots
\end{align*}}}
controls unfoldings of $\left( X, \mathcal{N}_{\mathcal{F}_0}^* \right)$ uner the trivial deformation of $X$. We will denote the $i$-th cohomology group of $\textnormal{Unf}_{\mathcal{N}_{\mathcal{F}_0}^*}^\bullet$ by $\mathbb{H}^i\left( X, \textnormal{Unf}_{\mathcal{N}_{\mathcal{F}_0}^*}^\bullet \right)$.

\begin{proposition}\label{uu10}
Let $\left(X, \mathcal{N}_{\mathcal{F}_0}^* \right)$ be a foliated nonsingular variety with $\mathcal{N}_{\mathcal{F}_0}^*$ locally free.
\begin{enumerate}
\item there is a $1-1$ correspondence
\begin{align*}
\textnormal{Unf}_{\mathcal{N}_{\mathcal{F}_0}^*}\left(k[\epsilon] \right)\cong \mathbb{H}^0\left(X, \textnormal{Unf}_{\mathcal{N}_{\mathcal{F}_0}^*}^\bullet  \right)
\end{align*}
\item Let $\eta=\left( \mathcal{X}, \mathcal{N}_\mathcal{F}^* \right)$ be an infinitesimal unfoldings over $A:=k[x]/ \left(x^m\right)$ of $\left( X, \mathcal{N}_{\mathcal{F}_0}^* \right)$ with $X$ fixed. To small extension $0\to (t)=\left( t^m\right)\to \tilde{A}:=k[x]/ \left( x^{m+1} \right) \to A : = k[x]/ \left( x^m \right) \to 0$, we can associate an element $o_\eta(e)\in \mathbb{H}^1\left( X, \textnormal{Unf}_{\mathcal{N}_{\mathcal{F}_0}^*}^\bullet \right)$ called the obstruction lifting of $\eta$ to $\tilde{A}$, which is $0$ if and only if a lifting of $\eta$ to $\tilde{A}$ exists.
\end{enumerate}
\end{proposition}

\begin{proof}
Given a first-order unfolding $\left( X\times \textnormal{Spec}\left(k[\epsilon] \right) , \tilde{\mathcal{N}}_{\mathcal{F}}^*\right)$ under the trivial deformation of $X$, we choose an affine open covering $\mathcal{U}=\left\{ U_i \right\}$ of $X$ such that $\Gamma\left(  U_i, \mathcal{N}_{\mathcal{F}_0}^* \right)$ is generated by $w_{0i}^1,..., w_{0i}^q \in \Gamma\left( U_i, \Omega_X^1 \right)$ with $w_{0i}^\alpha= \sum_{\beta=1}^q h_{0ij}^{\alpha\beta} w_{0j}^\beta$ for some $h_{0ij}^{\alpha\beta}\in \Gamma\left(U_{ij}, \mathcal{O}_X \right)$, and $\Gamma\left(  U_i, \tilde{\mathcal{N}}_{\mathcal{F}}^* \right)$ is generated by $W_i^\alpha= w_{0i}^\alpha + \epsilon  w_i'^\alpha + v_i^\alpha d\epsilon$ for some $w_i'^\alpha \in \Gamma\left( U_i, \Omega_X^1 \right)$ and $v_i^\alpha \in \Gamma\left( U_i, \mathcal{O}_X \right),\alpha=1,...,q$ with $w_{0i}^\alpha + \epsilon w_i'^\alpha +  v_i^\alpha d \epsilon = \sum_{\beta=1}^q \left( h_{0ij}^{\alpha\beta} + \epsilon g_{ij}^{\alpha\beta} \right)\left( w_{0j}^\alpha + \epsilon w_j'^\beta + v_j^\beta d \epsilon \right)$ for some $g_{ij}^{\alpha\beta}\in \Gamma\left( U_{ij}, \mathcal{O}_X \right)$ and $\left( w_{0i}^1 + \epsilon w_i'^1 + v_i^1 d\epsilon \right)\wedge \cdots \wedge \left( w_{0i}^q + \epsilon w_i'^q + \epsilon v_i^q d\epsilon \right) \wedge d \left( w_{0i}^\alpha + \epsilon w_i'^\alpha + v_i^\alpha d \epsilon \right)=0,\alpha=1,...,q$. Then we have
\begin{align}
&w_i'^\alpha= \sum_{\beta=1}^q h_{0ij}^{\alpha\beta} w_j'^\beta + \sum_{\beta=1}^q g_{ij}^{\alpha\beta} w_{0j}^\beta \label{uu6}\\
&v_i^\alpha = \sum_{\beta=1}^q h_{0ij}^{\alpha\beta} v_j^\beta \label{uu9} \\
\sum_{\beta=1}^q w_{0i}^1 \wedge \cdots \wedge w_i'^\beta &\wedge \cdots \wedge w_{0i}^q \wedge dw_{0i}^\alpha+ w_{0i}^1 \wedge \cdots \wedge w_{0i}^q \wedge dw_i'^\alpha =0 \label{uu7}\\
\sum_{\beta=1}^q (-1)^{q-\beta} v_i^\beta w_{0i}^1 \wedge \cdots \wedge \widehat{w}_{0i}^\beta \wedge \cdots \wedge w_{0i}^q \wedge & d w_{0i}^\alpha\wedge  d \epsilon  + w_{0i}^1 \wedge \cdots \wedge w_{0i}^q \wedge d v_i^\alpha \wedge d \epsilon - w_{0i}^1 \wedge \cdots \wedge w_{0i}^q \wedge w_i'^\alpha \wedge d\epsilon=0 \label{uu8}
\end{align}
We define $\overline{w_i'} \in \Gamma\left(U_i, \mathscr{H}om_{\mathcal{O}_X}\left( \mathcal{N}_{\mathcal{F}_0}^*, \frac{\Omega_X^1}{\mathcal{N}_{\mathcal{F}_0}^*} \right) \right)$ in the following way:
\begin{align}\label{uu25}
\overline{w_i'}:\Gamma\left( U_i, \mathcal{N}_{\mathcal{F}_0}^* \right) &\to \Gamma\left(  U_i, \frac{\Omega_X^1}{\mathcal{N}_{\mathcal{F}_0}^*}\right)\\
 w_{0i}^\alpha &\mapsto \overline{w_i'^\alpha} \notag
\end{align}
and linearly extends to $\Gamma\left( U_i , \mathcal{N}_{\mathcal{F}_0}^* \right)$. Here $\overline{w_i'^\alpha}$ is the image of $w_i'^\alpha$ in $\Gamma\left( U_i, \frac{\Omega_X^1}{\mathcal{N}_{\mathcal{F}_0}^*} \right)$. On the other hand, we define $v_i \in \Gamma\left( U_i, \mathscr{H}om_{\mathcal{O}_X}\left( \mathcal{N}_{\mathcal{F}_0}^*  , \mathcal{O}_X \right) \right)$ by
\begin{align}\label{uu26}
v_i:\Gamma\left( U_i, \mathcal{N}_{\mathcal{F}_0}^* \right) &\to \Gamma\left( U_i, \mathcal{O}_X  \right)\\
 w_{0i}^\alpha &\mapsto v_i^\alpha \notag
\end{align}
and linearly extends to $\Gamma\left( U_i, \mathcal{N}_{\mathcal{F}_0}^* \right)$. Then $\textnormal{(\ref{uu6})}, \textnormal{(\ref{uu9})}, \textnormal{(\ref{uu7})}$ and $\textnormal{(\ref{uu8})}$ imply that $\left( \overline{w'},v \right):= \left( \left\{ \overline{w_i'} \right\} , \left\{ v_i \right\} \right)\in \mathbb{H}^0\left( X, \textnormal{Unf}_{\mathcal{N}_{\mathcal{F}_0}}^\bullet  \right)$.

We define an inverse map. Given a $\left( \overline{w'}, v \right)\in \mathbb{H}^0\left( X, \textnormal{Unf}_{\mathcal{N}_{\mathcal{F}_0}^*}^\bullet \right)$, we have $\textnormal{(\ref{uu6})}, \textnormal{(\ref{uu9})}, \textnormal{(\ref{uu7})}$ and $\textnormal{(\ref{uu8})}$. By reversing the above process, we obtain a first-order unfolding of $\left( X , \mathcal{N}_{\mathcal{F}_0}^* \right)$ whose underlying flat deformation of $\left( X, \mathcal{N}_{\mathcal{F}_0}^* \right)$ is trivial.

Next we identify obstructions. Given a small extension $0\to (t)=(x^m)\to \tilde{A}=k[x]/(x^{m+1})\to A= k[x]/(x^m)\to 0$ and an infinitesimal unfolding $\eta:=\left(X\times \textnormal{Spec}(A), \tilde{\mathcal{N}}_{\mathcal{F}}^* \right)$ over $A$ of $\left( X ,  \mathcal{N}_{\mathcal{F}_0}^* \right)$ with $X$ fixed, we choose an affine open covering $\mathcal{U}=\left\{ U_i \right\}$ of $X$ such that $\Gamma\left(U_i, \tilde{\mathcal{N}}_{\mathcal{F}}^* \right)$ is generated by $ W_i^\alpha= w_i^\alpha + v_i^\alpha \in \left( \Gamma\left( U_i, \Omega_X^1 \right)\otimes A \right) \bigoplus \left( \Gamma\left(U_i, \mathcal{O}_X \right)\otimes \Omega_{A/k}^1 \right) , \alpha=1,...,q$, which induces $w_{0i}^\alpha$, and $w_i^\alpha +  v_i^\alpha = \sum_{\beta=1}^q h_{ij}^{\alpha\beta}\left( w_j^\beta + v_j^\beta \right)$ for some $h_{ij}^{\alpha\beta}\in \Gamma\left( U_{ij}, \mathcal{O}_X \right)\otimes A$ which induces $h_{0ij}^{\alpha\beta}$ and $\left( w_i^1 + v_i^1 \right) \wedge \cdots \wedge \left( w_i^q + v_i^q \right) \wedge d \left( w_i^\alpha + v_i^\alpha \right)=0,\alpha=1,...,q$. In order to see if there exits a lifting of $\left(X\times \textnormal{Spec}\left( \tilde{A} \right), \tilde{\mathcal{N}}_{\mathcal{F}}^* \right)$ to $\tilde{A}$, we choose an arbitrary lifting $\tilde{w}_i^\alpha \in \Gamma\left( U_i, \Omega_X^1 \right)\otimes \tilde{A}$ of $w_i^\alpha$ and arbitrary lifting of $\tilde{v}_i^\alpha\in \Gamma\left( \mathcal{U}_i, \mathcal{O}_X \right)\otimes \Omega_{\tilde{A}/k}^1$ and arbitrary lifting of $\tilde{h}_{ij}^{\alpha\beta}\in \Gamma\left( U_{ij}, \mathcal{O}_X \right) \otimes \tilde{A}$ of $h_{ij}^{\alpha\beta}$. Then we have
\begin{align*}
\left(\tilde{w}_i^\alpha + \tilde{v}_i^\alpha \right) - \sum_{\beta=1}^q \tilde{h}_{ij} \left(\tilde{w}_j^\beta + \tilde{v}_j^\beta \right)= tR_{ij}^{\alpha}+ G_{ij}^\alpha & dt\,\,\,\,\,\textnormal{for some}\,\,\,R_{ij}^\alpha \in \Gamma\left( U_{ij}, \Omega_X^1 \right),\,\,\,\, G_{ij}^\alpha \in \Gamma\left( U_{ij}, \mathcal{O}_X \right)\\
\left(\tilde{w}_i^1 + \tilde{v}_i^\alpha \right)\wedge \cdots \wedge \left(\tilde{w}_i^q+ \tilde{v}_i^q\right) \wedge d \left( \tilde{w}_i^\alpha  + \tilde{v}_i^\alpha \right)  = t B_i^\alpha+ &F_i^\alpha \wedge dt\,\,\,\,\,\textnormal{for some}\,\,\, B_i^\alpha \in \Gamma\left( U_i, \tilde{\mathcal{S}}^2 \right),\,\,\, F_i^\alpha \in \Gamma\left( U_i, \tilde{\mathcal{S}}^1 \right) \\
\tilde{h}_{ik}^{\alpha \gamma}- \sum_{\beta=1}^q \tilde{h}_{ij}^{\alpha\beta} \tilde{h}_{jk}^{\beta \gamma}= t G_{ijk}^{\alpha \gamma}\,\,\,\,\,\,&\textnormal{for some}\,\,\, G_{ijk}^{\alpha \gamma} \in \Gamma\left( U_{ijk}, \mathcal{O}_X \right)
\end{align*}
Here we note that we can show $B_i^\alpha \in \Gamma\left( U_i, \tilde{\mathcal{S}}^2\right)$, and $F_i^\alpha \in \Gamma\left( U_i, \tilde{\mathcal{S}}^1 \right)$ as in the proof of Proposition \ref{uu11}. Then we define $\overline{R_{ij}}\in \Gamma\left( U_{ij}, \mathscr{H}om_{\mathcal{O}_X}\left( \mathcal{N}_{\mathcal{F}_0}^*, \frac{\Omega_X^1}{\mathcal{N}_{\mathcal{F}_0}^*} \right) \right)$ in the following way:
\begin{align}\label{uu29}
\overline{R_{ij}}:\Gamma\left( U_{ij}, \mathcal{N}_{\mathcal{F}_0}^* \right) &\to \Gamma\left( U_{ij}, \frac{\Omega_X^1}{\mathcal{N}_{\mathcal{F}_0}^*} \right)\\
w_{0i}^\alpha &\mapsto \overline{ R_{ij}^\alpha } \notag
\end{align}
and linearly extends to $\Gamma\left(U_{ij}, \mathcal{N}_{\mathcal{F}_0}^* \right)$. Here $\overline{R_{ij}^\alpha }$ is the image of $R_{ij}^\alpha$ in $\Gamma\left( U_{ij}, \frac{\Omega_X^1}{\mathcal{N}_{\mathcal{F}_0}^*} \right)$, and we define $G_{ij}\in \Gamma\left( U_{ij}, \mathscr{H}om_{\mathcal{O}_X}\left( \mathcal{N}_{\mathcal{F}_0}^*, \mathcal{O}_X \right) \right)$ in the following way:
\begin{align}
G_{ij} : \Gamma\left( U_{ij}, \mathcal{N}_{\mathcal{F}_0}^* \right) &\to \Gamma\left( U_{ij}, \mathcal{O}_X \right) \label{uu30} \\
 w_{0i}^\alpha &\mapsto G_{ij}^\alpha \notag
\end{align}
and linearly extends to $\Gamma\left( U_i, \mathcal{N}_{\mathcal{F}_0}^*  \right)$, and we define $B_i\in \Gamma\left( U_i,  \mathscr{H}om_{\mathcal{O}_X}\left( \mathcal{N}_{\mathcal{F}_0}^*,  \tilde{\mathcal{S}}^2 \right) \right)$ in the following way:
\begin{align}
B_i:\Gamma\left(  U_i, \mathcal{N}_{\mathcal{F}_0}^* \right) &\to  \Gamma\left( U_i, \tilde{\mathcal{S}}^2 \right) \label{uu31} \\
 w_{0i}^\alpha &\mapsto B_i^\alpha \notag
\end{align}
and linearly extends to $\Gamma\left( U_i, \mathcal{N}_{\mathcal{F}_0}^* \right)$, and we define $F_i\in \Gamma\left(U_i, \mathscr{H}om_{\mathcal{O}_X}\left( \mathcal{N}_{\mathcal{F}_0}^* , \tilde{\mathcal{S}}^1 \right) \right)$ in the following way:
\begin{align}
F_i: \Gamma\left( U_i, \mathcal{N}_{\mathcal{F}_0}^* \right) &\to \Gamma\left( U_i, \tilde{\mathcal{S}}^1 \right) \label{uu32}\\
 w_{0i}^\alpha &\mapsto F_i^\alpha \notag
\end{align}
and linearly extends to $\Gamma\left( U_i, \mathcal{N}_{\mathcal{F}_0}^* \right)$. Then we claim that $\left( \left\{ B_i \right\},\left\{ F_i \right\},  \left\{ - \overline{R_{ij}} \right\} , \left\{ - G_{ij} \right\} \right)\in C^0 \left( \mathcal{U}, \mathscr{H}om_{\mathcal{O}_X}\left( \mathcal{N}_{\mathcal{F}_0}^*, \tilde{\mathcal{S}}^2 \right)   \right)$ $ \bigoplus C^0\left( \mathcal{U}, \mathscr{H}om_{\mathcal{O}_X}\left( \mathcal{N}_{\mathcal{F}_0}^*, \tilde{\mathcal{S}}^1 \right) \right)  \bigoplus C^1\left( \mathcal{U}, \mathscr{H}om_{\mathcal{O}_X}\left( \mathcal{N}_{\mathcal{F}_0}^*, \frac{\Omega_X^1}{\mathcal{N}_{\mathcal{F}_0}^*} \right) \right) \bigoplus C^1\left( \mathcal{U}, \mathscr{H}om_{\mathcal{O}_X}\left( \mathcal{N}_{\mathcal{F}_0}^*, \mathcal{O}_X \right) \right)$ defines a $1$-cocycle in the following \v Cech resolution of $\mathscr{H}om_{\mathcal{O}_X}\left( \mathcal{N}_{\mathcal{F}_0}^*, \frac{\Omega_X^1}{\mathcal{N}_{\mathcal{F}_0}^*} \right)^\bullet$:
{\tiny{\begin{center}
$\begin{CD}
\cdots \\
@AG_2AA\\
C^0\left( \mathcal{U}, \mathscr{H}om_{\mathcal{O}_X}\left( \mathcal{N}_{\mathcal{F}_0}^* , \tilde{\mathcal{S}}^2 \right) \right) \bigoplus C^0 \left( \mathcal{U}, \mathscr{H}om_{\mathcal{O}_X}\left( \mathcal{N}_{\mathcal{F}_0}^*, \tilde{\mathcal{S}}^1 \right) \right) @>-\delta>> \cdots \\
@AG_1'AA @AG_1' AA \\
C^0\left(\mathcal{U}, \mathscr{H}om_{\mathcal{O}_X}\left(  \mathcal{N}_{\mathcal{F}_0}^*  , \frac{\Omega_X^1}{\mathcal{N}_{\mathcal{F}_0}^*} \right) \right) \bigoplus C^0\left(\mathcal{U},  \mathscr{H}om_{\mathcal{O}_X}\left(  \mathcal{N}_{\mathcal{F}_0}^*, \mathcal{O}_X  \right) \right) @>\delta>> C^1\left( \mathcal{U}, \mathscr{H}om_{\mathcal{O}_X}\left(  \mathcal{N}_{\mathcal{F}_0}^* , \frac{\Omega_X^1}{\mathcal{N}_{\mathcal{F}_0}^*} \right) \right) \bigoplus C^1\left( \mathcal{U}, \mathscr{H}om_{\mathcal{O}_X}\left( \mathcal{N}_{\mathcal{F}_0}^*, \mathcal{O}_X \right)  \right) @>-\delta>> 
\end{CD}$
\end{center}}}
We note that
\begin{align}
&tR_{ij}^\alpha+ G_{ij}^\alpha dt - \left(tR_{ik}^\alpha+ G_{ik}^\alpha dt \right)+\sum_{\beta=1}^q h_{0ij}^{\alpha\beta} \left(tR_{jk}^\beta + G_{jk}^\beta dt \right) \\
&= \left(\tilde{w}_i^\alpha + \tilde{v}_i^\alpha \right) -\sum_{\beta=1}^q \tilde{h}_{ij}^{\alpha\beta} \left( \tilde{w}_j^\beta + \tilde{v}_j^\beta  \right) - \left(\tilde{w}_i^\alpha + \tilde{v}_i^\alpha \right)+\sum_{\beta=1}^q \tilde{h}_{ik}^{\alpha\beta} \left( \tilde{w}_k^\beta + \tilde{v}_i^\beta \right) +\sum_{\beta=1}^q \tilde{h}_{ij}^{\alpha\beta} \left( \tilde{w}_j^\beta + \tilde{v}_j^\beta \right) -\sum_{\beta,\gamma=1}^q \tilde{h}_{ij}^{\alpha\beta} \tilde{h}_{jk}^{\beta\gamma} \left( \tilde{w}_k^\gamma + \tilde{v}_k^\gamma \right) \notag \\
&=\sum_{\gamma=1}^q \left(\tilde{h}_{ik}^{\alpha\gamma}-\sum_{\beta=1}^q \tilde{h}_{ij}^{\alpha\beta}\tilde{h}_{jk}^{\beta\gamma} \right)\left( \tilde{w}_k^\gamma + \tilde{v}_i^\gamma \right) =t\sum_{\gamma=1}^q  G_{ijk}^{\alpha\gamma} w_{0k}^\gamma \notag
\end{align}
This implies that $-\delta\left( \left\{ - \overline{R_{ij}} \right\} \right) = 0$ and $-\delta\left( \left\{ - G_{ij} \right\} \right)=0$. As in the proof of Proposition \ref{uu11}, we have $G_2\left( \left\{ B_i \right\}, \left\{ F_i \right\} \right)=0$. On the other hand, by setting $h_{0ij}= \det\left( h_{0ij}^{\alpha\beta} \right)$, and $\tilde{h}_{ij}=\det\left( \tilde{h}_{ij}^{\alpha\beta} \right)$, we have
{\tiny{\begin{align}\label{uu34}
&t B_i^\alpha+ F_i^\alpha \wedge dt=\left(\tilde{w}_i^1+ \tilde{v}_i^1\right)\wedge \cdots \wedge \left(\tilde{w}_i^q+ \tilde{v}_i^q\right) \wedge d\left( \tilde{w}_i^\alpha+ \tilde{v}_i^\alpha\right)\\
     &=\left(tR_{ij}^1+G_{ij}^1 dt +\sum_{\beta_1=1}^q \tilde{h}_{ij}^{\alpha\beta_1}\left( \tilde{w}_j^{\beta_1}+ \tilde{v}_j^{\beta_1}\right) \right) \wedge \cdots \wedge \left( t R_{ij}^q+ G_{ij}^q dt +\sum_{\beta_q=1}^q \tilde{h}_{ij}^{\alpha\beta_q} \left(\tilde{w}_j^{\beta_q} + \tilde{v}_j^{\beta_q} \right)\right) \wedge d \left( tR_{ij}^\alpha+ G_{ij}^\alpha dt +\sum_{\beta=1}^q \tilde{h}_{ij}^{\alpha\beta}\left( \tilde{w}_j^\beta  + \tilde{v}_j^\beta\right) \right)\notag\\
          &=\left(tR_{ij}^1+G_{ij}^1 dt +\sum_{\beta_1=1}^q \tilde{h}_{ij}^{\alpha\beta_1}\left( \tilde{w}_j^{\beta_1}+ \tilde{v}_j^{\beta_1}\right) \right) \wedge \cdots \wedge \left( t R_{ij}^q+ G_{ij}^q dt +\sum_{\beta_q=1}^q \tilde{h}_{ij}^{\alpha\beta_q} \left(\tilde{w}_j^{\beta_q} + \tilde{v}_j^{\beta_q} \right)\right) \wedge  \left( dt\wedge R_{ij}^\alpha+tdR_{ij}^\alpha+ dG_{ij}^\alpha \wedge dt +d\left(\sum_{\beta=1}^q \tilde{h}_{ij}^{\alpha\beta}\left( \tilde{w}_j^\beta  + \tilde{v}_j^\beta\right)\right) \right)\notag\\
     &=t \sum_{\gamma=1}^q w_{0i}^1 \wedge \cdots \wedge R_{ij}^\gamma \wedge \cdots \wedge w_{0i}^q\wedge d w_{0i}^\alpha +t w_{0i}^1\wedge \cdots \wedge w_{0i}^q \wedge dR_{ij}^\alpha\notag\\
     &+ \left( \sum_{\gamma=1}^q(-1)^{q-\gamma} G_{ij}^\gamma w_{0i}^1\wedge \cdots \wedge \hat{w}_{0i}^\gamma \wedge \cdots \wedge w_{0i}^q\wedge  d w_{0i}^\alpha + w_{0i}^1 \wedge \cdots \wedge w_{0i}^q \wedge dG_{ij}^\alpha  - w_{0i}^1 \wedge \cdots \wedge w_{0i}^q \wedge R_{ij}^\alpha\right) \wedge dt \notag\\
     & + \tilde{h}_{ij} \sum_{\beta=1}^q \tilde{h}_{ij}^{\alpha\beta} \left(\tilde{w}_j^1+\tilde{v}_j^1\right) \wedge \cdots \wedge \left( \tilde{w}_j^q +\tilde{v}_j^q \right) \wedge d \left( \tilde{w}_j^\beta+\tilde{v}_j^\beta \right)\notag \\
     &= t \sum_{\gamma=1}^q w_{0i}^1 \wedge \cdots \wedge R_{ij}^\gamma \wedge \cdots \wedge w_{0i}^q\wedge d w_{0i}^\alpha +t w_{0i}^1\wedge \cdots \wedge w_{0i}^q \wedge dR_{ij}^\alpha \notag \\
     &+ \left( \sum_{\gamma=1}^q(-1)^{q-\gamma} G_{ij}^\gamma w_{0i}^1\wedge \cdots \wedge \hat{w}_{0i}^\gamma \wedge \cdots \wedge w_{0i}^q\wedge  d w_{0i}^\alpha +w_{0i}^1 \wedge \cdots \wedge w_{0i}^q \wedge dG_{ij}^\alpha  - w_{0i}^1 \wedge \cdots \wedge w_{0i}^q \wedge R_{ij}^\alpha\right) \wedge dt  + h_{0ij}\sum_{\beta=1}^q h_{0ij}^{\alpha\beta} \left( tB_j^\beta+ F_j^\beta \wedge dt \right) \notag
\end{align}}}
This implies that
\begin{align*}
B_i^\alpha-  h_{0ij}\sum_{\beta=1}^q h_{0ij}^{\alpha\beta} B_j^\beta&=\sum_{\gamma=1}^q w_{0i}^1 \wedge \cdots \wedge R_{ij}^\gamma \wedge \cdots \wedge w_{0i}^q \wedge dw_{0i}^\alpha+ w_{0i}^1 \wedge \cdots \wedge w_{0i}^q\wedge dR_{ij}^\gamma\\
F_i^\alpha- h_{0ij}\sum_{\beta=1}^q h_{0ij}^{\alpha\beta} F_j^\beta &= \sum_{\gamma=1}^q (-1)^{q-\gamma} G_{ij}^\gamma w_{0i}^1\wedge \cdots \wedge \hat{w}_{0i}^\gamma \wedge \cdots \wedge w_{0i}^q \wedge dw_{0i}^\alpha + w_{0i}^1\wedge \cdots \wedge w_{0i}^q\wedge dG_{ij}^\alpha-w_{0i}^1\wedge \cdots \wedge w_{0i}^q \wedge R_{ij}^\alpha
\end{align*}
This implies that $-\delta\left( \left\{ B_i\right\}, \left\{ F_i \right\} \right) + G_1\left( \left\{ - R_{ij}\right\}, \left\{ - G_{ij} \right\} \right)=0$. This proves the claim. We show that for another lifting $\tilde{w}_i'^\alpha$ of $w_i^\alpha$, and another lifting $\tilde{v}_i'^\alpha$ of $v_i^\alpha$, and another lifting $\tilde{h}_{ij}'^{\alpha\beta}$ of $h_{ij}^{\alpha\beta},\alpha,\beta=1,...,q$, the associated $1$-cocycle $\beta:=\left( \left\{ B_i'\right\}, \left\{ F_i'\right\} \left\{- \overline{R_{ij}'} \right\}, \left\{ - G_{ij}' \right\} \right)$ is cohomologous to $\alpha:=\left( \left\{ B_i\right\} , \left\{ F_i \right\} , \left\{ - \overline{R_{ij}} \right\}, \left\{ G_{ij}\right\} \right)$. Then $\tilde{w}_i'^\alpha= \tilde{w}_i^\alpha + t V_i^\alpha$ for some $V_i^\alpha \in \Gamma\left( U_i, \Omega_X^1 \right)$ and $\tilde{v}_i'^\alpha = \tilde{v}_i^\alpha + r_i^\alpha dt$ for some $r_i^\alpha \in \Gamma\left( U_i , \mathcal{O}_X \right)$, and $\tilde{h}_{ij}'^{\alpha\beta}= \tilde{h}_{ij}^{\alpha\beta} + t y_{ij}^{\alpha\beta}$ for some $y_{ij}^{\alpha\beta}\in \Gamma\left( U_{ij}, \mathcal{O}_X \right)$. We define $V_i\in \Gamma\left( U_i, \mathscr{H}om_{\mathcal{O}_X}\left( \mathcal{N}_{\mathcal{F}_0}^*, \frac{\Omega_X^1}{\mathcal{N}_{\mathcal{F}_0}^*} \right) \right)$ by $V_i: \Gamma\left( U_i, \mathcal{N}_{\mathcal{F}_0}^* \right) \to \Gamma\left( U_i, \frac{\Omega_X^1}{\mathcal{N}_{\mathcal{F}_0}^*} \right), w_{0i}^\alpha \mapsto \overline{V_i^\alpha}$, and we define $r_i\in \Gamma\left( U_i, \mathscr{H}om_{\mathcal{O}_X}\left( \mathcal{N}_{\mathcal{F}_0}^*, \mathcal{O}_X \right) \right)$ by $r_i:\Gamma\left( U_i, \mathcal{N}_{\mathcal{F}_0}^* \right)\to \Gamma\left( U_i, \mathcal{O}_X \right), w_{0i}^\alpha \mapsto r_i^\alpha$. Then as in $\textnormal{(\ref{uu21})}$, we have $G_1\left( \left\{ \overline{V_i} \right\}, \left\{ r_i \right\} \right)= \left( \left\{ B_i'- B_i \right\}, \left\{ F_i' - F_i \right\} \right)$. On the other hand,
\begin{align*}
&t\left(- R_{ij}'^\alpha + R_{ij}^\alpha \right)+ \left( - G_{ij}'^\alpha + G_{ij}^\alpha \right) dt \\
&=  - \tilde{w}_i^\alpha  - t V_i^\alpha - \tilde{v}_i^\alpha - r_i^\alpha dt + \sum_{\beta=1}^q \left( \tilde{h}_{ij}^{\alpha \beta} +  t y_{ij}^{\alpha\beta} \right)\left( \tilde{w}_j^\beta + t V_j^\beta + \tilde{v}_j^\beta +  r_j^\beta dt   \right) + \tilde{w}_i^\alpha + \tilde{v}_i^\alpha - \sum_{\beta=1}^q \tilde{h}_{ij}^{\alpha\beta} \left( \tilde{w}_j^\beta + \tilde{v}_j^\beta \right)\\
&= t \left( \sum_{\beta=1}^q h_{0ij}^{\alpha\beta} V_j^\beta - V_i^\alpha + \sum_{\beta=1}^q y_{ij}^{\alpha\beta} w_{0j}^\beta \right) + \left(\sum_{\beta=1}^q h_{0ij}^{\alpha\beta} r_j^\beta - r_i^\alpha \right) dt
\end{align*}
This implies that $\delta\left( \left\{ V_i \right\},\right)=\left\{ -\overline{R_{ij}'} - \left( - \overline{R_{ij}} \right) \right\}$ and $\delta\left(\left\{ r_i \right\} \right)= \left\{ -G_{ij}' - \left( - G_{ij} \right) \right\}$. This proves $\alpha$ is cohomologous to $\beta$. So given a small extension $e:0\to (t)=(x^m)\to \tilde{A}=k[x]/(x^{m+1})\to A=k[x]/(x^m)\to 0$ and an infinitesimal unfolding $\eta$ over $A$ of $\left(X, \mathcal{N}_{\mathcal{F}_0}^* \right)$ with $X$ fixed, we can associate an element $o_\eta(e):=$ the cohomology class of $\alpha\in \mathbb{H}^1\left(X, \textnormal{Unf}_{\mathcal{N}_{\mathcal{F}_0^*}}^\bullet \right)$. We note that $o_\eta(e)=0$ if and only if there exist $\left\{ \overline{ N_i} \right\} \in C^0\left( \mathcal{U}, \mathscr{H}om_{\mathcal{O}_X}\left( \mathcal{N}_{\mathcal{F}_0}^*, \frac{\Omega_X^1}{\mathcal{N}_{\mathcal{F}_0}^*} \right) \right)$ where $\overline{N}_i:\Gamma\left( U_i, \mathcal{N}_{\mathcal{F}_0}^* \right)\to \Gamma\left( U_i, \frac{\Omega_X^1}{\mathcal{N}_{\mathcal{F}_0}^*} \right), w_{0i}^\alpha \mapsto \overline{N_i^\alpha}$ and $ \left\{ M_i \right\} \in C^0\left( \mathcal{U}, \mathscr{H}om_{\mathcal{O}_X}\left( \mathcal{N}_{\mathcal{F}_0}^*, \mathcal{O}_X \right) \right)$ where $M_i:\Gamma\left( U_i, \mathcal{N}_{\mathcal{F}_0}^* \right)\to \Gamma\left( U_i, \mathcal{O}_X \right),w_{0i}^\alpha \mapsto M_i^\alpha$ such that $\overline{N_j}- \overline{N_i}=- \overline{R_{ij}}$ and $M_j- N_i = - G_{ij}$, so that
{\small{\begin{align}
&t\left(\sum_{\beta=1}^q h_{0ij}^{\alpha\beta} N_j^\beta - N_i^\alpha + t \sum_{\beta=1}^q u_{ij}^{\alpha\beta} w_{0j}^\beta\right) + \left( \sum_{\beta=1}^q h_{0ij}^{\alpha\beta} M_j^\beta - M_i^\alpha \right) dt= - \tilde{w}_i^\alpha - \tilde{v}_i^\alpha + \sum_{\beta=1}^q \tilde{h}_{ij}^{\alpha\beta} \left( \tilde{w}_j^\beta + \tilde{v}_j^\beta \right) \notag \\
& \Longrightarrow \tilde{w}_i^\alpha - tN_i^\alpha - M_i^\alpha dt = \sum_{\beta=1}^q \left(\tilde{h}_{0ij}^{\alpha\beta}- tu_{ij}^{\alpha\beta} \right)\left(\tilde{w}_j^\beta - t N_j^\beta -  M_j^\beta dt \right)\,\,\,\textnormal{for some} \,\,\, u_{ij}^{\alpha\beta}\in \Gamma\left( U_{ij}, \mathcal{O}_X \right) \end{align}}}
and as in $\textnormal{(\ref{uu22})}$, we have $\left(\tilde{w}_i^1 + \tilde{v}_i^1 - t N_i^1 - M_i^1 dt \right) \wedge \cdots \wedge \left( \tilde{w}_i^q + \tilde{v}_i^q - t N_i^q - M_i^q dt \right) \wedge d \left( \tilde{w}_i^\alpha + \tilde{v}_i^\alpha - t N_i^\alpha - M_i^\alpha dt \right)=0,\alpha=1,...,q$. Hence a lifting of $\eta$ to $\tilde{A}= k[x]/(x^{m+1})$ exists. This completes the proof of Proposition \ref{uu10}.

\end{proof}

\subsection{Relationships between $\textnormal{Unf}_{\left\{ X, \mathcal{N}_{\mathcal{F}_0}^* \right\}_{\textnormal{fixed}}}, \textnormal{Unf}_{\mathcal{N}_{\mathcal{F}_0}^*}$ and $\textnormal{FDef}_{\mathcal{N}_{\mathcal{F}_0}^*}$ }\

Let $\left( X, \mathcal{N}_{\mathcal{F}_0}^* \right)$ be a foliated nonsingular variety with $\mathcal{N}_{\mathcal{F}_0}^*$ locally free. Then we have the canonical morphism of deformation functors
\begin{align*}
\textnormal{Unf}_{\left\{ X , \mathcal{N}_{\mathcal{F}_0}^* \right\}_{\textnormal{fixed}}}\to \textnormal{Unf}_{ \mathcal{N}_{\mathcal{F}_0}^*}
\end{align*}
and the morphism of deformation functors from unfoldings to the underlying foliated deformations (Part II)
\begin{align*}
\textnormal{Unf}_{\mathcal{N}_{\mathcal{F}_0}^*} \to \textnormal{FDef}_{\mathcal{N}_{\mathcal{F}_0}^*}
\end{align*}
Then the relationship of those three functors $\textnormal{Unf}_{\left\{X, \mathcal{N}_{\mathcal{F}_0}^* \right\}_{\textnormal{fixed}}}, \textnormal{Unf}_{\mathcal{N}_{\mathcal{F}_0}^*}$ and $\textnormal{FDef}_{\mathcal{N}_{\mathcal{F}_0}^*}$ is expressed in the following exact sequence of complexes of sheaves
{\Small{\begin{center}
$\begin{CD}
@. \cdots @. \cdots @. \cdots \\
@. @AAA @AAA @AAA \\
0 @>>> \mathscr{H}om_{\mathcal{O}_X}\left( \mathcal{N}_{\mathcal{F}_0}^*, \tilde{\mathcal{S}}^2 \right) @>>> \mathscr{H}om_{\mathcal{O}_X}\left( \mathcal{N}_{\mathcal{F}_0}^*, \tilde{\mathcal{S}}^3  \right) \bigoplus \mathscr{H}om_{\mathcal{O}_X}\left( \mathcal{N}_{\mathcal{F}_0}^*, \tilde{\mathcal{S}}^2 \right) @>>> \mathscr{H}om_{\mathcal{O}_X}\left( \mathcal{N}_{\mathcal{F}_0}^*, \tilde{\mathcal{S}}^3 \right) @>>> 0 \\
@. @AAA @AAA @AAA \\
0 @>>> \mathscr{H}om_{\mathcal{O}_X}\left( \mathcal{N}_{\mathcal{F}_0}^*, \tilde{\mathcal{S}}^1 \right) @>>> \mathscr{H}om_{\mathcal{O}_X}\left( \mathcal{N}_{\mathcal{F}_0}^*, \tilde{\mathcal{S}}^2  \right) \bigoplus \mathscr{H}om_{\mathcal{O}_X}\left( \mathcal{N}_{\mathcal{F}_0}^*, \tilde{\mathcal{S}}^1 \right) @>>> \mathscr{H}om_{\mathcal{O}_X}\left( \mathcal{N}_{\mathcal{F}_0}^*, \tilde{\mathcal{S}}^2 \right) @>>> 0 \\
@. @AAA @AAA @AAA \\
0 @>>> \mathscr{H}om_{\mathcal{O}_X}\left( \mathcal{N}_{\mathcal{F}_0}^*, \mathcal{O}_X \right) @>>> \mathscr{H}om_{\mathcal{O}_X}\left( \mathcal{N}_{\mathcal{F}_0}^*, \frac{\Omega_X^1}{\mathcal{N}_{\mathcal{F}_0}^*} \right) \bigoplus \mathscr{H}om_{\mathcal{O}_X}\left( \mathcal{N}_{\mathcal{F}_0}^*, \mathcal{O}_X  \right) @>>> \mathscr{H}om_{\mathcal{O}_X}\left( \mathcal{N}_{\mathcal{F}_0}^*, \frac{\Omega_X^1}{\mathcal{N}_{\mathcal{F}_0}^*} \right) @>>> 0 
\end{CD}$
\end{center}}}
Then from the long exact sequence of cohomology groups, we have
\begin{align*}
 \mathbb{H}^0\left(X, \mathscr{H}om_{\mathcal{O}_X}\left(\mathcal{N}_{\mathcal{F}_0}^*, \mathcal{O}_X \right)^\bullet \right) \to \mathbb{H}^0\left( X, \textnormal{Unf}_{ \mathcal{N}_{\mathcal{F}_0}^* }^\bullet \right) \to H^0 \left(X, \mathscr{H}om_{\mathcal{O}_X}\left(\mathcal{N}_{\mathcal{F}_0}^*, \frac{\Omega_X^1}{\mathcal{N}_{\mathcal{F}_0}^*} \right)^\bullet \right)
\end{align*}
which gives a morphism of first-order deformations 
\begin{align*}
\textnormal{Unf}_{\left\{ X, \mathcal{N}_{\mathcal{F}_0}^* \right\}_{\textnormal{fixed}}}\left(k[\epsilon] \right)\to \textnormal{Unf}_{  \mathcal{N}_{\mathcal{F}_0}^*} \left( k[\epsilon] \right) \to \textnormal{FDef}_{ \mathcal{N}_{\mathcal{F}_0}^* } \left( k[\epsilon ] \right)
\end{align*}

and the long exact sequence also gives an exact sequence 
\begin{align*}
 \mathbb{H}^1\left(X, \mathscr{H}om_{\mathcal{O}_X}\left(\mathcal{N}_{\mathcal{F}_0}^*, \mathcal{O}_X \right)^\bullet \right) \to \mathbb{H}^1 \left( X, \textnormal{Unf}_{  \mathcal{N}_{\mathcal{F}_0}^* }^\bullet \right) \to H^1 \left(X, \mathscr{H}om_{\mathcal{O}_X}\left(\mathcal{N}_{\mathcal{F}_0}^*, \frac{\Omega_X^1}{\mathcal{N}_{\mathcal{F}_0}^*} \right)^\bullet \right)
\end{align*}
which gives obstruction maps of morphisms of deformation functors $\textnormal{Unf}_{\left\{ X, \mathcal{N}_{\mathcal{F}_0}^* \right\}_{\textnormal{fixed}}} \to \textnormal{Unf}_{  \mathcal{N}_{\mathcal{F}_0}^* } \to \textnormal{FDef}_{ \mathcal{N}_{\mathcal{F}_0}^* }$ for small extensions of the form $0\to (x^m)\to k[x]/(x^{m+1}) \to k[x]/(x^m)  \to 0$ in $\bold{Art}$. 

\subsection{Unfolidngs of $\left(X, \mathcal{N}_{\mathcal{F}_0}^* \right)$ with $\mathcal{N}_{\mathcal{F}_0}^*$ locally free under the trivial deformation of $X$ in terms of vector bundles}\

Let $\left( X, \mathcal{N}_{\mathcal{F}_0}^* \right)$ be a foliated nonsingular variety with $\mathcal{N}_{\mathcal{F}_0}^*$ locally free. The underlying deformations of locally free sheaf $\mathcal{N}_{\mathcal{F}_0}^*$ is implicit in the complex $\textnormal{Unf}_{\mathcal{N}_{\mathcal{F}_0}^*}^\bullet$. In order to make the underlying deformations of $\mathcal{N}_{\mathcal{F}_0}^*$ explicit, we define another complex of sheaves also controlling flat unfoldings of $\left(X, \mathcal{N}_{\mathcal{F}_0}^* \right)$ under the trivial deformation of $X$ and present unfoldings of $\left(X, \mathcal{N}_{\mathcal{F}_0}^* \right)$ with $X$ fixed by using the complex. As is known, deformations of a locally free sheaf $\mathcal{N}_{\mathcal{F}_0}^*$ on $X$ is controlled by $\mathscr{H}om_{\mathcal{O}_X}\left( \mathcal{N}_{\mathcal{F}_0}^* , \mathcal{N}_{\mathcal{F}_0}^* \right)$ (see \cite{Har10}). We note that we have an injection $j:\mathscr{H}om_{\mathcal{O}_X}\left( \mathcal{N}_{\mathcal{F}_0}^* , \mathcal{N}_{\mathcal{F}_0}^* \right) \to \mathscr{H}om_{\mathcal{O}_X}\left( \mathcal{N}_{\mathcal{F}_0}^* , \Omega_X^1 \right)$ induced by $\mathcal{N}_{\mathcal{F}_0}^* \hookrightarrow \Omega_X^1$. Then let us consider an exact sequence of complexes of sheaves
{\small{\begin{center}
$\begin{CD}
@. \cdots @. \cdots @. \cdots \\
@. @AAA @AG_3AA @AG_3AA \\
0@>>>0 @>>>  \left( \mathcal{N}_{\mathcal{F}_0}^*\right)^*\otimes \tilde{\mathcal{S}}^3 \bigoplus \left( \mathcal{N}_{\mathcal{F}_0}^* \right)^* \otimes \tilde{\mathcal{S}}^2 @>>>  \left( \mathcal{N}_{\mathcal{F}_0}^*\right)^*\otimes \tilde{\mathcal{S}}^3 \bigoplus \left( \mathcal{N}_{\mathcal{F}_0}^* \right)^* \otimes \tilde{\mathcal{S}}^2 @>>> 0 \\
@.@AAA @AG_2AA @AG_2AA \\
0@>>>0 @>>>  \left( \mathcal{N}_{\mathcal{F}_0}^* \right)^* \otimes \tilde{\mathcal{S}}^2  \bigoplus \left( \mathcal{N}_{\mathcal{F}_0}^* \right)^*\otimes \tilde{\mathcal{S}}^1 @>>>  \left( \mathcal{N}_{\mathcal{F}_0}^* \right)^* \otimes \tilde{\mathcal{S}}^2  \bigoplus \left( \mathcal{N}_{\mathcal{F}_0}^* \right)^*\otimes \tilde{\mathcal{S}}^1 @>>> 0\\
@. @AAA @AG_1''AA @AG_1AA \\
0@>>> \mathscr{H}om_{\mathcal{O}_X}\left( \mathcal{N}_{\mathcal{F}_0}^*, \mathcal{N}_{\mathcal{F}_0}^* \right) @>>> \mathscr{H}om_{\mathcal{O}_X}\left( \mathcal{N}_{\mathcal{F}_0}^* , \Omega_X^1 \right) \bigoplus \left( \mathcal{N}_{\mathcal{F}_0}^*\right)^* @>>> \mathscr{H}om_{\mathcal{O}_X}\left( \mathcal{N}_{\mathcal{F}_0}^* ,\frac{\Omega_X^1}{\mathcal{N}_{\mathcal{F}_0}^*} \right) \bigoplus \left( \mathcal{N}_{\mathcal{F}_0}^*\right)^*  @>>> 0 \\
@. @A\textnormal{Id}AA @AjAA @AAA \\
0@>>> \mathscr{H}om_{\mathcal{O}_X}\left(\mathcal{N}_{\mathcal{F}_0}^*, \mathcal{N}_{\mathcal{F}_0}^* \right) @>>> \mathscr{H}om_{\mathcal{O}_X}\left( \mathcal{N}_{\mathcal{F}_0}^*, \mathcal{N}_{\mathcal{F}_0}^* \right) @>>> 0@>>> 0
\end{CD}$
\end{center}}}
Let us denote the middle complex of sheaves by $\mathscr{H}om_{\mathcal{O}_X}\left(\mathcal{N}_{\mathcal{F}_0}^*, \mathcal{N}_{\mathcal{F}_0}^* \right)^\bullet$:
{\small{\begin{align}\label{au1}
\mathscr{H}om_{\mathcal{O}_X}\left( \mathcal{N}_{\mathcal{F}_0}^*, \mathcal{N}_{\mathcal{F}_0}^* \right)^\bullet: & \mathscr{H}om_{\mathcal{O}_X}\left( \mathcal{N}_{\mathcal{F}_0}^*, \mathcal{N}_{\mathcal{F}_0}^* \right) \xrightarrow{j}  \mathscr{H}om_{\mathcal{O}_X}\left( \mathcal{N}_{\mathcal{F}_0}^* ,  \Omega_X^1  \right) \bigoplus \mathscr{H}om_{\mathcal{O}_X}\left(  \mathcal{N}_{\mathcal{F}_0}^*, \mathcal{O}_X \right) \\
&\xrightarrow{ G_1'' } \mathscr{H}om_{\mathcal{O}_X}\left(  \mathcal{N}_{\mathcal{F}_0}^* , \tilde{\mathcal{S}}^2 \right) \bigoplus \mathscr{H}om_{\mathcal{O}_X}\left( \mathcal{N}_{\mathcal{F}_0}^*, \tilde{\mathcal{S}}^1 \right) \xrightarrow{G_2} \mathscr{H}om_{\mathcal{O}_X}\left( \mathcal{N}_{\mathcal{F}_0}^*, \tilde{\mathcal{S}}^3\right) \bigoplus \mathscr{H}om_{\mathcal{O}_X}\left( \mathcal{N}_{\mathcal{F}_0}^*, \tilde{\mathcal{S}}^2 \right) \xrightarrow{G_3} \cdots \notag
\end{align}}}
We will denote the $i$-th cohomology group by $\mathbb{H}_{\textnormal{Unf}}^i\left( X,  \mathscr{H}om_{\mathcal{O}_X}\left( \mathcal{N}_{\mathcal{F}_0}^*, \mathcal{N}_{\mathcal{F}_0}^* \right)  \right)$. From the exact sequence we see that
\begin{align*}
\mathbb{H}^i_{\textnormal{Unf}} \left(X, \mathscr{H}om_{\mathcal{O}_X}\left( \mathcal{N}_{\mathcal{F}_0}^* , \mathcal{N}_{\mathcal{F}_0}^* \right)^\bullet \right)\cong \mathbb{H}^{i-1}\left( X, \textnormal{Unf}_{ \mathcal{N}_{\mathcal{F}_0}^*}^\bullet\right),\,\,\,\,\,i\geq 0 . 
\end{align*}
This shows that $\mathscr{H}om_{\mathcal{O}_X}\left( \mathcal{N}_{\mathcal{F}_0}^* , \mathcal{N}_{\mathcal{F}_0}^* \right)^\bullet $ also controls flat unfoldings of $\left( X, \mathcal{N}_{\mathcal{F}_0}^* \right)$ with $X$ fixed. Moreover the complex $\textnormal{(\ref{au1})}$ explicitly contains $\mathscr{H}om_{\mathcal{O}_X}\left( \mathcal{N}_{\mathcal{F}_0}^* ,  \mathcal{N}_{\mathcal{F}_0}^* \right)$ as its component. We will reinterpret flat unfoldings of $\left(X, \mathcal{N}_{\mathcal{F}_0}^* \right)$ with $X$ fixed by using the complex $\textnormal{(\ref{au1})}$.

\begin{remark}
To emphasize deformations of the underlying locally free sheaf $\mathcal{N}_{\mathcal{F}_0}^*$, we can define infinitesimal unfoldings of $\left(X, \mathcal{N}_{\mathcal{F}_0}^* \right)$ with $X$ fixed equivalently in the following way: an infinitesimal unfolding $\left(\mathcal{X}=X\times \textnormal{Spec}(A), \mathcal{N}_{\mathcal{F}}^* \right)$ of $\left(X, \tilde{\mathcal{N}}_{\mathcal{F}_0}^* \right)$ under the trivial deformation of $X$ over $A\in \bold{Art}$ is a locally free sheaf $\mathcal{N}_{\mathcal{F}}^*$ and an exact sequence on $\mathcal{X}=X\times \textnormal{Spec}(A)$
\begin{align*}
0 \to \tilde{\mathcal{N}}_{\mathcal{F}}^* \xrightarrow{l} \Omega_{\mathcal{X}/k}^1 \to \frac{\Omega_{\mathcal{X}/k}^1}{l\left(\tilde{\mathcal{N}}_{\mathcal{F}}^* \right)}\to 0
\end{align*}
such that $d \left( l\left(\tilde{\mathcal{N}}_{\mathcal{F},x}^* \right) \right) \subset l\left( \tilde{\mathcal{N}}_{\mathcal{F}}^* \right)_x \wedge \Omega_{\mathcal{X}/k,x}^1$ for $x\in X-\textnormal{Sing}\left( \mathcal{F}_0 \right)$ and $\tilde{\mathcal{N}}_\mathcal{F}^* \cong \alpha \circ l\left( \tilde{\mathcal{N}}_\mathcal{F}^* \right)= \mathcal{N}_\mathcal{F}^*$  and $\left(\mathcal{X}, \mathcal{N}_\mathcal{F}^* \right)$ is a flat foliated deformation over $A$ of $\left( X , \mathcal{N}_{\mathcal{F}_0}^* \right)$ with $X$ fixed, and  induces $0\to \mathcal{N}_{\mathcal{F}_0}^* \to \Omega_X^1 \to \frac{\Omega_X^1}{\mathcal{N}_{\mathcal{F}_0}^*}\to 0$ over $\textnormal{Spec}(k)$. Two infinitesimal deformations $\left(\mathcal{X}, \tilde{\mathcal{N}}_{\mathcal{F}}^*\right)$ and $\left(\mathcal{X}, \tilde{ \mathcal{N} }_{\mathcal{F}'}^*\right)$ of $\left(X, \mathcal{N}_{\mathcal{F}_0}^* \right)$ with $X$ fixed are equivalent  if there is an isomorphism $\beta: \tilde{\mathcal{N}}_{\mathcal{F}}^* \to \tilde{\mathcal{N}}_{\mathcal{F}'}^* $ over $\mathcal{X}$ such that $l= l'\circ \beta$ where $0 \to \tilde{\mathcal{N}}_{\mathcal{F}'}^* \xrightarrow{l'} \Omega_{\mathcal{X}/k}^1 \to \frac{\Omega_{\mathcal{X}/k}^1}{\tilde{\mathcal{N}}_{\mathcal{F}'}^*}\to 0$. Then we have a functor of Artin rings
\begin{align*}
\textnormal{Unf}_{\mathcal{N}_{\mathcal{F}_0}}: \bold{Art} &\to \textnormal{(Sets)} \\
                             A &\to \left\{ \textnormal{Infinitesimal unfoldings of $\left( X, \mathcal{N}_{\mathcal{F}_0}^*\right)$ with $X$ fixed}   \right\}/\textnormal{isomorphisms}
\end{align*}
\end{remark}

By combining Proposition \ref{uu10} and the related Proposition in Part II, we obtain

\begin{proposition}
Let $\left(X, \mathcal{N}_{\mathcal{F}_0}^* \right)$ be a foliated nonsingular variety with $\mathcal{N}_{\mathcal{F}_0}^*$ locally free.
\begin{enumerate}
\item there is a $1-1$ correspondence
\begin{align*}
\textnormal{Unf}_{\mathcal{N}_{\mathcal{F}_0}^*}\left(k[\epsilon] \right)\cong \mathbb{H}^1_{\textnormal{Unf}}\left(X, \mathscr{H}om_{\mathcal{O}_X}\left(  \mathcal{N}_{\mathcal{F}_0}^* , \mathcal{N}_{\mathcal{F}_0}^* \right)^\bullet \right)
\end{align*}
\item Let $\eta=\left( \mathcal{X}, \tilde{\mathcal{N}}_\mathcal{F}^* \right)$ be an infinitesimal unfolding over $A:= k[x]/\left(x^m\right)$ of $\left( X, \mathcal{N}_{\mathcal{F}_0}^* \right)$ with $X$ fixed. To small extension $0\to (t)=\left(x^m\right)\to \tilde{A}: = k[x] / \left(x^{m+1} \right) \to A : = k[x] / \left(x^m\right) \to 0$, we can associate an element $o_\eta(e)\in \mathbb{H}_{\textnormal{Unf}}^2\left( X, \mathscr{H}om_{\mathcal{O}_X}\left( \mathcal{N}_{\mathcal{F}_0}^* , \mathcal{N}_{\mathcal{F}_0}^* \right)^\bullet \right)$ called the obstruction lifting of $\eta$ to $\tilde{A}$, which is $0$ if and only if a lifting of $\eta$ to $\tilde{A}$ exists.
\end{enumerate}
\end{proposition}

\section{Unfoldings of $\left(X, \mathcal{N}_{\mathcal{F}_0} \right)$ with $\mathcal{N}_{\mathcal{F}_0}^*$ locally free}

Let $\left( X, \mathcal{N}_{\mathcal{F}_0}^*\right)$ be a foliated nonsingular variety over $k$ with $\mathcal{N}_{\mathcal{F}_0}^*$ locally free. In this section we now vary the underlying nonsingular variety $X$. We recall  the definition of an infinitesimal unfolding of $\left(X, \mathcal{N}_{\mathcal{F}_0}^* \right)$ over $A\in \bold{Art}$ given in Appendix \ref{b2}. Then we have a functor of Artin rings
\begin{align*}
\textnormal{Unf}_{\left( X , \mathcal{N}_{\mathcal{F}_0}^* \right)}: \bold{Art} \to \left( \textnormal{Sets} \right)
\end{align*}

We will show that the unfolded dual leaf complex
{\small{\begin{align}
\tilde{\mathcal{N}}_{\mathcal{F}_0}^{*\bullet} :  T_X \xrightarrow{G_0} \mathscr{H}om_{\mathcal{O}_X}\left( \mathcal{N}_{\mathcal{F}_0}^* ,   \frac{\Omega_X^1}{\mathcal{N}_{\mathcal{F}_0}^*}  \right) \bigoplus \mathscr{H}om_{\mathcal{O}_X}\left(  \mathcal{N}_{\mathcal{F}_0}^*, \mathcal{O}_X \right) \xrightarrow{G_1} \mathscr{H}om_{\mathcal{O}_X}\left( \mathcal{N}_{\mathcal{F}_0}^* , \tilde{\mathcal{S}}^2 \right) \bigoplus \mathscr{H}om_{\mathcal{O}_X}\left( \mathcal{N}_{\mathcal{F}_0}^*, \tilde{\mathcal{S}}^1 \right) \label{au2}\\
 \xrightarrow{G_2} \mathscr{H}om_{\mathcal{O}_X}\left( \mathcal{N}_{\mathcal{F}_0}^*, \tilde{\mathcal{S}}^3 \right) \bigoplus \mathscr{H}om_{\mathcal{O}_X}\left( \mathcal{N}_{\mathcal{F}_0}^*, \tilde{\mathcal{S}}^2 \right) \xrightarrow{G_3} \cdots \notag
\end{align}}}
controls unfoldings of $\left( X, \mathcal{N}_{\mathcal{F}_0}^* \right)$.

\begin{lemma}[compare \cite{Ser06} Lemma 1.2.6 p.26]\label{uu23}
Let $B_0$ be a $k$-algebra with an element $w_0\in \Omega_{B_0/k}^1$, and $e:0  \to (t) \to \tilde{A} \to A \to 0$, where $(t)=(x^n), \tilde{A}=k[x]/(x^{n+1})$, and $A= k[x]/(x^n)$,  be a small extension in $\bold{Art}$. Let $\omega \in \Omega_{B_0\otimes A/k}^1\cong \left(\Omega_{B_0/k}^1\otimes_{B_0}(B_0\otimes A) \right) \bigoplus \left( \Omega_{A/k}^1\otimes_A (B_0\otimes A)\right) $, and let $\omega_1, \omega_2\in \Omega_{B_0\otimes \tilde{A}/k}^1\cong  \left(\Omega_{B_0/k}^1\otimes_k  \tilde{A} \right) \bigoplus \left(\Omega_{\tilde{A}/k}^1\otimes_k B_0 \right)$ which induces $\omega$, so that we can write $\omega_1-\omega_2= tw'+ v dt$ for some $w' \in  \Omega_{B_0/k}^1$ and $v\in B_0$. Then there is one to one correspondence
\begin{align*}
\left\{ \textnormal{isomorphisms on $\left(B_0\otimes_k \tilde{A}, \omega_1\right)$ and $\left(B_0\otimes_k \tilde{A}, \omega_2 \right)$ inducing the identity on $\left(B_0\otimes_k A, \omega \right)$ }   \right\} \\
\to \left\{ P\in \textnormal{Der}_k(B_0,B_0)=\textnormal{Hom}_{B_0}(\Omega_{B_0/k}^1, B_0) |w'+\mathcal{L}_P(w_0)=0, v+ i_P(w_0)=0  \right\}
\end{align*}
In particular, when $\omega_1= \omega_2$, there is a canonical isomorphism of groups
\begin{align*}
\left\{\textnormal{automorphisms on $\left(B_0\otimes_k \tilde{A}, \omega_1 \right)$ inducing the identity on $\left(B_0\otimes_k A, \omega \right)$}  \right\}\\
 \to \left\{P\in \textnormal{Der}_k(B_0,B_0)| \mathcal{L}_p(w_0)=0,  i_P(w_0)=0  \right\}
\end{align*}
\end{lemma}

\begin{proof}
Let $\theta:\left( B_0\otimes_k \tilde{A} , \omega_1\right)\to \left( B_0\otimes_k \tilde{A} , \omega_2 \right)$ be an isomorphism sending $\omega_1$ to $\omega_2$ inducing the identity on $\left( B_0\otimes_k A, T \right)$. Then $\theta(x)= x + P(dx)$ for some $P\in \textnormal{Der}_k\left( B_0, B_0 \right) = \textnormal{Hom}_{B_0}\left( \Omega_{B_0/k}^1, B_0 \right)$. We define the correspondence $\theta \mapsto P$. On the other hand, let $\omega_1 = \sum_{i=1}^r \tilde{g} d \tilde{f}_i$ for some $g_i, f_i \in B_0\otimes_k \tilde{A}$. Then $w_0=\sum_{i=1}^r g_i df_i$ where $g_i, f_i \in B_0$ which are the restriction of $\tilde{g}_i, \tilde{f}_i$ on $B_0$. Since $\omega_1$ corresponds to $\omega_2$ via $\theta= \textnormal{Id}+ t P$, we have

\begin{align*}
\omega_2&=\sum_{i=1}^r \left(\tilde{g}_i+tP\left(\tilde{g}_i\right)\right)d\left(\tilde{f}_i + tP\left(\tilde{f}_i\right)\right)=\sum_{i=1}^r \left(\tilde{g}_i+tP(\tilde{g}_i) \right) \left(d\tilde{f}_i+t dP\left(\tilde{f}_i\right)+P \left(\tilde{f}_i \right)dt \right)\\
       &=\sum_{i=1}^r \left( \tilde{g}_i d\tilde{f}_i+t P \left(g_i \right) df_i+ tg_i dP\left(f_i \right)+ g_iP(f_i)dt  \right) = \omega_1+ t \mathcal{L}_P\left(\sum_{i=1}^r g_i df_i \right)+ i_{P}\left(\sum_{i=1}^r g_idf_i\right)dt
\end{align*}
This implies that $w'+ \mathcal{L}_P \left(w_0 \right)=0$ and $v+i_P\left(w_0\right)=0$. Since $\theta$ is determined by $P$, the correspondence is one to one. Now assume that $ \omega_1=\omega_2$. Then $\theta$ corresponds to $P$ with $\mathcal{L}_P\left( w_0 \right)=0$ and $i_P \left( w_0 \right)=0$, and the correspondence $\theta \mapsto P$ is a group isomorphism.

\end{proof}

\begin{proposition}\label{uu35}
Let $\left( X, \mathcal{N}_{\mathcal{F}_0}^* \right)$ be a foliated nonsingular variety with $\mathcal{N}_{\mathcal{F}_0}^*$ locally free.
\begin{enumerate}
\item there is a $1-1$ correspondence
\begin{align*}
\textnormal{Unf}_{\left(X, \mathcal{N}_{\mathcal{F}_0}^*\right)}\left(k[\epsilon] \right)\cong \mathbb{H}^1\left(X, \tilde{\mathcal{N}}_{\mathcal{F}_0}^{* \bullet } \right)
\end{align*}
\item Let $\eta=\left( \mathcal{X}, \tilde{\mathcal{N}}_\mathcal{F}^* \right)$ be an infinitesimal unfolding of $\left( X, \mathcal{N}_{\mathcal{F}_0}^* \right)$ over $A$. To  small extension $0\to (t)=(x^m) \to \tilde{A}=k[x]/(x^{m+1}) \to A=k[x]/(x^m) \to 0$, we can associate an element $o_\eta(e)\in \mathbb{H}^2\left( X,  \tilde{\mathcal{N}}_{\mathcal{F}_0}^{* \bullet } \right)$ called the obstruction lifting of $\eta$ to $\tilde{A}$, which is $0$ if and only if a lifting of $\eta$ to $\tilde{A}$ exists.
\end{enumerate}
\end{proposition}

\begin{proof}
Given a first-order unfolding $\left( \mathcal{X}, \tilde{\mathcal{N}}_\mathcal{F}^* \right)$ of $\left( X, \mathcal{N}_{\mathcal{F}_0}^* \right)$ over $k[\epsilon]$, we choose an affine open covering $\mathcal{U}=\left\{ U_i \right\}$ of $X$ such that $\Gamma\left( U_i, \mathcal{N}_{\mathcal{F}_0}^* \right)$ is generated by $w_{0i}^1,.., w_{0i}^q\in \Gamma\left( U_i, \Omega_X^1 \right)$ with $w_{0i}^\alpha = \sum_{\beta=1}^q h_{0ij}^{\alpha\beta} w_{0j}^\beta$ for some $h_{0ij}^{\alpha\beta}\in \Gamma\left( U_{ij}, \mathcal{O}_X \right)$ and $w_{0i}^1 \wedge \cdots \wedge w_{0i}^q \wedge d w_{0i}^\alpha = 0, \alpha=1,...,q$, and we have $\mathcal{X}|_{U_i} \cong U_i\times \textnormal{Spec}\left( k[\epsilon] \right)$ is trivial for all $i$ with the induced unfolded structure $\tilde{\mathcal{N}}_{\mathcal{F}_i}^*$ on $U_i\times \textnormal{Spec}\left( k[\epsilon] \right)$, and $\Gamma\left( U_i, \tilde{\mathcal{N}}_{\mathcal{F}_i}^*\right)$ is generated by $w_{0i}^\alpha + \epsilon w_i'^\alpha + v_i^\alpha d \epsilon, \alpha=1,...,q$, where $w_i'^\alpha \in \Gamma\left( U_i, \Omega_X^1 \right)$ and $v_i^\alpha \in \Gamma\left( U_i, \mathcal{O}_X \right)$, and $\left( w_{0i}^1 + \epsilon w_i'^1 + v_i^1 d \epsilon \right) \wedge \cdots \wedge\left( w_{0i}^q + \epsilon w_i'^q + v_i^q d\epsilon \right) \wedge d \left( w_{0i}^\alpha + \epsilon w_i'^\alpha + v_i^\alpha d \epsilon \right)=0$. For all $i$, we have unfolded isomorphisms
\begin{align*}
\theta_i: \left( U_i \times \textnormal{Spec}\left( k[\epsilon] \right), \tilde{\mathcal{N}}_{\mathcal{F}_i}^* \right) \to \left( \mathcal{X}|_{U_i}, \tilde{\mathcal{N}}_{\mathcal{F}}^*|_{U_i} \right)
\end{align*}
Then for each $i,j, \theta_{ij}:=\theta_j^{-1}\theta_i: \left( U_{ij} \times \textnormal{Spec}\left( k[\epsilon] \right), \tilde{\mathcal{N}}_{\mathcal{F}_i}^*  \right) \to \left( U_{ij}\times \textnormal{Spec}\left( k[\epsilon] \right), \tilde{\mathcal{N}}_{\mathcal{F}_j}^* \right)$ is an unfolded isomorphism inducing the identity on $\left( U_{ij}, \mathcal{N}_{\mathcal{F}_0}^* \right)$ by modulo $\epsilon, d\epsilon$. Then by Lemma \ref{uu23}, $\theta_{ij}$ corresponds to $\textnormal{Id} + \epsilon p_{ij}:\left( \Gamma\left( U_{ij}, \mathcal{O}_X \right), \tilde{\mathcal{N}}_{\mathcal{F}_j}^* \right) \to \left( \Gamma\left( U_{ij},  \mathcal{O}_X \right), \tilde{\mathcal{N}}_{\mathcal{F}_i}^* \right)$ for some $p_{ij}\in \Gamma\left( U_{ij}, T_X \right)$ such that
\begin{align}\label{uu24}
\theta_{ji}^*\left( w_{0i}^\alpha + \epsilon w_i'^\alpha + v_i^\alpha d\epsilon \right) = w_{0i}^\alpha + \epsilon w_i'^\alpha + v_i^\alpha d \epsilon - \epsilon \mathcal{L}_{p_{ij}}\left( w_{0i}^\alpha \right) - i_{p_{ij}}\left( w_{0i}^\alpha \right) d \epsilon = \sum_{\beta=1}^q \left( h_{0ij}^{\alpha\beta} + \epsilon g_{ij}^{\alpha\beta} \right)\left( w_{0j}^\beta + \epsilon w_j'^\beta + v_j^\beta d \epsilon \right)
\end{align}
for some $g_{ij}^{\alpha\beta}\in \Gamma\left( U_{ij}, \mathcal{O}_X \right)$. Then $\textnormal{(\ref{uu24})}$ implies that
\begin{align}
w_i'^\alpha - \mathcal{L}_{p_{ij}}\left( w_{0i}^\alpha \right) &= \sum_{\beta=1}^q h_{0ij}^{\alpha\beta} w_j'^\beta + \sum_{\beta=1}^q g_{ij}^{\alpha\beta} w_{0j}^\beta \label{uu27} \\
v_i^\alpha - i_{p_{ij}}\left( w_{0i}^\alpha \right) &= \sum_{\beta=1}^q  h_{0ij}^{\alpha\beta} v_j^\beta \label{uu28}
\end{align}
We define $\overline{w_i'} \in \Gamma\left( U_i, \mathscr{H}om_{\mathcal{O}_X}\left( \mathcal{N}_{\mathcal{F}_0}^*,  \frac{\Omega_X^1}{\mathcal{N}_{\mathcal{F}_0}^*} \right) \right) $ as in $\textnormal{(\ref{uu25})}$ and $v_i \in \Gamma\left( U_i , \mathscr{H}om_{\mathcal{O}_X} \left( \mathcal{N}_{\mathcal{F}_0}^*, \mathcal{O}_X  \right) \right)$ as in $\textnormal{(\ref{uu26})}$. Then $\textnormal{(\ref{uu7})}$ and $\textnormal{(\ref{uu8})}$ and $\textnormal{(\ref{uu27})}$ and $\textnormal{(\ref{uu28})}$ imply that $\left( \left\{ \overline{ w_i' }  \right\} , \left\{ v_i \right\} ,  \left\{ p_{ij} \right\} \right)\in C^0\left( \mathcal{U}, \mathscr{H}om_{\mathcal{O}_X}\left( \mathcal{N}_{\mathcal{F}_0}^*, \frac{\Omega_X^1}{\mathcal{N}_{\mathcal{F}_0}^*} \right) \right) $ $ \bigoplus C^0\left( \mathcal{U}, \mathscr{H}om_{\mathcal{O}_X}\left( \mathcal{N}_{\mathcal{F}_0}^*, \mathcal{O}_X \right) \right) \bigoplus C^1\left(  \mathcal{U}, T_X \right)$ defines a $1$-cocycle in the following \v Cech resolution of $\tilde{\mathcal{N}}_{\mathcal{F}_0}^\bullet$:
\begin{center}
$\begin{CD}
\cdots \\
@AG_2 AA \\
C^0\left( \mathcal{U}, \left( \mathcal{N}_{\mathcal{F}_0}^* \right)^* \otimes \tilde{\mathcal{S}}^2 \bigoplus \left( \mathcal{N}_{\mathcal{F}_0}^* \right)^*\otimes \tilde{\mathcal{S}}^1 \right) @>-\delta >> \cdots \\
@AG_1AA @AG_1AA \\
C^0\left( \mathcal{U}, \left( \mathcal{N}_{\mathcal{F}_0}^* \right)^* \otimes  \frac{\Omega_X^1}{\mathcal{N}_{\mathcal{F}_0}^*} \bigoplus \left( \mathcal{N}_{\mathcal{F}_0}^* \right)^* \right) @>\delta >> C^1\left( \mathcal{U}, \left( \mathcal{N}_{\mathcal{F}_0}^* \right)^* \otimes \frac{\Omega_X^1}{\mathcal{N}_{\mathcal{F}_0}^*} \bigoplus \left( \mathcal{N}_{\mathcal{F}_0}^* \right)^* \right) @> - \delta >> \cdots \\
@AG_0AA @AG_0AA @AG_0AA \\
C^0\left(\mathcal{U}, T_X \right) @>- \delta>> C^1\left( \mathcal{U}, T_X \right) @>\delta >> C^2\left( \mathcal{U}, T_X \right) @>- \delta >> \cdots
\end{CD}$
\end{center}
We show that for two equivalent first-order unfoldings of $\left( X, \mathcal{N}_{\mathcal{F}_0}^* \right)$, the associated $1$-cocycles are equivalent. Let $\left( \mathcal{X}', \tilde{\mathcal{N}}_{\mathcal{F}'}^* \right)$ be another unfolding of $\left( X , \mathcal{N}_{\mathcal{F}_0}^* \right)$ over $k[\epsilon]$ which induces $\left( \left\{ \overline{w_i''} \right\} , \left\{ v_i' \right\}, \left\{ p_{ij}' \right\} \right)$ and $\Phi : \left( \mathcal{X}, \tilde{\mathcal{N}}_\mathcal{F}^* \right) \to \left( \mathcal{X}' , \tilde{\mathcal{N}}_{\mathcal{F}'}^* \right)$ is an unfolded isomorphism inducing the identity on $\left( X, \mathcal{N}_{\mathcal{F}_0}^* \right)$. Then there is an unfolded isomorphism $\alpha_i$ on $U_i\times \textnormal{Spec}\left( k[\epsilon] \right)$ given by
\begin{align*}
\alpha_i : \left( U_i\times \textnormal{Spec}\left( k[\epsilon] \right), \tilde{\mathcal{N}}_{\mathcal{F}_i}^*   \right) \xrightarrow{\theta_i} \left( \mathcal{X}|_{U_i} , \tilde{\mathcal{N}}_\mathcal{F}^* |_{U_i}    \right) \xrightarrow{\Phi|_{U_i}} \left( \mathcal{X}' |_{U_i} , \tilde{\mathcal{N}}_{\mathcal{F}'}^*|_{U_i}  \right) \xrightarrow{\theta_i'^{-1}} \left( U_i \times \textnormal{Spec}\left( k[\epsilon] \right) , \tilde{\mathcal{N}}_{\mathcal{F}_i'}^*   \right)
\end{align*}
which corresponds to $1+ \epsilon a_i : \left(\Gamma\left( U_i, \mathcal{O}_X \right)\otimes k[\epsilon] , \tilde{\mathcal{N}}_{\mathcal{F}_i'}^*  \right) \to \left( \Gamma\left( U_i, \mathcal{O}_X \right) \otimes k[\epsilon], \tilde{\mathcal{N}}_{\mathcal{F}_i}^* \right)$ for some $a_i\in \Gamma\left( U_i,  T_X \right) $ such that $w_{0i}^\alpha + \epsilon w_i'^\alpha + v_i^\alpha d \epsilon - \epsilon \mathcal{L}_{a_i}\left( w_{0i}^\alpha \right) - i_{a_i}\left( w_{0i}^\alpha \right) d \epsilon = \sum_{\beta=1}^q \left( \delta_\alpha^\beta + \epsilon b_i^{\alpha\beta} \right) \left( w_{0i}^\beta + \epsilon w_i''^\beta + v_i'^\beta d \epsilon \right) $ for some $b_i^{\alpha\beta}\in \Gamma\left( U_i, \mathcal{O}_X \right)$ by Lemma \ref{uu23}. Then we have $\mathcal{L}_{a_i}\left( w_{0i}^\alpha \right) = w_i'^\alpha - w_i''^\alpha - \sum_{\beta=1}^q b_i^{\alpha\beta} w_{0i}^\beta$ and $i_{a_i} \left( w_{0i}^\alpha \right) = v_i^\alpha - v_i'^\alpha$. On the other hand, since $\theta_i'\alpha_i= \Phi |_{U_i} \theta_i$, we have $\alpha_j^{-1}\theta_{ij}' \alpha_i= \left( \theta_j' \alpha_j \right)^{-1} \left( \theta_i' \alpha_i\right) = \theta_j^{-1}\Phi|_{U_{ij}}^{-1} \Phi|_{U_{ij}} \theta_i = \theta_{ij}$ which corresponds to $\left(\textnormal{Id} + \epsilon a_i \right)\left( \textnormal{Id} + \epsilon p_{ij}' \right)\left(  \textnormal{Id} - \epsilon a_j \right)= \textnormal{Id}+ t p_{ij}$, so that $a_i - a_j= p_{ij} - p_{ij}'$. This implies that $\left( \left\{ \overline{w_i'} \right\} , \left\{ v_i \right\} , \left\{ p_{ij} \right\}  \right)$ and $\left( \left\{ \overline{w_i''} \right\} , \left\{ v_i'\right\}, \left\{ p_{ij}' \right\} \right)$ are cohomologous.

We define an inverse map. Given an element in $\mathbb{H}^1\left( X, \tilde{\mathcal{N}}_{\mathcal{F}_0}^{*\bullet} \right)$, we represent it by a \v Cech $1$-cocycle $\left( \left\{ \overline{w_i'}  \right\} , \left\{ v_i  \right\}, \left\{ p_{ij}   \right\} \right)$ for an affine open covering $\mathcal{U}=\left\{ U_i \right\}$ of $X$ satisfying $\textnormal{(\ref{uu7})},\textnormal{(\ref{uu8})},\textnormal{(\ref{uu27})}$ and $\textnormal{(\ref{uu28})}$, and $p_{ij}- p_{ik} + p_{jk}=0$. By reversing the above process, we obtain a first-order unfolding of $\left( X, \mathcal{N}_{\mathcal{F}_0}^* \right)$.

Next we identify obstructions. Let $e: 0\to (t)=(x^m) \to \tilde{A}= k[x]/(x^{m+1}) \to A = k[x]/(x^m) \to 0$ a small extension and $\eta:=\left( \mathcal{X}, \tilde{\mathcal{N}}_\mathcal{F}^* \right)$ be an infinitesimal unfolding of $\left( X, \mathcal{N}_{\mathcal{F}_0}^* \right)$ over $A$. Let $\mathcal{U}=\left\{ U_i \right\}$ be an affine open covering of $X$ such that we have unfolded isomorphisms $\theta_i:\left( U_i \times \textnormal{Spec}(A) , \tilde{\mathcal{N}}_{\mathcal{F}_i}^* \right) \to \left( \mathcal{X}|_{U_i}, \tilde{\mathcal{N}}_{\mathcal{F}}^*|_{U_i} \right)$, where $\tilde{\mathcal{N}}_{\mathcal{F}_i}^* \subset \left( \Omega_X^1 |_{U_i} \otimes A \right) \bigoplus \left( \mathcal{O}_X|_{U_i}\otimes \Omega_{A/k}^1 \right)$ is the induced foliation on $U_i\times \textnormal{Spec}(A)$ from $\tilde{\mathcal{N}}_\mathcal{F}^*$ via $\theta_i$. Let $\Gamma\left( U_i, \tilde{\mathcal{N}}_{\mathcal{F}_i}^* \right) \subset \left( \Gamma\left( U_i, \Omega_X^1 \right) \otimes A \right) \bigoplus \left( \Gamma\left( U_i, \mathcal{O}_X \right) \otimes \Omega_{A/k}^1 \right)$ be generated by $w_i^\alpha + v_i^\alpha,\alpha=1,...,q$ which induces $w_{0i}^\alpha$, where $w_i^\alpha \in \Gamma\left( U_i, \Omega_X^1 \right)\otimes A $ and $v_i^\alpha\in \Gamma\left( U_i, \mathcal{O}_X \right)\otimes \Omega_{A/k}^1$ with $\left( w_i^1 + v_i^1 \right) \wedge \cdots \wedge \left( w_i^q + v_i^q \right) \wedge d \left( w_i^\alpha + v_i^\alpha \right)=0,\alpha=1,...,q$, and let $\theta_{ij}:=\theta_j^{-1} \theta_i: \left( U_{ij}\times \textnormal{Spec}(A), \tilde{\mathcal{N}}_{\mathcal{F}_j}^* \right) \to \left( U_{ij} \times \textnormal{Spec}(A), \tilde{\mathcal{N}}_{\mathcal{F}_i}^* \right)$ is an unfolded isomorphism with $\theta_{jk} \theta_{ij}= \theta_{ik}$ on $U_{ijk}\times \textnormal{Spec}(A)$ and $\theta_{ji}^*\left(  w_i^\alpha + v_i^\alpha \right)= \sum_{\beta=1}^q h_{ij}^{\alpha\beta}\left( w_j^\beta + v_j^\beta  \right)$ for some $h_{ij}^{\alpha\beta} \in \Gamma\left( U_{ij} , \mathcal{O}_X \right) \otimes A$ which induces $w_{0i}^\alpha =  \sum_{\beta=1}^q h_{0ij}^{\alpha\beta} w_{0j}^\beta$. We let $f_{ij}:\left( \Gamma\left( U_{ij}, \mathcal{O}_X \right)\otimes A ,  \tilde{\mathcal{N}}_{\mathcal{F}_j}^* \right) \to \left( \Gamma\left( U_{ij} , \mathcal{O}_X \right) \otimes A, \tilde{\mathcal{N}}_{\mathcal{F}_i}^* \right)$ be the ring homomorphism corresponding to $\theta_{ij}$. Then we have $f_{ij} f_{jk} = f_{ik}$. To give a lifting $\tilde{\eta}$ of $\eta$ to $\tilde{A}$ is equivalent to give a collection $\left\{ \tilde{\mathcal{N}}_{\tilde{\mathcal{F}}_i}^* \right\}$, where $\Gamma\left(U_i, \tilde{\mathcal{N}}_{\tilde{\mathcal{F}}_i}^* \right) \subset \left( \Gamma\left( U_i, \Omega_X^1 \right) \otimes \tilde{A} \right)$ $ \bigoplus  \left( \Gamma\left( U_i , \mathcal{O}_X \right)\otimes \Omega_{\tilde{A}/k}^1 \right)$ is generated by $\tilde{w}_i^\alpha + \tilde{v}_i^\alpha \in \left( \Gamma\left( U_i, \Omega_X^1 \right) \otimes \tilde{A} \right) \bigoplus \left( \Gamma\left( U_i, \mathcal{O}_X \right)\otimes \Omega_{\tilde{A}/k}^1\right),\alpha=1,...,q$ which induces $w_i^\alpha + v_i^\alpha$ such that $\left( \tilde{w}_i^1 + \tilde{v}_i^1 \right) \wedge \cdots \wedge \left( \tilde{w}_i^q + \tilde{v}_i^q \right) \wedge d \left( \tilde{w}_i^\alpha + \tilde{v}_i^\alpha \right)=0,\alpha=1,...,q$, and a collection of unfolded isomorphisms $\left\{ \tilde{\theta}_{ij} \right\}$ which induces $\theta_{ij}$, where $\tilde{\theta}_{ij} : \left( U_{ij} \times \textnormal{Spec}\left( \tilde{A} \right) ,  \tilde{\mathcal{N}}_{\tilde{\mathcal{F}}_i}^*   \right) \to \left( U_{ij} \times \textnormal{Spec}\left( \tilde{A} \right) , \tilde{\mathcal{N}}_{\tilde{\mathcal{F}}_j}^*   \right)$ such that $\tilde{\theta}_{ji}^* \left(  \tilde{w}_i^\alpha + \tilde{v}_i^\alpha \right)= \sum_{\beta=1}^q \tilde{h}_{ij}^{\alpha\beta} \left( \tilde{w}_j^\beta + \tilde{v}_j^\beta \right)$ for some $\tilde{h}_{ij}^{\alpha\beta}\in \Gamma\left( U_{ij} , \mathcal{O}_X \right)\otimes  \tilde{A}$ which induces $h_{ij}^{\alpha\beta}$, and $\tilde{\theta}_{jk} \tilde{\theta}_{ij} = \tilde{\theta}_{ik}$ (equivalently, $\tilde{f}_{ij} \tilde{f}_{jk} = \tilde{f}_{ik}$). From such data, we can glue together $\left( U_i \times \textnormal{Spec}\left( \tilde{A} \right)  , \tilde{\mathcal{N}}_{\tilde{\mathcal{F}}_i}^* \right)$ to make an infinitesimal unfolding $\left( \tilde{\mathcal{X}}, \tilde{\mathcal{N}}_{\tilde{\mathcal{F}}}^* \right)$ inducing $\left( \mathcal{X}, \tilde{\mathcal{N}}_\mathcal{F}^* \right)$. In order to see if there exists a lifting of $\left( \mathcal{X}, \tilde{\mathcal{N}}_{\mathcal{F}}^* \right)$ to $\tilde{A}$, we choose arbitrary automorphisms $\left\{ \tilde{\theta}_{ij} \right\}$ which restricts to $\theta_{ij}$ on $U_{ij}\times \textnormal{Spec}\left( A \right)$ and arbitrary lifting $\tilde{w}_i^\alpha + \tilde{v}_i^\alpha \in \left( \Gamma\left( U_i, \Omega_X^1 \right) \otimes \tilde{A} \right) \bigoplus \left( \Gamma\left( U_i, \mathcal{O}_X \right) \otimes \Omega_{\tilde{A}/k}^1 \right)$ of $w_i^\alpha + v_i^\alpha$ and arbitrary lifting $\tilde{h}_{ij}^{\alpha\beta}\in \Gamma\left( U_{ij}, \mathcal{O}_X \right)$ of $h_{ij}^{\alpha\beta},\alpha,\beta=1,...,q$. Since $\tilde{\theta}_{ijk}=\tilde{\theta}_{ki}\tilde{\theta}_{jk}\tilde{\theta}_{ij}$ (equivalently $\tilde{f}_{ijk}=\tilde{f}_{ij}\tilde{f}_{jk}\tilde{f}_{ki}$) is an automorphism on $U_{ijk}\times \textnormal{Spec}\left(\tilde{A}\right)$ inducing the identity on $U_{ijk}\times \textnormal{Spec}(A)$, $\tilde{\theta}_{ijk}$ corresponds to $\tilde{d}_{ijk}\in \Gamma\left(U_{ijk}, T_X\right)$ with $\tilde{d}_{jkl}-\tilde{d}_{ikl}+\tilde{d}_{ijl}-\tilde{d}_{jkl}=0$, so that $\delta(\{\tilde{d}_{ijk}\})=0$. On the other hand, we have
\begin{align}
\tilde{\theta}_{ji}^*\left( \tilde{w}_i^\alpha + \tilde{v}_i^\alpha \right)-\sum_{\beta=1}^q \tilde{h}_{ij}^{\alpha\beta} \tilde{w}_j^\beta = t R_{ij}^\alpha + G_{ij}^\alpha dt \,\,\,\,\,\,\,\,\,\textnormal{for some}\,\,\, R_{ij}^\alpha \in \Gamma\left( U_{ij}, \Omega_X^1 \right),\,\,\,\,\,G_{ij}^\alpha \in \Gamma\left( U_{ij}, \mathcal{O}_X \right)\label{uu33}\\
\left( \tilde{w}_i^1 + \tilde{v}_i^1 \right) \wedge \cdots \wedge \left( \tilde{w}_i^q + \tilde{v}_i^q \right) \wedge d \left( \tilde{w}_i^\alpha + \tilde{v}_i^\alpha \right)= t B_i^\alpha + F_i^\alpha \wedge dt\,\,\,\,\,\,\,\textnormal{for some}\,\,\, B_i^\alpha \in \Gamma\left( U_i, \tilde{\mathcal{S}}^2 \right),\,\,\,\,\,F_i^\alpha \in \Gamma\left( U_i, \tilde{\mathcal{S}}^1 \right)
\end{align}
We define $\overline{R_{ij}}\in \Gamma\left(  U_{ij}, \mathscr{H}om_{\mathcal{O}_X}\left( \mathcal{N}_{\mathcal{F}_0}^*, \frac{\Omega_X^1}{\mathcal{N}_{\mathcal{F}_0}^*} \right) \right)$ and $G_{ij}\in \Gamma\left( U_{ij} , \mathscr{H}om_{\mathcal{O}_X} \left( \mathcal{N}_{\mathcal{F}_0}^* , \mathcal{O}_X \right)  \right)$ as in a similar way to $\textnormal{(\ref{uu29})}$ and $\textnormal{(\ref{uu30})}$, and $B_i\in \Gamma\left( U_i, \mathscr{H}om_{\mathcal{O}_X}\left( \mathcal{N}_{\mathcal{F}_0}^*, \tilde{\mathcal{S}}^2 \right) \right)$ and $G_i \in \Gamma\left( U_i, \mathscr{H}om_{\mathcal{O}_X}\left( \mathcal{N}_{\mathcal{F}_0}^*, \tilde{\mathcal{S}}^1 \right) \right)$ as in a similar way to $\textnormal{(\ref{uu31})}$ and $\textnormal{(\ref{uu32})}$. We show that $\left\{ \overline{R_{ij}}\right\}\in C^1\left( \mathcal{U}, \mathscr{H}om_{\mathcal{O}_X}\left( \mathcal{N}_{\mathcal{F}_0}^*, \frac{\Omega_X^1}{\mathcal{N}_{\mathcal{F}_0}^*} \right) \right)$ and $\left\{ G_{ij} \right\} \in C^1\left( \mathcal{U}, \mathscr{H}om_{\mathcal{O}_X}\left( \mathcal{N}_{\mathcal{F}_0}^*, \mathcal{O}_X \right) \right)$. In fact, by applying $\tilde{\theta}_{ij}^*$ on $\textnormal{(\ref{uu33})}$, we have
\begin{align*}
&\tilde{w}_i^\alpha + \tilde{v}_i^\alpha - \sum_{\beta=1}^q \tilde{f}_{ij}\left( \tilde{h}_{ij}^{\alpha\beta} \right)\left( \sum_{\gamma=1}^q \tilde{h}_{ji}^{\beta \gamma} \left( \tilde{w}_i^\gamma + \tilde{v}_i^\gamma \right) + t R_{ji}^\beta + G_{ji}^\beta dt   \right) =  t R_{ij}^\alpha + G_{ij}^\alpha dt \\
&\iff \tilde{w}_i^\alpha + \tilde{v}_i^\alpha - \sum_{\gamma=1}^q \left( \delta_\alpha^\gamma +t c_{ij\alpha}^{ \gamma} \right)\left( \tilde{w}_i^\gamma + \tilde{v}_i^\gamma \right) - t \sum_{\beta=1}^q h_{0ij}^{\alpha\beta} R_{ji}^\beta - \sum_{\beta=1}^q h_{0ij}^{\alpha\beta} G_{ji}^\beta dt = t R_{ij}^\alpha + G_{ij}^\alpha dt\,\,\,\,\,\textnormal{for some}\,\,\, c_{ij\alpha}^{ \gamma}\in \Gamma\left( U_{ij}, \mathcal{O}_X \right)  \\
&\iff - t\sum_{\gamma=1}^q c_{ij\alpha}^\gamma w_{0i}^\gamma - t \sum_{\beta=1}^q h_{0ij}^{\alpha\beta} R_{ji}^\beta - \sum_{\beta=1}^q h_{0ij}^{\alpha\beta} G_{ji}^\beta dt = t R_{ij}^\alpha + G_{ij}^\alpha dt \\
&\Longrightarrow  R_{ij}^\alpha= - \sum_{\beta=1}^q h_{0ij}^{\alpha\beta} R_{ji}^\beta - \sum_{\gamma=1}^q c_{ij\alpha}^\gamma w_{0i}^\gamma,\,\,\,\,\,\,\,\,\, G_{ij}^\alpha = - \sum_{\beta=1}^q h_{0ij}^{\alpha\beta} G_{ji}^\beta \Longrightarrow \overline{R_{ij}} = - \overline{R_{ji}} ,\,\,\,\,\,\,\,G_{ij}= - G_{ji}
\end{align*}
We claim that 
{\small{\begin{align*}
&\left(\{(B_i ,F_i ) \}, \{\left( -\overline{R_{ij}}, - G_{ij} \right)\},\{  d_{ijk}\} \right)\\
&\in C^0\left( \mathcal{U},  \mathscr{H}om_{\mathcal{O}_X}\left(\mathcal{N}_{\mathcal{F}_0}^*, \tilde{\mathcal{S}}^2 \right)\bigoplus \mathscr{H}om_{\mathcal{O}_X}\left( \mathcal{N}_{\mathcal{F}_0}^*, \tilde{\mathcal{S}}^1 \right) \right)\bigoplus C^1\left(\mathcal{U}, \mathscr{H}om_{\mathcal{O}_X}\left( \mathcal{N}_{\mathcal{F}_0}^*, \frac{\Omega_X^1}{\mathcal{N}_{\mathcal{F}_0}^* } \right) \bigoplus \mathscr{H}om_{\mathcal{O}_X}\left(\mathcal{N}_{\mathcal{F}_0}^*, \mathcal{O}_X \right)\right) \bigoplus C^2(\mathcal{U}, T_X )
\end{align*}}}
in the above \v Cech resolution of $\tilde{\mathcal{N}}_{\mathcal{F}_0}^{*\bullet}$. In fact, we have $-\delta\left( \left\{ d_{ijk}  \right\}  \right)=0$ as above and $G_2\left( B_i, F_i \right)=0$ as in the proof of Proposition \ref{uu10}. We have $-\delta\left( \left\{ \left( B_i , F_i\right)\right\} \right)+ G_1\left( \left\{- \overline{R}_{ij}, - G_{ij} \right\} \right) = 0 $ which follows from
{\tiny{\begin{align*}
&t B_i^\alpha+ F_i^\alpha \wedge dt= \theta_{ji}^*\left( t B_i^\alpha+ F_i^\alpha \wedge dt\right)=\tilde{\theta}_{ji}^*\left(\tilde{w}_i^1+ \tilde{v}_i^1\right)\wedge \cdots \wedge \tilde{\theta}_{ji}^*\left(\tilde{w}_i^q+ \tilde{v}_i^q\right) \wedge d \tilde{\theta}_{ji}^*\left( \tilde{w}_i^\alpha+ \tilde{v}_i^\alpha\right)\\
     &=\left(tR_{ij}^1+G_{ij}^1 dt +\sum_{\beta_1=1}^q \tilde{h}_{ij}^{\alpha\beta_1}\left( \tilde{w}_j^{\beta_1}+ \tilde{v}_j^{\beta_1}\right) \right) \wedge \cdots \wedge \left( t R_{ij}^q+ G_{ij}^q dt +\sum_{\beta_q=1}^q \tilde{h}_{ij}^{\alpha\beta_q} \left(\tilde{w}_j^{\beta_q} + \tilde{v}_j^{\beta_q} \right)\right) \wedge d \left( tR_{ij}^\alpha+ G_{ij}^\alpha dt +\sum_{\beta=1}^q \tilde{h}_{ij}^{\alpha\beta}\left( \tilde{w}_j^\beta  + \tilde{v}_j^\beta\right) \right)\\
\end{align*}}}
and $\textnormal{(\ref{uu34})}$. On the other hand, we show that $-\delta\left( \left\{\left( -\overline{R}_{ij}, - G_{ij} \right) \right\} \right) + G_0\left( d_{ijk} \right)=0$. In fact, first we note that
\begin{align}
w_i^\alpha + v_i^\alpha  &=\theta_{ij}^* \theta_{jk}^* \theta_{ki}^*\left( w_i^\alpha +v_i^\alpha\right) = \theta_{ij}^*\theta_{jk}^*\left(\sum_{\beta=1}^q h_{ik}^{\alpha\beta} \left(w_k^\beta + v_k^\beta \right) \right)= \theta_{ij}^*\left( \sum_{\beta, \gamma=1}^q f_{jk}\left(h_{ik}^{\alpha\beta} \right)  h_{kj}^{\beta \gamma} \left( w_j^\gamma + v_j^\gamma \right) \right) \\
& = \sum_{\beta, \gamma, \eta=1}^q \left(f_{ij} f_{jk} \right)\left(h_{ik}^{\alpha\beta} \right)f_{ij}\left( h_{kj}^{\beta \gamma} \right) h_{ji}^{\gamma \eta} \left( w_i^\eta + v_i^\eta \right) \notag \\
&\Longrightarrow \sum_{\beta, \gamma=1}^q  \left(\tilde{f}_{ij} \tilde{f}_{jk} \right)\left(\tilde{h}_{ik}^{\alpha\beta} \right) \tilde{f}_{ij}\left( \tilde{h}_{kj}^{\beta \gamma} \right) \tilde{h}_{ji}^{\gamma \eta} = \delta_\alpha^\eta + t u_{ijk}^{\alpha \eta}\,\,\,\,\,\,\,\textnormal{for some}\,\,\,  u_{ijk}^{\alpha \eta} \in \Gamma\left( U_{ijk}, \mathcal{O}_X \right) \label{tn6}
\end{align}
We have an isomorphism $\textnormal{Id}+ td_{ijk}:\Gamma\left( U_{ijk}, \mathcal{O}_X \right)\otimes \tilde{A}\to \Gamma\left(U_{ijk}, \mathcal{O}_X \right)\otimes\tilde{A}$ corresponding to $\tilde{\theta}_{ki}\tilde{\theta}_{jk} \tilde{\theta}_{ij}$ and so $\tilde{f}_{ij}\tilde{f}_{jk}\tilde{f}_{ki}=\textnormal{Id}+t d_{ijk}$. Then by Lemma \ref{uu23}, 
\begin{align*}
&\tilde{w}_{i}^\alpha + \tilde{v}_i^\alpha +t \mathcal{L}_{ d_{ijk}}\left( w_{0i}^\alpha \right) + i_{d_{ijk}}\left( w_{0i}^\alpha \right) dt =\tilde{\theta}_{ij}^*\tilde{\theta}_{jk}^*\tilde{\theta}_{ki}^* \left(\tilde{w}_i^\alpha + \tilde{v}_i^\alpha \right)=\tilde{\theta}_{ij}^*\tilde{\theta}_{jk}^*\left(tR_{ik}^\alpha + G_{ik}^\alpha dt +\sum_{\beta=1}^q \tilde{h}_{ik}^{\alpha\beta} \left(\tilde{w}_k^\beta + \tilde{v}_k^\beta \right)\right) \\
&=\tilde{\theta}_{ij}^*\left(tR_{ik}^\alpha + G_{ik}^\alpha dt +\sum_{\beta=1}^q\tilde{h}_{jk}\left(\tilde{h}_{ik}^{\alpha\beta} \right)\left(tR_{kj}^\beta + G_{kj}^\beta dt +\sum_{\gamma=1}^q \tilde{h}_{kj}^{\beta\gamma} \left( \tilde{w}_j^\gamma + \tilde{v}_j^\gamma \right) \right) \right)\\
&=tR_{ik}^\alpha + G_{ik}^\alpha dt +t\sum_{\beta=1}^q h_{0ik}^{\alpha\beta} R_{kj}^\beta +\sum_{\beta=1}^q h_{0ik}^{\alpha \beta} G_{kj}^\beta dt  +\sum_{\beta,\gamma=1}^q \left(\tilde{f}_{ij}\tilde{f}_{jk}\right)\left( \tilde{h}_{ik}^{\alpha\beta}\right) \tilde{f}_{ij}\left( \tilde{h}_{kj}^{\beta\gamma}\right)\left(tR_{ji}^\gamma + G_{ji}^\gamma dt +\sum_{\eta=1}^q \tilde{h}_{ji}^{\gamma\eta} \left( \tilde{w}_i^\eta  + \tilde{v}_i^\eta \right) \right)\\
&=t\left(R_{ik}^\alpha+ \sum_{\beta=1}^q h_{0ik}^{\alpha\beta} R_{kj}^\beta+\sum_{\gamma=1}^p r_{0ij}^{\alpha\gamma} R_{ji}^\gamma \right) + \left( G_{ik}^\alpha + \sum_{\beta=1}^q h_{0ik}^{\alpha\beta} G_{kj}^\beta + \sum_{\beta=1}^q h_{0ij}^{\alpha \beta} G_{ji}^\beta    \right) dt + \tilde{w}_i^\alpha + \tilde{v}_i^\alpha + t\sum_{\gamma=1}^q u_{ijk}^{\alpha\eta} w_{0i}^\eta
\end{align*}
Then we have
\begin{align*}
\mathcal{L}_{d_{ijk}}\left( w_{0i}^\alpha \right) - \left(R_{ik}^\alpha+ \sum_{\beta=1}^q h_{0ik}^{\alpha\beta} R_{kj}^\beta+\sum_{\gamma=1}^p r_{0ij}^{\alpha\gamma} R_{ji}^\gamma \right) &= \sum_{\gamma=1}^q u_{ijk}^{\alpha \eta} w_{0i}^\eta\\
i_{d_{ijk}}\left( w_{0i}^\alpha \right) - \left( G_{ik}^\alpha + \sum_{\beta=1}^q h_{0ik}^{\alpha\beta} G_{kj}^\beta + \sum_{\beta=1}^q h_{0ij}^{\alpha \beta} G_{ji}^\beta    \right) & = 0
\end{align*}
This proves our claim. We show that for another lifting $\tilde{\theta}_{ij}'$ of $\theta_{ij}$, and another lifting $\tilde{w}_i'^\alpha$ of $w_i^\alpha$, and another lifting of $\tilde{v}_i'^\alpha$ of $v_i^\alpha$, and another lifting of $\tilde{h}_{ij}'^{\alpha\beta}$ of $h_{ij}^{\alpha\beta}$, the associated $2$-cocycle $\beta:= \left( \left\{ \left( B_i', F_i' \right) \right\}, \left\{ \left( - \overline{R_{ij}'} , - G_{ij}' \right) \right\}, \left\{ d_{ijk}'\right\} \right)$ is cohomologous to $\alpha:=\left( \left\{ \left( B_i, F_i \right) \right\}, \left\{ \left( - \overline{R_{ij}}, - G_{ij} \right) \right\}, \left\{ d_{ijk} \right\} \right)$. In fact, we have $\tilde{f}_{ij}'=\tilde{f}_{ij} + t p_{ij}$ for some $p_{ij}\in \Gamma\left( U_{ij}, T_X \right)$ and $\tilde{w}_i'^\alpha= \tilde{w}_i^\alpha + t V_i^\alpha$ for some $V_i^\alpha \in \Gamma\left( U_i, \Omega_X^1 \right)$, and $\tilde{v}_i'^\alpha = \tilde{v}_i^\alpha + r_i^\alpha dt$ for some $r_i^\alpha \in \Gamma\left( U_i, \mathcal{O}_X \right)$, and $\tilde{h}_{ij}'^{\alpha\beta} = \tilde{h}_{ij}^{\alpha\beta} + t y_{ij}^{\alpha\beta}$ for some $y_{ij}^{\alpha\beta}\in \Gamma\left( U_{ij}, \mathcal{O}_X \right)$. We define $\overline{V_i}\in \Gamma\left( U_i, \mathscr{H}om_{\mathcal{O}_X}\left( \mathcal{N}_{\mathcal{F}_0}^*, \frac{\Omega_X^1}{\mathcal{N}_{\mathcal{F}_0}^*} \right) \right), w_{0i}^\alpha \mapsto \overline{V_i^\alpha}$ and $r_i\in \Gamma\left( U_i, \mathscr{H}om_{\mathcal{O}_X}\left( \mathcal{N}_{\mathcal{F}_0}^*, \mathcal{O}_X \right) \right), w_{0i}^\alpha \mapsto r_i^\alpha$. Then as in $\textnormal{(\ref{uu21})}$, we have $G_1\left( \left\{ \overline{V_i} \right\}, \left\{ r_i \right\} \right)= \left( \left\{ B_i'- B_i \right\}, \left\{ F_i' - F_i \right\} \right)$. We note that
\begin{align}
&t\left(- R_{ij}'^\alpha + R_{ij}^\alpha \right) + \left(- G_{ij}^\alpha + G_{ij}^\alpha \right) dt =  - \tilde{\theta}_{ji}'^*\left( \tilde{w}_i'^\alpha + \tilde{v}_i'^\alpha \right) + \sum_{\beta=1}^q \tilde{h}_{ij}'^{\alpha\beta} \left(\tilde{w}_j'^\beta + \tilde{v}_j'^\beta \right)+ \tilde{\theta}_{ji}^*\left( \tilde{w}_i^\alpha + \tilde{v}_i^\alpha \right) - \sum_{\beta=1}^q \tilde{h}_{ij}^{\alpha\beta} \left(\tilde{w}_j^\beta + \tilde{v}_j^\beta \right) \\
&=\left(- \tilde{\theta}_{ji}'^* +\tilde{\theta}_{ji}^*\right) \left( \tilde{w}_i^\alpha + \tilde{v}_i^\alpha \right)- t V_{i}^\alpha - r_i^\alpha dt + \sum_{\beta=1}^q \left( \tilde{h}_{ij}^{\alpha\beta} +t y_{ij}^{\alpha\beta} \right)\left( \tilde{w}_j^\beta + \tilde{v}_j^\beta + t V_j^\beta + r_j^\beta dt \right) - \sum_{\beta=1}^q \tilde{h}_{ij}^{\alpha\beta} \left(\tilde{w}_j^\beta + \tilde{v}_j^\beta \right) \notag \\
&= - t \mathcal{L}_{p_{ji}}\left( w_{0i}^\alpha \right) - i_{p_{ji}}\left( w_{0i}^\alpha \right)dt + t  \left(  \sum_{\beta=1}^q h_{0ij}^{\alpha\beta} V_j^\beta - V_i^\alpha \right) + t \sum_{\beta=1}^q y_{ij}^{\alpha\beta} w_{0j}^\beta + \left( \sum_{\beta=1}^q h_{0ij}^{\alpha\beta} r_j^\beta - r_i^\alpha \right) dt  \notag \\
&= t\left( \mathcal{L}_{p_{ij}}\left( w_{0i}^\alpha \right) + \sum_{\beta=1}^q h_{0ij}^{\alpha\beta} V_j^\beta - V_i^\alpha + \sum_{\beta=1}^q y_{ij}^{\alpha\beta} w_{0j}^\beta \right) + \left( i_{p_{ij}}\left( w_{0i}^\alpha \right) + \sum_{\beta=1}^q h_{0ij}^{\alpha\beta} r_j^\beta - r_i^\alpha \right) dt \notag
\end{align}
This implies that $\delta\left(\left\{ \left( \overline{V_i} , r_i \right) \right\} \right) + G_0\left( \left\{ p_{ij} \right\}\right)= \left\{ \left( -\overline{R_{ij}' }+ \overline{R_{ij}}, - G_{ij}' + G_{ij} \right) \right\}$ and $\delta\left(\left\{ p_{ij} \right\} \right)= \left\{ d_{ijk}' - d_{ijk}\right\}$ (see \cite{Ser06}). This proves that $\alpha$ is cohomologous to $\beta$. So given a small extension $0\to (t)=(x^m) \to \tilde{A} \to k[x]/(x^{m+1})\to A= k[x]/(x^m) \to 0$ and an infinitesimal unfolding $\eta$ of $\left( X, \mathcal{N}_{\mathcal{F}_0}^* \right)$ over $A$, we can associate an element $o_\eta(e):=$ the cohomology class of $\alpha\in \mathbb{H}^2\left(X,  \tilde{\mathcal{N}}_{\mathcal{F}_0}^{*\bullet} \right)$. We note that $o_\eta(e)=0$ if and only if there exists $\left( \left\{\left( \overline{N_i}, M_i \right) \right\}, \left\{ d_{ij} \right\} \right)\in C^0\left(\mathcal{U}, \mathscr{H}om_{\mathcal{O}_X}\left( \mathcal{N}_{\mathcal{F}_0}^*, \frac{\Omega_X^1}{\mathcal{N}_{\mathcal{F}_0}^*} \right) \right) \bigoplus C^0\left( \mathcal{U}, \mathscr{H}om_{\mathcal{O}_X}\left( \mathcal{N}_{\mathcal{F}_0}^*, \mathcal{O}_X \right) \right) \bigoplus C^1\left( \mathcal{U}, T_X \right)$ such that $G_1\left( \overline{N_i}, M_i \right)=\left( B_i,  F_i \right)$ and $\left(\overline{N_j} , M_j \right)-\left(\overline{N_i}, M_i \right) + G_0\left( d_{ij} \right)= \left( \overline{- R_{ij}}, - G_{ij} \right)$ and $d_{ij}- d_{ik} + d_{jk}= d_{ijk}$. Then $\left(\tilde{f}_{ij}- t d_{ij} \right)\left(\tilde{f}_{jk}- t d_{jk} \right)= \tilde{f}_{ik} - td_{ik}$. Let $\Phi_{ij}:U_{ij}\times \textnormal{Spec}\left(\tilde{A}\right) \to  U_{ij} \times \textnormal{Spec}\left( \tilde{A}\right)$ be an isomorphism corresponding to $\tilde{f}_{ij} - t d_{ij}$. Then $\Phi_{jk}\Phi_{ij} = \Phi_{ik}$. We also note that
{\small{\begin{align}
&t\left( \mathcal{L}_{d_{ij}}\left(w_{0i}^\alpha \right) +\sum_{\beta=1}^q h_{0ij}^{\alpha\beta} N_j^\beta - N_i^\alpha +  \sum_{\beta=1}^q u_{ij}^{\alpha\beta} w_{0j}^\beta\right) + \left( i_{d_{ij}}\left( w_{0i}^\alpha \right) + \sum_{\beta=1}^q h_{0ij}^{\alpha\beta} M_j^\beta - M_i^\alpha \right) dt= - \tilde{\theta}_{ji}^* \left(\tilde{w}_i^\alpha + \tilde{v}_i^\alpha \right) + \sum_{\beta=1}^q \tilde{h}_{ij}^{\alpha\beta} \left( \tilde{w}_j^\beta + \tilde{v}_j^\beta \right) \notag \\
& \Longrightarrow \tilde{\theta}_{ji}^*\left( \tilde{w}_i^\alpha + \tilde{v}_i^\alpha \right) - tN_i^\alpha - M_i^\alpha dt + t \mathcal{L}_{d_{ij}}\left( w_{0i}^\alpha \right) + i_{d_{ij}}\left( w_{0i}^\alpha \right) dt  = \sum_{\beta=1}^q \left(\tilde{h}_{ij}^{\alpha\beta}- tu_{ij}^{\alpha\beta} \right)\left(\tilde{w}_j^\beta - t N_j^\beta - M_j^\beta dt \right)\,\,\,\textnormal{for some} \,\,\, u_{ij}^{\alpha\beta}\in \Gamma\left( U_{ij}, \mathcal{O}_X \right) \notag\\
&\Longrightarrow \Phi_{ji}^*\left( \tilde{w}_i^\alpha- tN_i^\alpha - M_i^\alpha dt \right) = \sum_{\beta=1}^q \left(\tilde{h}_{ij}^{\alpha\beta}- tu_{ij}^{\alpha\beta} \right)\left(\tilde{w}_j^\beta - t N_j^\beta - M_j^\beta dt \right)
\end{align}}}
and as in $\textnormal{(\ref{uu22})}$, we have $\left(\tilde{w}_i^1 + \tilde{v}_i^1 - t N_i^1 - M_i^1 dt \right) \wedge \cdots \wedge \left( \tilde{w}_i^q + \tilde{v}_i^q - t N_i^q - M_i^q dt \right) \wedge d \left( \tilde{w}_i^\alpha + \tilde{v}_i^\alpha - t N_i^\alpha - M_i^\alpha dt \right)=0,\alpha=1,...,q$. Hence a lifting of $\eta$ to $\tilde{A}= k[x]/(x^{m+1})$ exists. This completes the proof of Proposition \ref{uu35}.

\end{proof}

\subsection{Relationship between unfoldings of $\left( X, \mathcal{N}_{\mathcal{F}_0}^*\right)$ and foliated deformations of $\left(X, \mathcal{N}_{\mathcal{F}_0}^* \right)$} \

Let $\left( X, \mathcal{N}_{\mathcal{F}_0}^* \right)$ be a foliated nonsingular variety with $\mathcal{N}_{\mathcal{F}_0}^*$ locally free. Let $\textnormal{FDef}_{\left(X, \mathcal{N}_{\mathcal{F}_0}^* \right)}$ be the functor of Artin rings of flat foliated deformations of $\left(X, \mathcal{N}_{\mathcal{F}_0}^* \right)$ as in Part II. Then we have a canonical morphism of deformation functors
\begin{align*}
\textnormal{Unf}_{\left(X, \mathcal{N}_{\mathcal{F}_0}^* \right)} \to \textnormal{FDef}_{\left(X, \mathcal{N}_{\mathcal{F}_0}^* \right)}
\end{align*}
by mapping to the underlying flat foliated deformations, and we have a  canonical morphism of deformation functors
\begin{align*}
\textnormal{Unf}_{\left\{ X, \mathcal{N}_{\mathcal{F}_0}^*\right\}_{\textnormal{fixed}}} \to \textnormal{Unf}_{\left(X, \mathcal{N}_{\mathcal{F}_0}^* \right)}
\end{align*}
Then the relationship of three deformation functors $\textnormal{Unf}_{ \left\{ X, \mathcal{N}_{\mathcal{F}_0}^* \right\}} , \textnormal{Unf}_{\left( X, \mathcal{N}_{\mathcal{F}_0}^* \right)}$ and $\textnormal{FDef}_{\left(X,  \mathcal{N}_{\mathcal{F}_0}^* \right)}$ is expressed by the following exact sequence of complexes of sheaves
{\small{\begin{center}
$\begin{CD}
@. \cdots @. \cdots @. \cdots \\
@. @AAA @AG_3AA @AE_3AA \\
0@>>> \left( \mathcal{N}_{\mathcal{F}_0}^* \right)^* \otimes \tilde{\mathcal{S}}^2  @>>>  \left( \mathcal{N}_{\mathcal{F}_0}^*\right)^*\otimes \tilde{\mathcal{S}}^3 \bigoplus \left( \mathcal{N}_{\mathcal{F}_0}^* \right)^* \otimes \tilde{\mathcal{S}}^2 @>>>  \left( \mathcal{N}_{\mathcal{F}_0}^*\right)^*\otimes \tilde{\mathcal{S}}^3 @>>> 0 \\
@.@AAA @AG_2AA @AE_2AA \\
0@>>> \left( \mathcal{N}_{\mathcal{F}_0}^* \right)^*\otimes \tilde{\mathcal{S}}^1 @>>>  \left( \mathcal{N}_{\mathcal{F}_0}^* \right)^* \otimes \tilde{\mathcal{S}}^2  \bigoplus \left( \mathcal{N}_{\mathcal{F}_0}^* \right)^*\otimes \tilde{\mathcal{S}}^1 @>>>  \left( \mathcal{N}_{\mathcal{F}_0}^* \right)^* \otimes \tilde{\mathcal{S}}^2   @>>> 0\\
@. @AAA @AG_1AA @AE_1AA \\
0@>>> \left( \mathcal{N}_{\mathcal{F}_0}^* \right)^* @>>> \mathscr{H}om_{\mathcal{O}_X}\left( \mathcal{N}_{\mathcal{F}_0}^* , \frac{ \Omega_X^1 }{\mathcal{N}_{\mathcal{F}_0}^*} \right) \bigoplus \left( \mathcal{N}_{\mathcal{F}_0}^*\right)^* @>>> \mathscr{H}om_{\mathcal{O}_X}\left( \mathcal{N}_{\mathcal{F}_0}^* ,\frac{\Omega_X^1}{\mathcal{N}_{\mathcal{F}_0}^*} \right)@>>> 0 \\
@. @AAA @AG_0AA @AE_0AA \\
0@>>> 0 @>>> T_X @>>> T_X @>>> 0
\end{CD}$
\end{center}}}
Then the long exact sequence induces maps of first-order deformations
{\small{\begin{align*}
\textnormal{Unf}_{\left\{ X, \mathcal{N}_{\mathcal{F}_0}^*  \right\}}\left(k[\epsilon] \right) \cong  \mathbb{H}^0\left( X, \mathscr{H}om_{\mathcal{O}_X}\left( \mathcal{N}_{\mathcal{F}_0}^{*\bullet} , \mathcal{O}_X \right) \right) \to \textnormal{Unf}_{\left( X, \mathcal{N}_{\mathcal{F}_0}^* \right)} \left( k[\epsilon ] \right) \cong \mathbb{H}^1\left( X,  \tilde{\mathcal{N}}_{\mathcal{F}_0}^{*\bullet} \right)  \to \textnormal{FDef}_{\left( X, \mathcal{N}_{\mathcal{F}_0}^* \right)}\left( k[\epsilon] \right) \cong \mathbb{H}^1\left( X, \mathcal{N}_{\mathcal{F}_0}^{ * \bullet} \right)
\end{align*}}}
and obstruction maps for a small extension $e:0\to (t)=(x^m) \to \tilde{A}= k[x]/(x^{m+1)}\to A= k[x]/(x^m) \to 0$,
\begin{align*}
\mathbb{H}^1\left( X, \mathscr{H}om_{\mathcal{O}_X}\left( \mathcal{N}_{\mathcal{F}_0}^*, \mathcal{O}_X \right)^\bullet \right) \to \mathbb{H}^2\left(X, \tilde{\mathcal{N}}_{\mathcal{F}_0}^{*\bullet} \right) \to \mathbb{H}^2\left( X, \mathcal{N}_{\mathcal{F}_0}^{*\bullet} \right)
\end{align*}

\subsection{Unfoldings of $\left( X , \mathcal{N}_{\mathcal{F}_0}^*\right)$ with $\mathcal{N}_{\mathcal{F}_0}^*$ locally free in terms of vector bundles} \label{us5}\

Let $\left( X , \mathcal{N}_{\mathcal{F}_0}^* \right)$ be a foliated nonsingular variety with $\mathcal{N}_{\mathcal{F}_0}^*$ locally free.The underlying deformation of locally free sheaf $\mathcal{N}_{\mathcal{F}_0}^*$ is implicit in the unfolded dual leaf complex $\textnormal{(\ref{au2})}$. In order to make the underlying deformations of $\mathcal{N}_{\mathcal{F}_0}^*$ explicit, we define a complex of sheaves also controlling flat unfoldings of $\left(X, \mathcal{N}_{\mathcal{F}_0}^* \right)$ and present unfoldings of $\left(X, \mathcal{N}_{\mathcal{F}_0}^* \right)$ by using the complex. We shall construct

{\small{\begin{align*}
\tilde{\mathcal{E}}_{\mathcal{N}_{\mathcal{F}_0}^*}^\bullet : \mathcal{E}_{\mathcal{N}_{\mathcal{F}_0}^*} \xrightarrow{G_0''} \mathscr{H}om_{\mathcal{O}_X}\left( \mathcal{N}_{\mathcal{F}_0}^* , \Omega_X^1 \right) \bigoplus \mathscr{H}om_{\mathcal{O}_X}\left(  \mathcal{N}_{\mathcal{F}_0}^*, \mathcal{O}_X \right) \xrightarrow{ G_1''} \mathscr{H}om_{\mathcal{O}_X}\left( \mathcal{N}_{\mathcal{F}_0}^* , \tilde{\mathcal{S}}^2 \right) \bigoplus \mathscr{H}om_{\mathcal{O}_X}\left( \mathcal{N}_{\mathcal{F}_0}^*, \tilde{\mathcal{S}}^1 \right) \\
\xrightarrow{G_2} \mathscr{H}om_{\mathcal{O}_X}\left(  \mathcal{N}_{\mathcal{F}_0}^*, \tilde{\mathcal{S}}^3 \right) \bigoplus \mathscr{H}om_{\mathcal{O}_X} \left( \mathcal{N}_{\mathcal{F}_0}^*, \tilde{\mathcal{S}}^2 \right) \xrightarrow{G_3} \cdots
\end{align*}}}
where $\mathcal{E}_{\mathcal{N}_{\mathcal{F}_0}^*}$ is the Atiyah extension associated to $\mathcal{N}_{\mathcal{F}_0}^*$ (see Part II) and $G_0''$ is defined in the following way: let $\left\{ U_i \right\}$ be an affine open covering of $X$ such that $\Gamma\left(U_i,\mathcal{N}_{\mathcal{F}_0}^* \right)\cong \bigoplus^q \Gamma \left(U_i, \mathcal{O}_X\right)$ is trivialized in terms of the basis $\left\{w_{0i}^\alpha \right\}_{\alpha=1,...,q}$ with the relation $w_{0i}^\alpha= \sum_{\beta=1}^q h_{0ij}^{\alpha \beta} w_{0j}^\beta,\alpha=1,...,q$, and $\mathcal{E}_{\mathcal{N}_{\mathcal{F}_0}^*}|_{U_i}\cong \mathscr{H}om_{\mathcal{O}_X}\left( \mathcal{N}_{\mathcal{F}_0}^*, \mathcal{N}_{\mathcal{F}_0}^* \right)|_{U_i} \bigoplus T_X |_{U_i}$ as in Part II. Then on $U_i$, we define 
\begin{align*}
G_0''\left(\left( \phi_i, X\right) \right)\left(w_{0i}^\alpha \right)=\left( E_0'\left(\phi_i, X \right)\left(w_{0i}^\alpha \right), i_X\left( w_{0i}^\alpha \right) \right)
\end{align*}
and then linearly extends to $\Gamma\left( U_i, \mathcal{N}_{\mathcal{F}_0}^* \right)$ (for $E_0'$, see Part II). Then $G_0''$ is well-defined, and $G_1''\circ G_0''=0$. We will denote the $i$-th cohomology group by $\mathbb{H}^i\left( X, \tilde{\mathcal{E}}_{\mathcal{N}_{\mathcal{F}_0}^*}^{\bullet} \right)$. We have an exact sequence of complexes of sheaves
{\small{\begin{center}
$\begin{CD}
@. \cdots @. \cdots @. \cdots \\
@. @AAA @AG_3AA @AG_3AA \\
0@>>>0 @>>>  \left( \mathcal{N}_{\mathcal{F}_0}^*\right)^*\otimes \tilde{\mathcal{S}}^3 \bigoplus \left( \mathcal{N}_{\mathcal{F}_0}^* \right)^* \otimes \tilde{\mathcal{S}}^2 @>>>  \left( \mathcal{N}_{\mathcal{F}_0}^*\right)^*\otimes \tilde{\mathcal{S}}^3 \bigoplus \left( \mathcal{N}_{\mathcal{F}_0}^* \right)^* \otimes \tilde{\mathcal{S}}^2 @>>> 0 \\
@.@AAA @AG_2AA @AG_2AA \\
0@>>>0 @>>>  \left( \mathcal{N}_{\mathcal{F}_0}^* \right)^* \otimes \tilde{\mathcal{S}}^2  \bigoplus \left( \mathcal{N}_{\mathcal{F}_0}^* \right)^*\otimes \tilde{\mathcal{S}}^1 @>>>  \left( \mathcal{N}_{\mathcal{F}_0}^* \right)^* \otimes \tilde{\mathcal{S}}^2  \bigoplus \left( \mathcal{N}_{\mathcal{F}_0}^* \right)^*\otimes \tilde{\mathcal{S}}^1 @>>> 0\\
@. @AAA @AG_1''AA @AG_1AA \\
0@>>> \mathscr{H}om_{\mathcal{O}_X}\left( \mathcal{N}_{\mathcal{F}_0}^*, \mathcal{N}_{\mathcal{F}_0}^* \right) @>>> \mathscr{H}om_{\mathcal{O}_X}\left( \mathcal{N}_{\mathcal{F}_0}^* , \Omega_X^1 \right) \bigoplus \left( \mathcal{N}_{\mathcal{F}_0}^*\right)^* @>>> \mathscr{H}om_{\mathcal{O}_X}\left( \mathcal{N}_{\mathcal{F}_0}^* ,\frac{\Omega_X^1}{\mathcal{N}_{\mathcal{F}_0}^*} \right) \bigoplus \left( \mathcal{N}_{\mathcal{F}_0}^*\right)^*  @>>> 0 \\
@. @A-\textnormal{Id}AA @AG_0''AA @AAA \\
0@>>> \mathscr{H}om_{\mathcal{O}_X}\left(\mathcal{N}_{\mathcal{F}_0}^*, \mathcal{N}_{\mathcal{F}_0}^* \right) @>>> \mathcal{E}_{\mathcal{N}_{\mathcal{F}_0}^*} @>>> T_X @>>> 0
\end{CD}$
\end{center}}}
The above exact sequence implies that
\begin{align*}
\mathbb{H}^i\left(X, \tilde{\mathcal{E}}_{\mathcal{N}_{\mathcal{F}_0}^*}^\bullet \right) \cong \mathbb{H}^i\left( X, \tilde{\mathcal{N}}_{\mathcal{F}_0}^{*\bullet}\right)
\end{align*}
By combining Proposition $\textnormal{\ref{uu35}}$ and the related Proposition in Part II, we obtain
\begin{proposition}
Let $\left( X, \mathcal{N}_{\mathcal{F}_0}^* \right)$ be a foliated nonsingular variety with $\mathcal{N}_{\mathcal{F}_0}^*$ locally free.
\begin{enumerate}
\item there is a $1-1$ correspondence
\begin{align*}
\textnormal{Unf}_{\left(X, \mathcal{N}_{\mathcal{F}_0}^* \right)}\left(k[\epsilon] \right)\cong \mathbb{H}^1\left(X, \tilde{\mathcal{E}}_{\mathcal{N}_{\mathcal{F}_0}^*}^\bullet \right)
\end{align*}
\item Let $\eta=\left( \mathcal{X},  \tilde{\mathcal{N}}_\mathcal{F}^* \right)$ be an infinitesimal unfoldings of $\left( X, \mathcal{N}_{\mathcal{F}_0}^* \right)$ over $A:= k[x] / \left(x^m \right)$. To  small extension $0\to (t)=(x^m) \to \tilde{A} = k[x]/(x^{m+1}) \to A = k[x]/(x^m) \to 0$, we can associate an element $o_\eta(e)\in \mathbb{H}^2\left( X,  \tilde{\mathcal{E}}_{ \mathcal{N}_{\mathcal{F}_0}^*}^\bullet \right)$ called the obstruction lifting of $\eta$ to $\tilde{A}$, which is $0$ if and only if a lifting of $\eta$ to $\tilde{A}$ exists.
\end{enumerate}
\end{proposition}

\subsection{Relationship between unfoldings of $\left( X, \mathcal{N}_{\mathcal{F}_0}^* \right)$ in terms of vector bundles and foliated deformations of $\left( X, \mathcal{N}_{\mathcal{F}_0}^* \right)$ in terms of vector bundles}\

Let $\left( X, \mathcal{N}_{\mathcal{F}_0}^* \right)$ be a foliated nonsingular variety with $\mathcal{N}_{\mathcal{F}_0}^*$ locally free. Then the relationship between unfoldings of $\left(X, \mathcal{N}_{\mathcal{F}_0}^* \right)$ and the underlying foliated deformations of $\left( X, \mathcal{N}_{\mathcal{F}_0}^* \right)$ in terms of vector bundles in Part II is described by the following exact sequence of complexes of sheaves
{\Small{\begin{center}
$\begin{CD}
@. \cdots @. \cdots @. \cdots \\
@. @AAA @AAA @AAA \\
0 @>>> \mathscr{H}om_{\mathcal{O}_X}\left( \mathcal{N}_{\mathcal{F}_0}^*, \tilde{\mathcal{S}}^2 \right) @>>> \mathscr{H}om_{\mathcal{O}_X}\left( \mathcal{N}_{\mathcal{F}_0}^*, \tilde{\mathcal{S}}^3 \right) \bigoplus  \mathscr{H}om_{\mathcal{O}_X}\left( \mathcal{N}_{\mathcal{F}_0}^*, \tilde{\mathcal{S}}^2 \right) @>>>  \mathscr{H}om_{\mathcal{O}_X}\left( \mathcal{N}_{\mathcal{F}_0}^*, \tilde{\mathcal{S}}^3 \right) @>>> 0  \\
@. @AAA @AAA @AAA \\
0 @>>> \mathscr{H}om_{\mathcal{O}_X}\left( \mathcal{N}_{\mathcal{F}_0}^*, \tilde{\mathcal{S}}^1 \right) @>>> \mathscr{H}om_{\mathcal{O}_X}\left( \mathcal{N}_{\mathcal{F}_0}^*, \tilde{\mathcal{S}}^2 \right) \bigoplus  \mathscr{H}om_{\mathcal{O}_X}\left( \mathcal{N}_{\mathcal{F}_0}^*, \tilde{\mathcal{S}}^1 \right) @>>>  \mathscr{H}om_{\mathcal{O}_X}\left( \mathcal{N}_{\mathcal{F}_0}^*, \tilde{\mathcal{S}}^2 \right) @>>> 0  \\
@. @AAA @AAA @AAA \\
0 @>>> \mathscr{H}om_{\mathcal{O}_X} \left(\mathcal{N}_{\mathcal{F}_0}^*, \mathcal{O}_X \right) @>>>  \mathscr{H}om_{\mathcal{O}_X}\left( \mathcal{N}_{\mathcal{F}_0}^*, \Omega_X^1 \right) \bigoplus \mathscr{H}om_{\mathcal{O}_X}\left( \mathcal{N}_{\mathcal{F}_0}^*, \mathcal{O}_X \right) @>>> \mathscr{H}om_{\mathcal{O}_X}\left( \mathcal{N}_{\mathcal{F}_0}^*, \Omega_X^1 \right) @>>> 0 \\
@. @AAA @AAA @AAA \\
0 @>>> 0 @>>> \mathcal{E}_{\mathcal{N}_{\mathcal{F}_0}^*} @>>> \mathcal{E}_{\mathcal{N}_{\mathcal{F}_0}^*} @>>> 0 
\end{CD}$
\end{center}}}
The long exact sequence induces the morphism between first-order unfoldings and first-order foliated deformations in terms of cotangent sheaves
\begin{align*}
\mathbb{H}^1\left( X, \tilde{\mathcal{E}}_{\mathcal{N}_{\mathcal{F}_0}^*}^\bullet \right) \to \mathbb{H}^1\left( X, \mathcal{E}_{\mathcal{N}_{\mathcal{F}_0}^*}^\bullet \right)
\end{align*}
and the morphism of obstructions between unfoldings and the underlying foliated deformations for a small extension $0\to (t)= (x^m) \to k[x]/(x^{m+1})\to k[x]/(x^m) \to 0$
\begin{align*}
\mathbb{H}^2\left( X, \tilde{\mathcal{E}}_{\mathcal{N}_{\mathcal{F}_0}^*}^\bullet \right) \to \mathbb{H}^2\left( X, \mathcal{E}_{\mathcal{N}_{\mathcal{F}_0}^*}^\bullet \right)
\end{align*}

\section{Unfoldings of $\left(X, \mathcal{N}_{\mathcal{F}_0}^* \right)$ with $\mathcal{N}_{\mathcal{F}_0}^*$ in terms of deformations of sections}

Let $\left( X , \mathcal{N}_{\mathcal{F}_0}^* \right)$ be a foliated nonsingular variety with $\mathcal{N}_{\mathcal{F}_0}^*$ locally free of rank $q$. The exact sequence $0 \to \mathcal{N}_{\mathcal{F}_0}^* \to \Omega_X^1 \to \frac{\Omega_X^1}{\mathcal{N}_{\mathcal{F}_0}^*} \to 0$ induces a section $w_0\in  H^0\left( X, \bigwedge^2 \Omega_X^1\otimes \mathcal{L}_0 \right)$ where $\mathcal{L}_0= \left(\bigwedge^q \mathcal{N}_{\mathcal{F}_0}^*\right)^*$. An infinitesimal unfolding $\left( \mathcal{X}, \tilde{\mathcal{N}}_\mathcal{F}^* \right)$ over $A\in \bold{Art}$ defined by the exact sequence $0 \to \tilde{\mathcal{N}}_\mathcal{F}^* \to \Omega_{\mathcal{X}/k}^1 \to \frac{\Omega_{\mathcal{X}/k}^1}{\tilde{\mathcal{N}}_\mathcal{F}^*}\to 0$ induces a section $\tilde{w}\in H^0\left(\mathcal{X}, \Omega_{\mathcal{X}/k}^1\otimes \mathcal{L} \right)$, where $\mathcal{L}=\left(\bigwedge^q \tilde{\mathcal{N}}_\mathcal{F}^*\right)^* \cong \left( \bigwedge^q \mathcal{N}_\mathcal{F}^* \right)^*$ where $\mathcal{N}_\mathcal{F}^*$ is a subsheaf of $\Omega_{\mathcal{X}/A}^1$ such that $\mathcal{N}_\mathcal{F}^*=\alpha\left( \tilde{\mathcal{N}}_\mathcal{F}^* \right)\cong \tilde{\mathcal{N}}_{\mathcal{F}}^*$ under the canonical map $\alpha: \Omega_{\mathcal{X}/k}^1 \to \Omega_{\mathcal{X}/A}^1$. Then $\tilde{w}$ can be considered to be an infinitesimal unfolding of a section $w_0$ and an infinitesimal deformation of a pair of the nonsingular variety $X$ and the line bundle $\mathcal{L}_0$. We will describe infinitesimal unfoldings of $\left(X, \mathcal{N}_{\mathcal{F}_0}^* \right)$ in terms of unfoldings of the section $w_0$ and infinitesimal deformations of the pair of $\left(X, \mathcal{L}_0 \right)$. We shall construct a complex of sheaves
 {\small{\begin{align*}
\tilde{\mathcal{E}}_{\mathcal{L}_0}^\bullet : \mathcal{E}_{\mathcal{L}_0}\xrightarrow{\tilde{G}_0} \mathcal{S}^1\bigoplus \left( \mathcal{N}_{\mathcal{F}_0}^*,  \mathcal{O}_X \right) \xrightarrow{\tilde{G}_1} \mathscr{H}om_{\mathcal{O}_X}\left( \mathcal{N}_{\mathcal{F}_0}^*, \tilde{\mathcal{S}}^2 \right) \bigoplus \mathscr{H}om_{\mathcal{O}_X}\left( \mathcal{N}_{\mathcal{F}_0}^*, \tilde{\mathcal{S}}^1  \right) \\
\xrightarrow{\tilde{G}_2} \mathscr{H}om_{\mathcal{O}_X}\left( \mathcal{N}_{\mathcal{F}_0}^*, \tilde{\mathcal{S}}^3 \right) \bigoplus \mathscr{H}om_{\mathcal{O}_X}\left( \mathcal{N}_{\mathcal{F}_0}^*, \tilde{\mathcal{S}}^2  \right) \xrightarrow{\tilde{G}_3} \cdots
\end{align*}}}
where $\mathcal{E}_{\mathcal{L}_0}$ is the Atiyah extension associated to $\mathcal{L}_0$ (see Part II) and $\tilde{G}_0$ is defined in the following way: let $\left\{ U_i \right\}$ be an affine open covering of $X$ such that $\Gamma\left(U_i,\mathcal{N}_{\mathcal{F}_0}^* \right)\cong \bigoplus^q \Gamma \left(U_i, \mathcal{O}_X\right)$ is trivialized in terms of the basis $\left\{w_{0i}^\alpha \right\}_{\alpha=1,...,q}$ with the relation $w_{0i}^\alpha= \sum_{\beta=1}^q h_{0ij}^{\alpha \beta} w_{0j}^\beta,\alpha=1,...,q$, and $\mathcal{E}_{\mathcal{L}_0}|_{U_i}\cong \mathcal{O}_X|_{U_i} \bigoplus T_X |_{U_i}$ as in Part II. Then on $U_i$, we define
\begin{align*}
\tilde{G}_0\left( a_i, b_i \right)= \left( \tilde{E}_0\left( a_i, b_i \right), G_0'(b_i)  \right)
\end{align*}
(for $\tilde{E}_0$, see Part II). Then $\tilde{G}_0$ is well-defined. We define $\tilde{G}_1$ by the composition of $\tilde{\mathcal{S}}^1 \bigoplus \mathscr{H}om_{\mathcal{O}_X}\left( \mathcal{N}_{\mathcal{F}_0}^* , \mathcal{O}_X \right) \to \mathscr{H}om_{\mathcal{O}_X}\left( \mathcal{N}_{\mathcal{F}_0}^*, \frac{\Omega_X^1}{\mathcal{N}_{\mathcal{F}_0}^*} \right)\bigoplus \mathscr{H}om_{\mathcal{O}_X}\left( \mathcal{N}_{\mathcal{F}_0}^*, \mathcal{O}_X \right)$ with $G_1$. Then $\tilde{\mathcal{E}}_{\mathcal{L}_0}^\bullet$ defines a complex of sheaves. We will denote the $i$-th cohomology group by $\mathbb{H}^i\left( X , \tilde{\mathcal{E}}_{\mathcal{L}_0}^\bullet \right)$. We have an exact sequence of complexes of sheaves
{\Small{\begin{center}
$\begin{CD}
@. \cdots @. \cdots @. \cdots \\
@. @AAA @AG_3AA @AG_3 AA \\
0 @>>> 0 @>>> \mathscr{H}om_{\mathcal{O}_X}\left( \mathcal{N}_{\mathcal{F}_0}^*, \tilde{\mathcal{S}}^3 \right) \bigoplus  \mathscr{H}om_{\mathcal{O}_X}\left( \mathcal{N}_{\mathcal{F}_0}^*, \tilde{\mathcal{S}}^2 \right) @>>> \mathscr{H}om_{\mathcal{O}_X}\left(\mathcal{N}_{\mathcal{F}_0}^*, \tilde{\mathcal{S}}^3 \right) \bigoplus  \mathscr{H}om_{\mathcal{O}_X}\left( \mathcal{N}_{\mathcal{F}_0}^*, \tilde{\mathcal{S}}^2 \right) @>>> 0 \\
@. @AAA @AG_2AA @AG_2AA \\
0@>>> 0 @>>> \mathscr{H}om_{\mathcal{O}_X}\left(\mathcal{N}_{\mathcal{F}_0}^*, \tilde{\mathcal{S}}^2 \right) \bigoplus \mathscr{H}om_{\mathcal{O}_X}\left( \mathcal{N}_{\mathcal{F}_0}^*, \tilde{\mathcal{S}}^1 \right) @>>> \mathscr{H}om_{\mathcal{O}_X}\left( \mathcal{N}_{\mathcal{F}_0}^*, \tilde{\mathcal{S}}^2 \right) \bigoplus \mathscr{H}om_{\mathcal{O}_X}\left( \mathcal{N}_{\mathcal{F}_0}^* , \tilde{\mathcal{S}}^1 \right) @>>> 0\\
@. @AAA @A\tilde{G}_1AA @AG_1AA \\
0@>>> \mathcal{O}_X @>>> \mathcal{S}^1\bigoplus \mathscr{H}om_{\mathcal{O}_X}\left( \mathcal{N}_{\mathcal{F}_0}^*, \mathcal{O}_X \right) @>\alpha >> \mathscr{H}om_{\mathcal{O}_X}\left( \mathcal{N}_{\mathcal{F}_0}^*, \frac{\Omega_X^1}{ \mathcal{N}_{\mathcal{F}_0}^* } \right) \bigoplus \mathscr{H}om_{\mathcal{O}_X}\left( \mathcal{N}_{\mathcal{F}_0}^*, \mathcal{O}_X \right) @>>> 0 \\
@. @A\textnormal{Id} AA @A\tilde{G}_0 AA @AG_0AA \\
0 @>>> \mathcal{O}_X @>>> \mathcal{E}_{\mathcal{L}_0}  @>>> T_X @>>> 0 
\end{CD}$
\end{center}}}
This implies that we have
\begin{align*}
\mathbb{H}^i\left(X, \tilde{\mathcal{E}}_{\mathcal{L}_0}^\bullet \right) \cong \mathbb{H}^1\left( X, \tilde{\mathcal{N}}_{\mathcal{F}_0}^{*\bullet} \right)
\end{align*}

By combining Proposition \ref{uu35} and the related Proposition in Part II, we obtatin

\begin{proposition}
Let $\left( X, \mathcal{N}_{\mathcal{F}_0}^* \right)$ be a foliated nonsingular variety with $\mathcal{N}_{\mathcal{F}_0}^*$ locally free.
\begin{enumerate}
\item there is a $1-1$ correspondence
\begin{align*}
\textnormal{Unf}_{\left(X, \mathcal{N}_{\mathcal{F}_0}^*\right)}\left(k[\epsilon] \right)\cong \mathbb{H}^1\left(X, \tilde{\mathcal{E}}_{\mathcal{L}_0}^\bullet \right)
\end{align*}
\item Let $\eta=\left( \mathcal{X}, \tilde{\mathcal{N}}_\mathcal{F}^* \right)$ be an infinitesimal unfolding of $\left( X, \mathcal{N}_{\mathcal{F}_0}^* \right)$ over $A:= k[x] / \left(x^m\right) $. To  small extension $0\to (t)= \left(x^m\right) \to \tilde{A}=k[x]/ \left(x^{m+1}\right) \to A=k[x]/\left(x^m\right) \to 0$, we can associate an element $o_\eta(e)\in \mathbb{H}^2\left( X, \tilde{\mathcal{E}}_{\mathcal{L}_0}^\bullet \right)$ called the obstruction lifting of $\eta$ to $\tilde{A}$, which is $0$ if and only if a lifting of $\eta$ to $\tilde{A}$ exists.
\end{enumerate}
\end{proposition}

\subsection{Relationship between unfoldings of $\left( X, \mathcal{N}_{\mathcal{F}_0}^* \right)$ in terms of unfoldings of sections and foliated deformations of $\left(X, \mathcal{N}_{\mathcal{F}_0}^* \right)$ in terms of deformations of  sections}\

Let $\left( X, \mathcal{N}_{\mathcal{F}_0}^* \right)$ be a foliated nonsingular variety with $\mathcal{N}_{\mathcal{F}_0}^*$ locally free of rank $q$, which induces a section $w_0\in H^0\left( X, \bigwedge^q\Omega_X^1\otimes \mathcal{L}_0 \right)$. where $\mathcal{L}_0= \left( \bigwedge^q \mathcal{N}_{\mathcal{F}_0}^* \right)^*$. Let $\left( \mathcal{X}, \tilde{\mathcal{N}}_{\mathcal{F}}^* \right)$ be an unfolding of $\left( X, \mathcal{N}_{\mathcal{F}_0}^* \right)$ over $A\in \bold{Art}$, which induces a section $\tilde{w}\in H^0\left(\mathcal{X}, \Omega_{\mathcal{X}/k}^1\otimes \mathcal{L} \right)$, where $\mathcal{L}=\left( \bigwedge^q  \tilde{\mathcal{N}}_\mathcal{F}^* \right)^*$. We note that $\tilde{w}$ induces a section $w\in H^0\left( \mathcal{X}, \bigwedge^q \Omega_{\mathcal{X}/A}^1\otimes \mathcal{L} \right)$ by the canonical map $\alpha:\Omega_{\mathcal{X}/k}^1 \to \Omega_{\mathcal{X}/A}^1$, and $\tilde{\mathcal{N}}_\mathcal{F}^*$ is isomorphic to the image $\mathcal{N}_\mathcal{F}^*$ of $\tilde{\mathcal{N}}_\mathcal{F}^*$ under $\alpha$, so that $w\in H^0\left( \mathcal{X}, \bigwedge^q \Omega_{\mathcal{X}/A}^1 \otimes \mathcal{L} \right)$ defines a foliated deformation of $w_0$ over $A$ in terms of deformations of sections in Part II. Then the unfoldings $\tilde{w}$ of $w_0$ in terms of sections and the underlying foliated deformation $w$ of $w_0$ in terms of sections are described by the following exact sequence of complexes of sheaves:
{\Small{\begin{center}
$\begin{CD}
@. \cdots @. \cdots @. \cdots \\
@. @AAA @AG_3AA @AE_3AA \\
0 @>>> \mathscr{H}om_{\mathcal{O}_X}\left( \mathcal{N}_{\mathcal{F}_0}^*, \tilde{\mathcal{S}}^2 \right) @>>> \mathscr{H}om_{\mathcal{O}_X}\left( \mathcal{N}_{\mathcal{F}_0}^*, \tilde{\mathcal{S}}^3 \right) \bigoplus \mathscr{H}om_{\mathcal{O}_X}\left( \mathcal{N}_{\mathcal{F}_0}^*, \tilde{\mathcal{S}}^2 \right) @>>> \mathscr{H}om_{\mathcal{O}_X}\left( \mathcal{N}_{\mathcal{F}_0}^*, \tilde{\mathcal{S}}^3 \right) @>>> 0 \\
@. @AAA @AG_2AA @AE_2AA \\
0 @>>> \mathscr{H}om_{\mathcal{O}_X}\left( \mathcal{N}_{\mathcal{F}_0}^*, \tilde{\mathcal{S}}^1 \right) @>>> \mathscr{H}om_{\mathcal{O}_X}\left( \mathcal{N}_{\mathcal{F}_0}^*, \tilde{\mathcal{S}}^2 \right) \bigoplus \mathscr{H}om_{\mathcal{O}_X}\left( \mathcal{N}_{\mathcal{F}_0}^*, \tilde{\mathcal{S}}^1 \right) @>>> \mathscr{H}om_{\mathcal{O}_X}\left( \mathcal{N}_{\mathcal{F}_0}^*, \tilde{\mathcal{S}}^2 \right) @>>> 0 \\
@. @AAA @A\tilde{G}_1AA @A\tilde{E}_1AA \\
0 @>>>  \mathscr{H}om_{\mathcal{O}_X}\left( \mathcal{N}_{\mathcal{F}_0}^*, \mathcal{O}_X \right) @>>> \mathcal{S}^1\bigoplus \mathscr{H}om_{\mathcal{O}_X}\left( \mathcal{N}_{\mathcal{F}_0}^*, \mathcal{O}_X \right) @>>> \mathcal{S}^1 @>>> 0 \\
@. @AAA @A\tilde{G}_0AA @A\tilde{E}_0AA \\
0 @>>> 0 @>>> \mathcal{E}_{\mathcal{L}_0} @>>> \mathcal{E}_{\mathcal{L}_0} @>>> 0
\end{CD}$
\end{center}}}
The long exact sequence induces the morphism between first-order unfolding in terms of sections and first-order foliated deformations in terms of sections
\begin{align*}
\mathbb{H}^1\left( X, \tilde{\mathcal{E}}_{\mathcal{L}_0}^\bullet  \right) \to \mathbb{H}^1\left( X,  \mathcal{E}_{\mathcal{L}_0}^\bullet \right)
\end{align*}
and morphism of obstructions in terms of sections for a small extension $0\to (t)=(x^m) \to k[x]/(x^{m+1})\to k[x]/(x^m) \to 0$
\begin{align*}
\mathbb{H}^2\left( X, \tilde{\mathcal{E}}_{\mathcal{L}_0}^\bullet \right) \to \mathbb{H}^2\left( X, \mathcal{E}_{\mathcal{L}_0}^\bullet \right)
\end{align*}

\section{Rigid unfoldings and stability for one-parameter Artin rings}

\subsection{Rigid unfoldings over $k[x]/ \left(x^m\right)$ } \

We recall that a foliated nonsingular variety $\left(X, \mathcal{N}_{\mathcal{F}_0}^* \right)$ with $\mathcal{N}_{\mathcal{F}_0}^*$ locally free is rigid in unfoldings if every infinitesimal unfolding $\left( \mathcal{X}, \mathcal{N}_\mathcal{F}^* \right)$ of $\left( X, \mathcal{N}_{\mathcal{F}_0}^* \right)$ over $A$ is trivial for every $A$ in $\bold{Art}$, in other words, $\left( \mathcal{X}, \tilde{\mathcal{N}}_\mathcal{F}^* \right) \cong \left( X\times \textnormal{Spec}\left( A \right), \mathcal{N}_{\mathcal{F}_0}^* \right)$. In the light of Theorem \ref{uc1} of completeness for unfoldings of a compact foliated complex manifold $\left(X, \mathcal{N}_{\mathcal{F}_0}^* \right)$ which implies the rigidity under the assumption $\mathbb{H}^1\left( X , \tilde{\mathcal{N}}_{\mathcal{F}_0}^{*\bullet} \right)=0$ as a corollary (we will give a direct proof in Part IV-2), the author believes that the following Proposition \ref{a1} holds for any $A\in \bold{Art}$. However the author could not prove for arbitrary $A\in \bold{Art}$ by using the same method employed in foliated deformations of $\left(X, \mathcal{N}_{\mathcal{F}_0}^* \right)$ which is based on the transitive action of $\mathbb{H}^1\left( X, \mathcal{N}_{\mathcal{F}_0}^{\bullet} \right)$ on the liftings of a foliated deformations for a small extension $0\to (t)\to \tilde{A} \to A \to 0$. Our method in the proof of Proposition \ref{a1} is based on the transitive action of $\mathbb{H}^1\left( X , \tilde{\mathcal{N}}_{\mathcal{F}_0}^{* \bullet} \right)$ on the liftings of an unfolding for a small extension of the form $0\to (t)=(x^m)\to k[x]/(x^{m+1})\to k[x]/(x^m)\to 0$.

\begin{proposition}\label{a1}
A foliated nonsingular variety $\left(X, \mathcal{N}_{\mathcal{F}_0}^* \right)$ with $\mathcal{N}_{\mathcal{F}_0}^*$ is rigid in unfoldings over $k[x]/(x^m)$ for any $m\geq  1$, i.e. any unfolding of $\left( X, \mathcal{N}_{\mathcal{F}_0}^* \right)$ over $k[x]/(x^m)$ is trivial, if and only if $\mathbb{H}^1\left( X, \tilde{\mathcal{N}}_{\mathcal{F}_0}^{*\bullet} \right)=0$
\end{proposition}

\begin{proof}
Assume that $\left(X, \mathcal{N}_{\mathcal{F}_0}^* \right)$ is rigid in unfoldings over $k[x]/(x^m)$. Then in particular any first-order unfolding of $\left( X, \mathcal{N}_{\mathcal{F}_0}^* \right)$ is trivial, so that $\mathbb{H}^1\left( X, \tilde{\mathcal{N}}_{\mathcal{F}_0}^{*\bullet} \right)=0$ by Proposition $\textnormal{\ref{uu35}}$. Assume that $\mathbb{H}^1\left( \right)+0$. First we claim that given an infinitesimal unfolding $\eta$ of $\left(X , \mathcal{N}_{\mathcal{F}_0}^* \right)$ over $k[x]/(x^m)$ and a small extension $0\to (t)=(x^m)\to \tilde{A}=k[x]/(x^{m+1})\to A=k[x]/(x^m)\to 0$, any two liftings $\xi, \xi'$ of $\eta$ to $\tilde{A}=k[x]/(x^{m+1})$ are equivalent. Let $\mathcal{U}=\left\{ U_i \right\}$ be an affine open covering of $X$ such that $\xi=\left(\mathcal{X}, \tilde{\mathcal{N}}_\mathcal{F}^* \right)$ is represented by the following data: we have unfolded isomorphisms $\theta_i:\left( U_i\times \textnormal{Spec}\left( \tilde{A}\right) , \tilde{\mathcal{N}}_{\mathcal{F}_i}^* \right)\to \left( \mathcal{X}|_{U_i},  \tilde{\mathcal{N}}_\mathcal{F}^*|_{U_i} \right)$ where $\tilde{\mathcal{N}}_{\mathcal{F}_i}^*$ is the unfolded structure on $U_i\times \textnormal{Spec}\left( \tilde{A} \right)$ induced from $\tilde{\mathcal{N}}_\mathcal{F}^*$ via $\theta$, and $\Gamma\left( U_i, \tilde{\mathcal{N}}_{\mathcal{F}_i}^* \right)$ is generated by $\omega_i^1=w_i^1+ v_i^1,..., \omega_i^q=w_i^q+ v_i^q\in \left(\Gamma\left(U_i, \Omega_X^1 \right)\otimes \tilde{A} \right)\bigoplus \left(\Gamma\left( U_i, \mathcal{O}_X \right)\otimes \Omega_{\tilde{A}/k}^1\right)$ with $\omega_i^1 \wedge \cdots \wedge \omega_i^q\wedge d\omega_i^\alpha=0$, and letting $\theta_{ij}= \theta_j^{-1}\theta_i$, we have unfolded isomorphisms $\theta_{ij}:\left( U_{ij}\times \textnormal{Spec}\left( \tilde{A} \right), \tilde{\mathcal{N}}_{\mathcal{F}_i}^* \right)\to \left( U_{ij}\times \textnormal{Spec}\left( \tilde{A}\right), \tilde{\mathcal{N}}_{\mathcal{F}_j}^* \right)$ such that $\theta_{ji}^* \omega_i^\alpha= \sum_{\beta=1}^q h_{ij}^{\alpha\beta} \omega_j^\beta$ for $h_{ij}^{\alpha\beta}\in \Gamma\left( U_{ij}, \mathcal{O}_X \right)\otimes \tilde{A}$, and on the other hand, $\xi'=\left(\mathcal{X}', \tilde{\mathcal{N}}_{\mathcal{F}'}^*  \right)$ is represented by the following data: similarly $\theta_i', \theta_{ij}'=\theta_j'^{-1}\theta_i', \left\{ \tilde{\mathcal{N}}_{\mathcal{F}_i'}^*\right\}, \omega_i'^\alpha=w_i'^\alpha + v_i'^\alpha, h_{ij}^{\alpha\beta}$. We may assume that $\omega_i^\alpha$ and $\omega_i'^\alpha$ induce the same element in $\left(\Gamma\left(U_i, \Omega_X^1\right)\otimes A\right) \bigoplus \left(\Gamma\left( U_i, \mathcal{O}_X \right)\otimes \Omega_{A/k}^1\right)$. Let $f_{ij}:\left( \Gamma\left( U_{ij}, \mathcal{O}_X \right), \tilde{\mathcal{N}}_{\mathcal{F}_i}^* \right)\to \left( \Gamma\left( U_{ij}, \mathcal{O}_X \right)\otimes \tilde{A}, \tilde{\mathcal{N}}_{\mathcal{F}_i}^* \right)$ be the ring homomorphism corresponding to $\theta_{ij}$ and similarly $f_{ij}'$ corresponding to $\theta_{ij}'$. Since $\xi, \xi'$ induce $\eta$ over $A$, we have $f_{ij}'= f_{ij} + t p_{ij}$ for some $p_{ij}\in \Gamma\left( U_{ij}, T_X \right)$ and $\omega_i'^\alpha = \omega_i^\alpha + t B_i^\alpha + C_i^\alpha dt$ for some $B_i^\alpha\in \Gamma\left( U_i, \Omega_X^1\right)$ and some $C_i^\alpha\in \Gamma\left(U_i, \mathcal{O}_X \right),\alpha=1,...,q$. We define $\overline{B_i}:\Gamma\left( U_i, \mathcal{N}_{\mathcal{F}_0}^* \right)\to \Gamma\left( U_i, \frac{\Omega_X^1}{\mathcal{N}_{\mathcal{F}_0}^*} \right)$ by $w_{0i}^\alpha \mapsto \overline{B_i^\alpha}$, and we define $C_i:\Gamma\left( U_i, \mathcal{N}_{\mathcal{F}_0}^* \right)\to \Gamma\left( U_i, \mathcal{O}_X \right)$ by $w_{0i}^\alpha \mapsto C_i^\alpha$. Then we claim that $\left( \left\{ \overline{B_i}\right\}, \left\{ C_i \right\}, \left\{ p_{ij}\right\} \right)\in C^0\left( U_i, \mathscr{H}om_{\mathcal{O}_X}\left(  \mathcal{N}_{\mathcal{F}_0}^*, \frac{\Omega_X^1}{\mathcal{N}_{\mathcal{F}_0}^*} \right) \right) \bigoplus C^0\left( \mathcal{U}, \mathscr{H}om_{\mathcal{O}_X}\left( \mathcal{N}_{\mathcal{F}_0}^*, \mathcal{O}_X \right) \right) \bigoplus C^1\left( \mathcal{U}, T_X \right)$ defines a $1$-cocycle in the \v Cech resolution of $\tilde{\mathcal{N}}_{\mathcal{F}_0}^{*\bullet}$. In fact, we have $p_{ij}- p_{ik}+ p_{jk}=0$, and
\begin{align*}
0&=\omega_1'\wedge \cdots \wedge \omega_i^q \wedge d\omega_i^\alpha= \left(\omega_i^1+ t B_i^1 + C_i^1 dt \right) \wedge \cdots \wedge \left( \omega_i^q + t B_i^q + C_i^q dt\right) \wedge d \left( \omega_i^\alpha + t B_i^\alpha + C_i^\alpha dt \right)\\
&= \omega_i^1 \wedge \cdots \wedge \omega_i^q \wedge d \omega_i^\alpha + t \sum_{\beta=1}^q w_{0i}^1 \wedge \cdots \wedge B_i^\beta \wedge \cdots \wedge w_{0i}^q \wedge d w_{0i}^\alpha + t w_{0i}^1 \wedge \cdots \wedge w_{0i}^q \wedge d B_i^\alpha \\
&\,\,\,\,\,+\left( \sum_{\beta=1}^q (-1)^{q-\beta} C_i^\beta w_{0i}^1 \wedge \cdots  \wedge \widehat{w}_{0i}^\beta \wedge \cdots \wedge w_{0i}^q \wedge dw_{0i}^\alpha  + w_{0i}^q\wedge \cdots \wedge w_{0i}^q \wedge d C_i^\alpha - w_{0i}^1\wedge \cdots \wedge w_{0i}^q \wedge B_i^\alpha \right) \wedge dt
\end{align*}
This implies that $G_1\left( \overline{B_i}, C_i \right)=0$. On the other hand, we note that
\begin{align*}
\theta_{ji}'^*\left( \omega_i'^\alpha \right)&= \theta_{ji}'^*\left(\omega_i^\alpha + t B_i^\alpha + C_i^\alpha dt \right) = \theta_{ji}^* \omega_i^\alpha + t B_i^\alpha + C_i^\alpha dt - t\mathcal{L}_{p_{ij}}\left( w_{0i}^\alpha \right) - i_{p_{ij}}\left(w_{0i}^\alpha \right)dt \\
\theta_{ji}'^*\left( \omega_i'^\alpha \right)&= \sum_{\beta=1}^q h_{ij}'^{\alpha\beta} \omega_j'^\beta = \sum_{\beta=1}^q \left( h_{ij}^{\alpha\beta} + t g_{ij}^{\alpha\beta} \right) \left( \omega_j^\beta + t B_j^\beta + C_j^\beta dt \right) \\
& = \sum_{\beta=1}^q h_{ij}^{\alpha\beta} \omega_j^\beta +t \left(\sum_{\beta=1}^q h_{0ij}^{\alpha\beta} B_j^\beta + \sum_{\beta=1}^q g_{ij}^{\alpha\beta} w_{0j}^\beta \right) + \left(\sum_{\beta=1}^q h_{0ij}^{\alpha\beta} C_j^\beta\right) dt
\end{align*}
This implies that $\delta\left( \left\{ \overline{B_i} \right\}, \left\{ C_i \right\} \right) + G_0\left( \left\{p_{ij} \right\} \right)=0$. This proves the claim. Since $\mathbb{H}^1\left( X, \tilde{\mathcal{N}}_{\mathcal{F}_0}^{*\bullet} \right)=0$, there exists $\left\{ a_i \right\}\in C^0\left( \mathcal{U} , T_X \right)$ such that $a_i- a_j= p_{ij}$ and $\mathcal{L}_{a_i}\left( w_{0i}^\alpha \right)- B_i^\alpha = \sum_{\beta=1}^q e_i^{\alpha\beta} w_{0i}^\beta$ for some $e_i^{\alpha\beta} \in \Gamma\left( U_i, \mathcal{O}_X \right)$ and $i_{a_i}\left( w_{0i}^\alpha\right)= C_i^\alpha$. We explicitly construct an unfolded isomorphism $\left( \mathcal{X}, \tilde{\mathcal{N}}_\mathcal{F}^* \right) \cong \left( \mathcal{X}',  \tilde{\mathcal{N}}_{\mathcal{F}'}^* \right)$. We claim that $\textnormal{Id}- t a_i: \left( \Gamma\left(  U_i, \mathcal{O}_X\right)\otimes \tilde{A} , \tilde{\mathcal{N}}_{\mathcal{F}_i'}^* \right)\to \left(  \Gamma\left( U_i, \mathcal{O}_X \right) \otimes \tilde{A}, \tilde{\mathcal{N}}_{\mathcal{F}_i}^* \right)$ is an unfolded isomorphism. In fact, $\omega_i'^\alpha - t \mathcal{L}_{a_i}\left( w_{0i}^\alpha \right) - i_{a_i}\left( w_{0i}^\alpha \right) dt = \omega_i^\alpha +t B_i^\alpha + C_i^\alpha dt - t\mathcal{L}_{a_i}\left( w_{0i}^\alpha \right) - i_{a_i}\left( w_{0i}^\alpha \right) dt = \omega_i^\alpha - \sum_{\beta=1}^q e_i^{\alpha\beta} w_{0i}^\beta = \sum_{\beta=1}^q \left( \delta_\alpha^\beta -t e_i^{\alpha\beta} \right) \omega_i^\beta $. On the other hand, since $p_{ij}= a_i- a_j$, we have $\left( \textnormal{Id}- t a_i \right)\circ f_{ij}'= f_{ij} \circ \left( \textnormal{Id} - t a_j \right)$. Hence $\left\{ \textnormal{Id}- t a_i \right\}$ glues together to give an unfolded isomorphism $\left( \mathcal{X}, \tilde{\mathcal{N}}_\mathcal{F}^*  \right)\cong \left(  \mathcal{X}', \tilde{\mathcal{N}}_{\mathcal{F}'}^* \right)$.

Now we prove that $\mathbb{H}^1\left( X , \tilde{\mathcal{N}}_{\mathcal{F}_0}^{*\bullet} \right)=0$, then $\left( X, \mathcal{N}_{\mathcal{F}_0}^*  \right)$ is rigid in unfoldings over $k[x]/(x^m)$ for any $m\geq 1$. We prove by induction $m$. Assume that any infinitesimal unfolding of $\left(X, \mathcal{N}_{\mathcal{F}_0}^* \right)$ over $k[x]/(x^m)$ is trivial. Since any two liftings to $\tilde{A}=k[x]/(x^{m+1})$ are equivalent and the trivial unfolding is also a lifting, any infinitesimal unfolding of $\left( X, \mathcal{N}_{\mathcal{F}_0}^* \right)$ over $\tilde{A}=k[x]/(x^{m+1})$ is trivial. This completes the proof of Proposition \ref{a1}.
\end{proof}

As mentioned as above, the proof of Proposition \ref{a1} is based on the transitive action of $\mathbb{H}^1\left( X, \tilde{\mathcal{N}}_{\mathcal{F}_0}^* \right)$ on liftings of an unfolding over $k[x]/(x^m)$ for a small extension of the form $0\to (t)=(x^m)\to k[x]/(x^{m+1})\to k[x]/(x^m)\to 0$. More precisely,
\begin{proposition} \label{a6}
Let $\left( X, \mathcal{N}_{\mathcal{F}_0}^* \right)$ be a foliated nonsingular variety with $\mathcal{N}_{\mathcal{F}_0}^*$ locally free and $e: 0\to t=(x^m)\to \tilde{A}=k[x]/(x^{m+1})\xrightarrow{\mu} A=k[x]/(x^m)\to 0$ be a small extension in $\bold{Art}$. Let $p:=\textnormal{Unf}_{\left( X, \mathcal{N}_{\mathcal{F}_0}^* \right)}\left( \mu\right): \textnormal{Unf}_{\left( X, \mathcal{N}_{\mathcal{F}_0}^* \right)}\left( \tilde{A}\right)\to \textnormal{Unf}_{\left( X, \mathcal{N}_{\mathcal{F}_0}^* \right)}\left(A\right)$. Then given $\xi=\left( \mathcal{X}, \tilde{\mathcal{N}}_\mathcal{F}^* \right)$, there is a transitive action of $\mathbb{H}^1\left( X, \tilde{\mathcal{N}}_{\mathcal{F}_0}^{*\bullet} \right)$ on $p^{-1}(\xi)$. Moreover, if $\mathbb{H}^0\left( X, \tilde{\mathcal{N}}_{\mathcal{F}_0}^{*\bullet} \right)=0$, then the action is free.
\end{proposition}

\begin{proof}
We will define a map $G: \mathbb{H}^1\left( X, \tilde{\mathcal{N}}_{\mathcal{F}_0}^{*\bullet} \right)\times p^{-1}(\xi) \to p^{-1}(\xi)$ which is a group action of $\mathbb{H}^1\left( X , \tilde{\mathcal{N}}_{\mathcal{F}_0}^{* \bullet} \right)$ acting on $p^{-1}(\xi)$. Let $\tilde{\xi}=\left(\tilde{\mathcal{X}}, \tilde{\mathcal{N}}_{\tilde{\mathcal{F}}}^\bullet\right)\in p^{-1}(\xi)$ over $\tilde{A}=k[x]/(x^{m+1})$. Let $\left\{ U_i \right\}$ be an affine open cover of $X$ such that $\Gamma\left( U_i , \mathcal{N}_{\mathcal{F}_0}^*  \right)$ is generated by $w_{0i}^1,..., w_{0i}^q \in \Gamma\left(  U_i , \Omega_X^1 \right)$ with $w_{0i}^\alpha = \sum_{\beta=1}^q h_{0ij}^{\alpha\beta} w_{0j}^\beta$ and $w_{0i}^1 \wedge \cdots \wedge w_{0i}^q \wedge dw_{0i}^\alpha=0$, and $\tilde{\theta}_i:\left( U_i\times \textnormal{Spec}\left( \tilde{A} \right), \tilde{\mathcal{N}}_{\mathcal{F}_i}^* \right)\to \left( \mathcal{X}|_{U_i}, \tilde{\mathcal{N}}_\mathcal{F}^* |_{U_i } \right)$ is an unfolded isomorphism such that $\Gamma\left( U_i, \tilde{\mathcal{N}}_{\mathcal{F}_i}^* \right)$ is generated by $\tilde{\omega}_i^1=\tilde{w}_i^1 + \tilde{v}_i^1,..., \tilde{\omega}_i^q= \tilde{w}_i^q + \tilde{v}_i^q\in  \left( \Gamma\left( U_i , \Omega_X^1 \right)\otimes \tilde{A}  \right) \bigoplus \left( \Gamma\left( U_i, \mathcal{O}_X\right)\otimes \Omega_{\tilde{A}}^1 \right)$ which induces $w_{0i}^1,..., w_{0i}^q$, and $\tilde{\theta}_{ji}^* \tilde{\omega}_i^\alpha= \sum_{\beta=1}^q \tilde{h}_{ij}^{\alpha\beta} \tilde{\omega}_j^\beta$ where $\tilde{\theta}_{ij}=\tilde{\theta}_j^{-1}\tilde{\theta}_i$ and $\tilde{h}_{ij}^{\alpha\beta}$ induces $h_{0ij}^{\alpha\beta}$, and $\tilde{\omega}_i^1 \wedge \cdots \wedge \tilde{\omega}_i^q \wedge d \tilde{\omega}_i^\alpha = 0$. Let $\tilde{\theta}_{ij}$ corresponds to $\tilde{f}_{ij}: \left( \Gamma\left(U_{ij}, \mathcal{O}_X \right)\otimes \tilde{A}, \tilde{\mathcal{N}}_{\tilde{\mathcal{F}}_j}^* \right) \to \left( \Gamma\left( U_{ij}, \mathcal{O}_X \right) \otimes \tilde{A}, \tilde{\mathcal{N}}_{\tilde{\mathcal{F}}_i}^* \right)$. Let $\left( \left\{ \overline{w_i'} , \left\{ v_i \right\}, \left\{ p_{ij}\right\}\right\}, \right)\in \mathbb{H}^1\left( X, \tilde{\mathcal{N}}_{\mathcal{F}_0}^{ *\bullet} \right)$ corresponding to a first-order unfolding as in the proof of Proposition \ref{uu35}. We keep the notations in the proof of Proposition \ref{uu35}. Then $\left\{ \tilde{f}_{ij} + t p_{ij}:\Gamma\left( U_{ij}, \mathcal{O}_X \right)\otimes \tilde{A} \to \Gamma\left( U_{ij} , \mathcal{O}_X \right)\otimes \tilde{A} \right\}$ which corresponds to isomorphisms $\left\{ \tilde{\theta}_{ij}': U_{ij}\times \textnormal{Spec}\left( \tilde{A} \right) \to U_{ij} \times \textnormal{Spec}\left( \tilde{A} \right) \right\}$ satisfies the gluing condition, and since $\tilde{\theta}_{ji}'^* \left( \tilde{\omega}_i^\alpha + t w_i'^\alpha + v_i^\alpha dt \right)= \tilde{\theta}_{ji}^* \tilde{\omega}_i^\alpha +t w_i'^\alpha + v_i^\alpha dt - t\mathcal{L}_{p_{ij}}\left(w_{0i}^\alpha\right) - i_{p_{ij}}\left( w_{0i}^\alpha \right) dt = \sum_{\beta=1}^q \tilde{h}_{ij}^{\alpha\beta} \omega_j^\beta + t\sum_{\beta=1}^q h_{0ij}^{\alpha\beta} w_j'^\beta + t\sum_{\beta=1}^q g_{ij}^{\alpha\beta} w_{0j}^\beta + \sum_{\beta=1}^q h_{0ij}^{\alpha\beta} v_j^\beta dt  = \sum_{\beta=1}^q \left( \tilde{h}_{ij}^{\alpha\beta} + t g_{ij}^{\alpha\beta}\right)\left( \tilde{\omega}_j^\beta + t w_j'^\beta + v_j^\beta dt \right)$ is an unfolded isomorphism. On the other hand, $\left( \tilde{\omega}_i^1 +t w_i'^1+ v_i^1 dt \right) \wedge \cdots \wedge \left(\tilde{\omega}_i^q + t w_i'^q + v_i^q dt \right)\wedge d\left( \tilde{w}_i^\alpha + t w_i'^\alpha + v_i^\alpha dt \right)= \tilde{\omega}_i^1 \wedge \cdots \wedge \tilde{\omega}_i^q \wedge d \tilde{\omega}_i^q +t \sum_{\beta=1}^q w_{0i}^1 \wedge \cdots \wedge w_i'^\beta \wedge \cdots \wedge w_{0i}^q \wedge d w_{0i}^\alpha +t w_{0i}^1 \wedge \cdots \wedge w_{0i}^q \wedge dw_i'^\alpha + \sum_{\beta=1}^q (-1)^{q-\beta} v_i^\beta w_{0i}^1 \wedge \cdots \widehat{w}_{0i}^\beta \wedge \cdots \wedge w_{0i}^q \wedge dw_{0i}^\alpha \wedge dt + w_{0i}^1 \wedge \cdots \wedge w_{0i}^q \wedge  d v_i^\alpha \wedge dt  - w_{0i}^1 \wedge \cdots \wedge w_{0i}^q \wedge w_i'^\alpha \wedge dt=0$. This defines an infinitesimal unfolding $\tilde{\xi}'=\left( \tilde{\mathcal{X}}', \tilde{\mathcal{N}}_{\tilde{\mathcal{F}}'}^*  \right)$ of $\left( X, \mathcal{N}_{\mathcal{F}_0}^* \right)$ over $\tilde{A}$ which induces $\xi=\left( \mathcal{X}, \tilde{\mathcal{N}}_\mathcal{F}^* \right)$, so that $\left( \tilde{\mathcal{X}}, \tilde{\mathcal{N}}_{\tilde{\mathcal{F}}'}^*  \right)$ lies in $p^{-1}(\xi)$. Then we define a map $G$ by
\begin{align*}
G: \mathbb{H}^1\left( X, \tilde{\mathcal{N}}_{ \mathcal{F}_0}^{*\bullet} \right) \times p^{-1}(\xi) & \to p^{-1}(\xi) \\
\left( \left\{ \overline{w_i'} \right\}, \left\{ v_i \right\}, \left\{ p_{ij} \right\} \right)\times \tilde{\xi} & \mapsto \tilde{\xi}'
\end{align*}
$G$ is well-defined. In fact, let $\left( \left\{ \overline{w_i''} \right\}, \left\{ v_i'\right\}, \left\{ p_{ij}' \right\} \right)$ represent the same  cohomology class with $\left( \left\{ \overline{w_i'}\right\}, \left\{ v_i\right\}, \left\{ p_{ij} \right\} \right)$. Then we have to show that $\tilde{\xi}'=G\left( \left( \left\{ \overline{w_i'} \right\}, \left\{ v_i \right\}, \left\{ p_{ij} \right\} \right), \tilde{\xi}  \right)$ is equivalent to $\tilde{\xi}''= G\left( \left( \left\{ \overline{w_i''} \right\}, \left\{ v_i' \right\}, \left\{ p_{ij}' \right\} \right), \tilde{\xi} \right)$. Then there exists $\left\{ a_i \right\}\in C^0\left( \mathcal{U}, T_X \right)$ such that $a_i- a_j= p_{ij}- p_{ij}'$ and $\mathcal{L}_{a_i}\left( w_{0i}^\alpha \right)= w_i'^\alpha - w_i''^\alpha - \sum_{\beta=1}^q c_i^{\alpha\beta} w_{0i}^\beta$ for $c_i^{\alpha\beta}\in \Gamma\left( U_i, \mathcal{O}_X \right)$, and $i_{a_i}\left( w_{0i}^\alpha \right)= v_i^\alpha - v_i'^\alpha$. We claim that $\textnormal{Id} - t a_i :\left( \Gamma\left( U_i, \mathcal{O}_X \right)\otimes \tilde{A}, \tilde{\mathcal{N}}_{\tilde{\mathcal{F}}_i'}^* \right) \to \left( \Gamma\left( U_i, \mathcal{O}_X \right)\otimes \tilde{A}, \tilde{\mathcal{N}}_{\tilde{\mathcal{F}}_i''}^* \right)$ is an unfolded isomorphism. In fact, $\tilde{\omega}_i^\alpha + t w_i'^\alpha + v_i^\alpha dt - t \mathcal{L}_{a_i}\left( w_{0i}^\alpha \right)- i_{a_i}\left( w_{0i}^\alpha \right) dt= \tilde{\omega}_i^\alpha + t w_i''^\alpha + t \sum_{\beta=1}^q c_i^{\alpha\beta} w_{0i}^\beta  + v_i'^\alpha dt =\sum_{\beta=1}^q \left( \delta_\alpha^\beta + t c_i^{\alpha\beta} \right)\left( \tilde{\omega}_i^\beta + t w_i''^\beta + v_i'^\beta dt  \right)$. On the other hand, $\left(  \tilde{f}_{ij} + t p_{ij}' \right)\circ \left( \textnormal{Id}- ta_j \right)= \left( \textnormal{Id}- ta_i \right)\circ \left( \tilde{f}_{ij} + t p_{ij} \right)$. Hence $\left\{ \textnormal{Id}- ta_i\right\}$ glues together to give an unfolded isomorphism between $\tilde{\xi}'$ and $\tilde{\xi}''$.

Next we show that the group action $G$ is transitive. In fact, let $\tilde{\xi}\in p^{-1}(\xi)$ be a lifting of $\xi$ as above. Choose arbitrary lifting $\eta\in p^{-1}(\xi)$ of $\xi$. We have to find $\left( \left\{  \overline{w_i' } \right\}, \left\{ v_i \right\}, \left\{ p_{ij} \right\} \right)\in \mathbb{H}^1\left( X, \tilde{\mathcal{N}}_{\mathcal{F}_0}^{*\bullet} \right)$ such that $\left( \left( \left\{ \overline{w_i'} \right\}, \left\{ v_i \right\}, \left\{ p_{ij} \right\} \right), \tilde{\xi} \right)$ is mapped to $\eta$ under the action $G$. Let $\eta$ consist of the data: $f_{ij}':\Gamma\left( \Gamma\left( U_{ij}, \mathcal{O}_X\right) \otimes \tilde{A}, \tilde{\mathcal{N}}_{\tilde{\mathcal{F}_j'}}^* \right)\to \left( \Gamma\left( U_{ij}, \mathcal{O}_X \right) \otimes \tilde{A} , \tilde{\mathcal{N}}_{\tilde{\mathcal{F}}_i'}^* \right)$ such that $\tilde{f}_{ij}' $ and $\tilde{f}_{ij}$ induce the same homomorphism in $\xi$, and $\Gamma\left( U_i, \tilde{\mathcal{N}}_{\tilde{\mathcal{F}}_i'}^* \right)$ is generated by $\tilde{\omega}_i^\alpha\in \left( \Gamma\left(U_i, \Omega_X^1 \right)\otimes \tilde{A} \right) \bigoplus \left(\Gamma\left( U_i, \mathcal{O}_X \right)\otimes \Omega_{\tilde{A}}^1 \right) ,\alpha=1,...,q$ such that $\tilde{\omega}_i'^\alpha$ and $\tilde{\omega}_i^\alpha$ induces the same element in $\left(  \Gamma\left( U_i, \Omega_X^1 \right) \otimes A \right) \bigoplus \left( \Gamma\left( U_i, \mathcal{O}_X \right)\otimes \Omega_A^1  \right)$. Since $\tilde{\xi}$ and $\eta$ induce $\xi$, we have $\tilde{f}_{ij}'=\tilde{f}_{ij} + t p_{ij}$ and $\tilde{\omega}_i'^\alpha= \tilde{\omega}_i^\alpha + t B_i^\alpha + C_i^\alpha dt$. We define $\overline{B_i}:\Gamma\left(U_i, \mathcal{N}_{\mathcal{F}_0}^* \right)\to \Gamma\left( U_i, \frac{\Omega_X^1}{\mathcal{N}_{\mathcal{F}_0}^* } \right)$ by $w_{0i}^\alpha \mapsto \overline{B_i^\alpha}$ and $C_i:\Gamma\left( U_i, \mathcal{N}_{\mathcal{F}_0}^* \right) \to \Gamma\left( U_i, \mathcal{O}_X \right)$ by $w_{0i}^\alpha \mapsto C_i^\alpha$. Then as in the same argument in the proof of Proposition \ref{a1}, $\left( \left\{ \overline{B_i} \right\}, \left\{ C_i \right\}, \left\{ p_{ij} \right\} \right)$ defines a $1$-cocycle of the \v Cech resolution of $\tilde{\mathcal{N}}_{\mathcal{F}_0}^{*\bullet}$.

Now assume that $\mathbb{H}^0\left( X, \tilde{\mathcal{N}}_{\mathcal{F}_0}^{*\bullet} \right)=0$. We claim that the group action $G$ is free. Assume that for given $v=\left( \left\{ \overline{w_i'}\right\}, \left\{ v_i \right\}, \left\{ p_{ij} \right\} \right)\in \mathbb{H}^1\left( X, \tilde{\mathcal{N}}_{\mathcal{F}_)}^{*\bullet} \right)$ and $\tilde{\xi}\in p^{-1}(\xi)$, we have $G\left( v, \tilde{\xi} \right)=\tilde{\xi}$. We have to show that $v=0\in \mathbb{H}^1\left( X, \tilde{\mathcal{N}}_{\mathcal{F}_0}^{* \bullet} \right)$. Let $\tilde{\xi}$ and $\tilde{\xi}'= G\left( v, \tilde{\xi}\right)$ as above. Then $\tilde{\xi}'= G\left(v, \tilde{\xi}  \right)=\tilde{\xi}$ implies that we have unfolded isomorphisms $g_i:\left( \Gamma\left( U_i , \mathcal{O}_X \right)\otimes \tilde{A}, \tilde{\mathcal{N}}_{\tilde{\mathcal{F}}_i}^* \right) \to \left( \Gamma\left( U_i, \mathcal{O}_X \right)\otimes \tilde{A}, \tilde{\mathcal{N}}_{\tilde{\mathcal{F}}_i'}^* \right)$ such that $g_i \circ f_{ij}= \left( f_{ij} + t p_{ij} \right)\circ g_j$. Since $\left\{ g_i \right\}$ induces an unfolded automorphism on $\xi=\left( \mathcal{X}, \tilde{\mathcal{N}}_\mathcal{F}^* \right)$ by modulo $(t)$, which is necessarily the identity by the following Lemma \ref{a5}, $g_i$ is of the form $\textnormal{Id}+ t q_i$ for some $q_i \in \Gamma\left( U_i, T_X \right)$. Then $\tilde{\omega}_i^\alpha + t \mathcal{L}_{q_i}\left( w_{0i}^\alpha \right) + i_{q_i}\left( w_{0i}^\alpha \right) dt = \sum_{\beta=1}^q \tilde{c}_i^{\alpha\beta} \left( \tilde{w}_i^\beta + t w_i'^\beta + v_i^\beta  dt \right)$ for some $\tilde{c}_i^{\alpha\beta} \in \Gamma\left( U_i, \mathcal{O}_X \right)\otimes \tilde{A}$ which induces $c_i^{\alpha\beta}\in \Gamma\left( U_i, \mathcal{O}_X \right)\otimes A$. By restricting to $A$, we have $\omega_i^\alpha= \sum_{\beta=1}^q c_i^{\alpha\beta} \omega_i^\beta$, so that $c_i^{\alpha\beta}=\delta_\alpha^\beta$. Then $\tilde{c}_i^{\alpha\beta}:= \delta_\alpha^\beta + t e_i^{\alpha\beta}$ for some $e_i^{\alpha\beta}\in \Gamma\left(  U_i, \mathcal{O}_X \right)$. Then we have $\tilde{\omega}_i^\alpha + t \mathcal{L}_{q_i}\left( w_{0i}^\alpha \right) + i_{q_i}\left( w_{0i}^\alpha \right) dt= \sum_{\beta=1}^q \left(\delta_\alpha^\beta + t e_i^{\alpha\beta} \right) \left( \tilde{\omega}_i^\beta + t w_i'^\beta + v_i^\beta dt \right)$, so that $\mathcal{L}_{q_i}\left( w_{0i}^\alpha \right)= w_i'^\alpha + \sum_{\beta=1}^q e_i^{\alpha\beta} w_{0i}^\beta$ and $i_{q_i}\left( w_{0i}^\alpha \right)= v_i^\alpha$. On the other hand,$\left( \textnormal{Id}+ t q_i \right)\circ f_{ij} = \left( \right)\circ \left( f_{ij} + tp_{ij} \right)\circ \left( \textnormal{Id} + t q_i \right)$, so that $q_i- q_j= p_{ij}$. This implies that $v=0\in \mathbb{H}^1\left(X, \tilde{\mathcal{N}}_{\mathcal{F}_0}^{*\bullet} \right)$. This completes the proof of Proposition \ref{a6}.
\end{proof}

\begin{lemma}\label{a5}
Let $\left( X, \mathcal{N}_{\mathcal{F}_0}^* \right)$ be a foliated nonsingular variety with $\mathcal{N}_{\mathcal{F}_0}^*$ locally free. Assume that $\mathbb{H}^0\left(\mathcal{X}, \tilde{\mathcal{N}}_{\mathcal{F}_0}^{*\bullet} \right) = 0 $. The for any infinitesimal unfolding $\left( \mathcal{X}, \tilde{\mathcal{X}}_\mathcal{F}^* \right)$ of $\left( X, \mathcal{N}_{\mathcal{F}_0}^* \right)$ over $k[x]/(x^m)$,
\begin{align*}
\textnormal{Aut}\left(  \frac{\left( \mathcal{X}, \tilde{\mathcal{N}}_\mathcal{F}^* \right) }{\left( X,  \mathcal{N}_{\mathcal{F}_0}^* \right) }\right)=  \textnormal{Id}
\end{align*}
where $\textnormal{Aut}\left( \frac{\left( \mathcal{X}, \tilde{\mathcal{N}}_\mathcal{F}^*\right)}{\left( X, \mathcal{N}_{\mathcal{F}_0}^* \right) } \right)$ is the set of unfolded automorphisms of $\left( \mathcal{X}, \tilde{\mathcal{N}}_\mathcal{F}^* \right)$ restricting to the identity on $\left( X, \mathcal{N}_{\mathcal{F}_0}^* \right)$.
\end{lemma}

\begin{proof}
We prove by the induction on $m$. Let's assume that Lemma \ref{a5} holds for $m-1$. Let $\left( \mathcal{X}, \tilde{\mathcal{N}}_\mathcal{F}^* \right)$ be an infinitesimal unfolding over $A= k[x]/(x^m)$. Let $g:\left( \mathcal{X}, \tilde{\mathcal{N}}_\mathcal{F}^* \right)\to \left( \mathcal{X}, \tilde{\mathcal{N}}_\mathcal{F}^* \right)$ over $A$ be an unfolded automorphism with $g^* \tilde{\mathcal{N}}_\mathcal{F}^* =\tilde{\mathcal{N}}_\mathcal{F}^*$. Let $\mathcal{U}=\left\{ U_i \right\}$ be an affine open covering of $X$ such that $\Gamma\left( U_i, \mathcal{N}_{\mathcal{F}_0}^* \right)$ is generated by $w_{0i}^1,..., w_{0i}^q \in \Gamma\left( U_i, \Omega_X^1\right)$ with $w_{0i}^\alpha= \sum_{\beta=1}^q h_{0ij}^{\alpha\beta} w_{0j}^\beta $ and $w_{0i}^1 \wedge \cdots \wedge w_{0i}^q \wedge d w_{0i}^\alpha = 0$, and $\theta_i:\left( U_i \times \textnormal{Spec}(A) , \tilde{\mathcal{N}}_{\mathcal{F}_i}^* \right)\to \left( \mathcal{X}|_{U_i}, \tilde{\mathcal{N}}_\mathcal{F}^*|_{U_i} \right)$ is an unfolded isomorphism such that $\Gamma\left( U_i, \tilde{\mathcal{N}}_{\mathcal{F}_i}^* \right)$ is generated by $\omega_i^1,... , \omega_i^q \in \left( \Gamma\left( U_i, \Omega_X^1 \right)\otimes A  \right) \bigoplus \left( \Gamma\left( U_i, \mathcal{O}_X \right)\otimes \Omega_A^1 \right)$ which induces $w_{0i}^1,..., w_{0i}^q$. Let $\theta_{ij}= \theta_j^{-1}\theta_i$ corresponds to an unfolded isomorphism $f_{ij}:\left( \Gamma\left( U_{ij}, \mathcal{O}_X \right) \otimes A , \tilde{\mathcal{N}}_{\mathcal{F}_j}^* \right)\to \left( \Gamma\left( U_{ij} , \mathcal{O}_X \right)\otimes A , \tilde{\mathcal{N}}_{\mathcal{F}_i}^*  \right)$. Then $g$ is described by the data $\left\{ g_i \right\}$ where $g_i:\left( \Gamma\left( U_i, \mathcal{O}_X \right), \tilde{\mathcal{N}}_{\mathcal{F}_i}^*  \right)\to \left( \Gamma\left( U_i , \mathcal{O}_X \right)\otimes A, \tilde{\mathcal{N}}_{\mathcal{F}_i}^*  \right)$ which is an unfolded isomorphism and $g_i\circ f_{ij}= f_{ij} \circ g_j$. By the induction hypothesis, $g_i$ induces the identity on $\Gamma\left( U_i, \mathcal{O}_X \right)\otimes k[x]/(x^{m-1})$. Then $g_i$ is of the form $g_i=\textnormal{Id} +t d_i$ for some $d_i \in \Gamma\left( U_i, \mathcal{O}_X \right)$ where $t=x^m$ such that $d_i= d_j$ and $\omega_i^\alpha + t \mathcal{L}_{d_i}\left( w_{0i}^\alpha\right) + i_{d_i}\left( w_{0i}^\alpha \right) dt =\sum_{\beta=1}^q s_i^{\alpha\beta} \omega_i^\beta$ for some $s_i^{\alpha\beta}\in \Gamma\left( U_i, \mathcal{O}_X\right)\otimes A$. On the other hand, we note that by tensoring $0\to (t)\to A=k[x]/(x^m)\to \overline{A}=k[x]/(x^{m-1})\to 0$ with $\tilde{\mathcal{N}}_{\mathcal{F}}^*$, we have $0\to t \mathcal{N}_{\mathcal{F}_0}^*\to \tilde{\mathcal{N}}_{\mathcal{F}}^* \to \tilde{\mathcal{N}}_\mathcal{F}^*\otimes \overline{A} \to 0$ ($\tilde{\mathcal{N}}_\mathcal{F}^*$ is locally free and $\mathcal{X}$ is flat over $A$). Then since $\omega_i^\alpha - \sum_{\beta=1}^q s_i^{\alpha\beta} \omega_i^\beta$ induces $0$ in $\Gamma\left( U_i,\tilde{\mathcal{N}}_\mathcal{F}^*  \right)\otimes \overline{A}$, we can write $\omega_i^\alpha - \sum_{\beta=1}^q s_i^{\alpha\beta} \omega_i^\beta = t Q_i^\alpha$ for some $Q_i^\alpha \in \Gamma\left( U_i, \mathcal{N}_{\mathcal{F}_0}^* \right)$. This implies that $\mathcal{L}_i\left( w_{0i}^\alpha \right) \in \Gamma\left( U_i, \mathcal{N}_{\mathcal{F}_0}^* \right)$ and $i_{d_i}\left(w_{0i}^\alpha \right)=0$. Hence $d:=\left\{ d_i \right\}\in \mathbb{H}^0\left( X, \tilde{\mathcal{N}}_{\mathcal{F}_0}^{*\bullet} \right)$. Since $\mathbb{H}^0\left( X, \tilde{\mathcal{N}}_{\mathcal{F}_0}^{ * \bullet} \right)=0$, we have $d=0$, so that $g=\textnormal{Id}$.
\end{proof}

\subsection{Stability of unfoldings over $k[x]/(x^m)$ in $X$}\

\begin{proposition}\label{a2}
Let $\left( X, \mathcal{N}_{\mathcal{F}_0}^* \right)$ be a foliated nonsingular variety with $\mathcal{N}_{\mathcal{F}_0}^*$ locally free. Assume that $\mathbb{H}^1\left( X, \textnormal{Unf}_{\mathcal{N}_{\mathcal{F}_0}^*}^\bullet \right)=0$. Then for an infinitesimal deformation $\mathcal{X}$ of $X$ over $k[x]/(x^m)$, there exists an infinitesimal unfolding $\left( \mathcal{X}, \tilde{\mathcal{N}}_\mathcal{F}^* \right)$ over $k[x]/(x^m)$ which induces $\mathcal{X}$.
\end{proposition}

\begin{proof}
We prove by induction on $m$. Assume that Proposition \ref{a2} holds for $m$. Let $\tilde{\mathcal{X}}$ be a flat infinitesimal deformation of $X$ over $\tilde{A}=k[x]/(x^{m+1})$, and $\left( \mathcal{X}, \tilde{\mathcal{N}}_\mathcal{F}^* \right)$ be an infinitesimal unfolding over $A=k[x]/(x^m)$ such that $\tilde{\mathcal{X}}$ induces $\mathcal{X}$ by the induction hypothesis. Let $\mathcal{U}=\left\{ U_i \right\}$ be an affine open covering of $X$ such that we have isomorphisms $\tilde{\theta}_i: U_i \times \textnormal{Spec}\left(\tilde{A} \right)\to \tilde{\mathcal{X}}|_{U_i}$, and $\tilde{\theta}_{ij}= \tilde{\theta}_j^{-1}\tilde{\theta}_i: U_{ij}\times \textnormal{Spec}\left( \tilde{A} \right) \to U_{ij} \times \textnormal{Spec}\left( \tilde{A} \right)$ corresponds to $\tilde{f}_{ij} : \Gamma\left( U_{ij}, \mathcal{O}_X \right) \otimes \tilde{A} \to \Gamma\left( U_{ij}, \mathcal{O}_X  \right) \otimes \tilde{A}$, and we have an unfolded isomorphism $\theta_i:\left( U_i\times \textnormal{Spec}\left( A \right) , \tilde{\mathcal{N}}_{\mathcal{F}_i}^* \right)\to \left( \mathcal{X}|_{U_i}, \tilde{\mathcal{N}}_{\mathcal{F}}^* \right)$ where $\theta_i$ is induced from $\tilde{\theta}_i$ and $\tilde{\mathcal{N}}_{\mathcal{F}_i}^* \subset \left( \Omega_X^1|_{U_i}\otimes A  \right)  \bigoplus \left( \mathcal{O}_X|_{U_i}\otimes  \Omega_A^1 \right)$ is induced from $\tilde{\mathcal{N}}_{\mathcal{F}}^*$ via $\theta_i$, and $\theta_{ij}= \theta_j^{-1}\theta_i$ induced from $\tilde{\theta}_{ij}$ corresponds to $f_{ij}: \Gamma\left( U_{ij} , \mathcal{O}_X \right)\otimes A \to \Gamma\left( U_{ij} , \mathcal{O}_X \right)\otimes A$ induced from $\tilde{f}_{ij}$, and $\Gamma\left(U_i, \tilde{\mathcal{N}}_{\mathcal{F}_i}^* \right)$ is generated by $\omega_i^1=w_i^1+ v_i^1,..., \omega_i^q = w_i^q + v_i^q \in \left( \Gamma\left( U_i , \Omega_X^1 \right) \otimes A \right) \bigoplus \left( \Gamma\left( U_i , \mathcal{O}_X \right) \otimes \Omega_A^1 \right)$ which induces $w_{0i}^1,..., w_{0i}^q$, and $\theta_{ji}^* \omega_i^\alpha= \sum_{\beta=1}^q h_{ij}^{\alpha\beta} \omega_j^\beta$ for some $h_{ij}^{\alpha\beta}\in \Gamma\left( U_{ij} , \mathcal{O}_X \right) \otimes A $ which induces $h_{0ij}^{\alpha\beta}\in \Gamma\left(  U_{ij} , \mathcal{O}_X \right)$. Let $\tilde{w}_i^\alpha \in \Gamma\left( U_i, \Omega_X^1 \right)\otimes \tilde{A}$ be an arbitrary lifting of $w_i^\alpha$ and $\tilde{v}_i^\alpha \in \Gamma\left( U_i, \mathcal{O}_X\right) \otimes \Omega_{\tilde{A}}^1 $ be an arbitrary lifting of $v_i^\alpha$ and $\tilde{h}_{ij}^{\alpha\beta}\in \Gamma\left( U_{ij}, \mathcal{O}_X \right) \otimes \tilde{A} $ be an arbitrary lifting of $h_{ij}^{\alpha\beta}$. Then as in the proof of Proposition \ref{uu35} (we keep the notations in the proof of Proposition \ref{uu35}), the obstruction
{\Small{\begin{align*}
&\left( \left\{B_i \right\}, \left\{ F_i\right\}, \left\{ - \overline{R_{ij}} \right\}, \left\{ - G_{ij}\right\}, \left\{ d_{ijk}=0 \right\} \right)  \\
&\in C^0\left(\mathcal{U}, \mathscr{H}om_{\mathcal{O}_X}\left( \mathcal{N}_{\mathcal{F}_0}^*, \tilde{\mathcal{S}}^2 \right) \right)\bigoplus C^0\left(\mathcal{U}, \mathscr{H}om_{\mathcal{O}_X}\left( \mathcal{N}_{\mathcal{F}_0}^*, \tilde{\mathcal{S}}^1  \right) \right) \bigoplus C^1\left(\mathcal{U}, \mathscr{H}om_{\mathcal{O}_X}\left( \mathcal{N}_{\mathcal{F}_0}^*, \frac{\Omega_X^1}{\mathcal{N}_{\mathcal{F}_0}^*} \right) \right)\\
& \bigoplus C^1\left(\mathcal{U}, \mathscr{H}om_{\mathcal{O}_X}\left(\mathcal{N}_{\mathcal{F}_0}^*, \mathcal{O}_X \right) \right) \bigoplus C^2\left( \mathcal{U}, T_X \right)
\end{align*}}}
lies in $\mathbb{H}^1\left( X, \textnormal{Unf}_{\mathcal{N}_{\mathcal{F}_0^*}}^\bullet \right)=0$, so that a lifting $\left(\tilde{\mathcal{X}}, \tilde{\mathcal{N}}_{\tilde{\mathcal{F}}}^* \right)$ over $\tilde{A}$  of $\left( \mathcal{X}, \tilde{\mathcal{N}}_\mathcal{F}^* \right)$ exists. This completes the proof of Proposition \ref{a2}.
\end{proof}

\subsection{Stability of unfoldings over $k[x]/(x^m)$ in $\left( X, \mathcal{N}_{\mathcal{F}_0}^* \right)$}\

\begin{proposition} \label{a3}
Let $\left( X, \mathcal{N}_{\mathcal{F}_0}^* \right)$ be a foliated nonsingular variety with $\mathcal{N}_{\mathcal{F}_0}^*$ locally free. Assume that $\mathbb{H}^1\left( X, \mathscr{H}om_{\mathcal{O}_X } \left( \mathcal{N}_{\mathcal{F}_0}^*, \mathcal{O}_X \right)^\bullet \right)=0$. Then for any infinitesimal foliated deformation $\left( \mathcal{X}, \mathcal{N}_\mathcal{F}^* \right)$ of $\left( X, \mathcal{N}_{\mathcal{F}_0}^* \right)$ over $k[x]/(x^m)$, there exists an infinitesimal unfolding $\left( \mathcal{X}, \tilde{\mathcal{N}}_\mathcal{F}^* \right)$ of $\left( X , \mathcal{N}_{\mathcal{F}_0}^* \right)$ over $k[x]/(x^m)$ which induces $\left(  \mathcal{X}, \mathcal{N}_\mathcal{F}^*  \right)$.
\end{proposition}

\begin{proof}
We prove by induction on $m$. Assume that Proposition \ref{a3} holds for $m$. Let $\left( \tilde{X}, \mathcal{N}_{\tilde{\mathcal{F}} }^* \right)$ be an infinitesimal foliated deformation of $\left( X , \mathcal{N}_{\mathcal{F}_0}^* \right)$ over $\tilde{A}=k[x]/(x^{m+1})$, and $\left( \mathcal{X}, \tilde{\mathcal{N}}_\mathcal{F}^* \right)$ be an infinitesimal unfolding over $A=k[x]/(x^m)$ such that $\left( \tilde{\mathcal{X}}, \mathcal{N}_{\tilde{\mathcal{F}}}^*\right)$ induces $\left(\mathcal{X}, \mathcal{N}_\mathcal{F}^* \right)$ and $\left( \mathcal{X}, \mathcal{N}_\mathcal{F}^*\right)$ is the underlying foliated deformation of the unfolding $\left( \mathcal{X}, \tilde{\mathcal{N}}_\mathcal{F}^* \right)$  by the induction hypothesis. Let $\mathcal{U}=\left\{ U_i \right\}$ be an affine open covering of $X$ such that we have foliated isomorphisms $\tilde{\theta}_i: \left( U_i \times \textnormal{Spec}\left(\tilde{A} \right) , \mathcal{N}_{\mathcal{F}_i}^*  \right)\to \left(\tilde{\mathcal{X}}|_{U_i} , \mathcal{N}_\mathcal{F}^*|_{U_i} \right)$, and $\Gamma\left(U_i, \mathcal{N}_{\tilde{\mathcal{F}}_i}^* \right)$ is generated by $\tilde{w}_i^1,..., \tilde{w}_i^q \in \Gamma\left( U_i, \Omega_X^1 \right)\otimes \tilde{A}$, and let $\tilde{\theta}_{ij}= \tilde{\theta}_j^{-1}\tilde{\theta}_i: \left( U_{ij}\times \textnormal{Spec}\left( \tilde{A} \right) , \mathcal{N}_{\tilde{\mathcal{F}}_i}^*\right)\to U_{ij} \times \textnormal{Spec}\left( \tilde{A} \right)$ where $\theta_i$ corresponds to $\tilde{f}_{ij}: \left( \Gamma\left( U_{ij} , \mathcal{O}_X \right)\otimes \tilde{A} , \mathcal{N}_{\tilde{\mathcal{F}}_j}^*  \right)\to \left( \Gamma\left( U_{ij}, \mathcal{O}_X \right)\otimes \tilde{A} , \tilde{\mathcal{N}}_{\tilde{\mathcal{F}}_i}^* \right)$ and $\tilde{\theta}_{ji}^* \tilde{w}_i^\alpha= \sum_{\beta=1}^q \tilde{h}_{ij}^{\alpha\beta} \tilde{w}_j^\beta$ for some $\tilde{h}_{ij}^{\alpha\beta}\in \Gamma\left( U_{ij}, \mathcal{O}_X \right)\otimes \tilde{A}$, and on the other hand, we have unfolded isomorphisms $\theta_i: \left( U_i\times \textnormal{Spec}\left(A \right) , \tilde{\mathcal{N}}_{\mathcal{F}_i}^*\right)\to \left( \mathcal{X}|_{U_i}, \tilde{\mathcal{N}}_\mathcal{F}^*|_{U_i}\right)$ where $\theta_i$ is induced from $\tilde{\theta}_i$ and $\tilde{\mathcal{N}}_{\mathcal{F}_i}^*\subset \left(\Omega_X|_{U_i}\otimes A\right) \bigoplus \left( \mathcal{O}_X|_{U_i}\otimes \Omega_A^1\right)$ is induced from $\tilde{\mathcal{N}}_\mathcal{F}^*$ via $\theta_i$, and $\theta_{ij}=\theta_j^{-1}\theta_i$ induced from $\tilde{\theta}_{ij}$ corresponds to $f_{ij}:\Gamma\left( U_{ij}, \mathcal{O}_X\right)\otimes A \to \Gamma\left( U_{ij}, \mathcal{O}_X \right)\otimes A$ induced from $\tilde{f}_{ij}$, and $\Gamma\left( U_i, \tilde{\mathcal{N}}_{\mathcal{F}_i}^* \right)$ is generated by $\omega_i^1=w_i^1+ v_i^1,..., \omega_i^q= w_i^q + v_i^q\in \left( \Gamma\left( U_i, \Omega_X^1\right)\otimes A \right) \bigoplus \left( \Gamma\left( U_i, \mathcal{O}_X \right)\otimes \Omega_A^1 \right)$ such that $\tilde{w}_i^\alpha$ induces $w_i^\alpha$ which induces $w_{0i}^\alpha,\alpha=1,...,q$, and $\theta_{ji}^* \omega_i^\alpha= \sum_{\beta=1}^q h_{ij}^{\alpha\beta} \omega_j^\beta$ such that $\tilde{h}_{ij}^{\alpha\beta}$ induces $h_{ij}^{\alpha\beta}\in \Gamma\left( U_{ij} , \mathcal{O}_X \right)\otimes A$ which induces $h_{0ij}^{\alpha\beta}\in \Gamma\left( U_i, \mathcal{O}_X \right)$. Let $\tilde{v}_i^\alpha \in \Gamma\left(U_i, \mathcal{O}_X \right)\otimes \Omega_{\tilde{A}}^1$ be an arbitrary lifting of $v_i^\alpha,\alpha=1,...,q$. Then as in the proof of Proposition \ref{uu35} (we keep the notations in the proof of Proposition \ref{uu35}), the obstruction 
{\Small{\begin{align*}
&\left( \left\{B_i=0 \right\}, \left\{ F_i\right\}, \left\{ - \overline{R_{ij}}=0 \right\}, \left\{ - G_{ij}\right\}, \left\{ d_{ijk}=0 \right\} \right)  \\
&\in C^0\left(\mathcal{U}, \mathscr{H}om_{\mathcal{O}_X}\left( \mathcal{N}_{\mathcal{F}_0}^*, \tilde{\mathcal{S}}^2 \right) \right)\bigoplus C^0\left(\mathcal{U}, \mathscr{H}om_{\mathcal{O}_X}\left( \mathcal{N}_{\mathcal{F}_0}^*, \tilde{\mathcal{S}}^1  \right) \right) \bigoplus C^1\left(\mathcal{U}, \mathscr{H}om_{\mathcal{O}_X}\left( \mathcal{N}_{\mathcal{F}_0}^*, \frac{\Omega_X^1}{\mathcal{N}_{\mathcal{F}_0}^*} \right) \right)\\
& \bigoplus C^1\left(\mathcal{U}, \mathscr{H}om_{\mathcal{O}_X}\left(\mathcal{N}_{\mathcal{F}_0}^*, \mathcal{O}_X \right) \right) \bigoplus C^2\left( \mathcal{U}, T_X \right)
\end{align*}}}
lies in $\mathbb{H}^1\left( X, \mathscr{H}om_{\mathcal{O}_X}\left( \mathcal{N}_{\mathcal{F}_0}^*, \mathcal{O}_X \right)^\bullet \right)=0$, so that the lifting $\left( \tilde{\mathcal{X}} , \tilde{\mathcal{N}}_{\tilde{\mathcal{F}}}^* \right)$ over $\tilde{A}$ of $\left( \mathcal{X}, \tilde{\mathcal{N}}_{\mathcal{F}}^* \right)$ exists such that $\left( \tilde{\mathcal{X}}, \tilde{\mathcal{N}}_{\tilde{\mathcal{F}}}^*  \right)$ induces $\left( \tilde{\mathcal{X}}, \mathcal{N}_{\tilde{\mathcal{F}}}^*  \right)$. This completes the proof of Proposition \ref{a3}.
\end{proof}

\bibliographystyle{amsalpha}
\bibliography{References-Rev9}

\end{document}